БАКИНСКИЙ ГОСУДАРСТВЕННЫЙ УНИВЕРСИТЕТ

ПРЕПРИНТ № 1

САДЫГОВ М.А.

СУБДИФФЕРЕНЦИАЛ ВТОРОГО ПОРЯДКА.
ЭКСТРЕМАЛЬНЫЕ ЗАДАЧИ ДЛЯ ОПЕРАТОРНЫХ ВКЛЮЧЕНИЙ

БАКУ-2017

УДК 517.9

В работе изучается свойство субдифференциала второго порядка и получены условия оптимальности второго порядка для задачи минимизации. В работе также получены необходимые и достаточные условие экстремума для экстремальной задачи операторных включений.

# ВВЕДЕНИЕ

В работе изучается свойство субдифференциала второго порядка и получены условия оптимальности второго порядка для задачи минимизации. В работе также получены необходимые и достаточные условие экстремума для экстремальной задачи операторных включений.

Субдифференциал второго порядка изучен разными авторами, но в настоящее время однозначно принятое определение субдифференциала не имеется. Известно, что определение дифференциала второго порядка является аналогом дифференциала первого порядка. Но в определениях субдифференциала второго порядка такие аналогии в общем случае не имеются. Если рассмотреть субдифференциал Кларка, то можно получить разные аналогии субдифференциала второго порядка. Рассмотрим производную по направлению Кларка (см.[7],[29])

$$f^0(x_0;x) = \overline{\lim_{z \to x_0, \lambda \downarrow 0}} \frac{1}{\lambda}(f(z+\lambda x) - f(z)) \ .$$

Положив

$$f^{00}(x_0;x_1,x_2) = \overline{\lim_{z \to x_0, \lambda \downarrow 0}} \frac{1}{\lambda}(f^0(z+\lambda x_2;x_1) - f^0(z;x_1))$$

можно рассмотреть субдифференциал второго порядка следующего вида:

$$\partial_2^0 f(x_0) = \{b \in \overline{B}(X^2;R) : f^{00}(x_0;x_1,x_2) \geq b(x_1,x_2) \text{ при } (x_1,x_2) \in X^2\},$$

где через $\overline{B}(X^2;R)$ обозначено множество всех непрерывных симметричных билинейных функций из $X^2$ в $R$. Такое определение впервые введено автором в 1980 году, изучен ряд их свойств, обсуждены В.Ф.Демьяновым и А.М. Рубиновым и опубликовано только в 1988 году в более общем виде, т.е. рассмотрен субдифференциал произвольного порядка, где производная второго порядка по направлению $f^{00}(x_0;x_1,x_2)$ заменена через

$$f^{[2]}(x_0;x_1,x_2) = \overline{\lim_{\substack{z \to x_0, \\ \lambda_1 \downarrow 0, \lambda_2 \downarrow 0}}} \frac{1}{\lambda_1 \lambda_2}(f(z+\lambda_1 x_1 + \lambda_2 x_2) - f(z+\lambda_1 x_1) - f(z+\lambda_2 x_2) + f(z))$$

и рассмотрен субдифференциал второго порядка (см. [16],[17])

$$\partial_2 f(x_0) = \{b \in \overline{B}(X^2;R) : f^{[2]}(x_0;x_1,x_2) \geq b(x_1,x_2) \text{ при } (x_1,x_2) \in X^2\}.$$

Отметим, что

$$f^{00}(x_0;x_1,x_2) = \overline{\lim_{z \to x_0, t \downarrow 0}} \frac{1}{t}(f^0(z+tx_2;x_1) - f^0(z;x_1)) =$$

$$= \overline{\lim_{z \to x_0, t \downarrow 0}} \frac{1}{t}[\overline{\lim_{\upsilon \to z+tx_2, \lambda \downarrow 0}} \frac{1}{\lambda}(f(\upsilon+\lambda x_1) - f(\upsilon)) - \overline{\lim_{w \to z, \lambda \downarrow 0}} \frac{1}{\lambda}(f(w+\lambda x_1) - f(w))] =$$

$$= \overline{\lim_{z \to x_0, t \downarrow 0}} \frac{1}{t}[\overline{\lim_{w \to z,, \lambda \downarrow 0}} \frac{1}{\lambda}(f(w+tx_2+\lambda x_1) - f(w+tx_2)) - \overline{\lim_{w \to z, \lambda \downarrow 0}} \frac{1}{\lambda}(f(w+\lambda x_1) - f(w))] \leq$$

$$\leq \overline{\lim_{w \to x_0, t \downarrow 0, \lambda \downarrow 0}} \frac{1}{\lambda t}[f(w+tx_2+\lambda x_1) - f(w+tx_2) - f(w+\lambda x_1) + f(w)].$$

Близкое определение субдифференциала второго порядка также рассмотрено в [28] (см.также [30]). Отметим, что в [2] рассмотрено также идейно близкое определение субдифференциала второго порядка, но он различно и работать с ним трудно.

Хотя производные по направлению $f^{00}(x_0;x_1,x_2)$ и $f^{[2]}(x_0;x_1,x_2)$ близкие, но они в общем случае не совпадают. Производная по направлению $f^{[2]}(x_0;x_1,x_2)$ имеет ряд хороших свойств. При некоторых условиях $(x_1,x_2) \to f^{[2]}(x_0;x_1,x_2)$ бисублинейная функция (см.[16]).

Обозначим $U(x_0,\delta)=\{x\in X: \|x-x_0\|<\delta\}$. Если функция $f$ и производная Фреше $f'$ удовлетворяют условию Липшица на множестве $U(x_0,\delta)$, то

$$f^{00}(x_0;x_1,x_2)=\varlimsup_{z\to x_0,\lambda\downarrow 0}\frac{1}{\lambda}(f'(z+\lambda x_2)x_1-f'(z)x_1).$$

В работе [18],[19] рассмотрено также другое определение субдифференциала произвольного порядка. В частности отсюда следует, что

$$f^{\{2\}}(x_0;x_1,x_2)=\sup_{z_1,z_2\in X}\varlimsup_{\lambda_1\downarrow 0,\lambda_2\downarrow 0}\frac{1}{\lambda_1\lambda_2}(f(x_0+\lambda_1 z_1+\lambda_2 z_2+\lambda_1 x_1+\lambda_2 x_2)-$$
$$-f(x_0+\lambda_1 z_1+\lambda_2 z_2+\lambda_1 x_1)-f(x_0+\lambda_1 z_1+\lambda_2 z_2+\lambda_2 x_2)+f(x_0+\lambda_1 z_1+\lambda_2 z_2)),$$

и $\partial_{\{2\}}f(x_0)=\{b\in\overline{B}(X^2;R): f^{\{2\}}(x_0;x_1,x_2)\geq b(x_1,x_2)$ при $(x_1,x_2)\in X^2\}$.

Работа состоит из трех глав. В первой главе исследован ряд свойств субдифференциала второго порядка. В 1.1 изучен ряд свойств субдифференциала $\partial_{\{2\}}f(x_0)$.

В 1.2 рассматривается также следующее определение субдифференциала и изучен ряд их свойств. Положим

$$f^{\{2\}+}(x_0;x)=\sup_{z\in X}\varlimsup_{\lambda\downarrow 0}\frac{1}{\lambda^2}(f(x_0+\lambda z+2\lambda x)-2f(x_0+\lambda z+\lambda x)+f(x_0+\lambda z)),$$

$$f^{\{2\}-}(x_0;x)=\inf_{z\in X}\varliminf_{\lambda\downarrow 0}\frac{1}{\lambda^2}(f(x_0+\lambda z+2\lambda x)-2f(x_0+\lambda z+\lambda x)+f(x_0+\lambda z))$$

при $x\in X$.

Множество $D_2f(x_0)=\{Q\in B_0(X): f^{\{2\}-}(x_0;x)\leq Q(x)\leq f^{\{2\}+}(x_0;x), x\in X\}$ назовем бидифференциалом функции $f$ в точке $x_0$, где $B_0(X)$ множество всех непрерывных квадратичных функционалов.

В 1.3 получены условия оптимальности второго порядка. Отметим, что рассмотренный субдифференциал для получения условия оптимальности более эффективен. В негладком анализе существенную роль играет функция расстояния. Автор в работах [17],[19] показал, что при исследовании субдифференциала и условия оптимальности высокого порядка также существенную роль играет степень функции расстояния.

В гл.2 и гл.3 единая методика применяется к экстремальной задаче для операторных включений. Схема получения необходимых условий состоит из нескольких этапов. Сначала изучается непрерывная зависимость решения включения от возмущения. Затем исследуется выпуклая вариационная задача с оператором, заданным в соответствующем пространстве. Рассматривается возмущенная задача и изучается стабильность возмущенной задачи. Отметим, что можно рассматривать разные возмущенные задачи для экстремальной задачи операторных включений. Но естественная возмущенная задача упрощает исследование экстремальных задач. В работе изучаются субдифференцируемость интегрального, терминального и граничного функционала с оператором заданным в функциональном пространстве.

Далее рассматривается выпуклая экстремальная задача для операторного включения. Исследование выпуклый экстремальной задачи для операторных включений аналогично исследованию выпуклой вариационной задачи для интегральных и дифференциальных включений. При этом существенно используется теорема о непрерывной зависимости решения операторного включения от возмущения.

Также используя теорему о непрерывной зависимости решения операторного включения от возмущения исследуется невыпуклая экстремальная задача для включений.

Отметим, что в экстремальной задаче для включений трудность состоит в том, что в общем случае не понятны типы включений, для которых верна теорема о непрерывной зависимости от возмущения решений включения.

В работе изучены зависимости решений операторных включений от возмущения и исследована экстремальная задача для операторных включений. Используя методы выпуклого и негладкого анализа, получены необходимые и достаточные условия оптимальности. Хотя идея исследования оптимальности экстремальной задачи для операторных включений ближе к исследованию экстремальной задачи для интегральных и дифференциальных включений(см.[17]-[24]), но имеются также существенные отличия.

В зависимости от функционального пространства можно рассмотреть разные экстремальные задачи для операторных включений.

В главе 2, состоящей из восьми параграфов, исследована экстремальная задача для операторного включения. В 2.1 изучена непрерывная зависимость решения от возмущения операторного включения, которая ближе к включению типа Вольтерра. В 2.2 изучены непрерывная зависимость решения от возмущения операторного включения, которая ближе к включению типа Фредгольма. В 2.3 получены необходимые и достаточные условия экстремума выпуклых вариационных задач, которые зависят от оператора. В 2.4 изучена субдифференцируемость интегрального функционала, которая зависят от оператора. В 2.5 изучена субдифференцируемость терминального функционала, которая зависит от оператора. В 2.6 изучена экстремальная задача для выпуклых вариационных задач с оператором. Получены необходимые и достаточные условия экстремума выпуклых вариационных задач с оператором. В 2.7 получены необходимые и достаточные условия минимума для экстремальных задач выпуклых опраторных включений. В 2.8 изучена невыпуклая экстремальная задача для опраторных включений.

В главе 3 исследована экстремальная задача для многомерных операторных включений. Глава 3 состоит из пяти параграфов. В 3.1 изучена непрерывная зависимость от возмущения решения многомерного операторного включения. В 3.2 изучена субдифференцируемость интегрального функционала с оператором. В 3.3 изучена субдифференцируемость граничного функционала с оператором. В 3.4 рассматриваются вопросы минимизации многомерных вариационных задач с оператором. В 3.5 рассматривается невыпуклая экстремальная задача для операторных включений.

# Гл.1. Субдифференциал второго порядка и условия оптимальности

## 1.1. 2-субдифференциал и бисубдифференциал

В разделе 1.1 даны определения 2-субдифференциала, бисубдифференциала и изучены их свойства.

Пусть $f: X \to R$, где $X$ -банахово пространство. Положим

$$f^{\{2\}}(x_0;x_1,x_2) = \sup_{z_1,z_2 \in X} \overline{\lim_{\lambda_1 \downarrow 0, \lambda_2 \downarrow 0}} \frac{1}{\lambda_1 \lambda_2}(f(x_0 + \lambda_1 z_1 + \lambda_2 z_2 + \lambda_1 x_1 + \lambda_2 x_2) - f(x_0 + \lambda_1 z_1 + \lambda_2 z_2 + \lambda_1 x_1) -$$
$$- f(x_0 + \lambda_1 z_1 + \lambda_2 z_2 + \lambda_2 x_2) + f(x_0 + \lambda_1 z_1 + \lambda_2 z_2)),$$

$$f_{\{2\}}(x_0;x_1,x_2) = \inf_{z_1,z_2 \in X} \underline{\lim_{\lambda_1 \downarrow 0, \lambda_2 \downarrow 0}} \frac{1}{\lambda_1 \lambda_2}(f(x_0 + \lambda_1 z_1 + \lambda_2 z_2 + \lambda_1 x_1 + \lambda_2 x_2) - f(x_0 + \lambda_1 z_1 + \lambda_2 z_2 + \lambda_1 x_1) -$$
$$- f(x_0 + \lambda_1 z_1 + \lambda_2 z_2 + \lambda_2 x_2) + f(x_0 + \lambda_1 z_1 + \lambda_2 z_2)).$$

Из определения непосредственно следует, что $f^{\{2\}}(x_0;x_1,x_2) = f^{\{2\}}(x_0;x_2,x_1)$, $f_{\{2\}}(x_0;x_1,x_2) = f_{\{2\}}(x_0;x_2,x_1)$, $f^{\{2\}}(x_0;-x_1,-x_2) = f^{\{2\}}(x_0;x_1,x_2)$ и $f_{\{2\}}(x_0;-x_1,-x_2) = f_{\{2\}}(x_0;x_1,x_2)$ при $(x_1,x_2) \in X^2$, т.е. $f^{\{2\}}(x_0;x_1,x_2)$ и $f_{\{2\}}(x_0;x_1,x_2)$ четные и симметричные функции.

Используя определения $f_{\{2\}}(x_0;x_1,x_2)$ и заменив $z_1$ через $z_1 - x_1$ имеем

$$-f_{\{2\}}(x_0;x_1,x_2) = \sup_{z_1,z_2 \in X} \overline{\lim_{\lambda_1 \downarrow 0, \lambda_2 \downarrow 0}} \frac{1}{\lambda_1 \lambda_2}(-f(x_0 + \lambda_1 z_1 + \lambda_2 z_2 + \lambda_1 x_1 + \lambda_2 x_2) +$$
$$+ f(x_0 + \lambda_1 z_1 + \lambda_2 z_2 + \lambda_1 x_1) + f(x_0 + \lambda_1 z_1 + \lambda_2 z_2 + \lambda_2 x_2) - f(x_0 + \lambda_1 z_1 + \lambda_2 z_2)) =$$
$$= \sup_{z_1,z_2 \in X} \overline{\lim_{\lambda_1 \downarrow 0, \lambda_2 \downarrow 0}} \frac{1}{\lambda_1 \lambda_2}(f(x_0 + \lambda_1 z_1 + \lambda_2 z_2 - \lambda_1 x_1 + \lambda_2 x_2) - f(x_0 + \lambda_1 z_1 + \lambda_2 z_2 - \lambda_1 x_1) -$$
$$- f(x_0 + \lambda_1 z_1 + \lambda_2 z_2 + \lambda_2 x_2) + f(x_0 + \lambda_1 z_1 + \lambda_2 z_2)) = f^{\{2\}}(x_0;-x_1,x_2).$$

Отсюда следует, что $f_{\{2\}}(x_0;x_1,x_2) = -f^{\{2\}}(x_0;-x_1,x_2)$ при $(x_1,x_2) \in X^2$.

Положив $z_1 = y_1 - x_1$, $z_2 = y_2 - x_2$, где $y_1, y_2 \in X^2$, из определения $f^{\{2\}}(x_0;x_1,x_2)$ имеем, что $f^{\{2\}}(x_0;x_1,x_2) = f^{\{2\}}(x_0;-x_1,-x_2)$ при $(x_1,x_2) \in X^2$.

Множество всех непрерывных линейных операторов из $X$ в $X^*$ обозначим через $L(X,X^*)$, а множество всех непрерывных билинейных отображений из $X \times X$ в $R$, обозначим через $B(X^2,R)$. Множество всех билинейных симметричных непрерывных функций из $X \times X$ в $R$ обозначим через $\overline{B}(X^2,R)$. Соответствие между $B(X^2,R)$ и $L(X,X^*)$, определяемое равенством $x^*(x_1,x_2) = (Ax_1)x_2$ взаимно однозначно. Функционал $Q(x)$, заданный в пространстве X, называется квадратичным функционалом, если существует такой билинейный симметричный функционал $x^* \in \overline{B}(X^2,R)$, что $Q(x)=x^*(x,x)$. Отметим, что билинейный симметричный функционал $x^*$ определяется по $Q(x)$ однозначно и $x^*(x_1,x_2) = \frac{1}{2}(Q(x_1+x_2) - Q(x_1) - Q(x_2))$. Поэтому в дальнейшем симметричный билинейный функционал $x^*$, соответствующий симметричному оператору $A \in L(X, X^*)$ и квадратичному функционалу $Q(x)$ будем отождествлять.

Непрерывный билинейный симметричный функционал $x^*$, удовлетворяющий неравенству $f^{\{2\}}(x_0; x_1, x_2) \geq x^*(x_1, x_2)$ при $x_1, x_2 \in X$ назовем 2-субградиентом функции $f$ в точке $x_0$, а множество 2-субградиентов в точке $x_0$ назовем 2-субдифференциалом функции $f$ в точке $x_0$ и обозначим через $\partial_{\{2\}} f(x_0)$ (см. [19], [23]).

Так как $f_{\{2\}}(x_0; x_1, x_2) = -f^{\{2\}}(x_0; -x_1, x_2)$ при $(x_1, x_2) \in X^2$, то если $x^* \in \partial_{\{2\}} f(x_0)$, то $f_{\{2\}}(x_0; x_1, x_2) \leq x^*(x_1, x_2) \leq f^{\{2\}}(x_0; x_1, x_2)$.

Если функции $x_1 \to g(x_1, x_2)$, $x_2 \to g(x_1, x_2)$ выпуклые и положительно однородные, то функцию $g$ назовем бисублинейной. Если функции $x_1 \to g(x_1, x_2)$, $x_2 \to g(x_1, x_2)$ положительно однородные и $g(0, x_2) = g(x_1, 0) = 0$, то функцию $g$ назовем биположительно однородной.

Отметим, что $\|x_1\| \cdot \|x_2\|$, $g(x_1, x_2) = \sup_{i \in I} b_i(x_1, x_2)$, где $b_i \in B(X^2, R)$ при $i \in I$ и $g(x_1, x_2) < +\infty$ при $(x_1, x_2) \in X^2$, бисублинейные функции.

Если $f(x) = b(x, x)$, где $b \in \overline{B}(X^2, R)$, то $f^{\{2\}}(x_0; x_1, x_2) = 2b(x_1, x_2)$.

Функцию $f : X \to R$ назовем 2-липшицевой с постоянной $K$ в точке $x_0$, если $f$ для некоторого $\varepsilon > 0$ удовлетворяет условию
$$|f(x_0 + x_1 + x_2) - f(x_0 + x_1) - f(x_0 + x_2) + f(x_0)| \leq K \|x_1\| \|x_2\|$$
при $x_1, x_2 \in \varepsilon B$, где $B = \{x \in X : \|x\| \leq 1\}$.

Если для некоторого $\varepsilon > 0$ удовлетворяется неравенство
$$|f(x + x_1 + x_2) - f(x + x_1) - f(x + x_2) + f(x)| \leq K \|x_1\| \|x_2\|$$
при $x \in x_0 + \varepsilon B$ и $x_1, x_2 \in \varepsilon B$, то функцию $f : X \to R$ назовем 2-липшицевой с постоянной $K$ в окрестности точки $x_0$.

Легко проверяется, что если $f : X \to R$ 2-липшицевая с постоянной $K$ в окрестности точки $x_0$ функция, то функция $f^{\{2\}}(x_0; x_1, x_2)$ бисублинейна, раздельно непрерывна и удовлетворяется неравенство $|f^{\{2\}}(x_0; x_1, x_2)| \leq K \|x_1\| \|x_2\|$ при $(x_1, x_2) \in X \times X$.

Покажем, что вычислить $f^{\{2\}}(x_0; x_1, x_2)$ трудно. Пусть $\alpha, \beta \in R$. Из определения $f^{\{2\}}(x_0; x_1, x_2)$ следует, что

$$f^{\{2\}}(x_0; x_1, x_2) \geq \varlimsup_{\lambda_1 \downarrow 0, \lambda_2 \downarrow 0} \frac{1}{\lambda_1 \lambda_2} (f(x_0 + \lambda_1 \alpha x_1 + \lambda_2 \beta x_2 + \lambda_1 x_1 + \lambda_2 x_2) - f(x_0 + \lambda_1 \alpha x_1 + \lambda_2 \beta x_2 + \lambda_1 x_1) -$$
$$- f(x_0 + \lambda_1 \alpha x_1 + \lambda_2 \beta x_2 + \lambda_2 x_2) + f(x_0 + \lambda_1 \alpha x_1 + \lambda_2 \beta x_2)) = \varlimsup_{\lambda_1 \downarrow 0, \lambda_2 \downarrow 0} \frac{1}{\lambda_1 \lambda_2} (f(x_0 + \lambda_1 (1+\alpha) x_1 +$$
$$+ \lambda_2 (1+\beta) x_2) - f(x_0 + \lambda_1 (1+\alpha) x_1 + \lambda_2 \beta x_2) - f(x_0 + \lambda_1 \alpha x_1 + \lambda_2 (1+\beta) x_2) + f(x_0 + \lambda_1 \alpha x_1 + \lambda_2 \beta x_2)).$$

Обозначив
$$f_{\alpha, \beta}(x_0; x_1, x_2) = \varlimsup_{\lambda_1 \downarrow 0, \lambda_2 \downarrow 0} \frac{1}{\lambda_1 \lambda_2} (f(x_0 + \lambda_1 (1+\alpha) x_1 + \lambda_2 (1+\beta) x_2) - f(x_0 + \lambda_1 (1+\alpha) x_1 + \lambda_2 \beta x_2) -$$
$$- f(x_0 + \lambda_1 \alpha x_1 + \lambda_2 (1+\beta) x_2) + f(x_0 + \lambda_1 \alpha x_1 + \lambda_2 \beta x_2))$$
имеем, что $f^{\{2\}}(x_0; x_1, x_2) \geq f_{\alpha, \beta}(x_0; x_1, x_2)$ при $x_1, x_2 \in X$. Поэтому

$$f^{\{2\}}(x_0;x_1,x_2) \geq \sup_{\alpha,\beta \in R} f_{\alpha,\beta}(x_0;x_1,x_2)$$

при $x_1, x_2 \in X$.

Положим $\bar{\partial}_{\{2\}}f(x_0) = \{x^* \in B(X^2,R) : f^{\{2\}}(x_0;x_1,x_2) \geq x^*(x_1,x_2), x_1,x_2 \in X\}$.

Функцию $f$ назовем второго порядка s-дифференцируемой в точке $x_0$, если существует оператор $A \in L(X,X^*)$ такой, что

$$\lim_{\substack{z \to x_0, \\ \lambda_1 \downarrow 0, \lambda_2 \downarrow 0}} \frac{1}{\lambda_1 \lambda_2}(f(z+\lambda_1 x_1+\lambda_2 x_2) - f(z+\lambda_1 x_1) - f(z+\lambda_2 x_2) + f(z)) = \langle Ax_1, x_2 \rangle$$

при $x_1, x_2 \in X$. Ясно, что $A$ симметричный оператор.

Если функция $f$ второго порядка s-дифференцируема в точке $x_0$, то

$$\partial_2 f(x_0) = \partial_{\{2\}} f(x_0) = \{\langle Ax_1, x_2 \rangle\}.$$

Легко проверяется, что если $f:X \to R$ 2-липшицевая с постоянной $K$ в окрестности точки $x_0$ функция, то функция $f$ второго порядка s-дифференцируема в точке $x_0$ тогда и только тогда когда $\partial_2 f(x_0)$ состоит из единственного элемента.

**Следствие 1.** Если $f$ 2-липшицевая функция в окрестности точки $x_0$, то

$$f^{\{2\}}(x_0;x_1,x_2) = \max\{x^*(x_1,x_2) : x^* \in \bar{\partial}_{\{2\}}f(x_0)\} = \max\{x^*(x_1,x_2) : x^* \in \partial_{\{2\}}f(x_0)\}.$$

Справедливость следствия 1 следует из [23] (см.с.92 и с.89).

Непрерывный квадратичный функционал $Q(x)$, удовлетворяющий неравенству $f_{\{2\}}(x_0;x,x) \leq Q(x) \leq f^{\{2\}}(x_0;x,x)$ назовем бисубградиентом функции $f$ в точке $x_0$, а множество бисубградиентов в точке $x_0$ назовем бисубдифференциалом функции $f$ в точке $x_0$ и обозначим $d_{\{2\}}f(x_0)$.

Множество всех непрерывных квадратичных функционалов обозначим $B_0(X)$.

**Следствие 2**. Если $f$ 2-липшицевая функция в окрестности точки $x_0$, то

$$f^{\{2\}}(x_0;x;x) = \max\{Q(x) : Q \in d_{\{2\}}f(x_0)\}, \quad f_{\{2\}}(x_0;x;x) = \min\{Q(x) : Q \in d_{\{2\}}f(x_0)\}.$$

Если учесть, что $\partial_{\{2\}}f(x_0) \subset d_{\{2\}}f(x_0)$, то справедливость следствия 2 вытекает из следствия 1 (см. также [23], с.89).

Если $C \subset X$, то положим $d_C(y) = \inf\{\|y-z\| : z \in C\}$ и $d_2(y) = d_C^2(y)$.

Пусть $(x_1,x_2) \in X^2$, $\|x_2\| \leq \|x_1\|$ и $z \in C$. Тогда имеем

$$|d_2(z+x_1+x_2) - d_2(z+x_1) - d_2(z+x_2) + d_2(z)| \leq |d_2(z+x_1+x_2) - d_2(z+x_1)| +$$
$$+ |d_2(z+x_2) - d_2(z)| = |d(z+x_1+x_2) - d(z+x_1)||d(z+x_1+x_2) + d(z+x_1)| +$$
$$+ |d(z+x_2) - d(z)||d(z+x_2) + d(z)| \leq \|x_2\|(|d(z+x_1+x_2) - d(z)| + |d(z+x_1) - d(z)|) +$$
$$+ |d(z+x_2) - d(z)|^2 \leq \|x_2\|(\|x_1+x_2\| + \|x_1\|) + \|x_2\|^2 \leq \|x_2\|(\|x_1\| + \|x_2\| + \|x_1\|) + \|x_1\|\|x_2\| \leq 4\|x_1\|\|x_2\|.$$

Если $X$ гильбертово пространство, $C \subset X$ выпуклое замкнутое множество и $x_1,x_2,z \in X$, то (см. [16], стр. 9)

$$|d_2(z+x_1+x_2) - d_2(z+x_1) - d_2(z+x_2) + d_2(z)| \leq 6\|x_1\|\|x_2\|.$$

Положим

$$Q_C(x_0) = \{(x_1,x_2) \in X \times X : d_2^{\{2\}+}(x_0;x_1,x_2) \leq 0\},$$
$$\Omega_C(x_0) = \{b \in \bar{B}(X^2;R) : b(x_1,x_2) \leq 0 \text{ при } (x_1,x_2) \in Q_C(x_0)\}.$$

**Лемма 1.** Если $d_C^2(x)$ удовлетворяет 2-липшицеву условию в окрестности точки $x_0$, то $\bigcup_{\lambda>0} \lambda \partial_{\{2\}} d_2(x_0) \subset \Omega_C(x_0)$.

**Доказательство.** Пусть $b \in \bigcup_{\lambda>0} \lambda \partial_{\{2\}} d_2(x_0)$. Тогда существует $\lambda > 0$ такое, что $b \in \lambda \partial_{\{2\}} d_2(x_0)$. Поэтому $\frac{b}{\lambda} = b_1 \in \partial_{\{2\}} d_2(x_0)$ и $d_2^{\{2\}+}(x_0; x_1, x_2) \geq b_1(x_1, x_2)$ при $(x_1, x_2) \in X \times X$. Так как из $(x_1, x_2) \in Q_C(x_0)$ следует, что $d_2^{\{2\}+}(x_0; x_1, x_2) \leq 0$. Поэтому $b_1(x_1, x_2) \leq 0$ при $(x_1, x_2) \in Q_C(x_0)$. Тогда имеем, что $b(x_1, x_2) = \lambda b_1(x_1, x_2) \leq 0$ при $(x_1, x_2) \in Q_C(x_0)$. Отсюда следует, что $b \in \Omega_C(x_0)$, т.е. $\bigcup_{\lambda>0} \lambda \partial_{\{2\}} d_2(x_0) \subset \Omega_C(x_0)$. Лемма доказана.

В доказательстве леммы 1 условие, $d_C^2(x)$ удовлетворяет 2-липшицеву условию в окрестности точки $x_0$, не используется. Если $d_C^2(x)$ удовлетворяет 2-липшицеву условию в окрестности точки $x_0$, то $\partial_{\{2\}} d_2(x_0)$ непусто.

Отметим, что множество $K \subset X \times X$ называется биконусом (см.[23]), если для любого $(x, y) \in K$ множества $K_y = \{x : (x, y) \in K\}$ и $K_x = \{y : (x, y) \in K\}$ выпуклые конусы.

**Лемма 2.** Если $d_C^2(x)$ удовлетворяет 2-липшицеву условию в окрестности точки $x_0$, то $Q_C(x_0)$ является биконусом.

**Доказательство.** Пусть $(x, y) \in Q_C(x_0)$. Покажем, что $K_y = \{x \in X : (x, y) \in Q_C(x_0)\}$ является конусом. Если $x \in K_y$ и $\lambda \geq 0$, то имеем, что $d_2^{\{2\}+}(x_0; \lambda x, y) = \lambda d_2^{\{2\}+}(x_0; x, y) \leq 0$, т.е. $\lambda x \in K_y$. Если $x_1, x_2 \in K_y$, то имеем, что
$$d_2^{\{2\}+}(x_0; x_1 + x_2, y) \leq d_2^{\{2\}+}(x_0; x_1, y) + d_2^{\{2\}+}(x_0; x_2, y) \leq 0.$$
Отсюда следует, что $x_1 + x_2 \in K_y$. Получим, что $K_y$ является конусом. Аналогично проверяется, что $K_x = \{y \in X : (x, y) \in Q_C(x_0)\}$ является конусом. Лемма доказана.

Пусть $g$ функция из $X$ в $R$ и $Q \in B_0(X)$. Положим
$$g_+^*(Q) = \sup_{x \in X} \{Q(x) - g(x)\}, \quad g_-^*(Q) = \inf_{x \in X} \{Q(x) - g(x)\}.$$

**Теорема 1.** Если $f$ 2-липшицевая функция в окрестности точки $x_0$, то $Q \in d_{(2)}f(x_0)$ тогда и только тогда, когда $Q \in \text{dom} \varphi_+^* \cap \text{dom} \psi_-^*$, где $\varphi(x) = f^{\{2\}}(x_0; x, x)$, $\psi(x) = f_{\{2\}}(x_0; x, x)$.

**Доказательство.** Если $Q \in d_{\{2\}}f(x_0)$, то из следствия 2 получим
$$\min\{Q'(x) : Q' \in d_{\{2\}}f(x_0)\} \leq Q(x) \leq \max\{Q''(x) : Q'' \in d_{\{2\}}f(x_0)\}.$$
Поэтому $Q(x) - \max_{Q \in d_{\{2\}}f(x_0)} Q(x) \leq 0 \leq Q(x) - \min_{Q \in d_{\{2\}}f(x_0)} Q(x)\}$. Отсюда получим, что
$$\sup_x \{Q(x) - \max_{Q \in d_{\{2\}}f(x_0)} Q(x)\} = \inf_x \{Q(x) - \min_{Q \in d_{\{2\}}f(x_0)} Q(x)\} = 0,$$
т.е. $\varphi_+^*(Q) = \psi_-^*(Q) = 0$ при $Q \in d_{\{2\}}f(x_0)$. Поэтому $Q \in \text{dom} \varphi_+^* \cap \text{dom} \psi_-^*$. Так как
$$\varphi_+^*(Q) = \sup_x \{Q(x) - \varphi(x)\}, \quad \psi_-^*(Q) = \inf_x \{Q(x) - \psi(x)\},$$
то $\varphi_+^*(Q) \geq 0$, $\psi_-^*(Q) \leq 0$ при $Q \in B_0(X)$. Пусть существуют $x_1$ и $x_2$ такие, что $Q(x_1) - \varphi(x_1) > 0$ и $Q(x_2) - \psi(x_2) < 0$. Тогда
$$\varphi_+^*(Q) \geq \sup_{\lambda \geq 0} \{Q(\lambda x_1) - \varphi(\lambda x_1)\} = \sup_{\lambda \geq 0} \lambda^2 \{Q(x_1) - \varphi(x_1)\} = +\infty,$$
$$\psi_-^*(Q) \leq \inf_{\lambda \geq 0} \{Q(\lambda x_2) - \psi(\lambda x_2)\} = \inf_{\lambda \geq 0} \lambda^2 \{Q(x_2) - \psi(x_2)\} = -\infty.$$

Поэтому

$$\varphi_+^*(Q) = \begin{cases} 0, & Q \in \mathrm{dom}\varphi_+^* \\ +\infty, & Q \notin \mathrm{dom}\varphi_+^* \end{cases}, \quad \psi_-^*(Q) = \begin{cases} 0, & Q \in \mathrm{dom}\psi_-^* \\ -\infty, & Q \notin \mathrm{dom}\psi_-^* \end{cases}.$$

Отсюда следует, что если $Q \in \mathrm{dom}\varphi_+^* \cap \mathrm{dom}\psi_-^*$, то $\varphi_+^*(Q) = \psi_-^*(Q) = 0$. Поэтому $Q \in d_{\{2\}}f(x_0)$. Теорема доказана.

Если $f_1 : X \to R$ и $f_2 : X \to R$ удовлетворяют 2-липшицеву условию с постоянной K, в окрестности точки $x_0$, то непосредственно проверяется, что

$$(f_1 + f_2)^{\{2\}}(x_0; x_1, x_2) \le f_1^{\{2\}}(x_0; x_1, x_2) + f_2^{\{2\}}(x_0; x_1, x_2)$$

при $x_1, x_2 \in X$. Поэтому $\partial_{\{2\}}(f_1 + f_2)^{\{2\}}(x_0) \subset \partial_{\{2\}}(f_1^{\{2\}}(x_0; \cdot) + f_2^{\{2\}}(x_0; \cdot))$, где

$$\partial_{\{2\}}(f_1^{\{2\}}(x_0; \cdot) + f_2^{\{2\}}(x_0; \cdot)) = \{b \in \overline{B}(X^2, R) : f_1^{\{2\}}(x_0; x_1, x_2) + f_2^{\{2\}}(x_0; x_1, x_2) \ge b(x_1, x_2), x_1, x_2 \in X\}.$$

Пусть $A \subset \overline{B}(X^2, R)$. S–оболочкой множества A назовём множество

$$SA = \{b \in \overline{B}(X^2, R) : \sup\{\tilde{b}(x_1, x_2) : \tilde{b} \in A\} \ge b(x_1, x_2) \text{ при } x_1, x_2 \in X\}.$$

Используя S–оболочку множества имеем, что

$$\partial_{\{2\}}(f_1 + f_2)^{\{2\}}(x_0) \subset S(\partial_{\{2\}}f_1(x_0) + \partial_{\{2\}}f_2(x_0)).$$

**Лемма 3.** Если $f$ удовлетворяет 2-липшицеву условию с постоянной K в окрестности точки $x_0$, то множества $\partial_{\{2\}}f(x_0)$ и $d_{\{2\}}f(x_0)$ ограничены.

Если $x^* \in \partial_{\{2\}}f(x_0)$, то $-f^{\{2\}}(x_0; -x_1, x_2) \le x^*(x_1, x_2) \le f^{\{2\}}(x_0; x_1, x_2)$. Поэтому легко проверяется, что множество $\partial_{\{2\}}f(x_0)$ ограничено.

Пусть $Q \in d_{\{2\}}f(x_0)$ и $x^*(x_1, x_2) = \frac{1}{2}(Q(x_1 + x_2) - Q(x_1) - Q(x_2))$. Легко проверяется, что

$$-f^{\{2\}}(x_0; -x_1 - x_2, x_1 + x_2) - f^{\{2\}}(x_0; x_1, x_1) - f^{\{2\}}(x_0; x_2, x_2) \le 2x^*(x_1, x_2) \le$$
$$\le f^{\{2\}}(x_0; x_1 + x_2, x_1 + x_2) + f^{\{2\}}(x_0; -x_1, x_1) + f^{\{2\}}(x_0; -x_2, x_2)$$

при $x_1, x_2 \in X$. Поэтому множество $d_{\{2\}}f(x_0)$ ограничено. Лемма доказана.

**Лемма 4.** Пусть J конечное множество, $f_i$ ($i \in J$) непрерывны и удовлетворяют 2-липшицеву условию в окрестности точки $x_0$; $f(x) = \max_{i \in J} f_i(x)$. Тогда $f^{\{2\}}(x_0; x_1, x_2) \le \max_{i \in J(x_0)} f_i^{\{2\}}(x_0; x_1, x_2)$, где $J(x_0) = \{i \in J : f_i(x_0) = f(x_0)\}$.

**Доказательство.** Покажем, что существует $\alpha > 0$ такое, что $J(y) \subset J(x_0)$ при $\|y - x_0\| \le \alpha$. Случай $J(x_0) = J$ тривиален. Пусть $J(x_0) \ne J$. В самом деле, пусть $a = f(x_0) - \max_{i \in J \setminus J(x_0)} f_i(x_0) > 0$, $\varepsilon = a/3$ и $\alpha > 0$ такое, что для всякого $i \in J$ выполнено $|f_i(y) - f_i(x_0)| \le \varepsilon$ при $\|y - x_0\| \le \alpha$. Легко можно проверить, что $|f(y) - f(x_0)| \le \varepsilon$ при $\|y - x_0\| \le \alpha$. Если $j \in J(y)$, то
$f_j(x_0) \ge f_j(y) - \varepsilon = f(y) - \varepsilon \ge f(x_0) - 2\varepsilon = \max_{i \notin J(x_0)} f_i(x_0) + a - 2\varepsilon > \max_{i \notin J(x_0)} f_i(x_0)$.

Поэтому $j \in J(x_0)$. Возьмём $z_1, z_2 \in X$ и $x_1, x_2 \in X$. Пусть $\lambda_1 > 0, \lambda_2 > 0$ такие, что $\|\lambda_1 z_1 + \lambda_2 z_2\| \le \alpha$. Ясно, что $J(x_0 + \lambda_1 z_1 + \lambda_2 z_2) \subset J(x_0)$. Аналогично имеем, что существует $\delta > 0$ такое, что $J(y) \subset J(x_0 + \lambda_1 z_1 + \lambda_2 z_2)$ при $\|y - x_0 - \lambda_1 z_1 - \lambda_2 z_2\| \le \delta$. Пусть $\mu > 0$ достаточно малое число такое, что $\|\lambda_1 \mu x_1 + \lambda_2 \mu x_2\| \le \delta$.

$$f(x_0 + \lambda_1 z_1 + \lambda_2 z_2 + \lambda_1 \mu x_1 + \lambda_2 \mu x_2) - f(x_0 + \lambda_1 z_1 + \lambda_2 z_2 + \lambda_1 \mu x_1) -$$
$$- f(x_0 + \lambda_1 z_1 + \lambda_2 z_2 + \lambda_2 \mu x_2) - f(x_0 + \lambda_1 z_1 + \lambda_2 z_2) =$$
$$= \max_{i \in J(x_0 + \lambda_1 z_1 + \lambda_2 z_2)} (f(x_0 + \lambda_1 z_1 + \lambda_2 z_2 + \lambda_1 \mu x_1 + \lambda_2 \mu x_2) +$$
$$+ f_i(x_0 + \lambda_1 z_1 + \lambda_2 z_2)) - \max_{i \in J} f_i(x_0 + \lambda_1 z_1 + \lambda_2 z_2 + \lambda_1 \mu x_1) -$$
$$- \max_{i \in J} f_i(x_0 + \lambda_1 z_1 + \lambda_2 z_2 + \lambda_2 \mu x_2) \leq \max_{i \in J(x_0 + \lambda_1 z_1 + \lambda_2 z_2)} (f_i(x_0 + \lambda_1 z_1 + \lambda_2 z_2 + \lambda_1 \mu x_1 + \lambda_2 \mu x_2) +$$
$$+ f_i(x_0 + \lambda_1 z_1 + \lambda_2 z_2) - f_i(x_0 + \lambda_1 z_1 + \lambda_2 z_2 + \lambda_1 \mu x_1) - f_i(x_0 + \lambda_1 z_1 + \lambda_2 z_2 + \lambda_2 \mu x_2)) \leq$$
$$\leq \max_{i \in J(x_0)} (f_i(x_0 + \lambda_1 z_1 + \lambda_2 z_2 + \lambda_1 \mu x_1 + \lambda_2 \mu x_2) + f_i(x_0 + \lambda_1 z_1 + \lambda_2 z_2) -$$
$$- f_i(x_0 + \lambda_1 z_1 + \lambda_2 z_2 + \lambda_1 \mu x_1) - f_i(x_0 + \lambda_1 z_1 + \lambda_2 z_2 + \lambda_2 \mu x_2)).$$

Поэтому получим, что

$$f^{\{2\}}(x_0; \beta x_1, \beta x_2) \leq d = \sup_{z_1, z_2 \in X} \inf_{\substack{\alpha > 0, \\ \beta > 0}} \max_{i \in J(x_0)} \sup_{\substack{0 < \lambda_1 \leq \alpha, \\ 0 < \lambda_2 \leq \beta}} \frac{1}{\lambda_1 \lambda_2} (f_i(x_0 + \lambda_1 z_1 + \lambda_2 z_2 + \lambda_1 \mu x_1 + \lambda_2 \mu x_2) -$$
$$- f_i(x_0 + \lambda_1 z_1 + \lambda_2 z_2 + \lambda_1 \mu x_1) - f_i(x_0 + \lambda_1 z_1 + \lambda_2 z_2 + \lambda_2 \mu x_2) + f_i(x_0 + \lambda_1 z_1 + \lambda_2 z_2)).$$

Из определения супремума следует, что для $\varepsilon > 0$ существуют $z_1^\varepsilon, z_2^\varepsilon \in X$ такие, что

$$\inf_{\substack{\alpha > 0, \\ \beta > 0}} \max_{i \in J(x_0)} \sup_{\substack{0 < \lambda_1 \leq \alpha, \\ 0 < \lambda_2 \leq \beta}} \frac{1}{\lambda_1 \lambda_2} (f_i(x_0 + \lambda_1 z_1^\varepsilon + \lambda_2 z_2^\varepsilon + \lambda_1 \mu x_1 + \lambda_2 \mu x_2) - f_i(x_0 + \lambda_1 z_1^\varepsilon + \lambda_2 z_2^\varepsilon + \lambda_1 \mu x_1) -$$
$$- f_i(x_0 + \lambda_1 z_1^\varepsilon + \lambda_2 z_2^\varepsilon + \lambda_2 \mu x_2) + f_i(x_0 + \lambda_1 z_1^\varepsilon + \lambda_2 z_2^\varepsilon))) > d - \varepsilon.$$

Ясно, что

$$\inf_{\substack{\alpha > 0, \\ \beta > 0}} \sup_{\substack{0 < \lambda_1 \leq \alpha, \\ 0 < \lambda_2 \leq \beta}} \frac{1}{\lambda_1 \lambda_2} (f_i(x_0 + \lambda_1 z_1^\varepsilon + \lambda_2 z_2^\varepsilon + \lambda_1 \mu x_1 + \lambda_2 \mu x_2) - f_i(x_0 + \lambda_1 z_1^\varepsilon + \lambda_2 z_2^\varepsilon + \lambda_1 \mu x_1) -$$
$$- f_i(x_0 + \lambda_1 z_1^\varepsilon + \lambda_2 z_2^\varepsilon + \lambda_2 \mu x_2) + f_i(x_0 + \lambda_1 z_1^\varepsilon + \lambda_2 z_2^\varepsilon)) \leq f_i^{\{2\}}(x_0; \mu x_1, \mu x_2).$$

Поэтому для $\varepsilon > 0$ существуют $\alpha_i^\varepsilon > 0$, $\beta_i^\varepsilon > 0$ такие, что

$$\sup_{\substack{0 < \lambda_1 \leq \alpha_i^\varepsilon, \\ 0 < \lambda_2 \leq \beta_i^\varepsilon}} \frac{1}{\lambda_1 \lambda_2} (f_i(x_0 + \lambda_1 z_1^\varepsilon + \lambda_2 z_2^\varepsilon + \lambda_1 \mu x_1 + \lambda_2 \mu x_2) - f_i(x_0 + \lambda_1 z_1^\varepsilon + \lambda_2 z_2^\varepsilon + \lambda_1 \mu x_1) -$$
$$- f_i(x_0 + \lambda_1 z_1^\varepsilon + \lambda_2 z_2^\varepsilon + \lambda_2 \mu x_2) + f_i(x_0 + \lambda_1 z_1^\varepsilon + \lambda_2 z_2^\varepsilon)) \leq f_i^{\{2\}}(x_0; \mu x_1, \mu x_2) + \varepsilon.$$

Положив $\alpha^\varepsilon = \min_{i \in J(x_0)} \alpha_i^\varepsilon$, $\beta^\varepsilon = \min_{i \in J(x_0)} \beta_i^\varepsilon$ получим

$$\max_{i \in J(x_0)} \sup_{\substack{0 < \lambda_1 \leq \alpha^\varepsilon, \\ 0 < \lambda_2 \leq \beta^\varepsilon}} \frac{1}{\lambda_1 \lambda_2} (f_i(x_0 + \lambda_1 z_1^\varepsilon + \lambda_2 z_2^\varepsilon + \lambda_1 \mu x_1 + \lambda_2 \mu x_2) - f_i(x_0 + \lambda_1 z_1^\varepsilon + \lambda_2 z_2^\varepsilon + \lambda_1 \mu x_1) -$$
$$- f_i(x_0 + \lambda_1 z_1^\varepsilon + \lambda_2 z_2^\varepsilon + \lambda_2 \mu x_2) + f_i(x_0 + \lambda_1 z_1^\varepsilon + \lambda_2 z_2^\varepsilon)) \leq \max_{i \in J(x_0)} f_i^{\{2\}}(x_0; \mu x_1, \mu x_2) + \varepsilon.$$

Отсюда следует, что

$$f^{\{2\}}(x_0; \mu x_1, \mu x_2) \leq \max_{i \in J(x_0)} \sup_{\substack{0 < \lambda_1 \leq \alpha^\varepsilon, \\ 0 < \lambda_2 \leq \beta^\varepsilon}} \frac{1}{\lambda_1 \lambda_2} (f_i(x_0 + \lambda_1 z_1^\varepsilon + \lambda_2 z_2^\varepsilon + \lambda_1 \mu x_1 + \lambda_2 \mu x_2) -$$
$$- f_i(x_0 + \lambda_1 z_1^\varepsilon + \lambda_2 z_2^\varepsilon + \lambda_1 \mu x_1) - f_i(x_0 + \lambda_1 z_1^\varepsilon + \lambda_2 z_2^\varepsilon + \lambda_2 \mu x_2) +$$
$$+ f_i(x_0 + \lambda_1 z_1^\varepsilon + \lambda_2 z_2^\varepsilon)) + \varepsilon \leq \max_{i \in J(x_0)} f_i^{\{2\}}(x_0; \mu x_1, \mu x_2) + 2\varepsilon.$$

Отсюда при $\varepsilon \to 0$ получим $f^{\{2\}}(x_0; \mu x_1, \mu x_2) \le \max_{i \in J(x_0)} f_i^{\{2\}}(x_0; \mu x_1, \mu x_2)$. Поэтому $f^{\{2\}}(x_0; x_1, x_2) \le \max_{i \in J(x_0)} f_i^{\{2\}}(x_0; x_1, x_2)$. Лемма доказана.

**Следствие 3.** Пусть $f(x) = \max_{i \in J} x_i^*(x, x)$, где $x_i^*$ -симметричный непрерывный билинейный функционал, J-конечное множество. Тогда $f^{\{2\}}(x_0; x_1, x_2) \le 2 \max_{i \in J(x_0)} x_i^*(x_1, x_2)$, где $J(x_0) = \{i \in J : f(x_0) = x_i^*(x_0, x_0)\}$.

Рассмотрим некоторый аналог $\varepsilon$ - субдифференциала. Положим
$$f_\varepsilon^{\{2\}}(x_0; x_1, x_2) = f^{\{2\}}(x_0; x_1, x_2) + \varepsilon \|x_1\| \cdot \|x_2\|, \quad \varepsilon > 0,$$
$$\partial_{\{2\}}^\varepsilon f(x_0) = \{x^* \in \overline{B}(X^2, R) : f_\varepsilon^{\{2\}}(x_0; x_1, x_2) \ge x^*(x_1, x_2), x_1, x_2 \in X\}.$$

Из определения следует, что $\partial_{\{2\}} f(x_0) \subset \partial_{\{2\}}^\varepsilon f(x_0)$ для всех $\varepsilon > 0$, $\partial_{\{2\}} f(x_0) = \partial_{\{2\}}^0 f(x_0)$ и более того, если 2-липшицевая функция в окрестности точки $x_0$, то $\partial_{\{2\}} f(x_0) = \bigcap_{\varepsilon > 0} \partial_{\{2\}}^\varepsilon f(x_0)$.

Обозначим $B(x_0, \delta) = \{x \in X : \|x - x_0\| \le \delta\}$. Множество всех дважды непрерывно дифференцируемых функций из $B(x_0, \delta)$ в R обозначим через $C^2(B(x_0, \delta))$. По определению
$$f^{\{2\}}(x_0; x_1, x_2) = \sup_{z_1, z_2 \in X} \inf_{\alpha > 0} \sup_{\substack{0 < \lambda_1 \le \alpha, \\ 0 < \lambda_2 \le \alpha}} \frac{1}{\lambda_1 \lambda_2}(f(x_0 + \lambda_1 z_1 + \lambda_2 z_2 + \lambda_1 x_1 + \lambda_2 x_2) - f(x_0 + \lambda_1 z_1 + \lambda_2 z_2 + \lambda_1 x_1) -$$
$$- f(x_0 + \lambda_1 z_1 + \lambda_2 z_2 + \lambda_2 x_2) + f(x_0 + \lambda_1 z_1 + \lambda_2 z_2)).$$

Пусть $z_1, z_2 \in X$ фиксированы и $\alpha > 0$ такое, что $x_0 + \lambda_1 z_1 + \lambda_2 z_2 + \lambda_1 x_1 + \lambda_2 x_2 \in B(x_0, \delta)$ $x_0 + \lambda_1 z_1 + \lambda_2 z_2 + \lambda_1 x_1 \in B(x_0, \delta)$, $x_0 + \lambda_1 z_1 + \lambda_2 z_2 + \lambda_2 x_2 \in B(x_0, \delta)$ и $x_0 + \lambda_1 z_1 + \lambda_2 z_2 \in B(x_0, \delta)$ при $0 < \lambda_1 \le \alpha, \ 0 < \lambda_2 \le \alpha$. Если $f \in C^2(B(x_0, \delta))$, то имеем
$$f^{\{2\}}(x_0; x_1, x_2) = \sup_{z_1, z_2 \in X} \inf_{\alpha > 0} \sup_{\substack{0 < \lambda_1 \le \alpha, \\ 0 < \lambda_2 \le \alpha}} \frac{1}{\lambda_1 \lambda_2}(\int_0^1 (f'(x_0 + \lambda_1 z_1 + \lambda_2 z_2 + \lambda_1 x_1 + t\lambda_2 x_2)\lambda_2 x_2 -$$
$$- f'(x_0 + \lambda_1 z_1 + \lambda_2 z_2 + t\lambda_2 x_2)\lambda_2 x_2) dt =$$
$$= \sup_{z_1, z_2 \in X} \inf_{\alpha > 0} \sup_{\substack{0 < \lambda_1 \le \alpha, \\ 0 < \lambda_2 \le \alpha}} \frac{1}{\lambda_1 \lambda_2}(\int_0^1 (f''(x_0 + \lambda_1 z_1 + \lambda_2 z_2 + \theta \lambda_1 x_1 + t\lambda_2 x_2)[\lambda_1 x_1, \lambda_2 x_2] d\theta = f''(x_0)[x_1, x_2].$$

Рассмотрим модификации определения субдифференциала второго порядка (см.[21]).

Пусть X банахово пространство, $f : X \to \overline{R} = R \cup \{\pm \infty\}$, $\text{dom} f = \{x \in X : |f(x)| < +\infty\}$, $x_0 \in \text{dom} f$. Положим
$$f^{(2)+}(x_0; x_1, x_2) = \varlimsup_{t_1 \downarrow 0, t_2 \downarrow 0} \frac{1}{t_1 t_2}(f(x_0 + t_1 x_1 + t_2 x_2) - f(x_0 + t_1 x_1) - f(x_0 + t_2 x_2) + f(x_0)),$$
$$f^{(2)-}(x_0; x_1, x_2) = \varliminf_{t_1 \downarrow 0, t_2 \downarrow 0} \frac{1}{t_1 t_2}(f(x_0 + t_1 x_1 + t_2 x_2) - f(x_0 + t_1 x_1) - f(x_0 + t_2 x_2) + f(x_0)),$$
где считаем, что $f(x_0 + t_1 x_1 + t_2 x_2) - f(x_0 + t_1 x_1) - f(x_0 + t_2 x_2) + f(x_0) = +\infty$, если $f(x_0 + t_1 x_1 + t_2 x_2) - f(x_0 + t_1 x_1) - f(x_0 + t_2 x_2) + f(x_0)$ не определено.

Обозначим (см.[21])

$$f^{2\mp}(x_0;x_1,x_2) = \max\{f^{(2)+}(x_0;x_1,x_2), -f^{(2)-}(x_0;x_1,-x_2), -f^{(2)-}(x_0;-x_1,x_2)\}.$$

Положим

$$\partial_2^{\mp} f(x_0) = \{b \in \overline{B}(X^2, R) : f^{2\mp}(x_0;x_1,x_2) \geq b(x_1,x_2), (x_1,x_2) \in X \times X\}.$$

Следующие утверждения взяты из [21].

**Лемма 5.** Если $X = R$, $(x_1,x_2) \to f^{(2)+}(x_0;x_1,x_2)$ и $(x_1,x_2) \to f^{(2)-}(x_0;x_1,x_2)$ конечные функции, то $f^{2\mp}(x_0;x_1,x_2)$ бисублинейные функции.

**Лемма 6.** Если функция $f$ 2-липшицевая в точке $x_0$, то

$$(f+g)^{2\mp}(x_0;x_1,x_2) \leq f^{2\mp}(x_0;x_1,x_2) + g^{2\mp}(x_0;x_1,x_2)$$

при $(x_1,x_2) \in X \times X$.

**Лемма 7.** $f^{2\mp}(x_0;x_1,x_2) = f^{2\mp}(x_0;x_2,x_1)$.

Ясно, что $f^{\{2\}}(x_0;x_1,x_2) \geq f^{2\mp}(x_0;x_1,x_2)$ при $(x_1,x_2) \in X \times X$.

Множество всех бисублинейных непрерывных симметричных четных функций из $X \times X$ в $R$ обозначим через $BS(X^2)$.

Отметим, что если $h$ бисублинейная непрерывная симметричная четная функция, то $x \to h(x,x)$ второй степени положительно однородная четная непрерывная функция. Если $X$ гильбертово пространство, $\varphi$ четная второй степени положительно однородная полунепрерывная снизу функция, то из следствия 1.3.7 [23] следует, что

$$\varphi(x) = \sup\{Q(x) : Q \in d_2\varphi\},$$

где $d_2\varphi = \{Q \in B_0(X) : \varphi(x) \geq Q(x) \text{ при } x \in X\}$.

Квадратичная функция $Q$ из $d_2\varphi$ называется главной, если не существует другая квадритичная функция $Q_1$ из $d_2\varphi$ такая, что $Q_1(x) \geq Q(x)$ при $x \in X$. Множество главных квадритичных функций из $d_2\varphi$ назовем главным бидифференциалом функции $\varphi$ и обозначим через $d_2^m\varphi$. Ясно, что если $X$ гильбертово пространство, то $\varphi(x) = \sup_{Q \in d_2^m f(x_0)} Q(x)$ и

$$h(x_1,x_2) = \sup\{0{,}5(Q(x_1+x_2) - Q(x_1) - Q(x_2)) : Q \in d_2^m \varphi\}$$

бисублинейная симметричная четная функция.

Отметим, что вместо множества $d_2^m\varphi$ можно взять любое по включение наименьшее выпуклое подмножество $d_2^s\varphi$ множества $d_2\varphi$, которое удовлетворяет равенству $\varphi(x) = \sup_{Q \in d_2^s f(x_0)} Q(x)$.

Если $h \in BS(X^2)$ удовлетворяет условию

$$f^{\{2\}}(x_0;x_1,x_2) \geq h(x_1,x_2) \geq f^{2\mp}(x_0;x_2,x_1)$$

при $(x_1,x_2) \in X \times X$, то функцию $h$ назовем верхней бисублинейной аппроксимацией функции $f$ в точке $x_0$. Верхняя бисублинейная аппроксимация $h$ функции $f$ в точке $x_0$ называется главной верхней аппроксимацией, если не существует другая верхняя бисублинейная аппроксимация $h_1$ такая, что $h(x_1,x_2) \geq h_1(x_1,x_2)$ при $(x_1,x_2) \in X \times X$.

Если $h(x_1,x_2)$ верхняя бисублинейная аппроксимация функции $f$ в точке $x_0$, то

$\partial_2 h = \{b \in \overline{B}(X^2, R) : h(x_1, x_2) \geq b(x_1, x_2), (x_1, x_2) \in X \times X\}$ назовем верхним аппроксимативным бисубдифференциалом функции f в точке $x_0$ и обозначим через $\partial_2^a f(x_0)$. Если h(x) главная верхняя бисублинейная аппроксимация функции f в точке $x_0$, то $\partial_2 h = \partial_2 h(0)$ назовем главным верхним аппроксимативным бисубдифференциалом функции f в точке $x_0$ и обозначим через $\partial_2^m f(x_0)$. Ясно, что $\partial_2^{\mp} f(x_0) \subset \partial_2^m f(x_0)$. Если $X = R$, то $\partial_2^{\mp} f(x_0) = \partial_2^m f(x_0)$.

**Теорема 2.** Если функция f достигает локального минимума (максимума) в пространстве X в точке $x_0$, то $-f^{(2)-}(x_0; x, -x) \geq 0$ $(f^{(2)+}(x_0; x, -x) \geq 0)$ при $x \in X$.

**Следствие 4.** Если функция f достигает локального минимума в пространстве X в точке $x_0$, то $f^{2\mp}(x_0; x, x) \geq 0$ при $x \in X$.

**Следствие 5.** Если функция f достигает локального минимума в пространстве X в точке $x_0$ и h верхняя бисублинейная аппроксимация функции f в точке $x_0$, то $h(x, x) = \sup_{x^* \in \partial_2 h} x^*(x, x) \geq 0$ при $x \in X$.

Пусть X и Y банаховы пространства.

**Лемма 8.** Пусть I конечное множество, $f_i : X \times Y \to R, i \in I$, непрерывные функции в множестве $(x_0 + 2\delta B_X) \times (y_0 + \delta B_Y)$, где $\delta > 0$, и существует число $L > 0$ такое, что
$$f_i(u + x, \upsilon + y) - f_i(u + x, \upsilon) - f_i(u, \upsilon + y) + f_i(u, \upsilon) \leq L\|x\| \cdot \|y\|$$
при $u \in x_0 + \delta B_X$, $\upsilon \in y_0 + \delta B_Y$, $x \in \delta B_X$, $y \in \delta B_Y$ и $f(x, y) = \max_{i \in I} f_i(x, y)$. Тогда для каждого $u \in x_0 + \delta B_X$, $\upsilon \in y_0 + \delta B_Y$ существует $\delta_0 > 0$, где $\delta_0 \leq \delta$ такое, что
$$f(u + x, \upsilon + y) - f(u + x, \upsilon) - f(u, \upsilon + y) + f(u, \upsilon) \leq L\|x\| \cdot \|y\|$$
при $x \in \delta_0 B_X$, $y \in \delta_0 B_Y$.

**Доказательство.** Обозначим $I(u, \upsilon) = \{i \in I : f(u, \upsilon) = f_i(u, \upsilon)\}$. Покажем, что существует $\alpha > 0$ такое, что $I(x, y) \subset I(u, \upsilon)$ при $x \in u + \alpha B_X$, $y \in \upsilon + \alpha B_Y$. Случай $I(u, \upsilon) = I$ тривиален. Пусть $I(u, \upsilon) \neq I$. Положим $b = f(u, \upsilon) - \max_{i \in I/I(u,\upsilon)} f_i(u, \upsilon) > 0$ и $\varepsilon = \dfrac{b}{3}$. Тогда существует $\alpha > 0$ такое, что для каждого $i \in I$ выполнено $|f_i(x, y) - f_i(u, \upsilon)| \leq \varepsilon$ при $x \in u + \alpha B_X$, $y \in \upsilon + \alpha B_Y$. Легко можно проверить, что $|f(x, y) - f(u, \upsilon)| \leq \varepsilon$ при $x \in u + \alpha B_X$, $y \in \upsilon + \alpha B_Y$. Если $j \in I(x, y)$, то
$$f_j(u, \upsilon) \geq f_j(x, y) - \varepsilon = f(x, y) - \varepsilon \geq f(u, \upsilon) - 2\varepsilon =$$
$$= \max_{i \notin I(u,\upsilon)} f_i(u, \upsilon) + b - 2\varepsilon > \max_{i \notin I(u,\upsilon)} f_i(u, \upsilon).$$

Поэтому $j \in I(u, \upsilon)$, т.е. $I(x, y) \subset I(u, \upsilon)$. Обозначим $\delta_0 = \min\{\delta, \alpha\}$. Если $\|x\| \leq \delta_0$, $\|y\| \leq \delta_0$, то

$$f(u+x,\upsilon+y)-f(u+x,\upsilon)-f(u,\upsilon+y)+f(u,\upsilon) \leq \max_{i\in I(u,\upsilon)}(f_i(u+x,\upsilon+y)+f_i(u,\upsilon))-$$
$$-\max_{i\in I(u,\upsilon)}f_i(u+x,\upsilon)-\max_{i\in I(u,\upsilon)}f_i(u,\upsilon+y) \leq \max_{i\in I(u,\upsilon)}(f_i(u+x,\upsilon+y)-f_i(u+x,\upsilon)-$$
$$-f_i(u,\upsilon+y)+f_i(u,\upsilon)) \leq L\|x\|\|y\|.$$

Лемма доказана.

Функция $f: X \times Y \to R$ называется раздельно 2-липшицевой в окрестности точки $(x_0, y_0)$ с постоянной $L$, если существует $\delta > 0$ такое, что
$$|f(u+x,\upsilon+y)-f(u+x,\upsilon)-f(u,\upsilon+y)+f(u,\upsilon)| \leq L\|x\|\cdot\|y\|$$
при $u \in x_0 + \delta B_X$, $\upsilon \in y_0 + \delta B_Y$, $x \in \delta B_X$, $y \in \delta B_Y$.

Обозначим
$$f_{xy}((x_0,y_0);(x,y)) = \sup_{u\in X, \upsilon\in Y} \overline{\lim_{\lambda\downarrow 0, \mu\downarrow 0}} \frac{1}{\lambda\mu}(f(x_0+\lambda u+\lambda x, y_0+\mu\upsilon+\mu y)-f(x_0+\lambda u+\lambda x, y_0+\mu\upsilon)-$$
$$-f(x_0+\lambda u, y_0+\mu\upsilon+\mu y)+f(x_0+\lambda u, y_0+\mu\upsilon)).$$

**Лемма 9.** Если функция $f: X \times Y \to R$ раздельно 2-липшицевая в окрестности точки $(x_0, y_0)$, то $(x,y) \to f_{xy}((x_0,y_0);(x,y))$ бисублинейная функция.

**Доказательство.** Пусть $x_1, x_2 \in X$. Ясно, что
$$f_{xy}((x_0,y_0);x_1+x_2,y) = \sup_{u\in X, \upsilon\in Y} \overline{\lim_{\lambda\downarrow 0, \mu\downarrow 0}} \frac{1}{\lambda\mu}(f(x_0+\lambda u+\lambda x_1+\lambda x_2, y_0+\mu\upsilon+\mu y)-$$
$$-f(x_0+\lambda u+\lambda x_1+\lambda x_2, y_0+\mu\upsilon)-f(x_0+\lambda u, y_0+\mu\upsilon+\mu y)+f(x_0+\lambda u, y_0+\mu\upsilon)) =$$
$$\leq \sup_{u\in X, \upsilon\in Y} \overline{\lim_{\lambda\downarrow 0, \mu\downarrow 0}} \frac{1}{\lambda\mu}(f(x_0+\lambda u+\lambda x_1+\lambda x_2, y_0+\mu\upsilon+\mu y)-f(x_0+\lambda u+\lambda x_1+\lambda x_2, y_0+\mu\upsilon)-$$
$$-f(x_0+\lambda u+\lambda x_1, y_0+\mu\upsilon+\mu y)+f(x_0+\lambda u+\lambda x_1, y_0+\mu\upsilon))+$$
$$+\sup_{u\in X, \upsilon\in Y} \overline{\lim_{\lambda\downarrow 0, \mu\downarrow 0}} \frac{1}{\lambda\mu}(f(x_0+\lambda u+\lambda x_1, y_0+\mu\upsilon+\mu y)-f(x_0+\lambda u+\lambda x_1, y_0+\mu\upsilon)-$$
$$-f(x_0+\lambda u, y_0+\mu\upsilon+\mu y)+f(x_0+\lambda u, y_0+\mu\upsilon))\} = f_{xy}((x_0,y_0);x_1,y)+f_{xy}((x_0,y_0);x_2,y).$$

Кроме того, $f_{xy}((x_0,y_0);\mu x,y) = \mu f_{xy}((x_0,y_0);x,y)$ при $\mu \geq 0$. Получим, что $x \to f_{xy}((x_0,y_0);x,y)$ сублинейная функция. Аналогично имеем, что $y \to f_{xy}((x_0,y_0);x,y)$ сублинейная функция. Лемма доказана.

Если $f$ биположительно однородная функция, то
$$f_{xy}((0,0);(x,y)) = \sup_{u\in X, \upsilon\in Y} \overline{\lim_{\lambda\downarrow 0, \mu\downarrow 0}} \frac{1}{\lambda\mu}(f(\lambda u+\lambda x, \mu\upsilon+\mu y)-f(\lambda u+\lambda x, \mu\upsilon)-f(\lambda u, \mu\upsilon+\mu y)+$$
$$+f(\lambda u, \mu\upsilon)) = \sup_{u\in X, \upsilon\in Y}(f(u+x,\upsilon+y)-f(u+x,\upsilon)-f(u,\upsilon+y)+f(u,\upsilon)).$$

Множество
$$d_{xy}f(x_0,y_0) = \{x^* \in B(X\times Y, R) : f_{xy}((x_0,y_0);(x,y)) \geq x^*(x,y), \ (x,y) \in X \times Y\}$$
назовем частным бисубдифференциалом раздельно 2-липшицевых функций в точке

$(x_0, y_0)$.

Если $b: X \times Y \to R$ билинейная функция, то легко проверяется, что $b_{xy}((x_0, y_0); (x, y)) = b(x, y)$ при $(x, y) \in X \times Y$.

Если $P: X \times Y \to R$ бисублинейная функция, то положим
$$\partial_2 P = \{x^* \in B(X \times Y, R): P(x, y)) \geq x^*(x, y), \quad (x, y) \in X \times Y\}.$$

**Следствие 6**. $f^{\{2\}}((x_0, y_0); (\overline{x}, 0), (0, \overline{y})) = f_{xy}((x_0, y_0); (\overline{x}, \overline{y}))$.

**Лемма 10.** Если $P: X \times Y \to R$ бисублинейная функция, то $\partial_2 P \subset d_{xy} P(0, 0)$.

**Доказательство.** Положив в определении $P_{xy}((0,0); (x, y))$ $(z_1, z_2) = (0, 0)$ и $(z_1, z_2) = (-\lambda x, -\mu y)$ соответственно имеем, что $P_{xy}((0,0); (x, y)) \geq \max\{P(x, y), P(-x, -y)\}$. Поэтому из определения $\partial_2 P$ и $d_{xy} P(0,0)$ вытекает, что $\partial_2 P \subset d_{xy} P(0,0)$. Лемма доказана.

Отметим, что аналогично [19] и [23] можно рассмотреть геометрические аспекты бисубдифференциала и 2-субдифференциала.

## 1.2. О свойствах бидифференциала

В разделе 1.2 даны определения бидифференциала и двадифференциала, введены различные классы функций, связанные с этими понятиями и изучены их свойства.

Пусть $X$ банахово пространство, $f: X \to R$, $B = \{x \in X: \|x\| \leq 1\}$. Положим

$$f^{\{2\}+}(x_0; x) = \sup_{z \in X} \overline{\lim_{\lambda \downarrow 0}} \frac{1}{\lambda^2} (f(x_0 + \lambda z + 2\lambda x) - 2f(x_0 + \lambda z + \lambda x) + f(x_0 + \lambda z)),$$

$$f^{\{2\}-}(x_0; x) = \inf_{z \in X} \underline{\lim_{\lambda \downarrow 0}} \frac{1}{\lambda^2} (f(x_0 + \lambda z + 2\lambda x) - 2f(x_0 + \lambda z + \lambda x) + f(x_0 + \lambda z))$$

при $x \in X$. Положив $z = y - x$ отсюда имеем, что

$$f^{\{2\}+}(x_0; x) = \sup_{y \in X} \overline{\lim_{\lambda \downarrow 0}} \frac{1}{\lambda^2} (f(x_0 + \lambda y + \lambda x) - 2f(x_0 + \lambda y) + f(x_0 + \lambda y - \lambda x)),$$

$$f^{\{2\}-}(x_0; x) = \inf_{y \in X} \underline{\lim_{\lambda \downarrow 0}} \frac{1}{\lambda^2} (f(x_0 + \lambda y + \lambda x) - 2f(x_0 + \lambda y) + f(x_0 + \lambda z - \lambda x)).$$

Множество $D_2 f(x_0) = \{Q \in B_0(X): f^{\{2\}-}(x_0; x) \leq Q(x) \leq f^{\{2\}+}(x_0; x), x \in X\}$ назовем бидифференциалом функции $f$ в точке $x_0$.

Из определения следует, что $D_2(f(x_0)$ выпуклое множество.

Ясно, что $f_{\{2\}}(x_0; x, x) \leq f^{\{2\}-}(x_0; x) \leq f^{\{2\}+}(x_0; x) \leq f^{\{2\}}(x_0; x, x)$.

В частности, если $f$ положительно однородная степени 2 функция, то

$$f^{\{2\}+}(0; x) = \sup_{z \in X} \overline{\lim_{\lambda \downarrow 0}} \frac{1}{\lambda^2} (f(\lambda z + 2\lambda x) - 2f(\lambda z + \lambda x) + f(\lambda z)) =$$
$$= \sup_{z \in X} (f(z + 2x) - 2f(z + x) + f(z)),$$

$$f^{\{2\}-}(0; x) = \inf_{z \in X} \underline{\lim_{\lambda \downarrow 0}} \frac{1}{\lambda^2} (f(\lambda z + 2\lambda x) - 2f(\lambda z + \lambda x) + f(\lambda z)) =$$

$$= \inf_{z \in X} (f(z+2x) - 2f(z+x) + f(z)).$$

Из определения непосредственно следует следующая лемма.

**Лемма 1.** Справедливы следующие соотношения:
1) $f^{\{2\}+}(x_0;x) = f^{\{2\}+}(x_0;-x)$,
2) $f^{\{2\}-}(x_0;x) = f^{\{2\}-}(x_0;-x)$,
3) $f^{\{2\}+}(x_0;\alpha x) = \alpha^2 f^{\{2\}+}(x_0;x)$,
4) $f^{\{2\}-}(x_0;\alpha x) = \alpha^2 f^{\{2\}-}(x_0;x)$,
5) $D_2(\alpha f)(x_0) = \alpha D_2 f(x_0)$.

Функцию $f$ назовем $\{2\}$-липшицевой (двалипшицевой) с постоянной K в окрестности $x_0$, если для некоторого $\varepsilon > 0$ $f$ удовлетворяет условию

$$|f(z+2x) - 2f(z+x) + f(z)| \le K\|x\|^2 \quad (или\ |f(z+x) - 2f(z) + f(z-x)| \le K\|x\|^2)$$

при $x \in \varepsilon B$, $z \in x_0 + \varepsilon B$.

Если $f$ $\{2\}$-липшицевая (двалипшицевая) с постоянной K в окрестности $x_0$ функция, то легко проверяется, что $D_2 f(x_0)$ ограниченное множество.

Функцию $f$ назовем $\theta$-билипшицевой ($\theta$-двалипшицевой) с постоянной K в окрестности $x_0$, если для некоторого $\varepsilon > 0$ $f$ удовлетворяет условию

$$|f(z+2x) - 2f(z+x) - f(z+2y) + 2f(z+y)| \le K\|x-y\|^\theta (\|x\| + \|y\|)^{2-\theta},$$
$$(|f(z+x) + f(z-x) - f(z+y) - f(z-y)| \le K\|x-y\|^\theta (\|x\| + \|y\|)^{2-\theta}) \qquad (1)$$

при $x, y \in \varepsilon B$, $z \in x_0 + \varepsilon B$, $0 < \theta \le 2$.

Если при $z = x_0$ удовлетворяется соотношение (1), то функцию $f$ назовем $\theta$-билипшицевой ($\theta$-двалипшицевой) с постоянной K в точке $x_0$.

Функцию $f$ назовем $\{2\}$-липшицевой (двалипшицевой) с постоянной K в точке $x_0$, если для некоторого $\varepsilon > 0$ $f$ удовлетворяет условию

$$|f(x_0+2x) - 2f(x_0+x) + f(x_0)| \le K\|x\|^2, \quad x \in \varepsilon B,$$
$$(|f(x_0+x) - 2f(x_0) + f(x_0-x)| \le K\|x\|^2, \quad x \in \varepsilon B).$$

Функцию $f$ назовем слабо $\theta$-билипшицевой (слабо $\theta$-двалипшицевой) с постоянной K в окрестности $x_0$, если для некоторого $\varepsilon > 0$ $f$ удовлетворяет условию

$$|f(z+2x) - 2f(z+x) - f(z+2y) + 2f(z+y)| \le K\|x-y\|^\theta (\|x\|+\|y\|)^{2-\theta} + o((\|x\|+\|y\|)^2),$$
$$(|f(z+x) + f(z-x) - f(z+y) - f(z-y)| \le K\|x-y\|^\theta (\|x\|+\|y\|)^{2-\theta} + o((\|x\|+\|y\|)^2))$$

при $x, y \in \varepsilon B$, $z \in x_0 + \varepsilon B$, $0 < \theta \le 2$, где $\dfrac{o(\lambda)}{\lambda} \to 0$ при $\lambda \downarrow 0$.

Функцию $f$ назовем $(\theta, \delta)$-сильно билипшицевой с постоянной K в точке $x_0$, если $f$ удовлетворяет условию

$$|f(x_0+x) - f(x_0+y)| \le K\|x-y\|^\theta (\|x\|+\|y\|)^{2-\theta}, \quad x, y \in \delta B, \ \delta > 0, \ 0 < \theta \le 2.$$

Ясно, что если $f$ $\theta$-билипшицевая ($\theta$-двалипшицевая) функция в окрестности $x_0$, то $f$ является $\{2\}$-липшицевой (двалипшицевой) в окрестности $x_0$.

Функцию $f$ назовем второго порядка строго дифференцируемой в точке $x_0$, если существует симметричный оператор $A \in L(X, X^*)$ такой, что

$$\lim_{\lambda \downarrow 0,\ z \to x_0} \frac{1}{\lambda^2}(f(z+2\lambda x) - 2f(z+\lambda x) + f(z)) = \langle Ax, x \rangle$$

при $x \in X$, т.е. функция $f$ второго порядка строго дифференцируема в точке $x_0$, если для всякого $\varepsilon > 0$ существует симметричный оператор $A \in L(X, X^*)$ и $\delta > 0$ такие, что

$$\left|f(z+2x) - 2f(z+x) + f(z) - \langle Ax, x\rangle\right| \le \varepsilon\|x\|^2$$

при $z \in x_0 + \delta B$ и $x \in \delta B$.

Если функция $f$ второго порядка строго дифференцируема в точке $x_0$, то $D_2f(x_0) = \{\langle Ax, x\rangle\}$.

Ясно, что если функция $f$ второго порядка s-дифференцируема в точке $x_0$, то $f$ второго порядка строго дифференцируемо в точке $x_0$ и их можно отождествить.

Если $f(x) = Q(x) = \langle Ax, x\rangle$ квадратичная функция, где $A \in L(X, X^*)$ симметричный оператор, то $\lim\limits_{\lambda\downarrow 0,\ z\to x_0} \frac{1}{\lambda^2}(f(z+2\lambda x) - 2f(z+\lambda x) + f(z)) = 2\langle Ax, x\rangle$. Кроме того, имеем

$$f^{\{2\}+}(0; x) = \sup_{z\in X}\overline{\lim_{\lambda\downarrow 0}}\frac{1}{\lambda^2}(f(\lambda z + 2\lambda x) - 2f(\lambda z + \lambda x) + f(\lambda z)) =$$
$$= \sup_{z\in X}(\langle A(z+2x), z+2x\rangle - 2\langle A(z+x), z+x\rangle + \langle Az, z\rangle = 2\langle Ax, x\rangle,$$
$$f^{\{2\}-}(0; x) = \inf_{z\in X}(\langle A(z+2x), z+2x\rangle - 2\langle A(z+x), z+x\rangle + \langle Az, z\rangle = 2\langle Ax, x\rangle.$$

Поэтому $D_2f(x_0) = = \{2\langle Ax, x\rangle\}$.

Если $z, x, y \in X$ и $f$ выпуклая функция, то
$$f(z+2x) - 2f(z+x) + f(z) = 2\{\frac{1}{2}f(z+2x) - f(z+x) + \frac{1}{2}f(z)\} \ge 0.$$

Пусть $C \subset X$. Положим $d_C(y) = \inf\{\|y-z\| : z \in C\}$, $d_2(y) = d_C^2(y)$.

Если $c \in C$, $x \in X$, то
$$\left|d_2(c+2x) - 2d_2(c+x) + d_2(c)\right| = \left|d^2(c+2x) - 2d^2(c+x) + d^2(c)\right| \le$$
$$\le \left|(d(c+2x) - d(c+x))\right|\left|(d(c+2x) - d(c) + d(c+x) - d(c))\right| + (d(c+x) - d(c))^2 \le$$
$$\le \|x\|(2\|x\| + \|x\|) + \|x\|^2 = 4\|x\|^2.$$

Если $C \subset X$ выпуклое множество, то $d_C(y)$ выпуклая функция (см.[7], стр.56). Поэтому легко проверяется, что $d_2(y) = d_C^2(y)$ выпуклая функция. Тогда имеем, что
$$d_2(z+2x) - 2d_2(z+x) + d_2(z) = 2\{\frac{1}{2}d_2(z+2x) - d_2(z+x) + \frac{1}{2}d_2(z)\} \ge 0.$$

**Лемма 2.** Если $C \subset X$ непустое замкнутое выпуклое подмножество гильбертово пространства $X$, то
$$0 \le d_2(z+2x) - 2d_2(z+x) + d_2(z) \le 2\|x\|^2$$
при $z, x \in X$.

**Доказательство.** Так как $C$ замкнутое выпуклое подмножество гильбертово пространства $X$, то существует точка $c \in C$ такая, что $d_2(z+x) = \|z+x-c\|^2$. Поэтому
$$0 \le d_2(z+2x) - 2d_2(z+x) + d_2(z) \le \|z+2x-c\|^2 - \|z+x-c\|^2 + \|z-c\|^2 - \|z+x-c\|^2 =$$
$$= \langle 2z+3x-2c, x\rangle - \langle 2z+x-2c, x\rangle = \langle 2x, x\rangle = 2\|x\|^2$$
при $z, x \in X$. Лемма доказана.

**Лемма 3.** Если $C \subset X$ непустое подмножество гильбертово пространства $X$, то

$$d_2(z+2x) - 2d_2(z+x) + d_2(z) \le 2\|x\|^2$$

при $z, x \in X$.

**Доказательство.** Так как $d_2(z+x) = \inf\limits_{c \in C}\|z+x-c\|^2$, то для любого $\varepsilon > 0$ существует точка $c \in C$ такая, что $\|z+x-c\|^2 < d_2(z+x) + \varepsilon$. Поэтому

$$d_2(z+2x) - 2d_2(z+x) + d_2(z) \le \|z+2x-c\|^2 - \|z+x-c\|^2 + \|z-c\|^2 - \|z+x-c\|^2 + 2\varepsilon =$$
$$= \langle 2z+3x-2c, x\rangle - \langle 2z+x-2c, x\rangle + 2\varepsilon = \langle 2x, x\rangle + 2\varepsilon = 2\|x\|^2 + 2\varepsilon.$$

Так как $\varepsilon > 0$ любая, то имеем, что $d_2(z+2x) - 2d_2(z+x) + d_2(z) \le 2\|x\|^2$. Лемма доказана.

При условии леммы 3 имеем, что

$$2\|x\|^2 \ge \varlimsup_{t \downarrow 0,\, z \to x_0} \frac{1}{t^2}(d_2(z+2tx) - 2d_2(z+tx) + d_2(z)) \ge \varlimsup_{t \downarrow 0} \frac{1}{t^2}(d_2(x_0+2tx) - 2d_2(x_0+tx) +$$
$$+ d_2(x_0)) \ge \varlimsup_{t \downarrow 0} \frac{1}{t^2}(d_2(x_0+2tx) - d_2(x_0)) + \varliminf_{t \downarrow 0}\frac{1}{t^2}(-2d_2(x_0+tx) + 2d_2(x_0)) =$$
$$= 4\varlimsup_{t \downarrow 0}\frac{1}{t^2}(d_2(x_0+tx) - d_2(x_0)) - 2\varlimsup_{t \downarrow 0}\frac{1}{t^2}(d_2(x_0+tx) - d_2(x_0)) = 2\varlimsup_{t \downarrow 0}\frac{1}{t^2}d_2(x_0+tx) \ge 0.$$

По лемме 3 имеем, что если $C$ выпуклое множество, то лемма 2 также верна.

**Лемма 4.** Если $C \subset X$ непустое подмножество гильбертово пространства, то
$$d_2(z+x) - 2d_2(z) + d_2(z-x) \le 2\|x\|^2$$
при $z, x \in X$.

**Доказательство.** Так $d_2(z) = \inf\limits_{c \in C}\|z-c\|^2$, то для любого $\varepsilon > 0$ существует точка $c \in C$ такая, что $\|z-c\|^2 < d_2(z) + \varepsilon$. Поэтому

$$d_2(z+x) - 2d_2(z) + d_2(z-x) \le \|z+x-c\|^2 - \|z-c\|^2 + \|z-x-c\|^2 - \|z-c\|^2 + 2\varepsilon =$$
$$= \langle 2z+x-2c, x\rangle - \langle 2z-x-2c, x\rangle + 2\varepsilon = \langle 2x, x\rangle + 2\varepsilon = 2\|x\|^2 + 2\varepsilon.$$

Так как $\varepsilon > 0$ любое число, то имеем, что
$$d_2(z+x) - 2d_2(z) + d_2(z-x) \le 2\|x\|^2.$$

Лемма доказана.

**Следствие 1.** Если $C \subset X$ непустое выпуклое множество гильбертово пространства, то
$$0 \le d_2(z+x) - 2d_2(z) + d_2(z-x) \le 2\|x\|^2$$
при $z, x \in X$.

При условии леммы 4 имеем, что

$$2\|x\|^2 \ge \varlimsup_{t \downarrow 0,\, z \to x_0}\frac{1}{t^2}(d_2(z+tx) - 2d_2(z) + d_2(z-tx)) \ge \varlimsup_{t \downarrow 0}\frac{1}{t^2}(d_2(x_0+tx) - d_2(x_0)) +$$
$$+ \varliminf_{t \downarrow 0}\frac{1}{t^2}(d_2(x_0-tx) - d_2(x_0)) \ge 0.$$

**Лемма 5.** Если $C$ непустое замкнутое выпуклое подмножество гильбертого пространства $X$, то при любых $z, x, \upsilon \in X$ выполняется соотношение
$$|d_2(z+2x) - d_2(z+2\upsilon) - 2d_2(z+x) + 2d_2(z+\upsilon)| \le 10\,\|x-\upsilon\|(\|x\|+\|\upsilon\|).$$

Лемма 5 доказана в [19] (см. леммы 2).

Отметим, что если $x_0 \in C$, то $|d_2(x_0+x) - d_2(x_0)| \leq \|x\|^2$, $\|x\| d_C(x_0+x) \leq \|x\|^2$ и $|d_2(x_0+2x) - 2d_2(x_0+x) + d_2(x_0)| \leq 4\|x\|^2$ при $x \in X$.

Ясно, что $D_2 d_2(x_0) = \{Q \in B_0(X): d_2^{\{2\}-}(x_0;x) \leq Q(x) \leq d_2^{\{2\}+}(x_0;x), x \in X\}$ бидифференциал функции $d_2(x_0)$ в точке $x_0 \in C$.

Множество $\tilde{\Omega}_C(x_0) = \bigcup_{\lambda \geq 0} \lambda D_2 d_2(x_0)$ назовем бинормальным конусом к $C$ в точке $x_0 \in C$, а множество $\tilde{Q}_C(x_0) = \{x \in X: Q(x) \leq 0$ при $Q \in \tilde{\Omega}_C(x_0)\}$ назовем бикасательным конусом к $C$ в точке $x_0 \in C$. Если $x \in \tilde{Q}_C(x_0)$, то имеем, что $Q(x) \leq 0$ при $Q \in D_2 d_2(x_0)$. Кроме того, если $D_2 d_2(x_0)$ непусто, то $d_2^{\{2\}-}(x_0;x) \leq \inf_{Q \in D_2 d_2(x_0)} Q(x) \leq \sup_{Q \in D_2 d_2(x_0)} Q(x) \leq d_2^{\{2\}+}(x_0;x)$.

Поэтому отсюда следует, что если $x \in \tilde{Q}_C(x_0)$, то $d_2^{\{2\}-}(x_0;x) \leq 0$, а если $d_2^{\{2\}+}(x_0;x) \leq 0$, то $x \in \tilde{Q}_C(x_0)$.

**Теорема 1**. Если $f$ $\theta$-билипшицевая функция с постоянной $K$ в окрестности $u$, то
$$|f^{\{2\}+}(u;x) - f^{\{2\}+}(u;\upsilon)| \leq K\|x-\upsilon\|^\theta (\|x\|+\|\upsilon\|)^{2-\theta},$$
$$|f^{\{2\}-}(u;x) - f^{\{2\}-}(u,\upsilon)| \leq K\|x-\upsilon\|^\theta (\|x\|+\|\upsilon\|)^{2-\theta}$$

при $x, \upsilon \in X$.

**Доказательсьтво.** Если $x, \upsilon \in X$ и $y \in X$, то

$f(u+\lambda y+2\lambda x) - 2f(u+\lambda y+\lambda x) + f(u+\lambda y) \leq f(u+\lambda y+2\lambda\upsilon) - 2f(u+\lambda y+\lambda\upsilon) + f(u+\lambda y) +$
$+ K\lambda^2 \|x-\upsilon\|^\theta (\|x\|+\|\upsilon\|)^{2-\theta}$

где $\lambda > 0$ и достаточно мало. Деление на $\lambda^2$ и переход к верхним пределам при $\lambda \downarrow 0$ и далее получив супремум по $y \in X$, дает неравенство

$$f^{\{2\}+}(u;x) \leq f^{\{2\}+}(u;\upsilon) + K\|x-\upsilon\|^\theta (\|x\|+\|\upsilon\|)^{2-\theta}. \qquad (2)$$

В этом неравенстве меняя местами $x$ и $\upsilon$ получим

$$f^{\{2\}+}(u;\upsilon) \leq f^{\{2\}+}(u;x) + K\|x-\upsilon\|^\theta (\|x\|+\|\upsilon\|)^{2-\theta}. \qquad (3)$$

Из (2) и (3) получим

$$|f^{\{2\}+}(u;x) - f^{\{2\}+}(u;\upsilon)| \leq K\|x-\upsilon\|^\theta (\|x\|+\|\upsilon\|)^{2-\theta}.$$

Аналогично доказывается второе неравенство. Теорема доказана.

Если функция $f$ слабо $\theta$-билипшицевая (слабо $\theta$-двалипшицевая) с постоянной K в окрестности u функция, то теорема 1 также верна.

Положим
$$D^+_2 f(x_0) = \{Q \in B_0(X): f^{\{2\}+}(x_0;x) \geq Q(x), \quad x \in X\},$$
$$D^-_2 f(x_0) = \{Q \in B_0(X): f^{\{2\}-}(x_0;x) \leq Q(x), \quad x \in X\}.$$

Ясно, что $D_2(f(x_0)) = D^+_2 f(x_0) \cap D^-_2 f(x_0)$. Кроме того, если $D_2 f(x_0)$ непусто, то
$$f^{\{2\}-}(x_0;x) \leq \inf_{Q \in D_2 f(x_0)} Q(x) \leq \sup_{Q \in D_2 f(x_0)} Q(x) \leq f^{\{2\}+}(x_0;x)$$

при $x \in X$.

Пусть X — гильбертово пространство. Если $f^{\{2\}+}(x_0;x)$ — полунепрерывная снизу функция, то $f^{\{2\}+}(x_0;x) = \sup\limits_{Q \in D_2^+ f(x_0)} Q(x)$, $f^{\{2\}-}(x_0;x)$ — полунепрерывная сверху функция, то $f^{\{2\}-}(x_0;x) = \inf\limits_{Q \in D_2^- f(x_0)} Q(x)$ (см. [23], следствие 3.2.2). Тогда при условии теоремы 1 имеем, что $f^{\{2\}+}(x_0;x) = \sup\limits_{Q \in D_2^+ f(x_0)} Q(x)$ и $f^{\{2\}-}(x_0;x) = \inf\limits_{Q \in D_2^- f(x_0)} Q(x)$. Отсюда следует, что в общем случае функция $x \to f^{\{2\}+}(x_0;x)$ не является выпуклой.

Пусть функция $f$ удовлетворяет {2}-липшицеву условию с постоянной K в окрестности точки $x_0$. Положим $f_K^{\{2\}+}(x_0;x) = f^{\{2\}+}(x_0;x) + K\|x\|^2$ при $x \in X$. Ясно, что $f_K^{\{2\}+}(x_0;x) \geq 0$ при $x \in X$. Если $f^{\{2\}+}(x_0;x)$ выпуклая функция, то $f_K^{\{2\}+}(x_0;x)$ также выпуклая функция.

Покажем, что кроме того если $f$ выпуклая функция, то $x \to f^{\{2\}+}(x_0;x)$ выпуклая функция. Пусть $x_1, x_2 \in X$ и $\alpha \in [0,1]$. Тогда по определению имеем

$$f^{\{2\}+}(x_0;\alpha x_1 + (1-\alpha)x_2) = \sup_{y \in X} \overline{\lim_{\lambda \downarrow 0}} \frac{1}{\lambda^2}(f(x_0 + \lambda y + \lambda(\alpha x_1 + (1-\alpha)x_2)) - 2f(x_0 + \lambda y) +$$

$$+ f(x_0 + \lambda y - \lambda(\alpha x_1 + (1-\alpha)x_2)) = \sup_{y \in X} \overline{\lim_{\lambda \downarrow 0}} \frac{1}{\lambda^2}(f(\alpha(x_0 + \lambda y + \lambda x_1) + (1-\alpha)(x_0 + \lambda y + \lambda x_2)) -$$

$$- 2f(x_0 + \lambda y) + f(\alpha(x_0 + \lambda y - \lambda x_1) + (1-\alpha)(x_0 + \lambda y - \lambda x_2)) \leq \alpha \sup_{y \in X} \overline{\lim_{\lambda \downarrow 0}} \frac{1}{\lambda^2}(f(x_0 + \lambda y + \lambda x_1) -$$

$$- 2f(x_0 + \lambda y) + f(x_0 + \lambda y - \lambda x_1)) + (1-\alpha)\sup_{y \in X} \overline{\lim_{\lambda \downarrow 0}} \frac{1}{\lambda^2}(f((x_0 + \lambda y + \lambda x_2)) - 2f(x_0 + \lambda y) +$$

$$+ f((x_0 + \lambda y - \lambda x_2)) = \alpha f^{\{2\}+}(x_0;x_1) + (1-\alpha)f^{\{2\}+}(x_0;x_2)$$

Поэтому $x \to f^{\{2\}+}(x_0;x)$ выпуклая функция.

Квадратичная функция Q из $D^+_2 f(x_0)$ называется главной, если не существует другая квадритичная функция $Q_1$ из $D^+_2 f(x_0)$ такая, что $Q_1(x) \geq Q(x)$ при $x \in X$. Множество главных квадратных функций из $D^+_2 f(x_0)$ назовем главным бидифференциалом функции $f$ в точке $x_0$ и обозначим через $D_2^m f(x_0)$. Ясно, что если X — гильбертово пространство, $f^{\{2\}+}(x_0;x)$ полунепрерывная снизу функция, то $f^{\{2\}+}(x_0;x) = \sup\limits_{Q \in D_2^m f(x_0)} Q(x)$.

Если $f_{\{2\}}(x_0;x,x) = f^{\{2\}-}(x_0;x)$ и $f^{\{2\}+}(x_0;x) = f^{\{2\}}(x_0;x,x)$ при $x \in X$, то функцию $f$ назовем второго порядка регулярным в точке $x_0$. Если функция $f$ второго порядка регулярно в точке $x_0$, то $D_2 f(x_0) = d_{\{2\}} f(x_0)$.

Пусть $\psi_1 : X \times X \to R$ и $\psi_2 : X \times X \to R$ биположительно однородные симметричные функции. Функцию $f : X \to R$ назовем $(\theta, \psi_1, \psi_2)$-двалипшицевой в окрестности точки $x_0$, если $f$ для некоторого $\varepsilon > 0$ удовлетворяет условию

$$|f(x_0+z_1+z_2+x_1+x_2)-f(x_0+z_1+z_2+x_1)-f(x_0+z_1+z_2+x_2)+f(x_0+z_1+z_2)-$$
$$-f(x_0+y_1+y_2+\upsilon_1+\upsilon_2)+f(x_0+y_1+y_2+\upsilon_1)+f(x_0+y_1+y_2+\upsilon_2)-f(x_0+y_1+y_2)|\leq$$
$$\leq(|\psi(z_1+x_1,z_2+x_2)-\psi(y_1+\upsilon_1,y_2+\upsilon_2)|+|\psi(x_1,x_2)-\psi(\upsilon_1,\upsilon_2)|+$$
$$+|\psi(z_1,z_2)-\psi(y_1,y_2)|)^\theta(|\psi(x_1,x_2)+\psi(\upsilon_1,\upsilon_2)|+|\psi(z_1,z_2)+\psi(y_1,y_2)|+$$
$$+|\psi(z_1+x_1,z_2+x_2)+\psi(y_1+\upsilon_1,y_2+\upsilon_2)|)^{2-\theta}+o(\|x_1\|\|x_2\|+\|\upsilon_1\|\|\upsilon_2\|+\|z_1\|\|z_2\|+\|y_1\|\|y_2\|)$$

при $x_1,x_2,\upsilon_1,\upsilon_2,z_1,z_2,y_1,y_2 \in \varepsilon B$, где $0 < \theta \leq 2$, $B=\{x \in X: \|x\| \leq 1\}$.

Ясно, что

$$f^{\{2\}+}(x_0;x) = \sup_{y \in X} \overline{\lim_{\lambda \downarrow 0}} \frac{1}{\lambda^2}(f(x_0+\lambda y+2\lambda x)-2f(x_0+\lambda y+\lambda x)+f(x_0+\lambda y)) =$$
$$= \sup_{y_1,y_2 \in X} \overline{\lim_{\lambda \downarrow 0}} \frac{1}{\lambda^2}(f(x_0+\lambda(y_1+y_2)+\lambda 2x)-2f(x_0+\lambda(y_1+y_2)+\lambda x)+f(x_0+\lambda(y_1+y_2))).$$

Поэтому, если f удовлетворяет $(\theta,\psi_1,\psi_2)$-двалипшицеву условию в окрестности точки $x_0$, то аналогично лемме 1.5 [18] имеем, что $D_2f(x_0)=d_{\{2\}}f(x_0)$.

Положим

$$f^{(2)}(x_0;x) = \lim_{\lambda \downarrow 0} \frac{1}{\lambda^2}(f(x_0+2\lambda x)-2f(x_0+\lambda x)+f(x_0)),$$

$$f^{(2)+}(x_0;x) = \overline{\lim_{\lambda \downarrow 0}} \frac{1}{\lambda^2}(f(x_0+2\lambda x)-2f(x_0+\lambda x)+f(x_0)),$$

$$f^{(2)-}(x_0;x) = \underline{\lim_{\lambda \downarrow 0}} \frac{1}{\lambda^2}(f(x_0+2\lambda x)-2f(x_0+\lambda x)+f(x_0))$$

при $x \in X$. Отметим, что для задачи с ограничением использование (верхнего и нижнего) производных по направлению $f^{(2)}(x_0;x)$ ($f^{(2)+}(x_0;x)$, $f^{(2)-}(x_0;x)$) более удобно.

Если функция f удовлетворяет $\{2\}$-липшицеву условию с постоянной K в точке $x_0$, то $|f^{(2)+}(x_0;x)| \leq K\|x\|^2$ и $|f^{(2)-}(x_0;x)| \leq K\|x\|^2$ при $x \in X$.

Положим

$$f^2(x_0;x) = \lim_{\lambda \downarrow 0} \frac{1}{\lambda^2}(f(x_0+\lambda x)-2f(x_0)+f(x_0-\lambda x)),$$

$$f^{2+}(x_0;x) = \overline{\lim_{\lambda \downarrow 0}} \frac{1}{\lambda^2}(f(x_0+\lambda x)-2f(x_0)+f(x_0-\lambda x)),$$

$$f^{2-}(x_0;x) = \underline{\lim_{\lambda \downarrow 0}} \frac{1}{\lambda^2}(f(x_0+\lambda x)-2f(x_0)+f(x_0-\lambda x)).$$

Если функция f удовлетворяет двалипшицеву условию с постоянной K в точке $x_0$ и f выпуклая функция, то $x \to f^{2+}(x_0;x)$ выпуклая функция.

Покажем, что в общем случае $f^2(x_0;x) \neq f^{(2)}(x_0;x)$.

**Пример 1.** Пусть $X=R$ и $f(x) = \begin{cases} x^2 : x \geq 0, \\ -x^2 : x < 0. \end{cases}$

Ясно, что $f^2(0;x) = \lim_{\lambda \downarrow 0} \frac{1}{\lambda^2}(f(\lambda x)-2f(0)+f(-\lambda x)) = \lim_{\lambda \downarrow 0} \frac{1}{\lambda^2}(\lambda^2 x^2 - \lambda^2 x^2) = 0$ при $x \in R$.

Легко проверяется, что

$$f^{(2)}(0;x) = \lim_{\lambda\downarrow 0}\frac{1}{\lambda^2}(f(2\lambda x) - 2f(\lambda x) + f(0)) = \lim_{\lambda\downarrow 0}\frac{1}{\lambda^2}(4\lambda^2 x^2 - 2\lambda^2 x^2) = 2x^2$$

при $x \in R$, $x \geq 0$;

$$f^{(2)}(0;x) = \lim_{\lambda\downarrow 0}\frac{1}{\lambda^2}(f(2\lambda x) - 2f(\lambda x) + f(0)) = \lim_{\lambda\downarrow 0}\frac{1}{\lambda^2}(-4\lambda^2 x^2 + 2\lambda^2 x^2) = -2x^2$$

при $x \in R$, $x < 0$.

Отсюда следует, что $f^2(x_0;x) \neq f^{(2)}(x_0;x)$ при $x \in R$, $x \neq 0$.

Множество $d_2 f(x_0) = \{Q \in B_0(X): f^{(2)-}(x_0;x) \leq Q(x) \leq f^{(2)+}(x_0;x), x \in X\}$ назовем 2-дифференциалом функции $f$ в точке $x_0$.

Положим

$$d_2^+ f(x_0) = \{Q \in B_0(X): f^{(2)+}(x_0;x) \geq Q(x), \ x \in X\},$$
$$d_2^- f(x_0) = \{Q \in B_0(X): f^{(2)-}(x_0;x) \leq Q(x), \ x \in X\}.$$

Ясно, что $d_2 f(x_0) = d_2^+ f(x_0) \cap d_2^- f(x_0)$. Кроме того, если $d_2 f(x_0)$ непусто, то

$$f^{(2)-}(x_0;x) \leq \inf_{Q \in d_2 f(x_0)} Q(x) \leq \sup_{Q \in d_2 f(x_0)} Q(x) \leq f^{(2)+}(x_0;x).$$

Если $f$ удовлетворяет $\theta$-билипшицеву условию с постоянной $K$ в точке $x_0$, то аналогично теореме 1 имеем, что

$$\left|f^{(2)+}(x_0;x) - f^{(2)+}(x_0;y)\right| \leq K\|x-y\|^\theta (\|x\|+\|y\|)^{2-\theta},$$
$$\left|f^{(2)-}(x_0;x) - f^{(2)-}(x_0;y)\right| \leq K\|x-y\|^\theta (\|x\|+\|y\|)^{2-\theta}.$$

Поэтому, если $X$ гильбертово пространство и $f$ удовлетворяет $\theta$-билипшицеву условию с постоянной $K$ в точке $x_0$, то по следствию 3.2.2[23] имеем, что $d_2^+ f(x_0)$ и $d_2^- f(x_0)$ непусто, $f^{(2)+}(x_0;x) = \sup_{Q \in d_2^+ f(x_0)} Q(x)$ и $f^{(2)-}(x_0;x) = \inf_{Q \in d_2^- f(x_0)} Q(x)$.

Пусть $f_i : X \to R$, где $i \in J$, $f(x) = \sup_{i \in J} f_i(x)$, $J(x_0) = \{i \in J : f(x_0) = f_i(x_0)\}$.

**Лемма 6.** Если $X$ банахово пространство, $J$ конечное множество, $f_i$ ($i \in J$) непрерывные функции и удовлетворяют $\{2\}$-липшицеву условию в точке $x_0$; $f(x) = \max_{i \in J} f_i(x)$, то $f^{(2)+}(x_0;x) \leq \max_{i \in J(x_0)} f_i^{(2)+}(x_0;x)$ при $x \in X$.

**Доказательство.** Покажем, что существует $\alpha > 0$ такое, что $J(y) \subset J(x_0)$ при $\|y - x_0\| \leq \alpha$. Случай $J(x_0) = J$ тривиален. Пусть $J(x_0) \neq J$. В самом деле, пусть $a = f(x_0) - \max_{i \in J \setminus J(x_0)} f_i(x_0) > 0$, $\varepsilon = a/3$ и $\alpha > 0$ такое, что для всякого $i \in J$ выполнено $|f_i(y) - f_i(x_0)| \leq \varepsilon$ при $\|y - x_0\| \leq \alpha$. Легко можно проверить, что $|f(y) - f(x_0)| \leq \varepsilon$ при $\|y - x_0\| \leq \alpha$. Если $j \in J(y)$, то

$f_j(x_0) \geq f_j(y) - \varepsilon = f(y) - \varepsilon \geq f(x_0) - 2\varepsilon = \max_{i \in J(x_0)} f_i(x_0) + a - 2\varepsilon > \max_{i \notin J(x_0)} f_i(x_0)$.

Поэтому $j \in J(x_0)$, т.е. $J(y) \subset J(x_0)$ при $\|y - x_0\| \leq \alpha$. Если $i \in J(x_0)$, то при достаточно малых $t$ имеем

$$\frac{1}{t^2}(f(x_0 + 2tx) - 2f(x_0 + tx) + f(x_0)) \leq \frac{1}{t^2}(\max_{i \in J(x_0)}(f_i(x_0 + 2tx) + f_i(x_0)) - 2\max_{i \in J(x_0)} f_i(x_0 + tx) \leq$$

$$\leq \frac{1}{t^2} \max_{i \in J(x_0)}(f_i(x_0 + 2tx) + f_i(x_0) - 2f_i(x_0 + tx)).$$

Поэтому

$$\overline{\lim_{t\downarrow 0}}\frac{1}{t^2}(f(x_0+2tx)-2f(x_0+tx)+f(x_0))\leq \overline{\lim_{t\downarrow 0}}\frac{1}{t^2}\max_{i\in J(x_0)}(f_i(x_0+2tx)+f_i(x_0)-2f_i(x_0+tx)),$$

т.е.

$$f^{(2)+}(x_0;x)\leq \inf_{\alpha>0}\sup_{0<t\leq\alpha}\max_{i\in J(x_0)}\frac{1}{t^2}(f_i(x_0+2tx)+f_i(x_0)-2f_i(x_0+tx)).$$

Отсюда следует, что

$$f^{(2)+}(x_0;x)\leq \inf_{\alpha>0}\max_{i\in J(x_0)}\sup_{0<t\leq\alpha}\frac{1}{t^2}(f_i(x_0+2tx)-2f_i(x_0+tx)+f_i(x_0)).$$

По определению $f_i^{(2)+}(x_0;x)$ для любого $\varepsilon>0$ существует число $\alpha_i>0$ такое, что

$$\sup_{0<t\leq\alpha_i}\frac{1}{t^2}(f_i(x_0+2tx)-2f_i(x_0+tx)+f_i(x_0))<f_i^{(2)+}(x_0;x)+\varepsilon.$$

Обозначив $\alpha=\min_{i\in J(x_0)}\alpha_i$ имеем

$$\max_{i\in J(x_0)}\sup_{0<t\leq\alpha}\frac{1}{t^2}(f_i(x_0+2tx)-2f_i(x_0+tx)+f_i(x_0))<\max_{i\in J(x_0)}f_i^{(2)+}(x_0;x)+\varepsilon.$$

Поэтому

$$f^{(2)+}(x_0;x)\leq \max_{i\in J(x_0)}\sup_{0<t\leq\alpha}\frac{1}{t^2}(f_i(x_0+2tx)-2f_i(x_0+tx)+f_i(x_0))<\max_{i\in J(x_0)}f_i^{(2)+}(x_0;x)+\varepsilon.$$

Так как $\varepsilon>0$ произвольно, то отсюда получим $f^{(2)+}(x_0;x)\leq \max_{i\in J(x_0)}f_i^{(2)+}(x_0;x)$. Лемма доказана.

**Лемма 7.** Пусть $J$ конечное множество, $f_i$ ($i\in J$) непрерывные в окрестности $x_0$ и удовлетворяют $\{2\}$-липшицеву условию с постоянной $K$ в окрестности $x_0$; $f(x)=\max_{i\in J}f_i(x)$. Тогда $f^{\{2\}+}(x_0;x)\leq \max_{i\in J(x_0)}f_i^{\{2\}+}(x_0;x)$ при $x\in X$.

**Доказательство.** Покажем, что существует $\alpha>0$ такое, что $J(y)\subset J(x_0)$ при $\|x_0-y\|\leq\alpha$. В самом деле, пусть $a=f(x_0)-\max_{i\in J\setminus J(x_0)}f_i(x_0)>0$, $\varepsilon=a/3$ и $\alpha>0$ такое, что для всякого $i\in J$ выполнено $|f_i(y)-f_i(x_0)|\leq\varepsilon$ при $\|y-x_0\|\leq\alpha$. Легко проверить, что $|f(y)-f(x_0)|\leq\varepsilon$ при $\|y-x_0\|\leq\alpha$. Если $j\in J(y)$, то

$$f_j(x_0)\geq f_j(y)-\varepsilon = f(y)-\varepsilon \geq f(x_0)-2\varepsilon = \max_{i\in J(x_0)}f_i(x_0)+a-2\varepsilon > \max_{i\notin J(x_0)}f_i(x_0).$$

Поэтому $j\in J(x_0)$, т.е. $J(y)\subset J(x_0)$ при $\|x_0-y\|\leq\alpha$. Возьмем $z\in X$ и $x\in X$. Пусть $\|\lambda z\|\leq\alpha$. Аналогично имеем, что существует $\delta>0$ такое, что $J(y)\subset J(x_0+\lambda z)$ при $\|y-x_0-\lambda z\|\leq\delta$. Пусть $\mu>0$ достаточно малое число такое, что $\|2\lambda\mu x\|\leq\delta$.

$$f(x_0+\lambda z+2\lambda\mu x)-2f(x_0+\lambda z+\lambda\mu x)+f(x_0+\lambda z)=$$
$$=\max_{i\in J(x_0+\lambda z)}(f_i(x_0+\lambda z+2\lambda\mu x)+f_i(x_0+\lambda z))-2\max_{i\in J}f_i(x_0+\lambda z+\lambda\mu x)\leq$$
$$\leq \max_{i\in J(x_0+\lambda z)}(f_i(x_0+\lambda z+2\lambda\mu x)+f_i(x_0+\lambda z)-2f_i(x_0+\lambda z+\lambda\mu x))\leq$$
$$\leq \max_{i\in J(x_0)}(f_i(x_0+\lambda z+2\lambda\mu x)+f_i(x_0+\lambda z)-2f_i(x_0+\lambda z+\lambda\mu x)).$$

Поэтому получим, что

$$f^{\{2\}+}(x_0;x) \le b = \sup_{z\in X}\inf_{\alpha>0,\, i\in J(x_0)}\max \sup_{0<\lambda\le\alpha} \frac{1}{\lambda^2}(f_i(x_0+\lambda z+2\lambda\mu x)-2f_i(x_0+\lambda z+\lambda\mu x)+f_i(x_0+\lambda z)).$$

Из определения супремума следует, что для $\varepsilon>0$ существуют $z^\varepsilon\in X$ такие, что

$$\inf_{\alpha>0}\max_{i\in J(x_0)}\sup_{0<\lambda\le\alpha} \frac{1}{\lambda^2}(f_i(x_0+\lambda z^\varepsilon+2\lambda\mu x)-2f_i(x_0+\lambda z^\varepsilon+\lambda\mu x)+f_i(x_0+\lambda z^\varepsilon))>b-\varepsilon.$$

Ясно, что

$$\inf_{\alpha>0}\sup_{0<\lambda\le\alpha} \frac{1}{\lambda^2}(f_i(x_0+\lambda z^\varepsilon+2\lambda\mu x)-2f_i(x_0+\lambda z^\varepsilon+\lambda\mu x)+f_i(x_0+\lambda z^\varepsilon)) \le f_i^{\{2\}+}(x_0;\mu x).$$

Поэтому для $\varepsilon>0$ существуют $\alpha_i^\varepsilon>0$ такие, что

$$\sup_{0<\lambda\le\alpha_i^\varepsilon} \frac{1}{\lambda^2}(f_i(x_0+\lambda z^\varepsilon+2\lambda\mu x)-2f_i(x_0+\lambda z^\varepsilon+\lambda\mu x)+f_i(x_0+\lambda z^\varepsilon)) \le f_i^{\{2\}+}(x_0;\mu x)+\varepsilon.$$

Положив $\alpha^\varepsilon = \min_{i\in J(x_0)} \alpha_i^\varepsilon$ получим

$$\max_{i\in J(x_0)}\sup_{0<\lambda\le\alpha^\varepsilon} \frac{1}{\lambda^2}(f_i(x_0+\lambda_1 z^\varepsilon+2\lambda\mu x)-2f_i(x_0+\lambda z^\varepsilon+\lambda\mu x)+f_i(x_0+\lambda z^\varepsilon)) \le \max_{i\in J(x_0)} f_i^{\{2\}+}(x_0;\mu x)+\varepsilon.$$

Отсюда следует, что

$$f^{\{2\}+}(x_0;\mu x) \le \max_{i\in J(x_0)}\sup_{0<\lambda_1\le\alpha^\varepsilon} \frac{1}{\lambda^2}(f_i(x_0+\lambda z^\varepsilon+2\lambda\mu x)-2f_i(x_0+\lambda z^\varepsilon+\lambda\mu x)$$
$$+f_i(x_0+\lambda z^\varepsilon))+\varepsilon \le \max_{i\in J(x_0)} f_i^{\{2\}+}(x_0;\mu x)+2\varepsilon.$$

Отсюда при $\varepsilon\to 0$ получим $f^{\{2\}+}(x_0;\mu x) \le \max_{i\in J(x_0)} f_i^{\{2\}+}(x_0;\mu x)$. Поэтому $f^{\{2\}+}(x_0;x) \le \max_{i\in J(x_0)} f_i^{\{2\}+}(x_0;x)$ при $x\in X$. Лемма доказана.

**Рассмотрим несколько примеров.**

**Пример 2.** Если $f(x)=|x_1 x_2|$, где $x=(x_1,x_2),$ то легко проверяется, что

$$f^{(2)+}(0;x) = f^{(2)-}(0;x) = f^{2+}(0;x) = f^{2-}(0;x) = 2|x_1 x_2|.$$

Если $x_1\ne 0$ ($x_2\ne 0$), то положив $z_1=-x_1\ne 0$ ($z_2=-x_2\ne 0$) при $z_2\to\infty$ ($z_1\to\infty$) имеем, что

$$f^{\{2\}+}(0;x) = \begin{cases} 0; & x_1=x_2=0, \\ +\infty; & \text{в других случаях.} \end{cases}$$

Положив $z_1=-2x_1, z_2=0$ или $z_1=0, z_2=-2x_2$ имеем, что $f^{\{2\}-}(0;x)=-2|x_1 x_2|$.

Поэтому $d_2 f(0)=\emptyset$, $\left\{\begin{pmatrix}\beta_1 & \alpha \\ \alpha & \beta_2\end{pmatrix}: |\alpha|\le 1, \beta_1\ge 0, \beta_2\ge 0\right\}\subset D_2 f(0).$

**Пример 3.** Если $f(x)=\begin{cases} x^2; & x\ge 0 \\ 0; & x<0, \end{cases}$ то $f^{\{2\}+}(0;x)=2x^2$, $f^{\{2\}-}(0;x)=0$. Поэтому $D_2 f(0)=[0,2]$.

Ясно, что $f^{(2)+}(0;x)=f^{(2)-}(0;x)=\begin{cases} 2x^2: x\ge 0, \\ 0: x<0. \end{cases}$ Поэтому $d_2^+ f(x_0)=(-\infty,0]$, $d_2^- f(x_0)=[2,+\infty)$.

Легко проверяется, что

$$f^{\{2\}}(0; x_1, x_2) = \begin{cases} 2x_1 x_2 : x_1 x_2 \geq 0, \\ 0 \quad : x_1 x_2 < 0. \end{cases}$$

Отсюда следует, что $\partial_{\{2\}} f(0) = [0, 2]$.

Положим $f^{[2]+}(x_0; x) = \overline{\lim_{\lambda \downarrow 0, \, z \to x_0}} \dfrac{1}{\lambda^2} (f(z + 2\lambda x) - 2f(z + \lambda x) + f(z))$.

**Лемма 8.** Если $f$ $\{2\}$-липшицевая в окрестности точки $0$ и положительно однородная степени два функция, то

$$f^{[2]+}(0; x) = f^{\{2\}+}(0; x) = \sup_{z \in X}(f(z + 2x) - 2f(z + x) + f(z))$$

при $x \in X$.

**Доказательство.** По определению имеем, что

$$f^{\{2\}+}(0; x) = \sup_{z \in X} \overline{\lim_{\lambda \downarrow 0}} \dfrac{1}{\lambda^2}(f(\lambda z + 2\lambda x) - 2f(\lambda z + \lambda x) + f(\lambda z)) =$$
$$= \sup_{z \in X}(f(z + 2x) - 2f(z + x) + f(z)).$$

Так как $f$ положительно однородная степени два функция, то

$$\dfrac{1}{\lambda^2}(f(z + 2\lambda x) - 2f(z + \lambda x) + f(z)) = f(\dfrac{z}{\lambda} + 2x) - 2f(\dfrac{z}{\lambda} + x) + f(\dfrac{z}{\lambda}).$$

Ясно, что любой элемент $u \in X$ может быть записан в виде $u = \dfrac{z}{\lambda}$, где $\|z\| < \varepsilon, \, 0 < \lambda < \alpha$, $\varepsilon$ и $\alpha$ - произвольные положительные числа. Поэтому

$$\sup_{\|z\| < \varepsilon, \, 0 < \lambda < \alpha} \dfrac{1}{\lambda^2}(f(z + 2\lambda x) - 2f(z + \lambda x) + f(z)) = \sup_{\|z\| < \varepsilon, \, 0 < \lambda < \alpha} (f(\dfrac{z}{\lambda} + 2x) - 2f(\dfrac{z}{\lambda} + x) + f(\dfrac{z}{\lambda})) =$$
$$= \sup_{z \in X}(f(z + 2x) - 2f(z + x) + f(z)).$$

По определению

$$f^{[2]+}(0; x) = \inf_{\varepsilon > 0, \, \alpha > 0} \sup_{\|z\| < \varepsilon, \, 0 < \lambda < \alpha} \dfrac{1}{\lambda^2}(f(z + 2\lambda x) - 2f(z + \lambda x) + f(z)).$$

Тогда имеем

$$f^{[2]+}(0; x) = \sup_{z \in X}(f(z + 2x) - 2f(z + x) + f(z))$$

при $x \in X$. Лемма доказана.

**Следствие 2.** Если $f$ $\{2\}$-липшицевая в окрестности точки $0$ и положительно однородная степени два функция, то $f^{[2]+}(0; x) = f^{\{2\}+}(0; x) \geq 2f(x)$ и $f^{[2]+}(0; x) = f^{\{2\}+}(0; x) \geq f(x) + f(-x)$ при $x \in X$.

**Теорема 2.** Пусть $X, Y$ банаховы пространства, отображение $\varphi: X \to Y$ дифференцируемо по Фреше в точке $x_0$ и $g: Y \to R$ удовлетворяет $(\theta, \delta)$-сильно билипшицеву условию с постоянной $K$ в точке $\varphi(x_0)$. Тогда

$$(g \circ \varphi)^{(2)+}(x_0; x) = g^{(2)+}(\varphi(x_0); \varphi'(x_0)x)$$

при $x \in X$.

**Доказательство.** По определению

$$(g \circ \varphi)^{(2)+}(x_0; x) = \overline{\lim_{\lambda \downarrow 0}} \frac{1}{\lambda^2} (g(\varphi(x_0 + 2\lambda x)) - 2g(\varphi(x_0 + \lambda x)) + g(\varphi(x_0))) =$$

$$= \overline{\lim_{\lambda \downarrow 0}} \frac{1}{\lambda^2} (g(\varphi(x_0) + 2\lambda \varphi'(x_0)x + o(2\lambda x)) - 2g(\varphi(x_0) + \lambda \varphi'(x_0)x + o(\lambda x)) + g(\varphi(x_0))) \leq$$

$$\leq \overline{\lim_{\lambda \downarrow 0}} \frac{1}{\lambda^2} (g(\varphi(x_0) + 2\lambda \varphi'(x_0)x) + K\|o(2\lambda x)\|^\theta (4\lambda\|\varphi'(x_0)x\| + \|o(2\lambda x)\|)^{2-\theta} -$$

$$- 2g(\varphi(x_0) + \lambda \varphi'(x_0)x) + 2K\|o(\lambda x)\|^\theta (2\lambda\|\varphi'(x_0)x\| + \|o(\lambda x)\|)^{2-\theta} + g(\varphi(x_0))) =$$

$$= \overline{\lim_{\lambda \downarrow 0}} \frac{1}{\lambda^2} (g(\varphi(x_0) + 2\lambda \varphi'(x_0)x) - 2g(\varphi(x_0) + \lambda \varphi'(x_0)x) + g(\varphi(x_0))) = g^{(2)+}(\varphi(x_0); \varphi'(x_0)x).$$

где $\frac{\|o(h)\|_Y}{\|h\|_X} \to 0$ при $\|h\|_X \to 0$. Отметим, что при $x \in X$ существует $\delta > 0$ такое, что

$$\|g(\varphi(x_0) + 2\lambda \varphi'(x_0)x) - 2g(\varphi(x_0) + \lambda \varphi'(x_0)x) + g(\varphi(x_0))\| \leq$$

$$\leq \|g(\varphi(x_0) + 2\lambda \varphi'(x_0)x) - g(\varphi(x_0) + \lambda \varphi'(x_0)x)\| +$$

$$+ \|g(\varphi(x_0) + \lambda \varphi'(x_0)x) - g(\varphi(x_0))\| \leq K\|\lambda \varphi'(x_0)x\|^\theta \|3\lambda \varphi'(x_0)x\|^{2-\theta} +$$

$$+ K\|\lambda \varphi'(x_0)x\|^\theta \|\lambda \varphi'(x_0)x\|^{2-\theta} \leq K(3^{2-\theta} + 1)\lambda^2 \|\varphi'(x_0)x\|^\theta \|\varphi'(x_0)x\|^{2-\theta}$$

при $0 < \lambda \leq \delta$. Получим, что

$$(g \circ \varphi)^{(2)+}(x_0; x) \leq g^{(2)+}(\varphi(x_0); \varphi'(x_0)x)$$

при $x \in X$. Аналогично имеем, что

$$(g \circ \varphi)^{(2)+}(x_0; x) \geq g^{(2)+}(\varphi(x_0); \varphi'(x_0)x)$$

при $x \in X$. Поэтому

$$(g \circ \varphi)^{(2)+}(x_0; x) = g^{(2)+}(\varphi(x_0); \varphi'(x_0)x)$$

при $x \in X$. Теорема доказана.

**Теорема 3.** Пусть $X, Y$ банаховы пространства, отображение $\varphi: X \to Y$ дифференцируемо по Фреше в точке $x_0$, $\operatorname{Im} \varphi'(x_0) = Y$ и $g: Y \to R$ удовлетворяет $(\theta, \delta)$-сильно билипшицеву условию с постоянной $K$ в точке $\varphi(x_0)$ и удовлетворяет $\{2\}$-липшицеву условию с постоянной $K$ в окрестности точки $\varphi(x_0)$. Тогда

$$(g \circ \varphi)^{\{2\}+}(x_0; x) = g^{\{2\}+}(\varphi(x_0); \varphi'(x_0)x)$$

при $x \in X$.

**Доказательство.** По определению

$$(g \circ \varphi)^{\{2\}+}(x_0; x) = \sup_{z \in X} \overline{\lim_{\lambda \downarrow 0}} \frac{1}{\lambda^2} (g(\varphi(x_0 + \lambda z + 2\lambda x)) - 2g(\varphi(x_0 + \lambda z + \lambda x)) + g(\varphi(x_0 + \lambda z))) =$$

$$= \sup_{z \in X} \overline{\lim_{\lambda \downarrow 0}} \frac{1}{\lambda^2} (g(\varphi(x_0) + \lambda \varphi'(x_0)z + 2\lambda \varphi'(x_0)x + o(\lambda z + 2\lambda x)) - 2g(\varphi(x_0) + \lambda \varphi'(x_0)z + \lambda \varphi'(x_0)x +$$

$$+ o(\lambda z + \lambda x)) + g(\varphi(x_0) + \lambda \varphi'(x_0)z + o(\lambda z)) \leq$$

$$\leq \sup_{z \in X} \overline{\lim_{\lambda \downarrow 0}} \frac{1}{\lambda^2} (g(\varphi(x_0) + \lambda \varphi'(x_0)z + 2\lambda \varphi'(x_0)x) + K\|o(\lambda z + 2\lambda x)\|^\theta (\lambda\|2\varphi'(x_0)z + 4\varphi'(x_0)x\| +$$

$$+ \|o(\lambda z + 2\lambda x)\|)^{2-\theta} - 2g(\varphi(x_0) + \lambda \varphi'(x_0)z + \lambda \varphi'(x_0)x) + 2K\|o(\lambda z + \lambda x)\|^{\theta}(\lambda\|2\varphi'(x_0)z + 2\varphi'(x_0)x\| +$$

$$+ \|o(\lambda z + \lambda x)\|)^{2-\theta} + g(\varphi(x_0) + \lambda \varphi'(x_0)z) + K\|o(\lambda z)\|^{\theta}(2\lambda\|\varphi'(x_0)z\| + \|o(\lambda z)\|)^{2-\theta} =$$

$$= \sup_{y \in Y} \overline{\lim_{\lambda \downarrow 0}} \frac{1}{\lambda^2}(g(\varphi(x_0) + \lambda y + 2\lambda \varphi'(x_0)x) - 2g(\varphi(x_0) + \lambda y + \lambda \varphi'(x_0)x) + g(\varphi(x_0) + \lambda y)) =$$

$$= g^{\{2\}+}(\varphi(x_0); \varphi'(x_0)x),$$

где $\dfrac{\|o(h)\|_Y}{\|h\|_X} \to 0$ при $\|h\|_X \to 0$. Получим, что

$$(g \circ \varphi)^{\{2\}+}(x_0; x) \le g^{\{2\}+}(\varphi(x_0); \varphi'(x_0)x)$$

при $x \in X$. Аналогично имеем, что

$$(g \circ \varphi)^{\{2\}+}(x_0; x) \ge g^{\{2\}+}(\varphi(x_0); \varphi'(x_0)x)$$

при $x \in X$. Поэтому

$$(g \circ \varphi)^{\{2\}+}(x_0; x) = g^{\{2\}+}(\varphi(x_0); \varphi'(x_0)x)$$

при $x \in X$. Теорема доказана.

Пусть $A \subset B_0(X)$ ограниченное множество, т.е. существует число $L > 0$ такое, что $|Q(x)| \le L\|x\|^2$ при $Q \in A$ и $f(x) = \sup\limits_{Q \in A} Q(x)$. Тогда имеем

$$f(z+2x) - 2f(z+x) + f(z) = \sup_{Q \in A} Q(z+2x) - 2\sup_{Q \in A} Q(z+x) + \sup_{Q \in A} Q(z)) \ge -2\sup_{Q \in A} Q(z+x) +$$

$$+ \sup_{Q \in A}(Q(z+2x) + Q(z)) \ge \inf_{Q \in A}(Q(z+2x) - 2Q(z+x) + Q(z)) = 2\inf_{Q \in A} Q(x) \ge -2L\|x\|^2$$

при $x \in X$.

Пусть $f : X \to R$. Положим

$$\text{uep } f = \{(x, \alpha) \in X \times R_+ : f(x) \le \alpha^2\}, \qquad \text{lep } f = \{(x, \beta) \in X \times R_+ : f(x) \ge -\beta^2\}.$$

Ясно, что если $f(x) \ge 0$ при $x \in X$, то $f$ однозначно определяется через uep $f$, а если $f(x) \le 0$ при $x \in X$, то $f$ однозначно определяется через lep $f$. Геометрические аспекты бидифференциала и двадифференциала изучаются с помщью введенных множества uep $g$, lep $g$, uep $h$, lep $h$, uep $g_z$, lep $g_z$, uep $h_z$ и lep $h_z$, где

$$g(x) = f(x_0 + 2x) - 2f(x_0 + x) + f(x_0), \quad g_z(x) = f(z+2x) - 2f(z+x) + f(z),$$

$$h(x) = f(x_0 + x) - 2f(x_0) + f(x_0 - x), \quad h_z(x) = f(z+x) - 2f(z) + f(z-x).$$

Если $f$ удовлетворяет $\{2\}$–липшицеву условию в $\delta$–окрестности точки $x_0$, то $g_z(x) + K\|x\|^2 \ge 0$ и $g_z(x) - K\|x\|^2 \le 0$ при $x \in \delta B$.

Отметим, что аналогично [19] и [23] можно рассмотреть геометрические аспекты бидифференциала и двадифференциала.

## 1.3. О необходимых условиях экстремума второго порядка

Пусть $X$ банахово пространство, $C \subset X$ и $x_0 \in C$, $f : X \to R$, $\varphi : X \to R$, $\varphi_1 : X \to R$, $o : R_+ \to R_+$, где $\lim\limits_{t \downarrow 0} \dfrac{o(t)}{t} = 0$, $o(0) = 0$, $\varphi(0) = 0$, $\mu > 0$. Положим (см.[30])

$$K_\mu(x_0; C, \varphi) = \{x \in X : \exists \lambda_x > 0, \exists o_2(x, \lambda) : [0, \lambda_x] \to R_+, \text{где } \dfrac{o_2(x, \lambda)}{\lambda^\mu} \to 0 \text{ при } \lambda \downarrow 0, \text{что}$$

$x_0 + \lambda x \in C$, $\varphi(\lambda x) \leq o_2(x,\lambda)$ при $\lambda \in [0,\lambda_x]\}$.

Ясно, что

$K_2(x_0;X,\varphi) = \{x \in X: \exists \lambda_x > 0$, $\exists o_2(x,\lambda):[0,\lambda_x] \to R_+$, где $\dfrac{o_2(x,\lambda)}{\lambda^2} \to 0$ при $\lambda \downarrow 0$, что $\varphi(\lambda x) \leq o_2(x,\lambda)$ при $\lambda \in [0,\lambda_x]\}$.

**Теорема 1.** Если $f: X \to R$, $\varphi: X \to R$, $x_0$ является точкой локального минимума функции f на X и существуют $\delta > 0$ и $K > 0$ такие, что $|f(x_0+x) - \varphi(x) - f(x_0)| \leq K\|x\|^2$ при $x \in \delta B$, то

$$f_\varphi^{(2)+}(x_0;x) = \varlimsup_{\lambda \downarrow 0} \dfrac{1}{\lambda^2}(f(x_0+2\lambda x) - \varphi(2\lambda x) - 2f(x_0+\lambda x) + 2\varphi(\lambda x) + f(x_0)) \geq 0$$

при $x \in K_2(x_0;X,\varphi)$.

**Доказательство.** Если $x \in K_2(x_0;X,\varphi)$, то существует $\lambda_x > 0$ и $\exists o_2(x,\lambda):[0,\lambda_x] \to R_+$, где $\dfrac{o_2(x,\lambda)}{\lambda^2} \to 0$ при $\lambda \downarrow 0$, что $\varphi(\lambda x) \leq o_2(x,\lambda)$ при $\lambda \in [0,\lambda_x]$. Так как $x_0$ является точкой локального минимума функции f на X, то существует $\delta > 0$, где $\delta \leq \lambda_x$, такое, что $f(x_0+\lambda x) \geq f(x_0)$ при $\lambda \in [0,\delta]$. Так как $f(x_0+\lambda x) - f(x_0) \leq K\|\lambda x\|^2 + \varphi(\lambda x) \leq$
$\leq K\|x\|^2 + o_2(x,\lambda)$ при $\lambda \in [0,\lambda_x]$, то получим, что

$$f_\varphi^{(2)+}(x_0;x) = \varlimsup_{\lambda \downarrow 0} \dfrac{1}{\lambda^2}(f(x_0+2\lambda x) - \varphi(2\lambda x) - 2f(x_0+\lambda x) + 2\varphi(\lambda x) + f(x_0)) \geq$$

$$\geq \varlimsup_{\lambda \downarrow 0} \dfrac{1}{\lambda^2}(f(x_0+2\lambda x) - \varphi(2\lambda x) - f(x_0)) + 2\varlimsup_{\lambda \downarrow 0} \dfrac{1}{\lambda^2}(f(x_0) - f(x_0+\lambda x) + \varphi(\lambda x)) =$$

$$= \varlimsup_{\lambda \downarrow 0} \dfrac{4}{\lambda^2}(f(x_0+\lambda x) - \varphi(\lambda x) - f(x_0)) - 2\varlimsup_{\lambda \downarrow 0} \dfrac{1}{\lambda^2}(f(x_0+\lambda x) - \varphi(\lambda x) - f(x_0)) =$$

$$= 2\varlimsup_{\lambda \downarrow 0} \dfrac{1}{\lambda^2}(f(x_0+\lambda x) - \varphi(\lambda x) - f(x_0)) \geq 2\varlimsup_{\lambda \downarrow 0} \dfrac{1}{\lambda^2}(f(x_0+\lambda x) - o_2(x,\lambda) - f(x_0)) \geq$$

$$\geq 2\varlimsup_{\lambda \downarrow 0} \dfrac{1}{\lambda^2}(f(x_0+\lambda x) - f(x_0)) + 2\varlimsup_{\lambda \downarrow 0} \dfrac{1}{\lambda^2}(-o_2(x,\lambda)) \geq 0$$

при $x \in K_2(x_0;X,\varphi)$. Теорема доказана.

**Следствие 1.** Если X банахово пространство, функция $f: X \to R$ достигает локального минимума в точке $x_0$, существуют $\delta > 0$ и $K > 0$ такие, что $|f(x_0+x) - f(x_0)| \leq K\|x\|^2$ при $x \in \delta B$, то $f^{(2)+}(x_0;x) \geq 0$ при $x \in X$.

Если существуют $\delta > 0$ и $K > 0$ такие, что $f(x_0+x) - f(x_0) \leq K\|x\|^2$ при $x \in \delta B$, то следствие 1 также верно.

Отметим, что если f удовлетворяет условию $|f(x_0+x) - f(x_0)| \leq K\|x\|^2$ при $x \in \delta B$, то

$$|f(x_0+2x) - 2f(x_0+x) + f(x_0)| \leq |f(x_0+2x) - f(x_0)| +$$
$$+ 2|f(x_0+x) - f(x_0)| \leq 4K\|x\|^2 + 2K\|x\|^2 \leq 6K\|x\|^2$$

при $x \in 0{,}5\,\delta B$.

Если f удовлетворяет 2–липшицеву условию в окрестности точки $x_0$, то

$$f^{[2]}(x_0;x,x) \geq f^{\{2\}}(x_0;x,x) \geq f^{\{2\}+}(x_0;x) = \sup_{y \in X} \overline{\lim_{\lambda \downarrow 0}} \frac{1}{\lambda^2}(f(x_0+\lambda y+\lambda x) - 2f(x_0+\lambda y)+$$
$$+ f(x_0+\lambda y-\lambda x)) \geq \overline{\lim_{\lambda \downarrow 0}} \frac{1}{\lambda^2}(f(x_0+\lambda x) - 2f(x_0) + f(x_0-\lambda x)).$$

Поэтому, если функция $f: X \to R$ достигает локального минимума в точке $x_0$, то
$$f^{[2]}(x_0;x,x) \geq f^{\{2\}}(x_0;x,x) \geq f^{\{2\}+}(x_0;x) \geq f^{2+}(x_0;x) \geq 0$$
при $x \in X$.

Пусть $\omega: R_+ \to R_+$, где $\omega(0) = 0$.

**Теорема 2.** Пусть $x_0$ является минимумом функции $f$ на множестве $C$, $f$ в точке $x_0$ удовлетворяет $\varphi$-$(\alpha,\beta,\nu,\delta,\omega)$ локально липшицеву условию с постоянной $K$ в точке $x_0$, $C \subset B(x_0,\delta)$ и $D = \{x \in C: \varphi(x-x_0) \leq 0\}$. Тогда для любого $\lambda \geq K$ функция
$$S_\lambda(x) = f(x) - \varphi(x-x_0) + \lambda(d_D^\alpha(x) + \|x-x_0\|^{\beta-\alpha\nu} d_D^\nu(x)) + \omega(\|x-x_0\|^\beta)$$
достигает минимума на $B(x_0,\delta)$ в точке $x_0$ и если $\lambda > K$ и $D$ замкнуто, то любая точка, минимизирующая $S_\lambda(x)$ на множестве $B(x_0,\delta)$ принадлежит $D$.

Теорема 2 доказана в [23] (см. следствие 4.5.2)

**Теорема 3.** Если $x_0$ является минимумом функции $f$ на множестве $C$, $f$ в точке $x_0$ удовлетворяет $\varphi$-$(1,2,1,\delta,o(2))$ локально липшицеву условию с постоянной $K$ в точке $x_0$, $C \subset B(x_0,\delta)$ и $D = \{x \in C: \varphi(x-x_0) \leq 0\}$, то $S_\lambda^{(2)+}(x_0;x) \geq 0$ при $x \in X$.

**Доказательство.** По теореме 2 для любого $\lambda \geq K$ функция
$$S_\lambda(x) = f(x) - \varphi(x-x_0) + \lambda(d_D^2(x) + \|x-x_0\| d_D(x)) + o(\|x-x_0\|^2)$$
достигает минимума на $B(x_0,\delta)$ в точке $x_0$. Ясно, что
$$|S_\lambda(x_0+x) - S_\lambda(x_0)| =$$
$$= |f(x_0+x) - \varphi(x) + \lambda(d_D^2(x_0+x) + \|x\| d_D(x_0+x)) + o(\|x\|^2) - f(x_0) + \lambda d_D^2(x_0)| \leq$$
$$\leq |f(x_0+x) - \varphi(x) - f(x_0)| + \lambda|d_D^2(x_0+x) - d_D^2(x_0)| + \lambda\|x\| d_D(x_0+x) + o(\|x\|^2)$$
при $x \in \delta B$. Так как $f$ удовлетворяет $\varphi$-$(1,2,1,\delta,o(2))$ локально липшицеву условию с постоянной $K$ в точке $x_0$, то
$$|f(x_0+x) - \varphi(x) - f(x_0)| \leq K\|x\|^2$$
при $x \in \delta B$. Ясно, что
$$|d_D^2(x_0+x) - d_D^2(x_0)| \leq \|x\|^2, \qquad \|x\| d_D(x_0+x) \leq \|x\|^2$$
при $x \in X$. Поэтому имеем, что
$$|S_\lambda(x_0+x) - S_\lambda(x_0)| \leq K\|x\|^2 + \lambda\|x\|^2 + \lambda\|x\|^2 + o(\|x\|^2) \leq (K+2\lambda)\|x\|^2 + o(\|x\|^2)$$
при $x \in \delta B$. Тогда из следствия 1 следует, что $S_\lambda^{(2)+}(x_0;x) \geq 0$ при $x \in X$. Теорема доказана.

Положим $d_D^{(1)+}(x_0;x) = \overline{\lim_{t \downarrow 0}} \frac{1}{t} d_D(x_0 + tx).$

**Следствие 2.** Если выполняется условие теоремы 3, то
$$f_\varphi^{(2)+}(x_0;x) + 8\lambda\|x\| d_D^{(1)+}(x_0;x) \geq 0$$
при $x \in X$, где $\lambda \geq K$.

**Доказательство.** Из теоремы 3 следует, что
$$0 \leq S_\lambda^{(2)+}(x_0;x) = \overline{\lim_{t \downarrow 0}} \frac{1}{t^2}(S_\lambda(x_0+2tx) - 2S_\lambda(x_0+tx) + S_\lambda(x_0)) \leq$$

$$\leq \varlimsup_{t\downarrow 0}\frac{1}{t^2}(f(x_0+2tx)-\varphi(2tx)-2f(x_0+tx)+2\varphi(tx)+f(x_0))+$$

$$+\varlimsup_{t\downarrow 0}\frac{\lambda}{t^2}(d_D^2(x_0+2tx)-2d_D^2(x_0+tx)+d_D^2(x_0))+$$

$$+\varlimsup_{t\downarrow 0}\frac{\lambda}{t^2}(2t\|x\|d_D(x_0+2tx)-2t\|x\|d_D(x_0+tx))\leq f_\varphi^{(2)+}(x_0;x)+$$

$$+\varlimsup_{t\downarrow 0}\frac{\lambda}{t^2}d_D^2(x_0+2tx)+\varlimsup_{t\downarrow 0}(-\frac{2\lambda}{t^2}d_D^2(x_0+tx))+$$

$$+\varlimsup_{t\downarrow 0}\frac{2\lambda\|x\|}{t}d_D(x_0+2tx)+\varlimsup_{t\downarrow 0}\frac{-2\lambda\|x\|}{t}d_D(x_0+tx))\leq f_\varphi^{(2)+}(x_0;x)+$$

$$+2\lambda\|x\|\varlimsup_{t\downarrow 0}\frac{1}{t}d_D(x_0+2tx)+2\lambda\|x\|\varlimsup_{t\downarrow 0}\frac{1}{t}d_D(x_0+2tx)=f_\varphi^{(2)+}(x_0;x)+$$

$$+4\lambda\|x\|\varlimsup_{t\downarrow 0}\frac{1}{t}d_D(x_0+2tx)=f_\varphi^{(2)+}(x_0;x)+8\lambda\|x\|\varlimsup_{t\downarrow 0}\frac{1}{t}d_D(x_0+tx).$$

Поэтому
$$f_\varphi^{(2)+}(x_0;x)+8\lambda\|x\|d_D^{(1)+}(x_0;x)\geq 0$$
при $x\in X$, где $\lambda\geq K$. Следствие доказано.

**Следствие 3.** Если выполняется условие теоремы 3, то $f_\varphi^{(2)+}(x_0;x)\geq 0$ при $x\in K_D(x_0)=\{x\in X: d_D^{(1)+}(x_0;x)=0\}$.

Если $\varphi:X\to R$, где $\varphi(0)=0$, и $C\subset X$, то положим

$$T_{\alpha,\mu}(x_0;C,\varphi)=\{x\in X:\ \exists o(x,\lambda):[0,\lambda_x]\to R_+,\ \text{где}\ \frac{o(x,\lambda)}{\lambda^\mu}\to 0\ \text{при}\ \lambda\downarrow 0\ \text{и}\ \exists\lambda_i\downarrow 0,\ \exists\upsilon_i\subset X,$$

где $\frac{1}{\lambda_i^{\alpha-1}}\|\upsilon_i-x\|\to 0$ при $i\to+\infty$, что $x_0+\lambda_i\upsilon_i\in C$, $\varphi(\lambda_i\upsilon_i)\leq o(x,\lambda_i)$ при всех $i\}$.

**Теорема 4.** Если $X$ банахово пространство, $x_0$ является точкой минимума функции $f$ на множестве $C$, существуют $\alpha>0$, $0<\nu\leq 2$, $\mu>0$, $\sigma>0$, где $\mu\geq 2-\alpha\nu$, $\sigma\geq\frac{2-\alpha\nu}{\alpha}$, конечная положительно однородная степени $\mu$ функция $\varphi_1$, функция $o:R_+\to R_+$, где $\lim_{t\downarrow 0}\frac{o(t)}{t}=0$, числа $\delta>0$ и $K$ такие, что

$$|f(x_0+x+y)-f(x_0+x)-\varphi(x+y)+\varphi(x)|\leq K\|y\|^\nu(\varphi_1(x)+\|y\|^\sigma)+o(\|x\|^2)$$

для $x\in T_{\alpha,2}(x_0;C,\varphi), \|x\|\leq\delta$, $y\in X$, $\|y\|\leq\|x\|$, $x_0+x+y\in C$, то

$$\varlimsup_{\lambda\downarrow 0}\frac{1}{\lambda^2}(f(x_0+\lambda x)-\varphi(\lambda x)-f(x_0))\geq 0\ \text{при}\ x\in T_{\alpha,2}(x_0;C,\varphi).$$

Теорема 4 доказана в [31] (см. теорему 4.1).

**Следствие 4.** Если выполняется условие теоремы 4, то

$$f_\varphi^{(2)+}(x_0;x)=\varlimsup_{\lambda\downarrow 0}\frac{1}{\lambda^2}(f(x_0+2\lambda x)-\varphi(2\lambda x)-2f(x_0+\lambda x)+2\varphi(\lambda x)+f(x_0))\geq 0$$

при $x\in T_{\alpha,2}(x_0;C,\varphi)$.

**Доказательство.** Из теоремы 4 имеем, что
$$\varlimsup_{\lambda\downarrow 0}\frac{1}{\lambda^2}(f(x_0+\lambda x)-\varphi(\lambda x)-f(x_0))\geq 0$$

при $x\in T_{\alpha,2}(x_0;C,\varphi)$. Поэтому

$$f_\varphi^{(2)+}(x_0;x)=\varlimsup_{\lambda\downarrow 0}\frac{1}{\lambda^2}(f(x_0+2\lambda x)-\varphi(2\lambda x)-2f(x_0+\lambda x)+2\varphi(\lambda x)+f(x_0))\geq$$

$$\geq \overline{\lim_{\lambda \downarrow 0}} \frac{1}{\lambda^2}(f(x_0+2\lambda x)-\varphi(2\lambda x)-f(x_0))+\overline{\lim_{\lambda \downarrow 0}} \frac{1}{\lambda^2}(2f(x_0)-2f(x_0+\lambda x)+2\varphi(\lambda x))=$$

$$=\overline{\lim_{\lambda \downarrow 0}} \frac{4}{\lambda^2}(f(x_0+\lambda x)-\varphi(\lambda x)-f(x_0))-\overline{\lim_{\lambda \downarrow 0}} \frac{2}{\lambda^2}(f(x_0+\lambda x)-\varphi(\lambda x)-f(x_0))=$$

$$=\overline{\lim_{\lambda \downarrow 0}} \frac{2}{\lambda^2}(f(x_0+\lambda x)-\varphi(\lambda x)-f(x_0))\geq 0$$

при $x \in T_{\alpha,2}(x_0;C,\varphi)$. Теорема доказана.

Положим $g_k(x)=f(x)+kd_2(x)$, где $d_2(x)=d_C^2(x)$.

Если из $\overline{\lim_{t \downarrow 0}} \frac{1}{t}(d(x_0+tx)-d(x_0))=0$ следует, что $\overline{\lim_{t \downarrow 0}} \frac{1}{t^2}(d(x_0+tx)-d(x_0))=0$, то множество $C$ назовем регулярным второго порядка в точке $x_0$.

**Теорема 5.** Если X банахово пространство, $x_0$ является точкой минимума функции f на множестве C, f липшицевая функция в окрестности точки $x_0$ и двалипшицевая функция с постоянной L в точке $x_0$, множество C регулярно второго порядка в точке $x_0$ или f $(\theta,\delta)$ – сильно билипшицевая функция с постоянной L в точке $x_0$, то для любого $x \in X$ существует число $k > 0$ такое, что $g_k^{2+}(x_0;x) \geq 0$ при $x \in X$.

**Доказательство.** Так как $f(x)$ двалипшицевая функция в точке $x_0$, то функция $x \to f^{2+}(x_0;x)$ конечна для любого $x \in X$ и

$$g_k^{2+}(x_0;x)=\overline{\lim_{t \downarrow 0}} \frac{1}{t^2}(f(x_0+tx)-2f(x_0)+f(x_0-tx)+k(d_2(x_0+tx)-2d_2(x_0)+d_2(x_0-tx)))\geq$$

$$\geq \overline{\lim_{t \downarrow 0}} \frac{1}{t^2}(-L\|tx\|^2+k(d_2(x_0+tx)-2d_2(x_0)+d_2(x_0-tx)))=-L\|x\|^2+kd_2^{2+}(x_0;x)$$

при $x \in X$.

Из определения $d_2^{2+}(x_0;x)$ имеем, что

$$d_2^{2+}(x_0;x)=\overline{\lim_{t \downarrow 0}} \frac{1}{t^2}(d_2(x_0+tx)-2d_2(x_0)+d_2(x_0-tx))\geq 0.$$

Если $d_2^{2+}(x_0;x)>0$, то $g_k^{2+}(x_0;x)\geq 0$ при $k \geq \dfrac{L\|x\|^2}{d_2^{2+}(x_0;x)}$.

Пусть $d_2^{2+}(x_0;x)=0$, т.е.

$$\overline{\lim_{t \downarrow 0}} \frac{1}{t^2}(d_2(x_0+tx)-2d_2(x_0)+d_2(x_0-tx))=0.$$

Отсюда следует, что $\overline{\lim_{t \downarrow 0}} \frac{1}{t}d(x_0+tx)=0$ и $\overline{\lim_{t \downarrow 0}} \frac{1}{t}d(x_0-tx)=0$. Так как множество C регулярно второго порядка в точке $x_0$, то имеем, что $\overline{\lim_{t \downarrow 0}} \frac{1}{t^2}d(x_0+tx)=0$ и $\overline{\lim_{t \downarrow 0}} \frac{1}{t^2}d(x_0-tx)=0$. Поэтому существуют $t_x>0$, $o_1(x,t):[0,t_x]\to X$, $o_2(x,t):[0,t_x]\to R_+$, где $\dfrac{o_1(x,t)}{t^2}\to 0$, $\dfrac{o_2(x,t)}{t^2}\to 0$ при $t \downarrow 0$, такие, что $x_0+tx+o_1(x,t)\in C$ и $x_0-tx+o_2(x,t)\in C$. Так как функция f удовлетворяет условию Липшица в окрестности точки $x_0$, то

$$g_k^{2+}(x_0;x) \geq \overline{\lim_{t\downarrow 0}} \frac{1}{t^2}(f(x_0+tx)-2f(x_0)+f(x_0-tx)+k\overline{\lim_{t\downarrow 0}}\frac{1}{t^2}(d_2(x_0+tx)-$$
$$-2d_2(x_0)+d_2(x_0-tx))=\overline{\lim_{t\downarrow 0}}\frac{1}{t^2}(f(x_0+tx)-2f(x_0)+f(x_0-tx))=$$
$$=\overline{\lim_{t\downarrow 0}}\frac{1}{t^2}(f(x_0+tx+o_1(x,t))-2f(x_0)+f(x_0-tx+o_2(x,t))) \geq 0$$

при $k \geq 0$.

Пусть $d_2^{2+}(x_0;x)=0$ и $f$ $(\theta,\delta)-$ сильно билипшицевая функция с постоянной $L$ в точке $x_0$.

Если $d_2^{2+}(x_0;x)=0$, то $\overline{\lim_{t\downarrow 0}}\frac{1}{t}d(x_0+tx)=0$ и $\overline{\lim_{t\downarrow 0}}\frac{1}{t}d(x_0-tx)=0$. Поэтому существуют $t_x>0$, $o_1(x,t):[0,t_x]\to X$, $o_2(x,t):[0,t_x]\to R_+$, где $\frac{o_1(x,t)}{t}\to 0$, $\frac{o_2(x,t)}{t}\to 0$ при $t\downarrow 0$, такие, что $x_0+tx+o_1(x,t)\in C$ и $x_0-tx+o_2(x,t)\in C$. Так как функция $f$ удовлетворяет $(\theta,\delta)-$ сильно билипшицеву условию с постоянной $L$ в точке $x_0$, то аналогично имеем, что $g_k^{2+}(x_0;x)\geq 0$ при $x\in X$. Теорема доказана.

**Следствие 5.** Если $X$ банахово пространство, $x_0$ является точкой минимума функции $f$ на множестве $C$, $f$ липшицевая и $2-$липшицевая функция с постоянной $L$ в окрестности точки $x_0$, $d_2(x)$ $2-$липшицевая функция в окрестности точки $x_0$ и множество $C$ регулярно второго порядка в точке $x_0$ или $f$ $(\theta,\delta)-$ сильно билипшицевая функция с постоянной $L$ в точке $x_0$, то функции $g_k^{[2]}(x_0;x,x)$, $g_k^{\{2\}+}(x_0;x,x)$, $g_k^{\{2\}+}(x_0;x)$ принимают конечные значения и для любого $x\in X$ существует число $k>0$ такие, что $g_k^{[2]}(x_0;x,x)\geq g_k^{\{2\}+}(x_0;x,x)\geq g_k^{\{2\}+}(x_0;x)\geq g_k^{2+}(x_0;x)\geq 0$ при $x\in X$.

**Теорема 6.** Если $X$ гильбертово пространство, $x_0$ является точкой минимума функции $f$ на множестве $C$, $f$ липшицевая и $2-$липшицевая функция с постоянной $L$ в окрестности точки $x_0$, $C$ замкнутое выпуклое множество и регулярно второго порядка в точке $x_0$ или $f$ $(\theta,\delta)-$ сильно билипшицевая функция с постоянной $L$ в точке $x_0$, то $\max\{b(x,x):b\in \partial_{\{2\}}f(x_0)+\Omega_C(x_0)\}\geq 0$ при $x\in X$.

**Доказательство.** Из [16], стр.19 (см. раздел 1.1) следует, что $d_2(x)$ $2-$липшицевая функция в окрестности точки $x_0$. Тогда по следствию 5 для любого $\bar{x}\in X$ существует $k>0$ такое, что $g_k^{\{2\}+}(x_0;\bar{x},\bar{x})\geq 0$. Поэтому имеем, что
$$0\leq g_k^{\{2\}+}(x_0;\bar{x},\bar{x})\leq f^{\{2\}+}(x_0;\bar{x},\bar{x})+kd_2^{\{2\}+}(x_0;\bar{x},\bar{x}).$$

Так как $d_2(x)$ $2$-липшицевая функция в окрестности точки $x_0$, то используя из [23], стр.89 имеем, что
$$0\leq \max\{b(\bar{x},\bar{x}):b\in \partial_{\{2\}}f(x_0)\}+\max\{b(\bar{x},\bar{x}):b\in k\partial_{\{2\}}d_2(x_0)\}=$$
$$=\max\{b(\bar{x},\bar{x}):b\in \partial_{\{2\}}f(x_0)+k\partial_{\{2\}}d_2(x_0)\}\leq \max\{b(\bar{x},\bar{x}):b\in \partial_{\{2\}}f(x_0)+\Omega_C(x_0)\}.$$

Отсюда следует, что
$$\max\{b(x,x):b\in \partial_{\{2\}}f(x_0)+\Omega_C(x_0)\}\geq 0$$

при $x\in X$. Теорема доказана.

Положим $g_k(x)=f(x)-\langle f'(x_0),x-x_0\rangle+kd_2(x)$, где $D=\{x\in C:\langle f'(x_0),x-x_0\rangle\leq 0\}$, $d_2(x)=d_D^2(x)$.

**Теорема 7.** Если $X$ банахово пространство, $x_0$ является точкой минимума функции $f$ на множестве $C$, $f$ липшицевая функция с постоянной $L$ в окрестности точки $x_0$,

дифференцируема в точке $x_0$ и
$$|f(x_0 + x + y) - f(x_0 + x) - \langle f'(x_0), y \rangle| \leq L\|y\|^{\nu}(\|x\|^{2-\nu} + \|y\|^{2-\nu}) + o(\|x\|^2) \qquad (1)$$
при $x, y \in X, \|x\| \leq \delta, \|y\| \leq \delta,$ то для любого $x \in X$ существует число $k > 0$ такое, что $g_k^{2+}(x_0; x) \geq 0$ при $x \in X$.

**Доказательство.** Ясно, что $x_0 \in X$ является минимумом функции $f(x) - \langle f'(x_0), x - x_0 \rangle$ на множестве $D = \{x \in C : \langle f'(x_0, x - x_0) \rangle \leq 0\}$. Из (1) имеем, что
$$|f(x_0 + x) - \langle f'(x_0), x \rangle - 2f(x_0) + f(x_0 - x) - \langle f'(x_0), -x \rangle| \leq 2L\|x\|^2$$
при $x \in X, \|x\| \leq \delta$. Так как $f(x) - \langle f'(x_0), x - x_0 \rangle$ и $d_2(x)$ двалипшицевые функции в точке $x_0$, то функция $x \to g_k^{2+}(x_0; x)$ конечна для любого $x \in X$ и
$$g_k^{2+}(x_0; x) = \overline{\lim_{t \downarrow 0}} \frac{1}{t^2}(f(x_0 + tx) - \langle f'(x_0), tx \rangle - 2f(x_0) + f(x_0 - tx) - \langle f'(x_0), -tx \rangle +$$
$$+ k(d_2(x_0 + tx) - 2d_2(x_0) + d_2(x_0 - tx))) \geq \overline{\lim_{t \downarrow 0}} \frac{1}{t^2}(-L\|tx\|^2 +$$
$$+ k(d_2(x_0 + tx) - 2d_2(x_0) + d_2(x_0 - tx))) = -L\|x\|^2 + kd_2^{2+}(x_0; x)$$
при $x \in X$.

Ясно, что
$$d_2^{2+}(x_0; x) = \overline{\lim_{t \downarrow 0}} \frac{1}{t^2}(d_2(x_0 + tx) - 2d_2(x_0) + d_2(x_0 - tx)) \geq 0.$$

Если $d_2^{2+}(x_0; x) > 0$, то $g_k^{2+}(x_0; x) \geq 0$ при $k \geq \dfrac{L\|x\|^2}{d_2^{2+}(x_0; x)}$.

Если $d_2^{2+}(x_0; x) = 0$, то $\overline{\lim_{t \downarrow 0}} \dfrac{1}{t^2}(d_2(x_0 + tx) - 2d_2(x_0) + d_2(x_0 - tx)) = 0$.

Отсюда следует, что $\overline{\lim_{t \downarrow 0}} \dfrac{1}{t}d(x_0 + tx) = 0$ и $\overline{\lim_{t \downarrow 0}} \dfrac{1}{t}d(x_0 - tx) = 0$.

Поэтому существуют $t_x > 0$, $o_1(x, \lambda) : [0, t_x] \to X$, где $\dfrac{o_1(x, t)}{t} \to 0$, $\dfrac{o_2(x, t)}{t} \to 0$ при $t \downarrow 0$, такие, что $x_0 + tx + o_1(x, t) \in D$ и $x_0 - tx + o_2(x, t) \in D$. Так как функция $f$ удовлетворяет условию (1), то
$$g_k^{2+}(x_0; x) \geq k \cdot \lim_{t \downarrow 0} \frac{1}{t^2}((d_2(x_0 + tx) - 2d_2(x_0) + d_2(x_0 - tx)) +$$
$$+ \overline{\lim_{t \downarrow 0}} \frac{1}{t^2}(f(x_0 + tx) - \langle f'(x_0), tx \rangle - 2f(x_0) + f(x_0 - tx) - \langle f'(x_0), -tx \rangle) =$$
$$= \overline{\lim_{t \downarrow 0}} \frac{1}{t^2}(f(x_0 + tx) - \langle f'(x_0), tx \rangle - 2f(x_0) + f(x_0 - tx) - \langle f'(x_0), -tx \rangle) =$$
$$= \overline{\lim_{t \downarrow 0}} \frac{1}{t^2}(f(x_0 + tx + o_1(x, t)) - \langle f'(x_0), tx + o_1(x, t) \rangle - 2f(x_0) + f(x_0 - tx + o_2(x, t)) -$$
$$- \langle f'(x_0), -tx + o_2(x, t) \rangle) + \lim_{t \downarrow 0} \frac{1}{t^2}(f(x_0 + tx) - f(x_0 + tx + o_1(x, t)) - \langle f'(x_0), tx \rangle +$$
$$+ \langle f'(x_0), tx + o_1(x, t) \rangle) + f(x_0 - tx) - f(x_0 - tx + o_2(x, t)) - \langle f'(x_0), -tx \rangle +$$
$$+ \langle f'(x_0), -tx + o_2(x, t) \rangle) \geq \overline{\lim_{t \downarrow 0}} \frac{1}{t^2}(f(x_0 + tx + o_1(x, t)) - \langle f'(x_0), tx + o_1(x, t) \rangle - 2f(x_0) +$$
$$+ f(x_0 - tx + o_2(x, t))) - \langle f'(x_0), -tx + o_2(x, t) \rangle) +$$

$$+ \varliminf_{t\downarrow 0} \frac{1}{t^2}(f(x_0+tx)-f(x_0+tx+o_1(x,t))-\langle f'(x_0),tx\rangle+\langle f'(x_0),tx+o_1(x,t)\rangle)+$$

$$+ \varliminf_{t\downarrow 0} \frac{1}{t^2}(f(x_0-tx)-f(x_0-tx+o_2(x,t))-\langle f'(x_0),-tx\rangle+\langle f'(x_0),-tx+o_2(x,t)\rangle) \ge$$

$$\ge \varlimsup_{t\downarrow 0} \frac{1}{t^2}(f(x_0+tx+o_1(x,t))-\langle f'(x_0),tx+o_1(x,t)\rangle-2f(x_0)+f(x_0-tx+o_2(x,t))-$$

$$-\langle f'(x_0),-tx+o_2(x,t)\rangle)-\varlimsup_{t\downarrow 0} \frac{1}{t^2}(L\|o_1(x,t)\|^\nu(\|tx\|^{2-\nu}+\|o_1(x,t)\|^{2-\nu})+o(\|tx\|^2))+$$

$$-\varlimsup_{t\downarrow 0} \frac{1}{t^2}(L\|o_2(x,t)\|^\nu(\|tx\|^{2-\nu}+\|o_2(x,t)\|^{2-\nu})+o(\|tx\|^2)) \ge 0$$

при $k \ge 0$. Теорема доказана.

**Следствие 6.** Если выполняется условие теоремы 7, $f$ $2-$ липшицевая функция с постоянной $L$ в окрестности точки $x_0$, $d_2(x)$ $2-$ липшицевая функция в окрестности точки $x_0$, то функции $g_k^{[2]}(x_0;x,x)$, $g_k^{\{2\}+}(x_0;x,x)$, $g_k^{\{2\}+}(x_0;x)$ принимают конечные значения и для любого $x \in X$ существует число $k > 0$ такое, что $g_k^{[2]}(x_0;x,x) \ge g_k^{\{2\}+}(x_0;x,x) \ge g_k^{\{2\}+}(x_0;x) \ge 0$ при $x \in X$.

**Теорема 8.** Если $X$ гильбертово пространство, $x_0$ является точкой минимума функции $f$ на множестве $C$, выполняется условие теоремы 7, $f$ $2-$ липшицевая функция с постоянной $L$ в окрестности точки $x_0$, $C$ замкнутое выпуклое множество, то $\max\{b(x,x):b \in \partial_{\{2\}}f(x_0)+\Omega_D(x_0)\} \ge 0$ при $x \in X$.

**Доказательство.** Из [17], стр.19 следует, что $d_2(x)$ $2-$ липшицевая функция в окрестности точки $x_0$. Тогда по следствию 6 для любого $\bar{x} \in X$ существует $k > 0$ такое, что $g_k^{\{2\}+}(x_0;\bar{x},\bar{x}) \ge 0$. Поэтому имеем

$$0 \le g_k^{\{2\}+}(x_0;\bar{x},\bar{x}) \le f^{\{2\}+}(x_0;\bar{x},\bar{x})+kd_2^{\{2\}+}(x_0;\bar{x},\bar{x}).$$

Так как $d_2(x)$ $2-$ липшицевая функция в окрестности точки $x_0$, то используя из [23], стр.89 имеем, что

$$0 \le \max\{b(\bar{x},\bar{x}):b \in \partial_{\{2\}}f(x_0)\}+\max\{b(\bar{x},\bar{x}):b \in k\partial_{\{2\}}d_2(x_0)\}=$$
$$= \max\{b(\bar{x},\bar{x}):b \in \partial_{\{2\}}f(x_0)+k\partial_{\{2\}}d_2(x_0)\} \le$$
$$\le \max\{b(\bar{x},\bar{x}):b \in \partial_{\{2\}}f(x_0)+\Omega_D(x_0)\}.$$

Отсюда следует, что

$$\max\{b(x,x):b \in \partial_{\{2\}}f(x_0)+\Omega_D(x_0)\} \ge 0$$

при $x \in X$. Теорема доказана.

Если существует $\delta > 0$ такое, что $C \cap B(x_0,\delta)$ выпуклое множество, то множество $C$ называется локально выпуклым в точке $x_0$. Отметим, что если в теореме 6 и 8 множество $C$ локально выпукло и локально замкнуто в точке $x_0$, то теоремы 6 и 8 также верны.

Будем говорить, что нижняя производная $x \to f^{(2)-}(x_0;x)$ существует равномерно по направлениям, если функция $x \to f^{(2)-}(x_0;x)$ конечна при $x \in X$ и для любого $\varepsilon > 0$ существует число $\delta > 0$ такое, что

$$\inf_{0 < t \le \delta} \frac{1}{t^2}(f(x_0+2tx)-2f(x_0+tx)+f(x_0)) > f^{(2)-}(x_0;x)-\varepsilon$$

при $x \in X$, $\|x\|=1$.

**Теорема 9.** Если X банахово пространство, существуют числа $\delta > 0$ и $K > 0$ такие, что $|f(x_0 + x) - f(x_0)| \leq K\|x\|^2$ при $x \in \delta B$, $x \to f^{(2)-}(x_0; x)$ существует равномерно по направлениям и существует число $\alpha > 0$ такое, что $f^{(2)-}(x_0; x) \geq \alpha\|x\|^2$ при $x \in X$, то функция $f : X \to R$ достигает локального минимума в точке $x_0$.

**Доказательство.** Так как $x \to f^{(2)-}(x_0; x)$ существует равномерно по направлениям, то для любого $\varepsilon > 0$, где $\alpha > \varepsilon$, существует число $\delta > 0$ такое, что

$$\alpha - \varepsilon \leq \inf_{0 < \lambda \leq \delta} \frac{1}{\lambda^2}(f(x_0 + 2\lambda x) - 2f(x_0 + \lambda x) + f(x_0)) \leq$$

$$\leq \inf_{0 < \lambda \leq \delta} \frac{1}{\lambda^2}(f(x_0 + 2\lambda x) - f(x_0)) + \sup_{0 < \lambda \leq \delta} \frac{1}{\lambda^2}(2f(x_0) - 2f(x_0 + \lambda x)) \leq$$

$$\leq \inf_{0 < \lambda \leq 2\delta} \frac{4}{\lambda^2}(f(x_0 + \lambda x) - f(x_0)) - 2\inf_{0 < \lambda \leq \delta} \frac{1}{\lambda^2}(f(x_0 + \lambda x) - f(x_0)) \leq$$

$$\leq 2\inf_{0 < \lambda \leq \delta} \frac{1}{\lambda^2}(f(x_0 + \lambda x) - f(x_0))$$

при $x \in X$, $\|x\| = 1$, т.е.

$$\inf_{0 < \lambda \leq \delta} \frac{1}{\lambda^2}(f(x_0 + \lambda x) - f(x_0)) \geq \frac{1}{2}(\alpha - \varepsilon)$$

при $x \in X$, $\|x\| = 1$. Поэтому

$$f(x_0 + \lambda x) - f(x_0) \geq \frac{1}{2}(\alpha - \varepsilon)\|\lambda x\|^2$$

при $x \in X$, $\|x\| = 1$, $0 < \lambda \leq \delta$. Отсюда следует

$$f(x_0 + y) - f(x_0)) \geq \frac{1}{2}(\alpha - \varepsilon)\|y\|^2$$

при $y \in X$, $\|y\| \leq \delta$, т.е. $f(x_0 + y) > f(x_0)$ при $y \in X$, $\|y\| \leq \delta$, $y \neq 0$. Теорема доказана.

Пусть $\varphi : X \to R$, где $\varphi(0) = 0$. Положим (см.[20])

$$f^{(1)-}(x_0; x) = \lim_{t \downarrow 0} \frac{f(x_0 + tx) - f(x_0)}{t},$$

$$f_\varphi^{(2)-}(x_0; x) = \lim_{t \downarrow 0} \frac{f(x_0 + 2tx) - \varphi(2tx) - 2f(x_0 + tx) + 2\varphi(tx) + f(x_0)}{t^2},$$

$$T_s(x_0; C) = \{x \in X : \exists x_k \to x \text{ в } (X)_s, \exists t_k \downarrow 0, \text{что } x_0 + t_k x_k \in C\},$$

$$M_1 = \{x \in T_s(x_0; C) : f^{(1)-}(x_0; x) > 0\}, \quad M_2 = \{x \in T_s(x_0; C) : f^{(1)-}(x_0; x) \leq 0\}.$$

**Теорема 10.** Пусть существует конечная положительно однородная степени $2 - \nu$ функция $\psi : X \to R$, $(0 < \nu \leq 2)$ и $o(t) : R_+ \to R_+$, где $\frac{o(t)}{t} \to 0$ при $t \downarrow 0$, числа $\delta > 0$, $L_1 > 0$ и $L_2 > 0$ такие, что

$$|f(x_0 + x + y) - f(x_0 + x)| \leq L_1\|y\| + o(\|x\|)$$

при $x \in M_1, \|x\| \leq \delta, y \in X, \|y\| \leq \delta, x_0 + x + y \in C$,

$$|f(x_0 + x + y) - f(x_0 + x) - \varphi(x + y) + \varphi(x)| \leq L_2\|y\|^\nu (\psi(x) + \|y\|^{2-\nu}) + o(\|x\|^2)$$

при $x \in M_2, \|x\| \leq \delta, y \in X, \|y\| \leq \delta, x_0 + x + y \in C$.

Кроме того, $f_\varphi^{(2)-}(x_0; x) > 0$ при $x \in M_2$, $\|x\| = 1$ и существует число $\alpha > 0$ такое, что

$\varphi(x) \geq 0$ при $x \in (C - x_0) \bigcap \alpha B$. Тогда $x_0 \in C$ является точкой строгого слабого локального минимума функции $f$ на $C$.

**Доказательство.** Допустим противное. Тогда существует индекс $\tau_0 \in A$ такой, что $x_0$ не является точкой строгого локального минимума функции $f$ на множестве $C_{\tau_0}$. Поэтому для любого $\delta_k > 0$ найдется $y_k \neq 0$ такое, что $x_0 + y_k \in C_{\tau_0}$, $\|y_k\| \leq \delta_k$ и $f(x_0 + y_k) \leq f(x_0)$. Положим $t_k = \|y_k\|$, $x_k = \dfrac{y_k}{\|y_k\|}$. Не умаляя общности можно считать, что $x_k \to x$. Если $t_k \downarrow 0$, то имеем, что $x \in T_s(x_0; C)$. Поэтому $x \in M_1 \cup M_2$. Также имеем, что
$$f(x_0 + t_k x) - f(x_0) \leq f(x_0 + t_k x) - f(x_0 + t_k x_k).$$
Предположим, что $\delta_k < \dfrac{\delta}{2}$. Если $x \in M_1$, то по условию имеем, что $f(x_0 + t_k x) - f(x_0) \leq L_1 t_k \|x_k - x\| + o(t_k \|x\|)$. Отсюда вытекает, что
$$f^{(1)-}(x_0; x) \leq \varliminf_{k \to \infty} \frac{f(x_0 + t_k x) - f(x_0)}{t_k} \leq \lim_{k \to \infty}(L_1 \cdot \|x_k - x\| + \frac{o(t_k \|x\|)}{t_k}) = 0.$$
Получим, что $f^{(1)-}(x_0; x) \leq 0$ и $x \in M_1$. Получим противоречие. Если $x \in M_2$, то по условию имеем, что
$$f(x_0 + t_k x) - \varphi(t_k x) - f(x_0) \leq f(x_0 + t_k x) - f(x_0 + t_k x_k) + \varphi(t_k x_k) - \varphi(t_k x) \leq$$
$$\leq |f(x_0 + t_k x + (t_k x_k - t_k x)) - f(x_0 + t_k x) - \varphi(t_k x_k) + \varphi(t_k x)| \leq$$
$$\leq L_2 t_k^\nu \|x_k - x\|^\nu (\psi(t_k x) + \|x_k - x\|^{2-\nu} \cdot t_k^{2-\nu}) + o(t_k^2 \|x\|^2)$$
при $\|y_k\| = \|t_k x_k\| \leq \min\{\alpha; 0{,}5\delta\}$. Поэтому
$$f_\varphi^{(2)-}(x_0; x) = \varliminf_{k \to \infty} \frac{f(x_0 + 2t_k x) - \varphi(2t_k x) - 2f(x_0 + t_k x) + 2\varphi(t_k x) + f(x_0)}{t_k^2} \leq$$
$$\leq \varliminf_{k \to \infty} \frac{f(x_0 + 2t_k x) - \varphi(2t_k x) - f(x_0)}{t_k^2} + \varliminf_{k \to \infty} \frac{-2f(x_0 + t_k x) + 2\varphi(t_k x) + 2f(x_0)}{t_k^2} \leq$$
$$\leq 4\varliminf_{k \to \infty} \frac{f(x_0 + t_k x) - \varphi(t_k x) - f(x_0)}{t_k^2} - 2\varliminf_{k \to \infty} \frac{f(x_0 + t_k x) - \varphi(t_k x) - f(x_0)}{t_k^2} \leq$$
$$\leq 2\varliminf_{k \to \infty} \frac{f(x_0 + t_k x) - \varphi(t_k x) - f(x_0)}{t_k^2} \leq$$
$$\leq 2\lim_{k \to \infty}(L_2 \|x_k - x\|^\nu (\psi(x) + \|x_k - x\|^{2-\nu}) + \frac{o(t_k^2 \|x\|^2)}{t_k^2}) = 0.$$
Отсюда следует, что $f_\varphi^{(2)-}(x_0; x) \leq 0$ и $x \in M_2$. Получим противоречие. Теорема доказана.

# Гл.2. ЭКСТРЕМАЛЬНАЯ ЗАДАЧА ДЛЯ ОПЕРАТОРНОГО ВКЛЮЧЕНИЯ

## 2.1. Непрерывная зависимость решений операторного включения

Пусть $R^n$ $n$-мерное евклидово пространство. Совокупность всех непустых компактных (выпуклых компактных) подмножеств $R^n$ обозначим через $\text{comp} R^n (\text{conv} R^n)$.

В дальнейшем равенства и включения, связанные с измеримыми функциями или отображениями, понимаются как почти всюду.

Если $A, B \in \text{comp} R^n$, то положим $\rho_x(A,B) = \max\left\{\sup_{y \in B} d(y,A), \sup_{x \in A} d(x,B)\right\}$, где через $\rho_x(A,B)$ обозначена метрика Хаусдорфа, $d(y,A) = \inf\{|x-y| : x \in A\}$, $|x-y| = \left(\sum_{i=1}^{n}(x_i - y_i)^2\right)^{\frac{1}{2}}$.

Пусть $t_0 < T$, $F:[t_0,T] \times R^n \to \text{comp} R^n$ многозначное отображение, $A$ линейный непрерывный оператор из $L_1^n[t_0,T]$ в $C^n[t_0,T]$.

Рассмотрим задачу для включения
$$u(t) \in F(t,(Au)(t)), \quad t \in [t_0,T], \quad u(\cdot) \in L_1^n[t_0,T]. \tag{1}$$

Функцию $u(\cdot) \in L_1^n[t_0,T]$, удовлетворяющую (1) назовем решением задачи (1), где $L_1^n[t_0,T] = \{u(t) = (u_1(t),\ldots,u_n(t)) \in R^n : u_i(\cdot) \in L_1[t_0,T], i = \overline{1,n}\}$ пространство с нормой
$$\|u(\cdot)\|_{L_1^n} = \int_{t_0}^{T} |u(s)| ds.$$

Пусть $A:L_1^n[t_0,T] \to C^n[t_0,T]$ линейный непрерывный оператор. Легко проверяется, что множество $\{(u,Au) \in L_1^n[t_0,T] \times C^n[t_0,T] : u \in L_1^n[t_0,T]\}$ относительно нормы $\|(u,Au)\| = \|u\|_{L_1^n[t_0,T]} + \|Au\|_{C^n[t_0,T]}$ является подпространством в $L_1^n[t_0,T] \times C^n[t_0,T]$. Так как $\|u\|_{L_1^n[t_0,T]} \leq \|u\|_{L_1^n[t_0,T]} + \|Au\|_{C^n[t_0,T]} \leq (1 + \|A\|)\|u\|_{L_1^n[t_0,T]}$, то имеем, что в пространстве $\{(u,Au) \in L_1^n[t_0,T] \times C^n[t_0,T] : u \in L_1^n[t_0,T]\}$ нормы $\|(u,Au)\|$ и $\|u\|_{L_1^n[t_0,T]}$ эквивалентны.

Если оператор $A:L_1^n[t_0,T] \to C^n[t_0,T]$ удовлетворяет условию $|(Au)(t) - (A\upsilon)(t)| \leq L\int_{t_0}^{t}|u(s) - \upsilon(s)|ds$ при $u(\cdot), \upsilon(\cdot) \in L_1^n[t_0,T]$, то оператор $A$ назовем условно типа Вольтерра, где $L > 0$. Если $A:L_1^n[t_0,T] \to C^n[t_0,T]$ линейный оператор, то условие $|(Au)(t) - (A\upsilon)(t)| \leq L\int_{t_0}^{t}|u(s) - \upsilon(s)|ds$ при $u(\cdot), \upsilon(\cdot) \in L_1^n[t_0,T]$ эквивалентно условию $|(Au)(t)| \leq L\int_{t_0}^{t}|u(s)|ds$ при $u(\cdot) \in L_1^n[t_0,T]$.

**Теорема 1.** Пусть $A:L_1^n[t_0,T] \to C^n[t_0,T]$ линейный оператор, $|(Au)(t)| \leq L\int_{t_0}^{t}|u(s)|ds$ при $u(\cdot) \in L_1^n[t_0,T]$, где $L > 0$, $\overline{u}(\cdot) \in L_1^n[t_0,T]$ и $\overline{x}(t) = (A\overline{u})(t)$, $F:[t_0,T] \times R^n \to \text{comp} R^n$ многозначное отображение, $t \to F(t,x)$ измеримо по $t$, существует суммируемая функция $M(t) > 0$ такая, что
$$\rho_x(F(t,x),F(t,x_1)) \leq M(t)|x - x_1| \tag{2}$$
при $x, x_1 \in R^n$. Кроме того, пусть $\rho(\cdot) \in L_1[t_0,T]$, $\overline{u}(\cdot)$ и $\overline{x}(\cdot)$ такие, что

$$d(\overline{u}(t), F(t, \overline{x}(t))) \le \rho(t), \quad t \in [t_0, T].$$

Тогда существует такое решение $u(\cdot) \in L_i^n[t_0, T]$ задачи (1), что

$$|x(t) - \overline{x}(t)| \le L \int_{t_0}^{t} e^{m(t)-m(\tau)} \rho(\tau) d\tau, \qquad |u(t) - \overline{u}(t)| \le \rho(t) + LM(t) \int_{t_0}^{t} e^{m(t)-m(\tau)} \rho(\tau) d\tau$$

при $t \in [t_0, T]$, где $x(t) = (Au)(t) \in C^n[t_0, T]$, $m(t) = L \int_{t_0}^{t} M(s) ds$.

**Доказательство.** Построим последовательность $x_i(t)$ ($i = 0,1,2...$) с помощью рекурентного соотношения при $t \in [t_0, T]$

$$x_0(t) = \overline{x}(t), \qquad x_{i+1}(t) = A(\upsilon_i)(t), \quad i = 0,1,2,... \qquad (3)$$

где $\upsilon_i(s) \in F(s, x_i(s))$ при $i \ge 0$, $\upsilon_i(s)$ измеримы и
$$|\upsilon_0(s) - \overline{u}(s)| = d(\overline{u}(s), F(s, x_0(s))), \qquad |\upsilon_i(s) - \upsilon_{i-1}(s)| = d(\upsilon_{i-1}(s), F(s, x_i(s))), \quad i = 1,2,...$$

при $s \in [t_0, T]$. По лемме 2.1.4 [25] такая функция $\upsilon_i(t)$ существует. По условию (2) при $i \ge 1$ имеем
$$|\upsilon_{i+1}(s) - \upsilon_i(s)| = d(\upsilon_i(s), F(s, x_{i+1}(s))) \le \rho_x(F(s, x_i(s)), F(s, x_{i+1}(s))) \le M(s)|x_{i+1}(s) - x_i(s)|$$

при $s \in [t_0, T]$). Ясно, что

$$|x_{i+1}(t) - x_i(t)| = |A(\upsilon_i)(t) - A(\upsilon_{i-1})(t)| \le L \int_{t_0}^{t} |(\upsilon_i(s) - \upsilon_{i-1}(s))| ds.$$

По условию, имеем, что $|\upsilon_0(s) - \overline{u}(s)| \le \rho(s)$ при $s \in [t_0, T]$. Так как $x_1(t) = (A\upsilon_0)(t)$ и $\overline{x}(t) = (A\overline{u})(t)$, то имеем, что

$$|x_1(t) - x_0(t)| \le L \int_{t_0}^{t} |\upsilon_0(s) - \overline{u}(s)| ds \le L \int_{t_0}^{t} \rho(s) ds.$$

Тогда получим

$$|\upsilon_1(s) - \upsilon_0(s)| = d(\upsilon_0(s), F(s, x_1(s))) \le \rho_x(F(s, x_0(s)), F(s, x_1(s))) \le M(s)|x_0(s) - x_1(s)| \le$$
$$\le LM(s) \int_{t_0}^{s} |\upsilon_0(\tau) - \overline{u}(\tau)| d\tau \le LM(s) \int_{t_0}^{s} \rho(\tau) d\tau$$

при $s \in [t_0, T]$. Ясно, что

$$|x_2(t) - x_1(t)| = |A(\upsilon_1 - \upsilon_0)(t)| \le L \int_{t_0}^{t} |\upsilon_1(s) - \upsilon_0(s)| ds \le$$
$$\le L^2 \int_{t_0}^{t} M(s) \int_{t_0}^{s} \rho(v) dv ds = L^2 \int_{t_0}^{t} \int_{t_0}^{s} \rho(v) dv d \int_{t_0}^{s} M(\tau) d\tau = L^2 \int_{t_0}^{t} M(v) dv \int_{t_0}^{t} \rho(v) dv -$$
$$- L^2 \int_{t_0}^{t} \rho(s) \int_{t_0}^{s} M(v) dv ds = Lm(t) \int_{t_0}^{t} \rho(v) dv - L \int_{t_0}^{t} \rho(s) m(s) ds = L \int_{t_0}^{t} \rho(s)(m(t) - m(s)) ds$$

при $t \in [t_0, T]$. Поэтому

$$|\upsilon_2(s) - \upsilon_1(s)| \le M(s)|x_2(s) - x_1(s)| \le LM(s) \int_{t_0}^{s} \rho(\tau)(m(s) - m(\tau)) d\tau = \frac{d}{ds} \int_{t_0}^{s} \frac{1}{2} \rho(\tau)(m(s) - m(\tau))^2 d\tau$$

при $s \in [t_0, T]$. Также получим, что

$$|x_3(t) - x_2(t)| = |A(\upsilon_2 - \upsilon_1)(t)| \le L \int_{t_0}^{t} |\upsilon_2(s) - \upsilon_1(s)| ds \le$$

$$\leq L\int_{t_0}^{t} LM(s)\int_{t_0}^{s}\rho(\tau)(m(s)-m(\tau))d\tau ds = L\int_{t_0}^{t} d\int_{t_0}^{s}\frac{1}{2}\rho(\tau)(m(s)-m(\tau))^2 d\tau =$$

$$= L\int_{t_0}^{t}\frac{1}{2}\rho(\tau)(m(t)-m(\tau))^2 d\tau$$

при $t \in [t_0, T]$. Тогда

$$|\upsilon_3(s)-\upsilon_2(s)| \leq M(s)|x_3(s)-x_2(s)| \leq LM(s)\int_{t_0}^{s}\frac{1}{2}\rho(\tau)(m(s)-m(\tau))^2 d\tau =$$

$$= \frac{d}{ds}\int_{t_0}^{s}\frac{1}{3!}\rho(\tau)(m(s)-m(\tau))^3 d\tau$$

при $s \in [t_0, T]$. Поэтому

$$|x_4(t)-x_3(t)| = |A(\upsilon_3 - \upsilon_2)(t)| \leq L\int_{t_0}^{t}|\upsilon_3(s)-\upsilon_2(s)|ds \leq$$

$$\leq L\int_{t_0}^{t} LM(s)\int_{t_0}^{s}\frac{1}{2}\rho(\tau)(m(s)-m(\tau))^2 d\tau ds = L\int_{t_0}^{t} d\int_{t_0}^{s}\frac{1}{3!}\rho(\tau)(m(s)-m(\tau))^3 d\tau =$$

$$= L\int_{t_0}^{t}\frac{1}{3!}\rho(\tau)(m(t)-m(\tau))^3 d\tau$$

при $t \in [t_0, T]$. Тогда

$$|\upsilon_4(s)-\upsilon_3(s)| \leq M(s)|x_4(s)-x_3(s)| \leq LM(s)\int_{t_0}^{s}\frac{1}{3!}\rho(\tau)(m(s)-m(\tau))^3 d\tau =$$

$$= \frac{d}{ds}\int_{t_0}^{s}\frac{1}{4!}\rho(\tau)(m(s)-m(\tau))^4 d\tau$$

при $s \in [t_0, T]$. Поэтому

$$|x_5(t)-x_4(t)| \leq |A(\upsilon_4 - \upsilon_3)(t)| \leq L\int_{t_0}^{t}|\upsilon_4(s)-\upsilon_3(s)|ds \leq$$

$$\leq L\int_{t_0}^{t} LM(s)\int_{t_0}^{s}\frac{1}{3!}\rho(\tau)(m(s)-m(\tau))^3 d\tau ds = L\int_{t_0}^{t} d\int_{t_0}^{s}\frac{1}{4!}\rho(\tau)(m(s)-m(\tau))^4 d\tau =$$

$$= L\int_{t_0}^{t}\frac{1}{4!}\rho(\tau)(m(t)-m(\tau))^4 d\tau$$

при $t \in [t_0, T]$. Продолжая процесс получим

$$|x_{m+1}(t)-x_m(t)| \leq L\int_{t_0}^{t}\frac{1}{m!}\rho(\tau)(m(t)-m(\tau))^m d\tau \leq L\frac{1}{m!}(m(T))^m \int_{t_0}^{T}\rho(\tau)d\tau, \quad (4)$$

$$|\upsilon_{m+1}(s)-\upsilon_m(s)| \leq LM(s)\int_{t_0}^{s}\frac{1}{m!}\rho(\tau)(m(s)-m(\tau))^m d\tau \quad (5)$$

при $t \in [t_0, T]$. Пользуясь неравенством $1+\frac{z}{1!}+\frac{z^2}{2!}+\cdots+\frac{z^m}{m!} \leq e^z$, где $z \geq 0$, получим

$$|x_{m+1}(t)-\bar{x}(t)| \leq |x_1(t)-\bar{x}(t)| + |x_2(t)-x_1(t)| + |x_3(t)-x_2(t)| + \ldots + |x_{m+1}(t)-x_m(t)| \leq$$

$$\leq L\int_{t_0}^{t}\rho(\tau)d\tau + L\int_{t_0}^{t}\frac{1}{1!}\rho(\tau)(m(t)-m(\tau))d\tau + L\int_{t_0}^{t}\frac{1}{2!}\rho(\tau)(m(t)-m(\tau))^2 d\tau + \quad (6)$$

$$+ L\int_{t_0}^{t} \frac{1}{3!}\rho(\tau)(m(t)-m(\tau))^3 d\tau + \cdots + L\int_{t_0}^{t} \frac{1}{m!}\rho(\tau)(m(t)-m(\tau))^m d\tau =$$

$$= L\int_{t_0}^{t} (1 + \frac{1}{1!}(m(t)-m(\tau)) + \frac{1}{2!}(m(t)-m(\tau))^2 + \frac{1}{3!}(m(t)-m(\tau))^3 + \cdots +$$

$$+ \frac{1}{m!}(m(t)-m(\tau))^m)\rho(\tau)d\tau \le L\int_{t_0}^{t} e^{m(t)-m(\tau)}\rho(\tau)d\tau;$$

и

$$|\upsilon_{m+1}(t) - \overline{u}(t)| \le |\upsilon_0(t) - \overline{u}(t)| + |\upsilon_1(t) - \upsilon_0(t)| + |\upsilon_2(t) - \upsilon_1(t)| + |\upsilon_3(t) - \upsilon_2(t)| + \ldots +$$

$$+ |\upsilon_{m+1}(t) - \upsilon_m(t)| \le \rho(t) + LM(t)\int_{t_0}^{t}\rho(s)ds + LM(t)\int_{t_0}^{t}\frac{1}{1!}\rho(\tau)(m(t)-m(\tau))d\tau \ldots +$$

$$+ LM(t)\int_{t_0}^{t}\frac{1}{m!}\rho(\tau)(m(t)-m(\tau))^m d\tau \le \rho(t) + LM(t)\int_{t_0}^{t}(1 + \frac{1}{1!}(m(t)-m(\tau)) +$$

$$+ \frac{1}{2!}(m(t)-m(\tau))^2 + \cdots + \frac{1}{m!}(m(t)-m(\tau))^m)\rho(\tau)ds \le \rho(t) + LM(t)\int_{t_0}^{t}e^{m(t)-m(\tau)}\rho(\tau)d\tau \quad (7)$$

при $t \in [t_0, T]$. Из оценки (4)-(7) вытекает, что последовательности $x_m(t)$ и $\upsilon_m(t)$ сходятся соответственно к функциям $x(\cdot)$ и $u(\cdot) \in L_1^n[t_0, T]$. Из теоремы 1.2.19[4] и из (3) вытекает, что $x(t)$ непрерывно. Так как $\upsilon_i(t) \in F(t, x_i(t))$, то из теоремы 1.2.23 и 1.2.28 [3] имеем, что $u(t) \in F(t, x(t))$ при $t \in [t_0, T]$. Кроме того, из (6) и (7) получим, что

$$|x(t) - \overline{x}(t)| \le L\int_{t_0}^{t} e^{m(t)-m(\tau)}\rho(\tau)d\tau, \quad |u(t) - \overline{u}(t)| \le \rho(t) + LM(t)\int_{t_0}^{t} e^{m(t)-m(\tau)}\rho(\tau)d\tau.$$

Из (3) вытекает, что $x(t) = (Au)(t)$. Теорема доказана.

**Следствие 1.** Пусть $A: L_1^n[t_0, T] \to C^n[t_0, T]$ линейный оператор, $|(Au)(t)| \le L\int_{t_0}^{t}|u(s)|ds$ при $u(\cdot) \in L_1^n[t_0, T]$, $\overline{u}(\cdot) \in L_1^n[t_0, T]$ и $\overline{x}(t) = (A\overline{u})(t)$, $F:[t_0, T] \times R^n \to \text{comp}R^n$ многозначное отображение, $t \to F(t, x)$ измеримо по $t$, существует $M > 0$ такое, что $\rho_x(F(t, x), F(t, x_1)) \le M|x - x_1|$ при $x, x_1 \in R^n$. Кроме того, пусть $\rho(\cdot) \in L_1[t_0, T]$, $\overline{u}(\cdot)$ и $\overline{x}(\cdot)$ такие, что $d(\overline{u}(t), F(t, \overline{x}(t))) \le \rho(t)$ при $t \in [t_0, T]$.

Тогда существует такое решение $u(\cdot) \in L_1^n[t_0, T]$ задачи (1), что

$$|x(t) - \overline{x}(t)| \le L\int_{t_0}^{t} e^{LM(t-\tau)}\rho(\tau)d\tau, \quad |u(t) - \overline{u}(t)| \le \rho(t) + LM\int_{t_0}^{t} e^{LM(t-\tau)}\rho(\tau)d\tau,$$

где $x(t) = (Au)(t)$, $x(\cdot) \in C^n[t_0, T]$.

**Доказательство.** Построим последовательность $x_i(t)$ ($i = 0,1,2\ldots$) с помощью рекурентного соотношения при $t \in [t_0, T]$

$$x_0(t) = \overline{x}(t), \quad x_{i+1}(t) = (A\upsilon_i)(t), \quad i = 0,1,2,\ldots \quad (8)$$

где $\upsilon_i(s) \in F(s, x_i(s))$ при $i \ge 0$, $\upsilon_i(s)$ измеримы и

$$|\upsilon_0(s) - \overline{u}(s)| = d(\overline{u}(s), F(s, x_0(s))), \quad |\upsilon_i(s) - \upsilon_{i-1}(s)| = d(\upsilon_{i-1}(s), F(s, x_i(s))), \; i = 1,2,\ldots$$

при $s \in [t_0, T]$. По лемме 2.1.4[25] такая функция $\upsilon_i(t)$ существует. Если $i \ge 0$, то по условию получим

$$|\upsilon_{i+1}(s) - \upsilon_i(s)| = d(\upsilon_i(s), F(s, x_{i+1}(s))) \le \rho_x(F(s, x_i(s)), F(s, x_{i+1}(s))) \le M|x_{i+1}(s) - x_i(s)|$$

при $s \in [t_0, T]$. Из (8) вытекает, что

$$|x_{i+1}(t) - x_i(t)| = |A(\upsilon_i - \upsilon_{i-1})(t)|$$

при $t \in [t_0, T]$. По условию, имеем, что $|\upsilon_0(s) - \bar{u}(s)| \leq \rho(s)$ при $s \in [t_0, T]$. Так как $x_1(t) = (A\upsilon_0)(t)$ и $\bar{x}(t) = (A\bar{u})(t)$, то имеем, что

$$|x_1(t) - x_0(t)| = |(A\upsilon_0)(t) - (A\bar{u})(t)| \leq L\int_{t_0}^{t}|\upsilon_0(s) - \bar{u}(s)|ds \leq L\int_{t_0}^{t}\rho(s)ds.$$

Тогда получим
$$|\upsilon_1(s) - \upsilon_0(s)| = d(\upsilon_0(s), F(s, x_1(s))) \leq \rho_x(F(s, x_0(s)), F(s, x_1(s))) \leq$$
$$\leq M|x_0(t) - x_1(t)| \leq LM\int_{t_0}^{t}\rho(s)ds$$

при $s \in [t_0, T]$.

$$|x_2(t) - x_1(t)| = |(A\upsilon_1)(t) - (A\upsilon_0)(t)| \leq L\int_{t_0}^{t}|\upsilon_1(s) - \upsilon_0(s)|ds \leq L^2M\int_{t_0}^{t}\int_{t_0}^{s}\rho(v)dvds$$

при $t \in [t_0, T]$. Поэтому

$$|\upsilon_2(s) - \upsilon_1(s)| \leq M|x_2(s) - x_1(s)| \leq L^2M^2\int_{t_0}^{s}\int_{t_0}^{s_1}\rho(v)dvds_1$$

при $s \in [t_0, T]$. Также получим, что

$$|x_3(t) - x_2(t)| = |(A\upsilon_2)(t) - (A\upsilon_1)(t)| \leq L\int_{t_0}^{t}|\upsilon_2(s) - \upsilon_1(s)|ds \leq L^3M^2\int_{t_0}^{t}\int_{t_0}^{t_2}\int_{t_0}^{t_1}\rho(v)dvdt_1dt_2$$

при $t \in [t_0, T]$. Тогда

$$|\upsilon_3(s) - \upsilon_2(s)| \leq M|x_3(s) - x_2(s)| \leq L^3M^3\int_{t_0}^{s}\int_{t_0}^{s_2}\int_{t_0}^{s_1}\rho(v)dvds_1ds_2$$

при $s \in [t_0, T]$. Поэтому

$$|x_4(t) - x_3(t)| = |A(\upsilon_3 - \upsilon_2)(t)| \leq L\int_{t_0}^{t}|\upsilon_3(s) - \upsilon_2(s)|ds \leq L^4M^3\int_{t_0}^{t}\int_{t_0}^{t_3}\int_{t_0}^{t_2}\int_{t_0}^{t_1}\rho(v)dvdt_1dt_2dt_3$$

при $t \in [t_0, T]$. Тогда

$$|\upsilon_4(s) - \upsilon_3(s)| \leq M|x_4(s) - x_3(s)| \leq L^4M^4\int_{t_0}^{s}\int_{t_0}^{s_3}\int_{t_0}^{s_2}\int_{t_0}^{s_1}\rho(v)dvds_1ds_2ds_3$$

при $s \in [t_0, T]$. Поэтому

$$|x_5(t) - x_4(t)| = |A(\upsilon_4 - \upsilon_3)(t))| \leq L\int_{t_0}^{t}|\upsilon_4(s) - \upsilon_3(s)|ds \leq L^5M^4\int_{t_0}^{t}\int_{t_0}^{t_4}\int_{t_0}^{t_3}\int_{t_0}^{t_2}\int_{t_0}^{t_1}\rho(v)dvdt_1dt_2dt_3dt_4$$

при $t \in [t_0, T]$. Тогда получим

$$|\upsilon_5(s) - \upsilon_4(s)| \leq M|x_5(s) - x_4(s)| \leq L^5M^5\int_{t_0}^{s}\int_{t_0}^{s_4}\int_{t_0}^{s_3}\int_{t_0}^{s_2}\int_{t_0}^{s_1}\rho(v)dvds_1ds_2ds_3ds_4$$

при $s \in [t_0, T]$. Продолжая процесс получим

$$|x_{m+1}(t) - x_m(t)| \leq L^{m+1}M^m\int_{t_0}^{t}\int_{t_0}^{t_m}\int_{t_0}^{t_{m-1}}\ldots\int_{t_0}^{t_1}\rho(v)dvdt_1dt_2\ldots dt_m \leq \tfrac{1}{m!}L^{m+1}M^m\int_{t_0}^{t}(t-\tau)^m\rho(\tau)d\tau, \quad (9)$$

$$|\upsilon_{m+1}(s) - \upsilon_m(s)| \leq L^{m+1}M^{m+1}\int_{t_0}^{s}\int_{t_0}^{s_m}\int_{t_0}^{s_{m-1}}\ldots\int_{t_0}^{s_1}\rho(v)dvds_1ds_2\ldots ds_m \leq \tfrac{1}{m!}L^{m+1}M^{m+1}\int_{t_0}^{s}(s-v)^m\rho(v)dv. \quad (10)$$

Таким образом получим, что

$$|x_{m+1}(t)-\overline{x}(t)|\le |x_1(t)-\overline{x}(t)|+|x_2(t)-x_1(t)|+|x_3(t)-x_2(t)|+\ldots+|x_{m+1}(t)-x_m(t)|\le L\int_{t_0}^{t}\rho(\tau)d\tau+$$

$$+L^2M\int_{t_0}^{t}(t-\tau)\rho(\tau)d\tau+\tfrac{1}{2!}L^3M^2\int_{t_0}^{t}(t-\tau)^2\rho(\tau)d\tau+\ldots+\tfrac{1}{m!}L^{m+1}M^m\int_{t_0}^{t}(t-\tau)^m\rho(\tau)d\tau= \qquad (11)$$

$$=L\int_{t_0}^{t}(1+aM(t-\tau)+\tfrac{1}{2!}L^2M^2(t-\tau)^2+\ldots+\tfrac{1}{m!}L^mM^m(t-\tau)^m)\rho(\tau)d\tau\le L\int_{t_0}^{t}e^{\|A\|M(t-\tau)}\rho(\tau)d\tau;$$

и

$$|\upsilon_{m+1}(t)-\overline{u}(t)|\le |\upsilon_0(t)-\overline{u}(t)|+|\upsilon_1(t)-\upsilon_0(t)|+|\upsilon_2(t)-\upsilon_1(t)|+\ldots+|\upsilon_{m+1}(t)-\upsilon_m(t)|=$$

$$\le \rho(t)+LM\int_{t_0}^{t}\rho(\tau)d\tau+L^2M^2\int_{t_0}^{t}(t-\tau)\rho(\tau)d\tau+\tfrac{1}{2!}L^3M^3\int_{t_0}^{t}(t-\tau)^2\rho(\tau)d\tau+\ldots+ \qquad (12)$$

$$+\tfrac{1}{m!}L^{m+1}M^{m+1}\int_{t_0}^{t}(t-\tau)^m\rho(\tau)d\tau=\rho(t)+LM\int_{t_0}^{t}(1+LM(t-\tau)+\tfrac{1}{2!}L^2M^2(t-\tau)^2+\ldots+$$

$$+\tfrac{1}{m!}L^mM^m(t-\tau)^m)\rho(\tau)d\tau\le \rho(t)+LM\int_{t_0}^{t}e^{LM(t-\tau)}\rho(\tau)d\tau$$

при $t\in[t_0,T]$. Из оценки (9)-(12) вытекает, что последовательности $x_m(t)$ и $\upsilon_m(t)$ сходятся соответственно к функциям $x(\cdot)$ и $u(\cdot)\in L_1^n[t_0,T]$. Из теоремы 1.2.19[4] и из (9) вытекает, что $x(t)$ непрерывно. Так как $\upsilon_i(t)\in F(t,x_i(t))$, то из теоремы 1.2.23 и 1.2.28 [3] имеем, что $u(t)\in F(t,x(t))$ при $t\in[t_0,T]$. Кроме того, из (11) и (12) получим, что

$$|x(t)-\overline{x}(t)|\le L\int_{t_0}^{t}e^{LM(t-\tau)}\rho(\tau)d\tau, \qquad |u(t)-\overline{u}(t)|\le \rho(t)+LM\int_{t_0}^{t}e^{LM(t-\tau)}\rho(\tau)d\tau.$$

Из (8) вытекает, что $x(t)=(Au)(t)$. Следствие доказано.

Положив $\overline{u}(t)=0$ и $\overline{x}(t)=(A\overline{u})(t)=0$, из теоремы 1 имеем, что верно следующее следствие.

**Следствие 2.** Пусть $A:L_1^n[t_0,T]\to C^n[t_0,T]$ линейный оператор, $|(Au)(t)|\le L\int_{t_0}^{t}|u(s)|ds$ при $u(\cdot)\in L_1^n[t_0,T]$, $F:[t_0,T]\times R^n\to \text{comp}\,R^n$ многозначное отображение, $t\to F(t,x)$ измеримо по $t$, существует суммируемая функция $M(t)>0$ такая, что $\rho_x(F(t,x),F(t,x_1))\le M(t)|x-x_1|$ при $x,x_1\in R^n$ и $\|F(t,0)\|$ суммируема. Тогда существует решение $u(\cdot)\in L_1^n[t_0,T]$ задачи (1).

Аналогично теореме 1, доказываются следующие теоремы 2-4.

**Теорема 2.** Пусть $A:L_1^n[t_0,T]\to C^n[t_0,T]$ оператор удовлетворяет условию $|(Au)(t)-(A\upsilon)(t)|\le L\int_{t_0}^{t}|u(s)-\upsilon(s)|ds$ при $u(\cdot),\upsilon(\cdot)\in L_1^n[t_0,T]$, $\overline{u}(\cdot)\in L_1^n[t_0,T]$ и $\overline{x}(t)=(A\overline{u})(t)$, $F:[t_0,T]\times R^n\to \text{comp}\,R^n$ многозначное отображение, $t\to F(t,x)$ измеримо по $t$, существует суммируемая функция $M(t)>0$ такая, что

$$\rho_x(F(t,x),F(t,x_1))\le M(t)|x-x_1|$$

при $x,x_1\in R^n$. Кроме того, пусть $\rho(\cdot)\in L_1[t_0,T]$, $\overline{u}(\cdot)$ и $\overline{x}(\cdot)$ такие, что

$$d(\overline{u}(t),F(t,\overline{x}(t)))\le \rho(t), \quad t\in[t_0,T].$$

Тогда существует такое решение $(u(\cdot),x(\cdot))\in L_1^n[t_0,T]\times C^n[t_0,T]$ задачи (1), что

$$|x(t)-\overline{x}(t)|\leq L\int_{t_0}^{t}e^{m(t)-m(\tau)}\rho(\tau)d\tau, \qquad |u(t)-\overline{u}(t)|\leq \rho(t)+LM(t)\int_{t_0}^{t}e^{m(t)-m(\tau)}\rho(\tau)d\tau,$$

при $t\in[t_0,T]$, где $x(t)=(Au)(t)$, $m(t)=L\int_{t_0}^{t}M(s)ds$.

**Теорема 3.** Пусть $A:L_1^n[t_0,T]\to L_\infty^n[t_0,T]$ линейный оператор, $|(Au)(t)|\leq L\int_{t_0}^{t}|u(s)|ds$ при $u(\cdot)\in L_1^n[t_0,T]$, $\overline{u}(\cdot)\in L_1^n[t_0,T]$ и $\overline{x}(t)=(A\overline{u})(t)$, $F:[t_0,T]\times R^n\to \mathrm{comp}R^n$ многозначное отображение, $t\to F(t,x)$ измеримо по $t$, существует суммируемая функция $M(t)>0$ такая, что

$$\rho_x(F(t,x),F(t,x_1))\leq M(t)|x-x_1|$$

при $x,x_1\in R^n$. Кроме того, пусть $\rho(\cdot)\in L_1[t_0,T]$, $\overline{u}(\cdot)$ и $\overline{x}(\cdot)$ такие, что

$$d(\overline{u}(t), F(t,\overline{x}(t)))\leq \rho(t), \quad t\in[t_0,T].$$

Тогда существует такое решение $(u(\cdot),x(\cdot))\in L_1^n[t_0,T]\times L_\infty^n[t_0,T]$ задачи (1), что

$$|x(t)-\overline{x}(t)|\leq L\int_{t_0}^{t}e^{m(t)-m(\tau)}\rho(\tau)d\tau, \qquad |u(t)-\overline{u}(t)|\leq \rho(t)+LM(t)\int_{t_0}^{t}e^{m(t)-m(\tau)}\rho(\tau)d\tau,$$

при $t\in[t_0,T]$, где $x(t)=(Au)(t)$, $m(t)=L\int_{t_0}^{t}M(s)ds$.

**Теорема 4.** Пусть $A:L_1^n[t_0,T]\to L_\infty^n[t_0,T]$ оператор удовлетворяет условию $|(Au)(t)-(A\upsilon)(t)|\leq L\int_{t_0}^{t}|u(s)-\upsilon(s)|ds$ при $u(\cdot),\upsilon(\cdot)\in L_1^n[t_0,T]$, $\overline{u}(\cdot)\in L_1^n[t_0,T]$ и $\overline{x}(t)=(A\overline{u})(t)$, $F:[t_0,T]\times R^n\to \mathrm{comp}R^n$ многозначное отображение, $t\to F(t,x)$ измеримо по $t$, существует суммируемая функция $M(t)>0$ такая, что

$$\rho_x(F(t,x),F(t,x_1))\leq M(t)|x-x_1|$$

при $x,x_1\in R^n$. Кроме того, пусть $\rho(\cdot)\in L_1[t_0,T]$, $\overline{u}(\cdot)$ и $\overline{x}(\cdot)$ такие, что

$$d(\overline{u}(t), F(t,\overline{x}(t)))\leq \rho(t), \quad t\in[t_0,T].$$

Тогда существует такое решение $(u(\cdot),x(\cdot))\in L_1^n[t_0,T]\times L_\infty^n[t_0,T]$ задачи (1), что

$$|x(t)-\overline{x}(t)|\leq L\int_{t_0}^{t}e^{m(t)-m(\tau)}\rho(\tau)d\tau, \qquad |u(t)-\overline{u}(t)|\leq \rho(t)+LM(t)\int_{t_0}^{t}e^{m(t)-m(\tau)}\rho(\tau)d\tau,$$

при $t\in[t_0,T]$, где $x(t)=(Au)(t)$, $m(t)=L\int_{t_0}^{t}M(s)ds$.

Из теоремы 3 имеем, что верно следующее следствие.

**Следствие 3.** Пусть $A:L_1^n[t_0,T]\to L_\infty^n[t_0,T]$ линейный оператор, $|(Au)(t)|\leq L\int_{t_0}^{t}|u(s)|ds$ при $u(\cdot)\in L_1^n[t_0,T]$, $\overline{u}(\cdot)\in L_1^n[t_0,T]$ и $\overline{x}(t)=(A\overline{u})(t)$, $F:[t_0,T]\times R^n\to \mathrm{comp}R^n$ многозначное отображение, $t\to F(t,x)$ измеримо по $t$, существует $M>0$ такое, что $\rho_x(F(t,x),F(t,x_1))\leq M|x-x_1|$ при $x,x_1\in R^n$. Кроме того, пусть $\rho(\cdot)\in L_1[t_0,T]$, $\overline{u}(\cdot)$ и $\overline{x}(\cdot)$ такие, что $d(\overline{u}(t), F(t,\overline{x}(t)))\leq \rho(t)$ при $t\in[t_0,T]$.

Тогда существует такое решение $(u(\cdot),x(\cdot))\in L_1^n[t_0,T]\times L_\infty^n[t_0,T]$ задачи (1), что

$$|x(t)-\overline{x}(t)|\leq L\int_{t_0}^{t}e^{LM(t-s)}\rho(s)ds, \qquad |u(t)-\overline{u}(t)|\leq \rho(t)+ML\int_{t_0}^{t}e^{LM(t-s)}\rho(s)ds,$$

при $t\in[t_0,T]$, где $x(t)=(Au)(t)$.

**Замечание 1.** Из доказательства теоремы 1 следует, что в теореме 1-4 условие $F:[t_0,T]\times R^n\to compR^n$ можно заменить условием: $F:[t_0,T]\times(x_0(t)+\alpha B)\to compR^n$, где $\alpha\geq L\int_{t_0}^{t}e^{m(t)-m(\tau)}\rho(\tau)d\tau$ при $t\in[t_0,T]$, $B=\{z\in R^n:\|z\|\leq 1\}$ -единичный шар в $R^n$. Ясно, что из неравенства $\alpha\geq Le^{m(T)}\int_{t_0}^{T}\rho(\tau)d\tau$ вытекает, что $\alpha\geq L\int_{t_0}^{t}e^{m(t)-m(\tau)}\rho(\tau)d\tau$ при $t\in[t_0,T]$.

**Лемма 2.** Пусть $A:L_1^n[t_0,T]\to L_\infty^n[t_0,T]$ линейный оператор, $|(Au)(t)|\leq L\int_{t_0}^{t}|u(s)|ds$ при $u(\cdot)\in L_1^n[t_0,T]$, $F:[t_0,T]\times(x_0(t)+\alpha B)\to compR^n$, где $\alpha>0$, $u_0(\cdot)\in L_1^n[t_0,T]$, $x_0(t)=(Au_0)(t)$, функция $u_0(t)$ -решение задачи $u(t)\in F(t,(Au)(t))$ и $F(t,x)$ удовлетворяет условиям теоремы 3 в $x_0(t)+\alpha B$. Тогда существует такое $\delta>0$, что при $s(\cdot)\in L_1^n[t_0,T]$, $\|s(\cdot)\|_{L_1^n}\leq\delta$ найдется решение $u_s(t)$ задачи $u(t)\in F(t,(Au)(t))+s(t)$, что $|x_0(t)-x_s(t)|\leq\alpha$ при $t\in[t_0,T]$, где $x_s(t)=(Au_s)(t)$ и $\|x_s(\cdot)-x_0(\cdot)\|_{L_\infty^n}\to 0$ при $\|s(\cdot)\|_{L_1^n}\to 0$.

**Доказательство.** Ясно, что
$$\rho_x(F(t,x)+s(t),F(t,x')+s(t))\leq M(t)|x-x'|$$
при $|x-x_0(t)|\leq\alpha$, $|x'-x_0(t)|\leq\alpha$ и $d(u_0(t),F(t,x_0(t))+s(t))\leq|s(t)|$. По теореме 3 существует решение $u_s(t)$ задачи $u_s(t)\in F(t,x_s(t))+s(t)$, такое, что

$$|x_0(t)-x_s(t)|\leq L\int_{t_0}^{t}e^{m(t)-m(\tau)}s(\tau)d\tau\leq Le^{M(T)}\|s(\cdot)\|_{L_1^n}, \qquad (13)$$

Если определять $\delta$ из неравенства $Le^{M(T)}\delta\leq\alpha$, то получим, что верна первая часть утверждения леммы.

Из (13) получим, что $\|x_s(\cdot)-x_0(\cdot)\|_{L_\infty^n[t_0,T]}\leq Le^{M(T)}\|s(\cdot)\|_{L_1^n}$, т.е. $\|x_s(\cdot)-x_0(\cdot)\|_{L_\infty^n[t_0,T]}\to 0$ при $\|s(\cdot)\|_{L_1^n}\to 0$. Лемма доказана.

**Лемма 3.** Пусть $A:L_1^n[t_0,T]\to L_\infty^n[t_0,T]$ линейный оператор, $|(Au)(t)|\leq L\int_{t_0}^{t}|u(s)|ds$ при $u(\cdot)\in L_1^n[t_0,T]$, $F:[t_0,T]\times R^n\to compR^n\cup\{\varnothing\}$, отображение $t\to F(t,x)$ измеримо на $[t_0,T]$, отображение $x\to F(t,x)$ замкнуто и выпукло почти при всех $t\in[t_0,T]$, т.е. $grF_t=\{(x,y):y\in F(t,x)\}$ замкнуто и выпукло почти при всех $t\in[t_0,T]$. Пусть существует такая суммируемая функция $\lambda(t)$, что
$$\|F(t,x)\|\leq\lambda(t)(1+|x|), \qquad (14)$$
при $x\in R^n$, где $\|F(t,x)\|=\sup\{|y|:y\in F(t,x)\}$, $\|\varnothing\|=0$ и существует решение $u_0(t)$ задачи $u_0(t)\in F(t,x_0(t))$, что $x_0(t)=(Au_0)(t)$ принадлежит $dom F_t=\{x:F(t,x)\neq\varnothing\}$ вместе с некоторой $\varepsilon$ трубкой, т.е. $\{x:|x_0(t)-x|\leq\varepsilon\}\subset dom F_t$. Тогда существуют такие $\delta>0$ и решение $u_s(t)$ задачи $u(t)\in F(t,(Au)(t))+s(t)$, где $s(\cdot)\in L_1^n[t_0,T]$, $\|s(\cdot)\|_{L_1^n}\leq\delta$, что $x_s(t)=(Au_s)(t)$, $|x_0(t)-x_s(t)|\leq\varepsilon$ при $t\in[t_0,T]$ и $\|x_s(\cdot)-x_0(\cdot)\|_{L_\infty^n[t_0,T]}\to 0$ при $\|s(\cdot)\|_{L_1^n[t_0,T]}\to 0$.

**Доказательство.** Рассмотрим множество $S_t=\{x:|x_0(t)-x|\leq\dfrac{\varepsilon}{2}\}$ и обозначим $\delta=\dfrac{\varepsilon}{3}$. Пусть

$x \in S_t$, тогда $x + \delta B \subset \text{int dom } a_t$, где $B$ единичный шар в $R^n$ с центром в нуле. Так как множество $\text{gr } a_t$ выпукло, то аналогично доказательству теоремы 3.1.1 [13] можно показать, что

$$\rho_x(F(t,y), F(t,x)) \le \frac{2\lambda(t)(1 + |x| + \delta)}{\delta}|y - x| \qquad (15)$$

для $y \in x + \delta B$ п.в. $t$. Ясно, что $S_t + \delta B \subset \text{int dom } F_t$. Поэтому (15) верно для всех $x, y \in S_t$ и почти всех $t$, т.е.

$$\rho_x(F(t,y), F(t,x)) \le \frac{2\lambda(t)\left(1 + |x_0(t)| + \frac{\varepsilon}{2} + \delta\right)}{\delta}|y - x|$$

при $x, y \in S_t$. Поэтому доказательство леммы 3 вытекает из леммы 2. Лемма доказана.

**Замечание 2.** Из замкнутости $x \to F(t,x)$ почти всех $t$ и из (14) вытекает, что отображение $x \to F(t,x)$ полунепрерывно сверху почти при всех $t$.

**Лемма 4.** Пусть $A: L_1^n[t_0, T] \to L_\infty^n[t_0, T]$ линейный оператор, $|(Au)(t)| \le L\int_{t_0}^t |u(s)|ds$ при $u(\cdot) \in L_1^n[t_0, T]$, многозначное отображение $F:[t_0,T] \times R^n \to \text{comp} R^n$ удовлетворяет условиям Каратеодори, существуют суммируемая неотрицательная функция $\alpha(t)$ и число $\beta > 0$ такие, что $\|F(t,x)\| \le \alpha(t) + \beta|x|$ при $x \in R^n$. Тогда множество решений $S$ задачи $u(t) \in F(t,(Au)(t))$ ограниченное множество в $L_1^n[t_0, T]$.

**Доказательство.** Если $x \in S$, то имеем, что

$$|u(t)| \le \|F(t,(Au)(t))\| \le \alpha(t) + \beta|(Au)(t)| \le \alpha(t) + L\beta\int_{t_0}^t |u(s)|ds.$$

Применяя неравенство Грануола (см. [4], стр. 219) имеем, что

$$|u(t)| \le \alpha(t) + L\beta e^{L\beta(T-t_0)}\int_{t_0}^t \alpha(s)ds \le \alpha(t) + L\beta e^{L\beta(T-t_0)}\int_{t_0}^T \alpha(s)ds.$$

Поэтому

$$\int_{t_0}^T |u(t)|dt \le \int_{t_0}^T \alpha(t)dt + L\beta(T-t_0)e^{L\beta(T-t_0)}\int_{t_0}^T \alpha(s)ds.$$

Лемма доказана.

## 2.2. Непрерывная зависимость решений операторного включения условного типа Фредгольма

Пусть $R^n$ $n$-мерное евклидово пространство. Совокупность всех непустых компактных (выпуклых компактных) подмножеств $R^n$ обозначим через $\text{comp} R^n (\text{conv} R^n)$.

Пусть $A: L_p^n[t_0, T] \to C^n[t_0, T]$ линейный непрерывный оператор, где $L_p^n[t_0, T] = \{u(t) = (u^1(t), \ldots, u^n(t)) \in R^n : u^i(\cdot) \in L_p[t_0, T], i = \overline{1,n}\}$ пространство с нормой $\|u(\cdot)\|_{L_p^n} = (\int_{t_0}^T |u(s)|^p ds)^{\frac{1}{p}}$. Множество $\{(u, Au) \in L_p^n[t_0, T] \times C^n[t_0, T] : u \in L_p^n[t_0, T]\}$ относительно нормы $\|(u, Au)\| = \|u\|_{L_p^n[t_0,T]} + \|Au\|_{C^n[t_0,T]}$ является подпространством в $L_p^n[t_0, T] \times C^n[t_0, T]$.

Так как $\|u\|_{L_p^n[t_0,T]} \le \|u\|_{L_p^n[t_0,T]} + \|Au\|_{C^n[t_0,T]} \le (1 + \|A\|)\|u\|_{L_p^n[t_0,T]}$, то имеем, что в пространстве

$\{(u, Au) \in L_p^n[t_0,T] \times C^n[t_0,T] : u \in L_p^n[t_0,T]\}$ нормы $\|(u, Au)\|$ и $\|u\|_{L_p^n[t_0,T]}$ эквивалентны.

Если оператор $A: L_p^n[t_0,T] \to C^n[t_0,T]$ удовлетворяет условию $|(Au)(t) - (A\upsilon)(t)| \leq L(\int_{t_0}^T |u(s) - \upsilon(s)|^p ds)^{\frac{1}{p}}$ при $u(\cdot), \upsilon(\cdot) \in L_p^n[t_0,T]$, то оператор A назовем условного типа Фредгольма, где $L > 0$. Если $A: L_p^n[t_0,T] \to C^n[t_0,T]$ линейный оператор, то условие $|(Au)(t) - (A\upsilon)(t)| \leq L(\int_{t_0}^T |u(s) - \upsilon(s)|^p ds)^{\frac{1}{p}}$ при $u(\cdot), \upsilon(\cdot) \in L_p^n[t_0,T]$ эквивалентно условию $|(Au)(t)| \leq L(\int_{t_0}^T |u(s)|^p ds)^{\frac{1}{p}}$ при $u(\cdot) \in L_p^n[t_0,T]$.

Если $A: L_p^n[t_0,T] \to L_\infty^n[t_0,T]$, то аналогично определяется оператор условного типа Фредгольма.

Пусть $t_0 < T$, $F:[t_0,T] \times R^n \to \text{comp} R^n$ многозначное отображение, $1 \leq p < +\infty$. Рассмотрим задачу для включения

$$u(t) \in F(t,(Au)(t)), \quad t \in [t_0,T], \quad u(\cdot) \in L_p^n[t_0,T] \qquad (1)$$

Функцию $u(\cdot) \in L_p^n[t_0,T]$, удовлетворяющую (1) назовем решением задачи (1).

**Теорема 1.** Пусть $F:[t_0,T] \times R^n \to \text{comp} R^n$ многозначное отображение, $t \to F(t,x)$ измеримо по $t$, существует функция $M(\cdot) \in L_p[t_0,T]$, где $M(t) > 0$ такая, что

$$\rho_x(F(t,x), F(t,x_1))) \leq M(t)|x - x_1|$$

при $x, x_1 \in R^n$. Кроме того, пусть $A: L_p^n[t_0,T] \to C^n[t_0,T]$ линейный непрерывный оператор и $\|A\|(\int_{t_0}^T M^p(s)ds)^{\frac{1}{p}} < 1$; $\rho(\cdot) \in L_p[t_0,T]$ и $\bar{u}(\cdot)$ такие, что

$$d(\bar{u}(t), F(t,(A\bar{u})(t))) \leq \rho(t), \quad t \in [t_0,T].$$

Тогда существует такое решение $u(\cdot) \in L_p^n[t_0,T]$ задачи (1), что

$$|u(t) - \bar{u}(t)| \leq \rho(t) + \frac{\|A\|(\int_{t_0}^T \rho^p(s)ds)^{\frac{1}{p}} M(t)}{1 - \|A\|(\int_{t_0}^T M^p(s)ds)^{\frac{1}{p}}}, \quad |(Au)(t) - (A\bar{u})(t)| \leq \frac{\|A\|(\int_{t_0}^T \rho^p(s)ds)^{\frac{1}{p}}}{1 - \|A\|(\int_{t_0}^T M^p(s)ds)^{\frac{1}{p}}}.$$

**Доказательство.** Построим последовательность $u_i(t)$ ($i = 0,1,2...$) с помощью рекурентного соотношения при $t \in [t_0,T]$

$$x_0(t) = \bar{x}(t) = (A\bar{u})(t), \quad x_{i+1}(t) = (Au_i)(t), \quad i = 0,1,2,... \qquad (2)$$

где $u_i(s) \in F(s,x_i(s))$ при $i \geq 0$, $u_i(s)$ измеримы и

$$|u_0(s) - \bar{u}(s)| = d(\bar{u}(s), F(s,x_0(s))), \quad u_0(s) \in F(s,x_0(s)),$$

$$|u_i(s) - u_{i-1}(s)| = d(u_{i-1}(s), F(s,(Au_{i-1})(s))), \quad i = 1,2,...$$

при $s \in [t_0,T]$. По лемме 2.1.4 [25] такая функция $u_i(t)$ существует. По условию, при $i \geq 1$ имеем

$$|u_{i+1}(s) - u_i(s)| = d(u_i(s), F(s,x_{i+1}(s))) \leq \rho_x(F(s,x_i(s)), F(s,x_{i+1}(s))) \leq M(s)|x_{i+1}(s) - x_i(s)|$$

при $s \in [t_0,T]$. Ясно, что

$$|x_{i+1}(s) - x_i(s)| = |(Au_i)(s) - (Au_{i-1})(s))| \leq \|A\|(\int_{t_0}^T |u_i(s) - u_{i-1}(s)|^p ds)^{\frac{1}{p}}.$$

По условию, имеем, что $|u_0(s) - \bar{u}(s)| \leq d(\bar{u}(s), F(s,(A\bar{u})(s))) \leq \rho(s)$ при $s \in [t_0, T]$. Поэтому имеем, что

$$|x_1(t) - x_0(t)| = |(Au_0)(t) - (A\bar{u})(t)| = |A(u_0 - \bar{u})(t)| \leq \|A\|(\int_{t_0}^T |u_0(s) - \bar{u}(s)|^p ds)^{\frac{1}{p}} \leq \|A\|(\int_{t_0}^T |\rho(s)|^p ds)^{\frac{1}{p}}.$$

Тогда получим

$$|u_1(s) - u_0(s)| = d(u_0(s), F(s, x_1(s))) \leq \rho_x(F(s, x_0(s)), F(s, x_1(s))) \leq$$
$$\leq M(s)|x_0(s) - x_1(s)| \leq \|A\| M(s)(\int_{t_0}^T \rho^p(\tau) d\tau)^{\frac{1}{p}}.$$

при $s \in [t_0, T]$. Ясно, что

$$|x_2(t) - x_1(t)| \leq \|A\|(\int_{t_0}^T |u_1(s) - u_0(s)|^p ds)^{\frac{1}{p}} \leq \|A\|(\int_{t_0}^T M^p(s)\|A\|^p \int_{t_0}^T \rho^p(v) dv ds)^{\frac{1}{p}} =$$
$$= \|A\|^2 (\int_{t_0}^T \rho^p(s) ds)^{\frac{1}{p}} (\int_{t_0}^T M^p(s) ds)^{\frac{1}{p}}.$$

при $t \in [t_0, T]$. Поэтому

$$|u_2(s) - u_1(s)| \leq M(s)\|A\|^2 (\int_{t_0}^T \rho^p(s) ds)^{\frac{1}{p}} (\int_{t_0}^T M^p(s) ds)^{\frac{1}{p}}$$

при $s \in [t_0, T]$. Отсюда получим, что

$$|x_3(t) - x_2(t)| \leq \|A\|^3 (\int_{t_0}^T \rho^p(s) ds)^{\frac{1}{p}} (\int_{t_0}^T M^p(s) ds)^{\frac{2}{p}}$$

при $t \in [t_0, T]$. Тогда

$$|u_3(s) - u_2(s)| \leq M(s)\|A\|^3 (\int_{t_0}^T \rho^p(s) ds)^{\frac{1}{p}} (\int_{t_0}^T M^p(s) ds)^{\frac{2}{p}}$$

при $s \in [t_0, T]$. Поэтому

$$|x_4(t) - x_3(t)| \leq \|A\|^4 (\int_{t_0}^T \rho^p(s) ds)^{\frac{1}{p}} (\int_{t_0}^T M^p(s) ds)^{\frac{3}{p}}$$

при $t \in [t_0, T]$. Также имеем

$$|u_4(s) - u_3(s)| \leq M(s)\|A\|^4 (\int_{t_0}^T \rho^p(s) ds)^{\frac{1}{p}} (\int_{t_0}^T M^p(s) ds)^{\frac{3}{p}}$$

при $s \in [t_0, T]$. Поэтому

$$|x_5(t) - x_4(t)| \leq \|A\|^5 (\int_{t_0}^T \rho^p(s) ds)^{\frac{1}{p}} (\int_{t_0}^T M^p(s) ds)^{\frac{4}{p}}$$

при $t \in [t_0, T]$. Продолжая процесс получим

$$|u_{m+1}(s) - u_m(s)| \leq M(s)\|A\|^{m+1} (\int_{t_0}^T \rho^p(s) ds)^{\frac{1}{p}} (\int_{t_0}^T M^p(s) ds)^{\frac{m}{p}} \qquad (3)$$

$$|x_{m+1}(s) - x_m(s)| \leq \|A\|^{m+1} (\int_{t_0}^T \rho^p(s) ds)^{\frac{1}{p}} (\int_{t_0}^T M^p(s) ds)^{\frac{m}{p}} \qquad (4)$$

при $s \in [t_0, T]$. Тогда получим

$$|u_{m+1}(t) - \bar{u}(t)| \leq |u_0(t) - \bar{u}(t)| + |u_1(t) - u_0(t)| + |u_2(t) - u_1(t)| + |u_3(t) - u_2(t)| + \ldots +$$

$$+\left|u_{m+1}(t)-u_m(t)\right| \le \rho(t)+M(t)(\int_{t_0}^T \rho^p(s)ds)^{\frac{1}{p}}\sum_{i=0}^m \|A\|^{i+1}(\int_{t_0}^T M^p(s)ds)^{\frac{i}{p}}$$

при $s \in [t_0, T]$.

Если $\|A\|(\int_{t_0}^T M^p(s)ds)^{\frac{1}{p}} < 1$, то имеем, что

$$\left|u_{m+1}(t)-\overline{u}(t)\right| \le \rho(t)+M(t)(\int_{t_0}^T \rho^p(s)ds)^{\frac{1}{p}} \frac{\|A\|}{1-\|A\|(\int_{t_0}^T M^p(s)ds)^{\frac{1}{p}}}. \tag{5}$$

Также имеем, что

$$\left|x_{m+1}(t)-\overline{x}(t)\right| \le \left|x_1(t)-\overline{x}(t)\right|+\left|x_2(t)-x_1(t)\right|+\left|x_3(t)-x_2(t)\right|+...+\left|x_{m+1}(t)-x_m(t)\right|=$$
$$\le (\int_{t_0}^T \rho^p(s)ds)^{\frac{1}{p}}\sum_{i=0}^m \|A\|^{i+1}(\int_{t_0}^T M^p(s)ds)^{\frac{i}{p}} \le (\int_{t_0}^T \rho^p(s)ds)^{\frac{1}{p}}\frac{\|A\|}{1-\|A\|(\int_{t_0}^T M^p(s)ds)^{\frac{1}{p}}}; \tag{6}$$

при $t \in [t_0, T]$. Из оценки (3)-(6) вытекает, что последовательности $x_m(t)$ и $u_m(t)$ сходятся соответственно к функциям $x(\cdot)$ и $u(\cdot) \in L_p^n[t_0, T]$. Из теоремы 1.2.19[4] и из (2) следует, что $x(t)$ непрерывно. Так как $u_i(t) \in F(t, x_i(t))$, то из теоремы 1.2.23 и 1.2.28 [3] имеем, что $u(t) \in F(t, x(t))$ при $t \in [t_0, T]$, т.е. из (2) следует, что $u(t) \in F(t,(Au)(t))$ при $t \in [t_0, T]$. Кроме того, из (5) и (6) получим, что

$$\left|u(t)-\overline{u}(t)\right| \le \rho(t)+\frac{\|A\|M(t)(\int_{t_0}^T \rho^p(s)ds)^{\frac{1}{p}}}{1-\|A\|(\int_{t_0}^T M^p(s)ds)^{\frac{1}{p}}}, \qquad \left|x(t)-\overline{x}(t)\right| \le \frac{\|A\|(\int_{t_0}^T \rho^p(s)ds)^{\frac{1}{p}}}{1-\|A\|(\int_{t_0}^T M^p(s)ds)^{\frac{1}{p}}}.$$

Теорема доказана.

В теореме 1 заменив неравенство $\left|(Au)(s)-(Au)(s))\right| \le \|A\|(\int_{t_0}^T \left|u(s)-u(s)\right|^p ds)^{\frac{1}{p}}$ неравенством $\left|(Au)(t)-(A\upsilon)(t)\right| \le L\|u-\upsilon\|_{L_p^n}$ аналогично теореме 1, доказывается следующая теорема.

**Теорема 2.** Пусть $F:[t_0, T] \times R^n \to \text{comp} R^n$ многозначное отображение, $t \to F(t, x)$ измеримо по $t$, существует функция $M(\cdot) \in L_p[t_0, T]$, где $M(t) > 0$ такая, что

$$\rho_x(F(t, x), F(t, x_1))) \le M(t)|x-x_1|$$

при $x, x_1 \in R^n$. Кроме того, пусть $A: L_p^n[t_0, T] \to C^n[t_0, T]$ оператор и существует $L > 0$ такое, что $\left|(Au)(t)-(A\upsilon)(t)\right| \le L\|u-\upsilon\|_{L_p^n}$ при $u(\cdot), \upsilon(\cdot) \in L_p^n[t_0, T]$; пусть $\rho(\cdot) \in L_p[t_0, T]$ и $\overline{u}(\cdot)$ такие, что

$$d(\overline{u}(t), F(t,(A\overline{u})(t))) \le \rho(t), \quad t \in [t_0, T],$$

где $L(\int_{t_0}^T M^p(s)ds)^{\frac{1}{p}} < 1$. Тогда существует такое решение $u(\cdot) \in L_p^n[t_0, T]$ задачи (1), что

$$\left|u(t)-\overline{u}(t)\right| \le \rho(t)+\frac{L(\int_{t_0}^T \rho^p(s)ds)^{\frac{1}{p}}M(t)}{1-L(\int_{t_0}^T M^p(s)ds)^{\frac{1}{p}}}, \qquad \left|(Au)(t)-(A\overline{u})(t)\right| \le \frac{L(\int_{t_0}^T \rho^p(s)ds)^{\frac{1}{p}}}{1-L(\int_{t_0}^T M^p(s)ds)^{\frac{1}{p}}}.$$

Аналогично теореме 1, доказываются следующие теоремы.

**Теорема 3.** Пусть $F:[t_0,T]\times R^n \to \text{comp}R^n$ многозначное отображение, $t \to F(t,x)$ измеримо по $t$, существует функция $M(\cdot) \in L_p[t_0,T]$, где $M(t) > 0$ такая, что
$$\rho_x(F(t,x),F(t,x_1))) \leq M(t)|x-x_1|$$
при $x, x_1 \in R^n$. Кроме того, пусть $A: L_p^n[t_0,T] \to L_\infty^n[t_0,T]$ линейный непрерывный оператор; $\rho(\cdot) \in L_p[t_0,T]$ и $\bar{u}(\cdot)$ такие, что
$$d(\bar{u}(t), F(t,(A\bar{u})(t))) \leq \rho(t), \quad t \in [t_0,T],$$
где $\|A\|(\int_{t_0}^T M^p(s)ds)^{\frac{1}{p}} < 1$. Тогда существует такое решение $u(\cdot) \in L_p^n[t_0,T]$ задачи (1), что

$$|u(t)-\bar{u}(t)| \leq \rho(t) + \frac{\|A\|(\int_{t_0}^T \rho^p(s)ds)^{\frac{1}{p}} M(t)}{1-\|A\|(\int_{t_0}^T M^p(s)ds)^{\frac{1}{p}}}, \qquad |(Au)(t)-(A\bar{u})(t)| \leq \frac{\|A\|(\int_{t_0}^T \rho^p(s)ds)^{\frac{1}{p}}}{1-\|A\|(\int_{t_0}^T M^p(s)ds)^{\frac{1}{p}}}.$$

**Теорема 4.** Пусть $F:[t_0,T]\times R^n \to \text{comp}R^n$ многозначное отображение, $t \to F(t,x)$ измеримо по $t$, существует функция $M(\cdot) \in L_p[t_0,T]$, где $M(t) > 0$ такая, что
$$\rho_x(F(t,x),F(t,x_1))) \leq M(t)|x-x_1|$$
при $x, x_1 \in R^n$. Кроме того, пусть $A: L_p^n[t_0,T] \to L_\infty^n[t_0,T]$ оператор и существует $L > 0$ такое, что $|(Au)(t)-(A\upsilon)(t)| \leq L\|u-\upsilon\|_{L_p^n}$ при $u(\cdot), \upsilon(\cdot) \in L_p^n[t_0,T]$; $\rho(\cdot) \in L_p[t_0,T]$ и $\bar{u}(\cdot)$ такие, что
$$d(\bar{u}(t), F(t,(A\bar{u})(t))) \leq \rho(t), \quad t \in [t_0,T],$$
где $L(\int_{t_0}^T M^p(s)ds)^{\frac{1}{p}} < 1$. Тогда существует такое решение $u(\cdot) \in L_p^n[t_0,T]$ задачи (1), что

$$|u(t)-\bar{u}(t)| \leq \rho(t) + \frac{L(\int_{t_0}^T \rho^p(s)ds)^{\frac{1}{p}} M(t)}{1-L(\int_{t_0}^T M^p(s)ds)^{\frac{1}{p}}}, \qquad |(Au)(t)-(A\bar{u})(t)| \leq \frac{L(\int_{t_0}^T \rho^p(s)ds)^{\frac{1}{p}}}{1-L(\int_{t_0}^T M^p(s)ds)^{\frac{1}{p}}}.$$

Отметим, что в теоремах 2 и 4 оператор $A$ в общем случае не является линейным. Положив $\bar{u}(t) = 0$ и $(A\bar{u})(t) = 0$, из теоремы 1 имеем, что верно следующее следствие.

**Следствие 1.** Пусть $A: L_p^n[t_0,T] \to C^n[t_0,T]$ линейный непрерывный оператор, $F:[t_0,T]\times R^n \to \text{comp}R^n$ многозначное отображение, $t \to F(t,x)$ измеримо по $t$, существует функция $M(\cdot) \in L_p[t_0,T]$, где $M(t) > 0$ такая, что $\rho_x(F(t,x),F(t,x_1)) \leq M(t)|x-x_1|$ при $x,x_1 \in R^n$, $\|A\|(\int_{t_0}^T M^p(s)ds)^{\frac{1}{p}} < 1$ и $\|F(t,0)\| \in L_p[t_0,T]$. Тогда существует решение $u(\cdot) \in L_p^n[t_0,T]$ задачи (1).

**Замечание 1.** Из доказательства теоремы 1 следует, что в теореме 1-4 условие $F:[t_0,T]\times R^n \to \text{comp}R^n$ можно заменить условием: $F:[t_0,T]\times(\bar{x}(t)+\alpha B) \to \text{comp}R^n$ при $t \in [t_0,T]$, где $\bar{x}(t) = (A\bar{u})(t)$, $\alpha \geq \dfrac{\|A\|(\int_{t_0}^T \rho^p(s)ds)^{\frac{1}{p}}}{1-\|A\|(\int_{t_0}^T M^p(s)ds)^{\frac{1}{p}}}$ или $\alpha \geq \dfrac{L(\int_{t_0}^T \rho^p(s)ds)^{\frac{1}{p}}}{1-L(\int_{t_0}^T M^p(s)ds)^{\frac{1}{p}}}$

соответственно, $B = \{z \in R^n : \|z\| \leq 1\}$ -единичный шар в $R^n$.

**Лемма 2.** Пусть $A: L_p^n[t_0, T] \to L_\infty^n[t_0, T]$ линейный непрерывный оператор, $F:[t_0,T] \times (x_0(t) + \alpha B) \to \text{comp} R^n$, где $\alpha > 0$, $u_0(\cdot) \in L_p^n[t_0,T]$, $x_0(t) = (Au_0)(t)$ и функция $u_0(t)$ - решение задачи $u(t) \in F(t,(Au)(t))$ и $F(t,x)$ удовлетворяет условиям теоремы 3 в $x_0(t) + \alpha B$, $\|A\|(\int_{t_0}^T M^p(s)ds)^{\frac{1}{p}} < 1$. Тогда существует такое $\delta > 0$, что при $s(\cdot) \in L_p^n[t_0,T]$, $\|s(\cdot)\|_{L_p^n} \leq \delta$ найдется решение $u_s(t)$ задачи $u(t) \in F(t,(Au)(t)) + s(t)$, что $|x_0(t) - x_s(t)| \leq \alpha$ при $t \in [t_0, T]$, где $x_s(t) = (Au_s)(t)$ и $\|x_s(\cdot) - x_0(\cdot)\|_{L_\infty^n} \to 0$ при $\|s(\cdot)\|_{L_p^n} \to 0$.

**Доказательство.** Ясно, что
$$\rho_x(F(t,x) + s(t), F(t,x') + s(t)) \leq M(t)|x - x'|$$
при $|x - x_0(t)| \leq \alpha$, $|x' - x_0(t)| \leq \alpha$ и $d(u_0(t), F(t,x_0(t)) + s(t)) \leq |s(t)|$. По теореме 3 существует решение $u_s(t)$ задачи $u_s(t) \in F(t, x_s(t)) + s(t)$, такое, что

$$|x_0(t) - x_s(t)| \leq \frac{\|A\|(\int_{t_0}^T s^p(t)dt)^{\frac{1}{p}}}{1 - \|A\|(\int_{t_0}^T M^p(s)ds)^{\frac{1}{p}}}. \qquad (7)$$

Если определять $\delta$ из неравенства $\dfrac{\|A\|\delta}{1 - \|A\|(\int_{t_0}^T M^p(s)ds)^{\frac{1}{p}}} \leq \alpha$, то получим, что верна первая часть утверждения леммы.

Из (7) получим, что $\|x_s(\cdot) - x_0(\cdot)\|_{L_\infty^n[t_0,T]} \leq \dfrac{\|A\|(\int_{t_0}^T s^p(t)dt)^{\frac{1}{p}}}{1 - \|A\|(\int_{t_0}^T M^p(s)ds)^{\frac{1}{p}}}$, т.е. $\|x_s(\cdot) - x_0(\cdot)\|_{L_\infty^n[t_0,T]} \to 0$ при $\|s(\cdot)\|_{L_p^n} \to 0$. Лемма доказана.

**Лемма 3.** Пусть $A: L_p^n[t_0,T] \to L_\infty^n[t_0,T]$ линейный непрерывный оператор, $F:[t_0,T] \times R^n \to \text{comp} R^n \cup \{\varnothing\}$, отображение $t \to F(t,x)$ измеримо на $[t_0,T]$, отображение $x \to F(t,x)$ замкнуто и выпукло почти при всех $t \in [t_0,T]$, т.е. $\text{gr} F_t = \{(x,y) : y \in F(t,x)\}$ замкнуто и выпукло почти при всех $t \in [t_0,T]$. Пусть существует функция $\lambda(\cdot) \in L_p[t_0,T]$, такая, что

$$\|F(t,x)\| \leq \lambda(t)(1 + |x|) \qquad (8)$$

при $x \in R^n$, где $\|F(t,x)\| = \sup\{|y| : y \in F(t,x)\}$, $\|\varnothing\| = 0$ и существует решение $u_0(t)$ задачи $u_0(t) \in F(t, x_0(t))$, что $x_0(t) = (Au_0)(t)$ принадлежит $\text{dom } F_t = \{x : F(t,x) \neq \varnothing\}$ вместе с некоторой $\varepsilon$ трубкой, т.е. $\{x : |x_0(t) - x| \leq \varepsilon\} \subset \text{dom } F_t$, $\|A\|(\int_{t_0}^T M^p(s)ds)^{\frac{1}{p}} < 1$, $M(t) = \dfrac{6\lambda(t)\left(1 + |x_0(t)| + \dfrac{5\varepsilon}{6}\right)}{\varepsilon}$.

Тогда существуют такие $\delta > 0$ и решение $u_s(t)$ задачи $u(t) \in F(t,(Au)(t)) + s(t)$, где $s(\cdot) \in L_p^n[t_0,T]$, $\|s(\cdot)\|_{L_p^n} \leq \delta$, что $x_s(t) = (Au_s)(t)$, $|x_0(t) - x_s(t)| \leq \varepsilon$ при $t \in [t_0, T]$ и $\|x_s(\cdot) - x_0(\cdot)\|_{L_\infty^n[t_0,T]} \to 0$ при $\|s(\cdot)\|_{L_p^n[t_0,T]} \to 0$.

**Доказательство.** Рассмотрим множество $S_t = \{x : |x_0(t) - x| \leq \dfrac{\varepsilon}{2}\}$ и обозначим $\delta = \dfrac{\varepsilon}{3}$. Пусть

$x \in S_t$, тогда $x + \delta B \subset \text{int dom} F_t$, где B единичный шар в $R^n$ с центром в нуле. Так как множество $grF_t$ выпукло, то аналогично доказательству теоремы 3.1.1 [13] можно показать, что

$$\rho_x(F(t,y), F(t,x)) \le \frac{2\lambda(t)(1+|x|+\delta)}{\delta}|y-x| \qquad (9)$$

для $y \in x + \delta B$ п.в. t. Ясно, что $S_t + \delta B \subset \text{int dom} F_t$. Поэтому (9) верно для всех $x, y \in S_t$ и почти всех t, т.е.

$$\rho_x(F(t,y), F(t,x)) \le \frac{2\lambda(t)\left(1+|x_0(t)|+\frac{\varepsilon}{2}+\delta\right)}{\delta}|y-x| \quad .$$

при $x, y \in S_t$. Поэтому доказательство леммы 3 вытекает из леммы 2. Лемма доказана.

Отметим, что используя теоремы 1 или 2 можно получить ряд аналог леммы 2 и 3.

### 2.3. О субдифференцируемости об одного функционала

Пусть $f:[t_0,T] \times R^n \times R^n \to R_{+\infty} = R \cup \{+\infty\}$ нормальный выпуклый интегрант, т.е. f - такая функция из $[t_0,T] \times R^n \times R^n$ в $R_{+\infty}$, что $f_t = f(t,\cdot)$ выпукла на $R^n \times R^n$ для каждого $t \in [t_0,T]$, $\text{epf}_t = \{(z,\alpha) \in R^n \times R^n \times R : f(t,z) \le \alpha\}$ -замкнуто и $t \to \text{epf}_t$ измеримо, $\varphi: R^n \times R^n \to R$ выпуклая функция, $A: L_p^n[t_0,T] \to C^n[t_0,T]$ линейный оператор.

Рассмотрим функционал

$$J(u(\cdot)) = \int_{t_0}^{T} f(t,(Au)(t),u(t))dt + \varphi((Au)(t_0),(Au)(T)) : L_p^n[t_0,T] \to \overline{R},$$

где $1 \le p < +\infty$. Отметим, что $J: L_p^n[t_0,T] \to \overline{R} = R \cup \{\pm\infty\}$ выпуклый функционал. Субградиентами J в точке $u_0(\cdot) \in L_p^n[t_0,T]$ являются по определению, элементы $\upsilon^* \in L_q^n[t_0,T]$, $\frac{1}{p} + \frac{1}{q} = 1$, для которых

$$J(u(\cdot)) - J(u_0(\cdot)) \ge \langle \upsilon^*, u(\cdot) - u_0(\cdot) \rangle$$

при всех $u(\cdot) \in L_p^n[t_0,T]$ или, что равносильно,

$$J^*(\upsilon^*) + J(u_0(\cdot)) = \langle \upsilon^*, u_0(\cdot) \rangle,$$

где $\langle \upsilon^*, u_0(\cdot) \rangle = \int_{t_0}^{T} (\upsilon^*(t)|u_0(t))dt$. Множество всех таких субградиентов обозначается через $\partial J(u_0(\cdot))$ и называется субдифференциалом функционала J в точке $u_0(\cdot)$.

В этом параграфе устанавливается связь между $\partial J(u_0(\cdot))$ и $\partial f_t((Au_0)(t), u_0(t))$, $\varphi((Au_0)(t_0),(Au_0)(T))$.

Множество всех финитных функций из $C^k[t_0,T]$ обозначим через $C_0^k[t_0,T]$ и положим $C_0^\infty[t_0,T] = \bigcap_{k=1}^{\infty} C_0^k[t_0,T]$ и $(C_0^\infty[t_0,T])^n = \prod_{k=1}^{n} C_0^\infty[t_0,T]$ (см.[9]). Отметим, что $(C_0^\infty[t_0,T])^n$ плотно в $L_p^n[t_0,T]$. Если $A: L_p^n[t_0,T] \to C^n[t_0,T]$ линейный суръективный непрерывный оператор, то $(C_0^\infty[t_0,T])^n \subset \{Au: L_p^n[t_0,T], (Au)(t_0) = (Au)(T) = 0\}$.

**Теорема 1.** Пусть $f:[t_0,T] \times R^n \times R^n \to R_{+\infty}$ нормальный выпуклый интегрант, существует $a(\cdot) \in L_1[t_0,T]$ и $c > 0$ такие, что $|f(t,y)| \le a(t) + c|y|^p$ при $y \in R^{2n}$, $\varphi: R^n \times R^n \to R$ выпуклая непрерывная функция, $A: L_p^n[t_0,T] \to C^n[t_0,T]$ линейный непрерывный оператор и $(C_0^\infty[t_0,T])^n \subset \{Au: L_p^n[t_0,T], (Au)(t_0) = (Au)(T) = 0\}$ и в пространстве $(C_0^\infty[t_0,T])^n$ оператор $A^{-1}$ существует. Если $u^*(\cdot) \in \partial J(\bar{u}(\cdot))$, то существуют $\bar{u}_1^* \in L_q^n[t_0,T]$ и числа $a_1, a_2 \in R^n$ такие, что

1) $\int_{t_0}^T ((\bar{u}_1^*)(t)|(Au)(t))dt - \int_{t_0}^T ((A^*\bar{u}_1^*)(t)|u(t))dt + (a_1|(Au)(t_0)) + (a_2|(Au)(T)) = 0$

при $u(\cdot) \in L_p^n[t_0,T]$,

2) $(-\bar{u}_1^*(t), u^*(t) + (A\bar{u}_1^*)(t)) \in \partial f(t,(A\bar{u})(t),\bar{u}(t))$,

3) $(-a_1,-a_2) \in \partial \varphi((A\bar{u})(t_0),(A\bar{u})(T))$.

**Доказательство.** Положим
$$C = \{((Au)(t), u(t), (Au)(t_0), (Au)(T)) : u \in L_p^n[t_0,T]\}.$$

Рассмотрим функционал $I(y_1(t), y(t), c_1, c_2) = \int_{t_0}^T f(t, y_1(t), y(t))dt + \varphi(c_1, c_2)$ в пространстве $L_p^n[t_0,T] \times L_p^n[t_0,T] \times R^n \times R^n$. По условию $\text{dom } I - C = L_p^n[t_0,T] \times L_p^n[t_0,T] \times R^n \times R^n$. Пусть $u^* \in L_q^n[t_0,T]$. Тогда по следствию 4.12[12] существуют $\bar{u}^* \in L_q^n[t_0,T]$, $\bar{u}_1^* \in L_q^n[t_0,T]$ и $a_1, a_2 \in R^n$ такие, что

$$J^*(u^*) = \sup_{u \in L_p^n[t_0,T]} \{\int_{t_0}^T (u^*(t)|u(t))dt - \int_{t_0}^T f(t,(Au)(t),u(t))dt - \varphi((Au)(t_0),(Au)(T))\} =$$

$$= \sup_{y_1, y \in L_p^n[t_0,T],(c_1,c_2) \in R^{2n}} \{\int_{t_0}^T (u^*(t)|y(t))dt - \int_{t_0}^T f(t,y_1(t),y(t))dt - \varphi(c_1,c_2) - \delta_C(y(\cdot),y_1(\cdot),c_1,c_2)\} =$$

$$= \sup_{y_1, y \in L_p^n[t_0,T],(c_1,c_2) \in R^{2n}} \{\int_{t_0}^T (u^*(t) - \bar{u}^*(t)|y(t))dt - \int_{t_0}^T (\bar{u}_1^*(t)|y_1(t))dt - (a_1|c_1) - (a_2|c_2) -$$

$$- \int_{t_0}^T f(t, y_1(t), y(t))dt - \varphi(c_1,c_2)\} + \sup_{y_1, y \in L_p^n[t_0,T],\ (c_1,c_2) \in R^{2n}} \{\int_{t_0}^T (\bar{u}^*(t)|y(t))dt + \int_{t_0}^T (\bar{u}_1^*(t)|y_1(t))dt + (a_1|c_1) +$$

$$+ (a_2|c_2) - \delta_C(y(\cdot), y_1(\cdot), c_1, c_2)\} = \int_{t_0}^T f^*(t, -\bar{u}_1^*(t), u^*(t) - \bar{u}^*(t))dt + \varphi^*(-a_1,-a_2).$$

Отсюда также имеем, что
$$\int_{t_0}^T (\bar{u}^*(t)|u(t))dt + \int_{t_0}^T (\bar{u}_1^*(t)|(Au)(t))dt + (a_1|(Au)(t_0)) + (a_2|(Au)(T)) = 0$$

при $u(\cdot) \in L_p^n[t_0,T]$. По условию имеем, что
$$(C_0^\infty[t_0,T])^n \subset \{Au : L_p^n[t_0,T], (Au)(t_0) = (Au)(T) = 0\}.$$

Поэтому $\int_{t_0}^T (\bar{u}^*(t) + (A^*\bar{u}_1^*)(t)|A^{-1}\upsilon(t))dt = 0$ при $\upsilon(\cdot) \in (C_0^\infty[t_0,T])^n$. Тогда получим, что $(A^*)^{-1}(\bar{u}^*(t) + (A^*\bar{u}_1^*)(t)) = 0$. Поэтому $\bar{u}^*(t) + (A^*\bar{u}_1^*)(t) = 0$ и

$$-\int_{t_0}^T ((A^*\bar{u}_1^*)(t)|u(t))dt + \int_{t_0}^T (\bar{u}_1^*(t)|(Au)(t))dt + (a_1|(Au)(t_0)) + (a_2|(Au)(T)) = 0$$

при $u(\cdot) \in L_p^n[t_0,T]$. Кроме того, получим

$$J^*(u^*(\cdot)) = \int_{t_0}^T f^*(t, -\bar{u}_1^*(t), (A^*\bar{u}_1^*)(t) + u^*(t))dt + \varphi^*(-a_1,-a_2).$$

Если $u^* \in \partial J(\bar{u})$, то имеем, что $J(\bar{u}(\cdot)) + J^*(u^*(\cdot)) = \langle u^*, \bar{u} \rangle$, т.е.

$$\int_{t_0}^T f(t,(A\bar{u})(t),\bar{u}(t))dt + \varphi((A\bar{u})(t_0),(A\bar{u})(T)) +$$

$$+\int_{t_0}^{T}f^*(t,-\overline{u}_1^*(t),(A^*\overline{u}_1^*)(t)+u^*(t))dt+\varphi^*(-a_1,-a_2)=\int_{t_0}^{T}(u^*(t)|\overline{u}(t))dt.$$

Тогда получим

$$\int_{t_0}^{T}f(t,(A\overline{u})(t),\overline{u}(t))dt+\varphi((A\overline{u})(t_0),(A\overline{u})(T))+$$

$$+\int_{t_0}^{T}f^*(t,-\overline{u}_1^*(t),(A^*\overline{u}_1^*)(t)+u^*(t))dt+\varphi^*(-a_1,-a_2)=\int_{t_0}^{T}(u^*(t)|\overline{u}(t))dt+$$

$$+\int_{t_0}^{T}((A^*\overline{u}_1^*)(t)|\overline{u}(t))dt-\int_{t_0}^{T}(\overline{u}_1^*(t)|(A\overline{u})(t))dt-(a_1|(A\overline{u})(t_0))-(a_2|(A\overline{u})(T)).$$

Отсюда по неравенство Фенхеля имеем

$$\int_{t_0}^{T}f(t,(A\overline{u})(t),\overline{u}(t))dt+\int_{t_0}^{T}f^*(t,-\overline{u}_1^*(t),(A^*\overline{u}_1^*)(t)+u^*(t))dt=\int_{t_0}^{T}(u^*(t)|\overline{u}(t))dt+$$

$$+\int_{t_0}^{T}((A^*\overline{u}_1^*)(t)|\overline{u}(t))dt-\int_{t_0}^{T}(\overline{u}_1^*(t)|(A\overline{u})(t))dt,$$

$$\varphi((A\overline{u})(t_0),(A\overline{u})(T))+\varphi^*(-a_1,-a_2)=-(a_1|(A\overline{u})(t_0))-(a_2|(A\overline{u})(T)).$$

Поэтому

$$f(t,(A\overline{u})(t),\overline{u}(t))+f^*(t,-\overline{u}_1^*(t),(A^*\overline{u}_1^*)(t)+\overline{u}^*(t))=(\overline{u}^*(t)|\overline{u}(t))+$$

$$+((A^*\overline{u}_1^*)(t)|\overline{u}(t))-(\overline{u}_1^*(t)|(A\overline{u})(t)).$$

Тогда получим, что

$$(-\overline{u}_1^*(t),\ u^*(t)+(A\overline{u}_1^*)(t))\in\partial f(t,(A\overline{u})(t),\overline{u}(t)),\quad (-a_1,-a_2)\in\partial\varphi((A\overline{u})(t_0),(A\overline{u})(T)),$$

и

$$\int_{t_0}^{T}((\overline{u}_1^*)(t)|(Au)(t))dt-\int_{t_0}^{T}((A^*\overline{u}_1^*)(t)|u(t))dt+(a_1|(Au)(t_0))+(a_2|(Au)(T))=0$$

при $u(\cdot)\in L_p^n[t_0,T]$. Отсюда следует справедливость теоремы. Теорема доказана.

**Следствие 1.** Пусть выполняется условие теоремы 1 и $\overline{u}(\cdot)\in L_p^n[t_0,T]$ минимизирует функционал $J(\cdot)$ в $L_p^n[t_0,T]$. Тогда существуют $\overline{u}^*\in L_q^n[t_0,T]$ и числа $a_1,a_2\in R^n$ такие, что

1) $\int_{t_0}^{T}((\overline{u}^*)(t)|(Au)(t))dt-\int_{t_0}^{T}((A^*\overline{u}^*)(t)|u(t))dt+(a_1|(Au)(t_0))+(a_2|(Au)(T))=0$

   при $u(\cdot)\in L_p^n[t_0,T]$,

2) $(-\overline{u}^*(t),(A\overline{u}^*)(t))\in\partial f(t,(A\overline{u})(t),\overline{u}(t))$,

3) $(-a_1,-a_2)\in\partial\varphi((A\overline{u})(t_0),(A\overline{u})(T))$.

## 2.4. О субдифференцируемости интегрального функционала с оператором

Пусть $f:[t_0,T]\times R^n\to R_{+\infty}$ нормальный выпуклый интегрант, т.е. $f$ -такая функция из $[t_0,T]\times R^n$ в $R_{+\infty}$, что $\mathrm{ep}f_t=\{(z,\alpha)\in R^n\times R:f(t,z)\leq\alpha\}$ -замкнуто и выпкло и $t\to\mathrm{ep}f_t$ измеримо.

Рассмотрим функционал

$$J(u(\cdot))=\int_{t_0}^{T}f(t,(Au)(t))dt,$$

где $u\in L_p^n[t_0,T],\ 1\leq p<+\infty.$

**Лемма 1.** Пусть $f:[t_0,T]\times R^n \to R_{+\infty}$ нормальный выпуклый интегрант, функция $x \to f(t,x)$ не убывает, $A:L_p^n[t_0,T] \to L_\infty^n[t_0,T]$ выпуклый оператор, то $J:L_p^n[t_0,T] \to \overline{R}$ выпуклый функционал.

**Доказательство.** Пусть $u_1(\cdot), u_2(\cdot) \in L_p^n[t_0,T]$ и $\alpha \in [0,1]$. Тогда получим, что

$$J(\alpha u_1 + (1-\alpha)u_2) = \int_{t_0}^T f(t, A(\alpha u_1 + (1-\alpha)u_2)(t))dt \le$$

$$\le \int_{t_0}^T f(t, \alpha A(u_1)(t) + (1-\alpha)A(u_2)(t))dt \le$$

$$\le \alpha \int_{t_0}^T f(t, A(u_1)(t))dt + (1-\alpha)\int_{t_0}^T f(t, A(u_2)(t))dt.$$

Лемма доказана.

Пусть $f:[t_0,T]\times R^n \to R_{+\infty}$ нормальный интегрант, $x \to f(t,x)$ не убывающая сублинейная функция в $R^n$, $A:L_p^n[t_0,T] \to L_\infty^n[t_0,T]$ сублинейный непрерывный оператор, для всякого $x \in R^n$ функция $f(t,x)$ суммируема. По условию имеем, что $J(u(\cdot)) = \int_{t_0}^T f(t,(Au)(t))dt$ сублинейная функция в $L_p^n[t_0,T]$. Так как $x \to f(t,x)$ не убывающая сублинейная функция в $R^n$, то имеем, что $\partial f(t,0) \subset R_+^n$. Так как $A:L_p^n[t_0,T] \to L_\infty^n[t_0,T]$ сублинейный непрерывный оператор, то $Au = \max\{Lu : L \in \partial A\}$ при $u \in L_p^n[t_0,T]$, где $L:L_p^n[t_0,T] \to L_\infty^n[t_0,T]$ линейный непрерывный оператор и $(Au)(t) \ge (Lu)(t)$ при $t \in [t_0,T]$, т.е. $L \in \partial A$. Покажем, что если $y^*(\cdot) = (L^*u^*)(t)$, где $L \in \partial A$, $u^*(t) \in \partial f(t,0)$, $u^*(\cdot) \in L_1^n[t_0,T]$, то $y^*(\cdot) \in \partial J(0)$. Ясно, что

$$\int_{t_0}^T (y^*(t)|u(t))dt = \int_{t_0}^T ((L^*u^*)(t)|u(t))dt = \int_{t_0}^T (u^*(t)|(Lu)(t))dt \le \int_{t_0}^T \sup_{L \in \partial A,\ u^*(t) \in \partial f(t,0)} (u^*(t)|(Lu)(t))dt =$$

$$= \int_{t_0}^T \sup_{u^*(t) \in \partial f(t,0)} (u^*(t)|\sup_{L \in \partial A}(Lu)(t))dt = \int_{t_0}^T \sup_{u^*(t) \in \partial f(t,0)} (u^*(t)|(Au)(t))dt = \int_{t_0}^T f(t,(Au)(t))dt.$$

Отсюда следует, что $y^*(\cdot) \in \partial J(0)$.

Пусть $A:L_p^n[t_0,T] \to C^n[t_0,T]$ линейный непрерывный оператор. Рассмотрим субдифференцируемость интегрального функционала

$$J(u(\cdot)) = \int_{t_0}^T f(t,(Au)(t))dt$$

в $L_p^n[t_0,T]$, где $1 \le p < +\infty$.

Отметим, что если $f:[t_0,T]\times R^n \to R_{+\infty}$ нормальный выпуклый интегрант, $A:L_p^n[t_0,T] \to C^n[t_0,T]$ линейный оператор, то $J:L_p^n[t_0,T] \to \overline{R}$ выпуклый функционал.

Субградиентами $J$ в точке $u_0(\cdot) \in L_p^n[t_0,T]$ являются по определению, элементы $\upsilon^* \in L_q^n[t_0,T]$, $\frac{1}{p}+\frac{1}{q}=1$, для которых

$$J(u(\cdot)) - J(u_0(\cdot)) \ge \langle \upsilon^*, u(\cdot) - u_0(\cdot) \rangle$$

при всех $u(\cdot) \in L_p^n[t_0, T]$ или, что равносильно,

$$J^*(\upsilon^*) + J(u_0(\cdot)) = \langle \upsilon^*, u_0(\cdot) \rangle.$$

Множество всех таких субградиентов обозначается через $\partial J(u_0(\cdot))$ и называется субдифференциалом функционала $J$ в точке $u_0(\cdot)$.

В этом параграфе устанавливается связь между $\partial J(u_0(\cdot))$ и $\partial f_t((Au_0)(t))$.

**Лемма 2.** Пусть $f:[t_0, T] \times R^n \to R_{+\infty}$ нормальный интегрант, $x \to f(t,x)$ сублинейная функция в $R^n$, $A: L_p^n[t_0, T] \to C^n[t_0, T]$ линейный непрерывный компактный оператор, для всякого $x \in R^n$ функция $f(t,x)$ суммируема и $y^*(\cdot) \in L_q^n[t_0, T]$, где $\dfrac{1}{p} + \dfrac{1}{q} = 1$. Тогда

$$J^*(y^*(\cdot)) = \begin{cases} 0, & y^*(\cdot) = (A^*u^*)(t), \text{ где } u^*(t) \in \partial f(t,0), \ u^*(\cdot) \in L_1^n[t_0, T], \\ +\infty, & \text{в других случаях.} \end{cases}$$

**Доказательство.** По условию леммы имеем, что $J(u(\cdot)) = \int_{t_0}^T f(t, (Au)(t)) dt$ сублинейная функция в $L_p^n[t_0, T]$. По условию для всякого $x \in R^n$ функция $f(t,x)$ суммируема. Тогда имеем, что (см.[14]) функционал $x(\cdot) \to \int_{t_0}^T f(t, x(t)) dt$ непрерывен в $C^n[t_0, T]$.

По условию оператор $A: L_p^n[t_0, T] \to C^n[t_0, T]$ непрерывен. Так как $x(\cdot) \to \int_{t_0}^T f(t, x(t)) dt$ непрерывен в $C^n[t_0, T]$, то существует $\delta > 0$ такое, что функционал $\int_{t_0}^T f(t, x(t)) dt$ ограничен в множестве $\{x(\cdot) \in C^n[t_0, T] : \|x(\cdot)\| \le \delta\}$. Так как $A: L_p^n[t_0, T] \to C^n[t_0, T]$ непрерывен, то существует $\alpha > 0$ такое, что $\|Au\|_{C^n} \le \delta$ при $\|u\|_{L_\infty^n} \le \alpha$. Поэтому получим, что функционал $J(u(\cdot)) = \int_{t_0}^T f(t, (Au)(t)) dt$ ограничен в множестве $\{u(\cdot) \in L_p^n[t_0, T] : \|u(\cdot)\|_{L_p^n} \le \alpha\}$. Так как $J(u)$ выпуклый функционал, то имеем, что $u \to J(u)$ непрерывный функционал в $L_p^n[t_0, T]$.

Ясно, что $I_0(y) = \int_{t_0}^T f(t, y(t)) dt$ сублинейная функция в $L_\infty^n[t_0, T]$ и по следствию 2A[14] непрерывный функционал в $L_\infty^n[t_0, T]$. Поэтому по предложению 3 ([6], стр.210) и по теореме 3 ([6], стр.362) $\partial I_0(0)$ слабо$^*$ компактно и $\partial I_0(0) \subset L_1^n[t_0, T]$. Из теоремы 6.4.8[8] следует, что

$$I_0(y(\cdot)) = \max_{u^* \in \partial I_0(0)} \langle u^*, y \rangle = \max_{u^* \in \partial I_0(0)} \int_{t_0}^T (u^*(t) | y(t)) dt$$

при $y(\cdot) \in L_\infty^n[t_0, T]$. Известно, что $u^* \in \partial I_0(0)$ в том и только в том случае, когда $u^*(t) \in \partial f(t,0)$ и $u^*(\cdot) \in L_1^n[t_0, T]$.

Ясно, что $J(u(\cdot)) = \max\limits_{u^* \in \partial I_0(0)} \int_{t_0}^T (u^*(t) | (Au)(t)) dt$ при $u(\cdot) \in L_p^n[t_0, T]$. Так как $A: L_p^n[t_0, T] \to C^n[t_0, T]$ линейный непрерывный компактный оператор и $y^*(\cdot) \in L_q^n[t_0, T]$, где $\dfrac{1}{p} + \dfrac{1}{q} = 1$, то применяя теоремы 6.2.7[12] имеем

$$J^*(y^*(\cdot)) = \sup_{u(\cdot)\in L_p^n[t_0,T]} \{\int_{t_0}^T (y^*(t)|u(t))dt - \int_{t_0}^T f(t,(Au)(t))dt\} =$$

$$= \sup_{u(\cdot)\in L_p^n[t_0,T]} \{\int_{t_0}^T (y^*(t)|u(t))dt - \sup_{u^*\in\partial I_0(0)} \int_{t_0}^T (u^*(t)|(Au)(t))dt\} =$$

$$= \sup_{u(\cdot)\in L_p^n[t_0,T]} \inf_{u^*\in\partial I_0(0)} \{\int_{t_0}^T (y^*(t)-(A^*u^*)(t)|u(t))dt = \inf_{u^*\in\partial I_0(0)} \sup_{u(\cdot)\in L_p^n[t_0,T]} \{\int_{t_0}^T (y^*(t)-(A^*u^*)(t)|u(t))dt =$$

$$= \begin{cases} 0, & y^*(\cdot) = (A^*u^*)(t),\ u^*(t)\in\partial f(t,0),\ u^*(\cdot)\in L_1^n[t_0,T], \\ +\infty, & \text{в других случаях}. \end{cases}$$

Лемма доказана.

**Следствие 1.** Если удовлетворяется условие леммы 1, то $y^*(\cdot)\in\partial J(0)$ в том и только в том случае, когда существует функция $u^*(t)\in\partial f(t,0)$, $u^*(\cdot)\in L_1^n[t_0,T]$, такая, что $y^*(\cdot)=(A^*u^*)(t)$.

**Теорема 1.** Пусть $f:[t_0,T]\times R^n\to R_{+\infty}$ нормальный выпуклый интегрант, существует $\delta>0$ такое, что $f(t,\overline{x}(t)+x)$ суммируема при $x\in R^n$, $|x|\le\delta$, где $\overline{x}(t)=(A\overline{u})(t)$, $\overline{u}(\cdot)\in L_p^n[t_0,T]$, $A:L_p^n[t_0,T]\to C^n[t_0,T]$ линейный непрерывный компактный оператор. Тогда $\partial J(\overline{u})$ непусто и $y^*(\cdot)\in L_q^n[t_0,T]$ принадлежит $\partial J(\overline{u})$ в том и только в том случае, когда существует функция $u^*(t)\in\partial f(t,\overline{x}(t))$, $u^*(\cdot)\in L_1^n[t_0,T]$, такая, что $y^*(\cdot)=(A^*u^*)(t)$.

**Доказательство.** Ясно, что

$$J'(\overline{u}(\cdot);u(\cdot)) = \lim_{\lambda\downarrow 0}\frac{1}{\lambda}(J(\overline{u}(\cdot)+\lambda u(\cdot))-J(\overline{u}(\cdot))=$$

$$= \lim_{\lambda\downarrow 0}\frac{1}{\lambda}(\int_{t_0}^T f(t,(A\overline{u})(t)+\lambda(Au)(t))-f(t,(A\overline{u})(t))dt.$$

Из доказательства предложения 4.1[11] имеем

$$f(t,(A\overline{u})(t)-f(t,(A\overline{u})(t)-(Au)(t))\le\frac{1}{\lambda}(f(t,(A\overline{u})(t)+\lambda(Au)(t))-f(t,(A\overline{u})(t)\le$$

$$\le (f(t,(A\overline{u})(t)+(Au)(t))-f(t,(A\overline{u})(t)).$$

при $\lambda\in(0,1)$. Пусть $\alpha>0$ такое, что $\|Au\|_{C^n}\le\delta$ при $\|u\|_{L_p^n}\le\alpha$. Так как $f(t,(A\overline{u})(t)-(Au)(t))$ и $f(t,(A\overline{u})(t)+(Au)(t))$ суммируемы при $u(\cdot)\in L_p^n[t_0,T]$, $\|u\|_{L_p^n}\le\alpha$, то используя теорему Лебега получим

$$J'(\overline{u}(\cdot);u(\cdot)) = \int_{t_0}^T f'(t,(A\overline{u})(t);(Au)(t))dt.$$

Если учесть, что $\partial J(\overline{u}(\cdot)=\partial J'(\overline{u}(\cdot);0)$, $\partial f'(t,(A\overline{u})(t))=\partial f'(t,(A\overline{u})(t);0)$, то из следствия 1 имеем, что $y^*(\cdot)\in L_q^n[t_0,T]$ принадлежит $\partial J(\overline{u})$ в том и только в том случае, когда существует функция $u^*(t)\in\partial f(t,\overline{x}(t))$, $u^*(\cdot)\in L_1^n[t_0,T]$, такая, что $y^*(\cdot)=(A^*u^*)(t)$. Теорема доказана.

Аналогично лемме 2 доказывается следующая лемма.

**Лемма 3.** Пусть $f:[t_0,T]\times R^n\to R_{+\infty}$ нормальный интегрант, $x\to f(t,x)$ не убывающая сублинейная функция в $R^n$, $A:L_p^n[t_0,T]\to L_\infty^n[t_0,T]$ сублинейный непрерывный оператор, для всякого $x\in R^n$ функция $f(t,x)$ суммируема, множество $\{(L^*u^*)(t):L\in\partial A,\ u^*(t)\in\partial f(t,0),\ u^*(\cdot)\in L_1^n[t_0,T]\}$ компактно в $L_q^n[t_0,T]$ и $y^*(\cdot)\in L_q^n[t_0,T]$. Тогда

$$J^*(y^*(\cdot)) = \begin{cases} 0, & y^*(\cdot) = (L^*u^*)(t), \text{ где } L \in \partial A, \ u^*(t) \in \partial f(t,0), \ u^*(\cdot) \in L_1^n[t_0,T], \\ +\infty, & \text{в других случаях.} \end{cases}$$

**Следствие 2.** Если удовлетворяется условие леммы 3, то $y^*(\cdot) \in \partial J(0)$ в том и только в том случае, когда существует функция $u^*(t) \in \partial f(t,0)$, $u^*(\cdot) \in L_1^n[t_0,T]$ и $L \in \partial A$ такая, что $y^*(\cdot) = (L^*u^*)(t)$.

## 2.5. О субдифференцируемости терминального функционала

Пусть $\varphi : R^{(k+1)n} \to (-\infty, +\infty]$ выпуклая функция, $A : L_p^n[t_0,T] \to C^n[t_0,T]$ линейный непрерывный оператор, $t_0 < t_1 < t_2 < \ldots < t_k = T$.

Рассмотрим субдифференцируемость функционала
$$F(u(\cdot)) = \varphi((Au)(t_0), (Au)(t_1), \ldots, (Au)(t_k))$$
в $L_p^n[t_0,T]$, где $1 \le p < +\infty$. По определению $F^* : L_q^n[t_0,T] \to \overline{R}$ и
$$F^*(\upsilon(\cdot)) = \sup_{u \in L_p^n} \{ \int_{t_0}^T (\upsilon(t)|u(t))dt - \varphi((Au)(t_0), (Au)(t_1), \ldots, (Au)(t_k)) \},$$
где $\upsilon(\cdot) \in L_q^n[t_0,T]$.

Пусть $\alpha \in [0,1]$. Ясно, что
$$F(\alpha u_1(\cdot) + (1-\alpha)u_2(\cdot)) = \varphi((\alpha((Au_1)(t_0), (Au_1)(t_1), \ldots, (Au_1)(t_k)) +$$
$$+ (1-\alpha)((Au_2)(t_0), (Au_2)(t_1), \ldots, (Au_2)(t_k))) \le$$
$$\le \alpha \varphi((Au_1)(t_0), (Au_1)(t_1), \ldots, (Au_1)(t_k)) +$$
$$+ (1-\alpha)\varphi((Au_2)(t_0), (Au_2)(t_1), \ldots, (Au_2)(t_k)) = \alpha F(u_1(\cdot)) + (1-\alpha)F(u_2(\cdot)).$$

Получим, что $F : L_p^n[t_0,T] \to R_{+\infty}$ выпуклый функционал.

**Теорема 1.** Если $\varphi : R^{(k+1)n} \to R_{+\infty}$ собственная выпуклая функция и существует точка $u_0(\cdot) \in L_p^n[t_0,T]$ такая, что функция $\varphi$ непрерывна в точке $((Au_0)(t_0), (Au_0)(t_1), \ldots, (Au_0)(t_k))$, то $\upsilon(\cdot) \in \partial F(\overline{u}(\cdot))$, где $\upsilon(\cdot) \in L_q^n[t_0,T]$, в том и только в том случае, когда существует вектор $(-c_0, -c_1, \ldots, -c_k) \in \partial \varphi((A\overline{u})(t_0), (A\overline{u})(t_1), \ldots, (A\overline{u})(t_k))$, где $(c_0, c_1, \ldots, c_k) \in R^{(k+1)n}$, такой, что
$$\int_{t_0}^T (\upsilon(t)|u(t))dt + \sum_{i=0}^k (c_i|(Au)(t_i)) = 0 \quad \text{при} \quad u(\cdot) \in L_p^n[t_0,T].$$

**Доказательство.** Так как $\varphi$ – непрерывная в точке $((Au_0)(t_0), (Au_0)(t_1), \ldots, (Au_0)(t_k))$, то для $\varepsilon > 0$ существует $\delta > 0$ такое, что
$$|\varphi((Au_0)(t_0), \ldots, (Au_0)(t_k)) + (x_0, x_1, \ldots, x_k) - \varphi((Au_0)(t_0), \ldots, (Au_0)(t_k))| \le \varepsilon$$
при $(x_0, x_1, \ldots, x_k) \in R^{(k+1)n}$, $|(x_0, x_1, \ldots, x_k)| \le \delta$. Так как $A : L_p^n[t_0,T] \to C^n[t_0,T]$ непрерывный оператор, то существует $\alpha > 0$ такое, что $\|Au - Au_0\|_{C^n[t_0,T]} \le \delta$ при $u(\cdot) \in L_p^n[t_0,T]$, $\|u - u_0\|_{L_p^n} \le \alpha$. Поэтому
$$|F(u(\cdot) - F(u_0(\cdot))| = |\varphi((Au)(t_0), \ldots, (Au)(t_k)) - \varphi((Au_0)(t_0), \ldots, (Au_0)(t_k))| \le \varepsilon$$
при $u(\cdot) \in L_p^n[t_0,T]$, $\|u - u_0\|_{L_p^n} \le \alpha$, т.е. $F(u(\cdot))$ непрерывен в точке $u_0(\cdot)$. Положим
$$C = \{(u(\cdot), (Au)(t_0), \ldots, (Au)(t_k)) : u(\cdot) \in L_p^n[t_0,T]\}.$$
Тогда по теореме 3.4.1 [6] имеем, что

$$J_0^*(\upsilon(\cdot)) = \sup_{u(\cdot) \in L_p^n[t_0,T]} \{\int_{t_0}^T (\upsilon(t)|u(t))dt - \varphi((Au)(t_0),\ldots,(Au)(t_k))\} =$$

$$= \sup_{(y(\cdot),a_0,\ldots,a_k) \in L_p^n[t_0,T] \times R^{(k+1)n}} \int_{t_0}^T (\upsilon(t)|y(t))dt - \varphi(a_0,\ldots,a_k) - \delta_C(y(\cdot),a_0,\ldots,a_k)\} =$$

$$= \sup_{(y(\cdot),a_0,\ldots,a_k) \in L_p^n[t_0,T] \times R^{(k+1)n}} \{\int_{t_0}^T (\upsilon(t)|y(t))dt + \sum_{i=0}^k (c_i|a_i) - \delta_C(y(\cdot),a_0,\ldots,a_k)\} +$$

$$+ \sup_{(y(\cdot),a_0,\ldots,a_k) \in L_p^n[t_0,T] \times R^{(k+1)n}} \{-\sum_{i=0}^k (c_i|a_i) - \varphi(a_0,\ldots,a_k)\} =$$

$$= \sup_{u \in L_p^n[t_0,T]} \{\int_{t_0}^T (\upsilon(t)|u(t))dt + \sum_{i=0}^k (c_i|(Au)(t_i))\} + \sup_{(a_0,\ldots,a_k) \in R^{(k+1)n}} \{\sum_{i=0}^k (-c_i|a_i) - \varphi(a_0,\ldots,a_k)\} =$$

$$= \begin{cases} \varphi^*(-c_0,-c_1,\ldots,-c_k) : \int_{t_0}^T (\upsilon(t)|u(t))dt + \sum_{i=0}^k (c_i|(Au)(t_i)) = 0 \text{ при } u(\cdot) \in L_p^n[t_0,T], \\ +\infty : \text{в других случаях.} \end{cases}$$

По определению $\upsilon(\cdot) \in \partial F(\overline{u}(\cdot))$ в том и только в том случае, когда $F(\overline{u}(\cdot)) + F^*(\upsilon(\cdot)) = \langle \upsilon(\cdot), \overline{u}(\cdot) \rangle$. Так как $\int_{t_0}^T (\upsilon(t)|u(t))dt + \sum_{i=0}^k (c_i|(Au)(t_i)) = 0$ при $u(\cdot) \in L_p^n[t_0,T]$, то имеем, что $\int_{t_0}^T (\upsilon(t)|u(t))dt = -\sum_{i=0}^k (c_i|(Au)(t_i))$. Отсюда следует, что

$$\varphi((A\overline{u})(t_0),\ldots,(A\overline{u})(t_k)) + \varphi^*(-c_0,-c_1,\ldots,-c_k) = \sum_{i=0}^k (-c_i|(A\overline{u})(t_i)).$$

Поэтому $(-c_0,-c_1,\ldots,-c_k) \in \partial \varphi((A\overline{u})(t_0),(A\overline{u})(t_0),\ldots,(A\overline{u})(t_k))$. Теорема доказана.

### 2.6. О минимизации вариационных задач с оператором

В работе получены необходимые и достаточные условия экстремума для вариационных задач с оператором.

Пусть $f:[t_0,T] \times R^n \times R^n \to R_{+\infty}$ нормальный выпуклый интегрант, $\varphi: R^n \times R^n \to R_{+\infty}$ выпуклая функция, $A: L_p^n[t_0,T] \to C^n[t_0,T]$ линейный непрерывный оператор.

Рассматривается задача минимизации функционала

$$J(u) = \varphi((Au)(t_0),(Au)(T)) + \int_{t_0}^T f(t,(Au)(t),u(t))dt, \qquad (1)$$

в пространстве $L_p^n[t_0,T]$, где $1 \le p < +\infty$.

Будем говорить, что кривая $\overline{u}(t)$, $t_0 \le t \le T$, является решением задачи (1), если $|J(\overline{u})| < +\infty$ и справедливо неравенство $J(u) \ge J(\overline{u})$ при любой $u(\cdot) \in L_p^n[t_0,T]$.

Рассмотрим функционал

$$\Phi(u,\upsilon) = \varphi((Au)(t_0),(Au)(T)) + \int_{t_0}^T f(t,(Au)(t),u(t)+\upsilon(t))dt,$$

где $\upsilon(\cdot) \in L_p^n[t_0,T]$. Ясно, что $\Phi: L_p^n[t_0,T] \times L_p^n[t_0,T] \to R \cup \{+\infty\}$ и $\Phi(u,0) = J(u)$. Положим $h(\upsilon) = \inf_{u \in L_p^n[t_0,T]} \Phi(u,\upsilon)$. Из предложения 2.5 [11] вытекает, что $h$ выпуклая функция.

Задача (1) называется стабильной, если $h(0)$ конечно и $h$ субдифференцируема в нуле.

**Лемма 1.** Допустим, что $A: L_p^n[t_0,T] \to C^n[t_0,T]$ линейный непрерывный оператор, $\inf_{u \in L_p^n[t_0,T]} J(u)$ конечен и существует такая функция $u_0(\cdot) \in L_p^n[t_0,T]$, что функция $f(t,(Au)(t)+y,u_0(t))$ суммируема при $y \in R^n$, $|y| < r$ для некоторого $r > 0$, а функция $\varphi$ непрерывна в точке $((Au_0)(t_0),(Au_0)(T))$. Тогда функция $h$ субдифференцируема в нуле, т.е. задача (1) стабильна (см. [26]).

**Доказательство.** Обозначим $x_0(t) = (Au_0)(t)$. Так как
$$I(x) = \int_{t_0}^{T} f(t, x_0(t) + x(t), u_0(t)) dt$$
непрерывен в точке нуль в пространстве $C^n[t_0,T]$, то существуют $\alpha_1 > 0$ и $L_1$, что $I(x) \leq L_1$ при $\|x(\cdot)\| \leq \alpha_1$, $x(\cdot) \in C^n[t_0,T]$. Те же свойства $\varphi(\cdot,\cdot)$ получим, что существуют $\alpha_2 > 0$ и $L_2$, что $\varphi(a,b) \leq L_2$ при $|a - x_0(t_0)| \leq \alpha_2, |b - x_0(T)| \leq \alpha_2, (a,b) \in R^n \times R^n$. Обозначим $\alpha = \min\{\alpha_1, \alpha_2\}$. Так как оператор $A: L_p^n[t_0,T] \to C^n[t_0,T]$ непрерывный, то существует $\delta > 0$ такое, что $\|Au - Au_0\| \leq \alpha$ при $u(\cdot) \in L_p^n[t_0,T]$ и $\|u - u_0\|_{L_p^n} \leq \delta$. Обозначив $\upsilon = u_0 - u_\upsilon$, где $\|u_\upsilon - u_0\|_{L_p^n} \leq \delta$, имеем
$$h(\upsilon) = \inf_{u \in L_p^n[t_0,T]} \Phi(u,\upsilon) \leq \Phi(u_\upsilon,\upsilon) = \varphi((A(u_0 - \upsilon))(t_0),(A(u_0-\upsilon))(T)) +$$
$$+ \int_{t_0}^{T} f(t,(A(u_0 - \upsilon))(t),u_0(t)) dt \leq L_1 + L_2$$
при $\upsilon(\cdot) \in L_p^n[t_0,T]$ и $\|\upsilon\|_{L_p^n} \leq \delta$. Так как $h(0) = \inf_{u \in L_p^n[t_0,T]} J(u)$ конечен, тогда из предложения 1.5.2 [26] вытекает, что $h$ субдифференцируема в точке нуль. Лемма доказана.

Пусть $\upsilon, w \in R^n$. Положим $f^0(t,z,\upsilon) = \inf_{w}\{(w|\upsilon) + f(t,z,w)\}$.

**Лемма 2.** Если $f$ нормальный выпуклый интегрант на $[t_0,T] \times (R^n \times R^n)$, $u_0(\cdot) \in L_p^n[t_0,T]$ и функция $f(t,(Au)(t)+y,u_0(t))$ суммируема при $y \in R^n$, $|y| < r$ для некоторого $r > 0$, $A: L_p^n[t_0,T] \to C^n[t_0,T]$ линейный непрерывный оператор, $\upsilon^*(\cdot) \in L_q^n[t_0,T]$ и $I_1(u) = \int_{t_0}^{T} f^0(t,(Au)(t),\upsilon^*(t)) dt$ конечен в точке $u_0(\cdot)$, то $I_1(u)$ непрерывен в точке $u_0(\cdot) \in L_p^n[t_0,T]$.

**Доказательство.** Если $\upsilon^*(\cdot) \in L_q^n[t_0,T]$, то ясно, что
$$f^0(t,(Au_0)(t)+y,\upsilon^*(t)) = \inf_{\omega \in R^n}\{(\omega|\upsilon^*(t)) + f(t,(Au_0)(t)+y,\omega)\} \leq$$
$$\leq (u_0(t)|\upsilon^*(t)) + f(t,(Au_0)(t)+y,u_0(t)).$$
Так как функция $f(t,(Au_0)(t)+y,u_0(t))$ суммируема при $y \in R^n$, $|y| < r$ для некоторого $r > 0$, то отсюда следует, что в пространстве $L_\infty^n[t_0,T]$ функционал $y(\cdot) \to \int_{t_0}^{T}((u_0(t)|\upsilon^*(t)) + f(t,(Au_0)(t)+y(t),u_0(t)) dt$ ограничено сверху в окрестности точки нуль.

Поэтому функционал $u(\cdot) \to \int_{t_0}^{T} f^0(t,(Au)(t),\upsilon^*(t)) dt$ также ограничено сверху в окрестности точки $u_0(\cdot)$ в $L_p^n[t_0,T]$. Тогда из предложения 1.2.5 [26] вытекает, функционал $I_1(u)$ непрерывен в точке $u_0(\cdot)$. Лемма доказана.

Если $I_1(y) = \int_{t_0}^{T} f^0(t,(Au_0)(t)+y(t),\upsilon^*(t)) dt$ собственный функционал в $L_\infty^n[t_0,T]$, то из доказательства леммы 2 имеем, что $f^0(t,(Au_0)(t)+y,\upsilon^*(t))$ суммируема при $y \in R^n$, $|y| < r$. Тогда из предложения 8.3.4[6] вытекает, что функционал $I_1(u)$ непрерывен в точке $u_0(\cdot) \in L_p^n[t_0,T]$.

**Теорема 1.** Пусть $f:[t_0,T] \times (R^n \times R^n) \to R_{+\infty}$ нормальный выпуклый интегрант, $\varphi: R^n \times R^n \to R_{+\infty}$ выпуклая функция, $A: L_p^n[t_0,T] \to C^n[t_0,T]$ линейный непрерывный оператор. Для того, чтобы функция $\bar{u}(t)$ среди всех функций $u \in L_p^n[t_0,T]$ минимизировала функционал (1) достаточно, чтобы нашлись $\upsilon^*(\cdot) \in L_q^n[t_0,T]$, $u^*(\cdot) \in L_1^n[t_0,T]$ и $c_1, c_2 \in R^n$ такие, что

1) $u^*(t) \in \partial f^0(t,(A\bar{u})(t),\upsilon^*(t))$,

2) $(-c_1,-c_2) \in \partial \varphi((A\bar{u})(t_0),(A\bar{u})(T))$,

3) $f^0(t,(A\bar{u})(t),\upsilon^*(t)) = (\bar{u}(t)|\upsilon^*(t)) + f(t,(A\bar{u})(t),\bar{u}(t))$,

4) $\int\limits_{t_0}^{T}(\upsilon^*(t) - (A^*u^*)(t)|u(t))dt + ((Au)(t_0)|c_1) + ((Au)(T)|c_2)) = 0$ при $u(\cdot) \in L_p^n[t_0,T]$,

а если выполнено условие леммы 1 при $u_0(\cdot) = \bar{u}(\cdot)$ и $A: L_p^n[t_0,T] \to C^n[t_0,T]$ линейный непрерывный компактный оператор, то условия 1)-4) и являются необходимыми.

**Доказательство.** Достаточность теоремы непосредственно проверяется.

**Необходимость.** Из леммы 1 вытекает, что $h$ субдифференцируема в точке нуле.

Поэтому из замечания 3.2.3 и предложения 3.2.4 [26] вытекает, что все решения $\bar{u}(\cdot)$ задачи $\inf\{J(u): u \in L_p^n[t_0,T]\}$ и все решения $\upsilon^*(\cdot)$ задачи $\sup\limits_{\upsilon^* \in L_q^n[t_0,T]}\{-\Phi^*(0,\upsilon^*)\}$ связаны экстремальным соотношением

$$\Phi(\bar{u},0) + \Phi^*(0,-\upsilon^*) = 0 \qquad (2)$$

По определению

$$\Phi^*(0,-\upsilon^*) = \sup_{\substack{u \in L_p^n[t_0,T] \\ \upsilon \in L_p^n[t_0,T]}}\left\{-\int\limits_{t_0}^{T}(\upsilon(t)|\upsilon^*(t))dt - \int\limits_{t_0}^{T}f(t,(Au)(t),u(t)+\upsilon(t))dt - \varphi((Au)(t_0),(Au)(T))\right\} =$$

$$= \sup_{u,\upsilon \in L_p^n[t_0,T]}\left\{\int\limits_{t_0}^{T}(u(t)|\upsilon^*(t))dt - \int\limits_{t_0}^{T}(u(t)+\upsilon(t)|\upsilon^*(t))dt - \int\limits_{t_0}^{T}f(t,(Au)(t),u(t)+\upsilon(t))dt - \varphi((Au)(t_0),(Au)(T))\right\} \quad (3)$$

$$= \sup_{u \in L_p^n[t_0,T]}\left\{\int\limits_{t_0}^{T}(u(t)|\upsilon^*(t))dt - \int\limits_{t_0}^{T}f^0(t,(Au)(t),\upsilon^*(t))dt - \varphi((Au)(t_0),(Au)(T))\right\}.$$

Обозначим $J_1(u) = \int\limits_{t_0}^{T}f^0(t,(Au)(t),\upsilon^*(t))dt$, $J_2(u) = \varphi((Au)(t_0),(Au)(T))$. Из (2) и (3) вытекает, что $J_1$ и $J_2$ собственные функционалы. При условии теоремы 1 из леммы 2 вытекает, что функционал $J_1$ непрерывен в точке $\bar{u}(\cdot)$. Также из доказательства теоремы 2.5.1 имеем, что $J_2(u(\cdot))$ непрерывен в точке $\bar{u}(\cdot)$.

По соотношению (2) имеем

$$\int\limits_{t_0}^{T}f(t,(A\bar{u})(t),\bar{u}(t))dt + \varphi((A\bar{u})(t_0),(A\bar{u})(T)) + \Phi^*(0,-\upsilon^*) = 0. \qquad (4)$$

Положив

$$S(u) = \int\limits_{t_0}^{T}f^0(t,(Au)(t),\upsilon^*(t))dt + \varphi((Au)(t_0),(Au)(T))$$

имеем, что $\Phi^*(0,-\upsilon^*) = S^*(\upsilon^*)$. Так как

$$S^*(\upsilon^*) \geq \int\limits_{t_0}^{T}(\bar{u}(t)|\upsilon^*(t))dt - S(\bar{u}), \qquad S(\bar{u}) \leq \int\limits_{t_0}^{T}(\bar{u}(t)|\upsilon^*(t))dt + J(\bar{u}).$$

то отсюда получим, что $S^*(\upsilon^*) \geq \int_{t_0}^{T} (\overline{u}(t)|\upsilon^*(t))dt - S(\overline{u}) \geq -J(\overline{u})$. Поэтому из соотношения (4)

вытекает, что

$$S^*(\upsilon^*) = \int_{t_0}^{T} (\overline{u}(t)|\upsilon^*(t))dt - S(\overline{u}), \qquad S(\overline{u}) = \int_{t_0}^{T} (\overline{u}(t)|\upsilon^*(t))dt + J(\overline{u}).$$

Из второго соотношения вытекает, что

$$\int_{t_0}^{T} f^0(t,(A\overline{u})(t),\upsilon^*(t))dt - \int_{t_0}^{T} f(t,(A\overline{u})(t),\overline{u}(t))dt = \int_{t_0}^{T} (\overline{u}(t)|\upsilon^*(t))dt.$$

Отсюда, используя неравенство Юнга-Фенхеля получим

$$f^0(t,(A\overline{u})(t),\upsilon^*(t)) - f(t,(A\overline{u})(t),\overline{u}(t)) = (\overline{u}(t)|\upsilon^*(t)).$$

Поэтому

$$f^0(t,(A\overline{u})(t),\upsilon^*(t)) = (\overline{u}(t)|\upsilon^*(t)) + f(t,(A\overline{u})(t),\overline{u}(t)).$$

Из равенства $S^*(\upsilon^*) = \int_{t_0}^{T} (\overline{u}(t)|\upsilon^*(t))dt - S(\overline{u})$ вытекает, что $\upsilon^* \in \partial S(\overline{u})$. Из теоремы Моро-Рокафеллара имеем, что $\partial S(\overline{u}) = \partial J_1(\overline{u}) + \partial J_1(\overline{u})$. Поэтому существуют функции $\upsilon_1^* \in \partial J_1(\overline{u})$, $\upsilon_2^* \in \partial J_2(\overline{u})$ такие, что $\upsilon^* = \upsilon_1^* + \upsilon_2^*$. По теореме 1.3.1 $\upsilon_1^* \in \partial J_1(\overline{u})$, где $\upsilon_1^*(\cdot) \in L_q^n[t_0,T]$, принадлежит $\partial J_1(\overline{u})$ в том и только в том случае, когда существует $u^*(t) \in \partial f^0(t,\overline{x}(t),\upsilon^*(t))$, $u^*(\cdot) \in L_1^n[t_0,T]$, такое, что $\upsilon_1^*(\cdot) = (A^*u^*)(t)$.

По теореме 1.3.2 $\upsilon_2^* \in \partial J_2(\overline{u})$, где $\upsilon_2^*(\cdot) \in L_q^n[t_0,T]$, в том и только в том случае, когда существует $(-c_1,-c_2) \in \partial \varphi((A\overline{u})(t_0),(A\overline{u})(T)))$ такое, что

$$\int_{t_0}^{T} (\upsilon_2^*(t)|u(t))dt + ((Au)(t_0)|c_1) + ((Au)(T)|c_2)) = 0$$

при $u(\cdot) \in L_p^n[t_0,T]$, где $\upsilon_2^* = \upsilon^* - \upsilon_1^* = \upsilon^* - (A^*u^*)(\cdot)$. Теорема доказана.

**Следствие 1.** Пусть $f:[t_0,T] \times (R^n \times R^n) \to R_{+\infty}$ нормальный выпуклый интегрант, $\varphi: R^n \times R^n \to R_{+\infty}$ выпуклая функция, $A: L_p^n[t_0,T] \to C^n[t_0,T]$ линейный непрерывный оператор. Для того, чтобы функция $\overline{u}(t)$ среди всех функций $u \in L_p^n[t_0,T]$ минимизировала функционал (1) достаточно, чтобы нашлись $\upsilon^*(\cdot) \in L_q^n[t_0,T]$, $u^*(\cdot) \in L_1^n[t_0,T]$ и $c_1,c_2 \in R^n$ такие, что

1) $(u^*(t),-\upsilon^*(t)) \in \partial f(t,(A\overline{u})(t),\overline{u}(t))$,
2) $(-c_1,-c_2) \in \partial \varphi((A\overline{u})(t_0),(A\overline{u})(T))$,
3) $\int_{t_0}^{T} (\upsilon^*(t) - (A^*u^*)(t)|u(t))dt + ((Au)(t_0)|c_1) + ((Au)(T)|c_2)) = 0$ при $u(\cdot) \in L_p^n[t_0,T]$,

а если выполнено условие леммы 1 при $u_0(\cdot) = \overline{u}(\cdot)$, $A: L_p^n[t_0,T] \to C^n[t_0,T]$ линейный непрерывный компактный оператор, то условия 1)-3) и являются необходимыми.

## 2.7. О необходимых и достаточных условиях минимума для операторных включений

Пусть $f:[t_0,T] \times R^n \to (-\infty,+\infty)$ нормальный выпуклый интегрант, $\varphi: R^n \to (-\infty,+\infty)$ выпуклая функция, $A: L_p^n[t_0,T] \to C^n[t_0,T]$ линейный непрерывный оператор, $1 \leq p < +\infty$. Пусть $t_0 < t_1 < t_2 < \ldots < t_k = T$, $F:[t_0,T] \times R^n \to \text{comp}R^n$ многозначное отображение.

Рассматривается задача минимизации функционала

$$J_0(u) = \varphi((Au)(t_0),(Au)(t_1),\ldots,(Au)(t_k)) + \int_{t_0}^{T} f(t,(Au)(t))dt, \qquad (1)$$

при следующих ограничениях

$$u(t) \in F(t,(Au)(t)), \quad t \in [t_0,T], \quad u(\cdot) \in L_p^n[t_0,T]. \tag{2}$$

Функцию $u(\cdot) \in L_p^n[t_0,T]$, удовлетворяющую (2) назовем решением задачи (2). Пусть $\bar{u}(\cdot) \in L_p^n[t_0,T]$ является решением задачи (2). Будем говорить, что кривая $\bar{u}(t)$, $t_0 \le t \le T$, является решением задачи (1), (2), если $|J_0(\bar{u})| < +\infty$ и справедливо неравенство $J_0(u) \ge J_0(\bar{u})$ среди всех решений $u(\cdot) \in L_p^n[t_0,T]$ задачи (2).

Обозначив $\omega(t,x,z) = \begin{cases} 0, & z \in F(t,x), \\ +\infty, & z \notin F(t,x) \end{cases}$ имеем, что задачи (1),(2) эквивалентны минимизации функционала

$$J(u) = \varphi((Au)(t_0),(Au)(t_1),\ldots,(Au)(t_k)) + \int_{t_0}^{T} f(t,(Au)(t))dt + \int_{t_0}^{T} \omega(t,(Au)(t),u(t))dt$$

среди всех функций $u(\cdot) \in L_p^n[t_0,T]$.

**1. Оператор типа Фредгольма.** Пусть $f:[t_0,T] \times R^n \to R_{+\infty}$ выпуклый нормальный интегрант, $\varphi: R^{(k+2)n} \to R_{+\infty}$ выпуклая функция, отображение $t \to \mathrm{gr}\, F_t(x)$ измеримо на $[t_0,T]$, множество $\mathrm{gr}\, F_t$ замкнуто и выпукло почти для всех $t \in [t_0,T]$, а $F(t,x)$ компактно при всех $t, x$. Отсюда вытекает, что $\omega(t,x,z)$ выпуклый нормальный интегрант на $[t_0,T] \times (R^n \times R^n)$. Положим $\omega^0(t,x,y) = \min\{(z|y) : z \in F(t,x)\}$.

Рассмотрим функционал

$$\Phi(u,\upsilon) = \varphi((Au)(t_0),(Au)(t_1),\ldots,(Au)(t_k)) + \int_{t_0}^{T} f(t,(Au)(t))dt + \int_{t_0}^{T} \omega(t,(Au)(t),u(t)+\upsilon(t))dt,$$

где $\upsilon(\cdot) \in L_p^n[t_0,T]$. Ясно, что $\Phi: L_p^n[t_0,T] \times L_p^n[t_0,T] \to R \cup \{+\infty\}$ и $\Phi(u,0) = J(u)$. Для любого $\upsilon(\cdot) \in L_p^n[t_0,T]$ рассмотрим задачу минимизации

$$\inf_{u \in L_p^n[t_0,T]} \Phi(u,\upsilon).$$

Положим $h(\upsilon) = \inf\limits_{u \in L_p^n[t_0,T]} \Phi(u,\upsilon)$. Из предложения 2.5 [11] вытекает, что $h$ выпуклая функция. Задача (1),(2) называется стабильной, если $h(0)$ конечно и функция $h$ субдифференцируема в нуле.

Положим $\tilde{p} = \begin{cases} p-1, & \text{если } p \in N, \\ p, & \text{если } p \notin N. \end{cases}$ Отметим, что если $a \ge 0, b \ge 0$, то по методу индукции легко проверяется, что $(a+b)^p \le 3^{\tilde{p}}(a^p + b^p)$.

**Лемма 1.** Пусть $F: [t_0,T] \times R^n \to \mathrm{comp}\, R^n$ многозначное отображение, $F(t,x)$ измеримо по $t$, существует такое число $\lambda > 0$, что $\|F(t,x)\| \le \lambda(1+|x|)$ при $x \in R^n$, $\mathrm{gr}\, F_t$ выпукло и замкнуто, $f:[t_0,T] \times R^n \to (-\infty,+\infty]$ нормальный выпуклый интегрант, $\varphi: R^{(k+1)n} \to (-\infty,+\infty]$ выпуклая функция, $A: L_p^n[t_0,T] \to C^n[t_0,T]$ линейный непрерывный оператор, $2\lambda \|A\|(T-t_0)^{\frac{1}{p}} < 1$, $\inf\limits_{u \in L_p^n[t_0,T]} J(u)$ конечен и существует решение $u_0(\cdot) \in L_p^n[t_0,T]$ задачи (2) такое, что функция $f(t,(Au_0)(t)+y)$ суммируема при $y \in R^n$, $|y| < r$, где $r > 0$, а функция $\varphi((Au)(t_0),(Au)(t_1),\ldots,(Au)(t_k))$ непрерывна в точке $u_0(\cdot) \in L_p^n[t_0,T]$. Тогда функция $h$ субдифференцируема в нуле, т.е. задачи (1),(2) стабильны.

**Доказательство.** По условию имеем, что функционал

$$I(u) = \int_{t_0}^{T} f(t,(Au)(t))dt$$

непрерывен в точке $u_0(\cdot) \in L_p^n[t_0,T]$. Из непрерывности функционала $I(u)$ в точке $u_0(\cdot) \in L_p^n[t_0,T]$ получим, что существуют $\alpha_1 > 0$ и $M_1$ такие, что $I(u) \leq M_1$ при $\|u(\cdot)-u_0(\cdot)\|_{L_p^n} < \alpha_1$. Аналогично, существуют $\alpha_2 > 0$ и $M_2$ такие, что $\varphi((Au)(t_0),(Au)(t_1),\ldots,(Au)(t_k)) \leq M_2$ при $\|u(\cdot)-u_0(\cdot)\|_{L_p^n} < \alpha_2$. Обозначим $\alpha = \min\{\alpha_1,\alpha_2\}$.

Так как множество $grF_t$ выпукло, то аналогично доказательству теоремы 3.1.1 [13] (см. раздел 2.2.1) можно показать, что $\rho_x(F(t,y),F(t,x)) \leq 2\lambda|y-x|$ для $y,x \in R^n$ п.в. $t \in [t_0,T]$.

Если $\upsilon(\cdot) \in L_p^n[t_0,T]$, то обозначив $M = 2\lambda$ имеем, что
$$\rho_x(F(t,x)+\upsilon(t),F(t,y)+\upsilon(t)) \leq M|x-y|$$
при $x,y \in R^n$ и $d(u_0(t),F(t,(Au_0)(t))+\upsilon(t)) \leq |\upsilon(t)|$. Тогда по теореме 2.2.1 существует решение $u_\upsilon$ задачи
$$u(t) \in F(t,(Au)(t))-\upsilon(t), \quad t \in [t_0,T], \quad u(\cdot) \in L_p^n[t_0,T],$$
такое, что
$$|u_\upsilon(t)-u_0(t)| \leq |\upsilon(t)|+\frac{\|A\|(\int_{t_0}^T|\upsilon(s)|^p ds)^{\frac{1}{p}}M}{1-\|A\|M(T-t_0)^{\frac{1}{p}}}, \quad |(Au_\upsilon)(t)-(A\overline{u})(t)| \leq \frac{\|A\|(\int_{t_0}^T|\upsilon(s)|^p ds)^{\frac{1}{p}}}{1-\|A\|M(T-t_0)^{\frac{1}{p}}}.$$

Поэтому
$$(\int_{t_0}^T|u_\upsilon(t)-u_0(t)|^p dt)^{\frac{1}{p}} \leq (\int_{t_0}^T(|\upsilon(t)|+\frac{\|A\|(\int_{t_0}^T|\upsilon(s)|^p ds)^{\frac{1}{p}}M}{1-\|A\|M(T-t_0)^{\frac{1}{p}}})^p dt)^{\frac{1}{p}} \leq$$
$$\leq ((\int_{t_0}^T(|\upsilon(t)|)^p dt)^{\frac{1}{p}}+\frac{\|A\|(\int_{t_0}^T|\upsilon(s)|^p ds)^{\frac{1}{p}}M(T-t_0)^{\frac{1}{p}}}{1-\|A\|M(T-t_0)^{\frac{1}{p}}})3^{\frac{\tilde{p}+1}{p}} = (\int_{t_0}^T(|\upsilon(t)|)^p dt)^{\frac{1}{p}}(1+\frac{\|A\|M(T-t_0)^{\frac{1}{p}}}{1-\|A\|M(T-t_0)^{\frac{1}{p}}})3^{\frac{\tilde{p}+1}{p}},$$

где $\tilde{p} = \begin{cases} p-1, \text{ если } p \in N, \\ p, \text{ если } p \notin N. \end{cases}$ Если $\|\upsilon(\cdot)\|_{L_p^n} \leq \alpha(1-\|A\|M(T-t_0)^{\frac{1}{p}})3^{-\frac{\tilde{p}+1}{p}}$, то имеем, что $\|u_\upsilon(t)-u_0(t)\|_{L_p^n} \leq \alpha$. Тогда получим, что
$$h(\upsilon) = \inf_{u \in L_p^n[t_0,T]}\Phi(u,\upsilon) \leq \Phi(u_\upsilon,\upsilon) \leq M_1+M_2$$

при $\upsilon(\cdot) \in L_p^n[t_0,T]$, $\|\upsilon(\cdot)\|_{L_p^n} \leq \alpha(1-\|A\|M)3^{-\frac{\tilde{p}+1}{p}}$. Так как $h(0) = \inf_{u \in L_p^n[t_0,T]}J(u)$ конечно, тогда из предложения 1.2.5 [26] вытекает, что функция $h$ непрерывна в нуле. Тогда из предложения 1.5.2 [26] имеем, что $h$ субдифференцируема в точке нуле. Лемма доказана.

Отметим, что если функция $\varphi: R^{(k+1)n} \to R_{+\infty}$ непрерывная в точке $((Au_0)(t_0),(Au_0)(t_1),\ldots,(Au_0)(t_k))$, то функция $\varphi((Au)(t_0),(Au)(t_1),\ldots,(Au)(t_k))$ непрерывна в точке $u_0(\cdot) \in L_p^n[t_0,T]$.

**Лемма 2.** Пусть $F:[t_0,T]\times R^n \to 2^{R^n}$, отображение $t \to F(t,x)$ измеримо на $[t_0,T]$, отображение $x \to F(t,x)$ замкнуто и выпукло почти при всех $t \in [t_0,T]$, т.е. $grF_t = \{(x,y): y \in F(t,x)\}$ замкнуто и выпукло почти при всех $t \in [t_0,T]$, существует такая функция $\lambda(\cdot) \in L_p^n[t_0,T]$, что $\|F(t,x)\| \leq \lambda(t)(1+|x|)$ при $x \in R^n$, где $\|F(t,x)\| = \sup\{|y|: y \in F(t,x)\}$, $\|\varnothing\| = 0$ и существует решение $u_0(t)$ задачи $u_0(t) \in F(t,(Au_0)(t))$ такое, что $x_0(t) = (Au_0)(t)$ принадлежит

$\mathrm{dom}\, F_t = \{x : F(t,x) \neq \emptyset\}$ вместе с некоторой $\varepsilon$ трубкой, т.е. $\{x : |x_0(t) - x| \leq \varepsilon\} \subset \mathrm{dom}\, F_t$, $f : [t_0, T] \times R^n \to (-\infty, +\infty]$ нормальный выпуклый интегрант, $\varphi : R^{(k+1)n} \to R_{+\infty}$ выпуклая функция, $A : L_p^n[t_0, T] \to C^n[t_0, T]$ линейный непрерывный оператор, $\|A\|(\int_{t_0}^T M^p(s)ds)^{\frac{1}{p}} < 1$, где

$$M(t) = \frac{6\lambda(t)\left(1 + |x_0(t)| + \frac{5\varepsilon}{6}\right)}{\varepsilon},\quad \inf_{u \in L_p^n[t_0,T]} J(u)$$

конечен, функция $f(t, (Au_0)(t) + y)$ суммируема при $y \in R^n$, $|y| < r$, где $r > 0$, а функция $\varphi((Au)(t_0), (Au)(t_1), \ldots, (Au)(t_k))$ непрерывна в точке $u_0(\cdot) \in L_p^n[t_0, T]$. Тогда функция $h$ субдифференцируема в нуле, т.е. задача (1), (2) стабильна.

**Доказательство.** По условию имеем, что функционал

$$I(u) = \int_{t_0}^T f(t, (Au)(t))dt$$

непрерывен в точке $u_0(\cdot) \in L_p^n[t_0, T]$. Из непрерывности функционала $I(u)$ в точке $u_0(\cdot) \in L_p^n[t_0, T]$ получим, что существуют $\alpha_1 > 0$ и $M_1$ такие, что $I(u) \leq M_1$ при $\|u(\cdot) - u_0(\cdot)\|_{L_p^n} < \alpha_1$. Аналогично, существуют $\alpha_2 > 0$ и $M_2$ такие, что $\varphi((Au)(t_0), (Au)(t_1), \ldots, (Au)(t_k)) \leq M_2$ при $\|u(\cdot) - u_0(\cdot)\|_{L_p^n} < \alpha_2$. Обозначим $\alpha = \min\{\alpha_1, \alpha_2\}$.

Рассмотрим множество $S_t = \{x : |x_0(t) - x| \leq \frac{\varepsilon}{2}\}$. Тогда из доказательства леммы 2.1.3 вытекает, что

$$\rho_x(F(t,y), F(t,x)) \leq \frac{6\lambda(t)\left(1 + |x_0(t)| + \frac{5\varepsilon}{6}\right)}{\varepsilon}|y - x|$$

при $x, y \in S_t$. Если $\upsilon(\cdot) \in L_p^n[t_0, T]$, то обозначив $M(t) = \dfrac{6\lambda(t)\left(1 + |x_0(t)| + \frac{5\varepsilon}{6}\right)}{\varepsilon}$ имеем, что

$$\rho_x(F(t,x) + \upsilon(t), F(t,y) + \upsilon(t)) \leq M(t)|x - y|$$

при $x, y \in S_t$. Тогда по теореме 2.2.1 и замечанию 2.2.1 существует решение $u_\upsilon(\cdot) \in L_p^n[t_0, T]$ задачи

$$u(t) \in F(t, (Au)(t)) - \upsilon(t),\quad t \in [t_0, T],\quad u(\cdot) \in L_p^n[t_0, T],$$

такое, что

$$|u_\upsilon(t) - u_0(t)| \leq |\upsilon(t)| + \frac{\|A\|(\int_{t_0}^T |\upsilon(s)|^p ds)^{\frac{1}{p}} M(t)}{1 - \|A\|(\int_{t_0}^T M^p(s)ds)^{\frac{1}{p}}},\quad |(Au_\upsilon)(t) - (Au_0)(t)| \leq \frac{\|A\|(\int_{t_0}^T |\upsilon(s)|^p ds)^{\frac{1}{p}}}{1 - \|A\|(\int_{t_0}^T M^p(s)ds)^{\frac{1}{p}}}.$$

Поэтому

$$\left(\int_{t_0}^T |u_\upsilon(t) - u_0(t)|^p dt\right)^{\frac{1}{p}} \leq \left(\int_{t_0}^T \left(|\upsilon(t)| + \frac{\|A\|(\int_{t_0}^T |\upsilon(s)|^p ds)^{\frac{1}{p}} M(t)}{1 - \|A\|(\int_{t_0}^T M^p(s)ds)^{\frac{1}{p}}}\right)^p dt\right)^{\frac{1}{p}} \leq$$

$$\leq \left(\left(\int_{t_0}^T (|\upsilon(t)|)^p dt\right)^{\frac{1}{p}} + \frac{\|A\|(\int_{t_0}^T |\upsilon(s)|^p ds)^{\frac{1}{p}}(\int_{t_0}^T M^p(s)ds)^{\frac{1}{p}}}{1 - \|A\|(\int_{t_0}^T M^p(s)ds)^{\frac{1}{p}}}\right) 3^{\frac{\tilde{p}+1}{p}} = \left(\int_{t_0}^T (|\upsilon(t)|)^p dt\right)^{\frac{1}{p}} \left(1 + \frac{\|A\|(\int_{t_0}^T M^p(s)ds)^{\frac{1}{p}}}{1 - \|A\|(\int_{t_0}^T M^p(s)ds)^{\frac{1}{p}}}\right) 3^{\frac{\tilde{p}+1}{p}},$$

где $\tilde{p} = \begin{cases} p-1, & \text{если } p \in N, \\ p, & \text{если } p \notin N. \end{cases}$

Если $\|\upsilon(\cdot)\|_{L_p^n} \leq 3^{-\frac{\tilde{p}+1}{p}} \alpha (1 - \|A\| (\int_{t_0}^T M^p(s)ds)^{\frac{1}{p}})$, то имеем, что $\|u_\upsilon(t) - u_0(t)\|_{L_p^n} \leq \alpha$.

Тогда получим, что
$$h(\upsilon) = \inf_{u \in L_p^n[t_0,T]} \Phi(u,\upsilon) \leq \Phi(u_\upsilon, \upsilon) \leq M_1 + M_2$$

при $\upsilon(\cdot) \in L_p^n[t_0,T]$, $\|\upsilon(\cdot)\|_{L_p^n} \leq 3^{-\frac{\tilde{p}+1}{p}} \alpha (1 - \|A\|(\int_{t_0}^T M^p(s)ds)^{\frac{1}{p}})$. Так как $h(0) = \inf_{u \in L_p^n[t_0,T]} J(u)$ конечен, тогда из предложения 1.2.5 [26] следует, что функция $h$ непрерывна в нуле. Тогда из предложения 1.5.2 [26] имеем, что $h$ субдифференцируема в точке нуле. Лемма доказана.

Пусть $\upsilon, w \in R^n$. Положим $\omega^0(t,z,\upsilon) = \inf_w \{(w|\upsilon) + \omega(t,z,w)\}$.

**Теорема 1.** Пусть $f$ нормальный выпуклый интегрант на $[t_0,T] \times R^n$, $\varphi: R^{(k+1)n}$ выпуклая функция, $A: L_p^n[t_0,T] \to C^n[t_0,T]$ линейный непрерывный оператор. Для того, чтобы функция $\bar{u}(\cdot) \in L_p^n[t_0,T]$ среди всех решений задачи (2) минимизировала функционал (1) достаточно, чтобы нашлись $u_1^*(\cdot), u_2^*(\cdot) \in L_1^n[t_0,T]$ и $(c_0, c_1, \ldots c_k) \in R^{(k+1)n}$ такие, что

1) $u_1^*(t) \in \partial f(t,(A\bar{u})(t))$,
2) $u_2^*(t) \in \partial \omega^0(t,(A\bar{u})(t),\upsilon^*(t))$,
3) $(-c_0,-c_1,\ldots,-c_k) \in \partial \varphi((A\bar{u})(t_0),(A\bar{u})(t_1),\ldots,(A\bar{u})(t_k))$,
4) $\omega^0(t,(A\bar{u})(t),\upsilon^*(t)) = (\bar{u}(t)|\upsilon^*(t)) + \omega(t,(A\bar{u})(t),\bar{u}(t))$,
5) $\int_{t_0}^T (-(A^*u_1^*)(t) - (A^*u_2^*)(t)|u(t))dt + \sum_{i=0}^k ((Au)(t_i)|c_i) = 0$ при $u(\cdot) \in L_p^n[t_0,T]$,

а если выполнено условие леммы 1 при $u_0(\cdot) = \bar{u}(\cdot)$ и $A: L_p^n[t_0,T] \to C^n[t_0,T]$ линейный непрерывный компактный оператор, $2\lambda \|A\|(T-t_0)^{\frac{1}{p}} < 1$, то условия 1)-5) и являются необходимыми.

**Доказательство.** Достаточность теоремы непосредственно проверяется.

**Необходимость.** Из леммы 1 вытекает, что $h$ субдифференцируема в точке нуле. Поэтому из замечания 3.2.3 и из предложения 3.2.4 [26] вытекает, что все решения $\bar{u}(\cdot)$ задачи $\inf\{J(u): u \in L_p^n[t_0,T]\}$ и все решения $\bar{u}(\cdot)$ задачи $\sup_{\upsilon^* \in L_q^n[t_0,T]} \{-\Phi^*(0,\upsilon^*)\}$ связаны экстремальным соотношением

$$\Phi(\bar{u},0) + \Phi^*(0,-\upsilon^*) = 0 \qquad (3)$$

По определению
$$\Phi^*(0,-\upsilon^*) = \sup_{u,\upsilon \in L_p^n[t_0,T]} \{-\int_{t_0}^T (\upsilon(t)|\upsilon^*(t))dt - \int_{t_0}^T f(t,(Au)(t))dt - \varphi((Au)(t_0),(Au)(t_1),\ldots,(Au)(t_k)) -$$
$$-\int_{t_0}^T \omega(t,(Au)(t),u(t)+\upsilon(t))dt\} = \sup_{u,\upsilon \in L_p^n[t_0,T]} \{\int_{t_0}^T (u(t)|\upsilon^*(t))dt -$$
$$-\int_{t_0}^T (u(t)+\upsilon(t)|\upsilon^*(t))dt - \int_{t_0}^T f(t,(Au)(t))dt - \varphi((Au)(t_0),(Au)(t_1),\ldots,(Au)(t_k)) - \qquad (4)$$
$$-\int_{t_0}^T \omega(t,(Au)(t),u(t)+\upsilon(t))dt\} = \sup_{u \in L_p^n[t_0,T]} \{\int_{t_0}^T (u(t)|\upsilon^*(t))dt - \int_{t_0}^T \omega^0(t,(Au)(t),\upsilon^*(t))dt -$$
$$-\int_{t_0}^T f(t,(Au)(t))dt - \varphi((Au)(t_0),(Au)(t_1),\ldots,(Au)(t_k))\},$$

где

$$\omega^0(t,(Au)(t),\upsilon^*(t)) = \inf_{y \in R^n}\{(y|\upsilon^*(t)) + \omega(t,(Au)(t),y)\} = \inf_{y \in R^n}\{(y|\upsilon^*(t)) : y \in F(t,(Au)(t))\}.$$

Обозначая

$$J_1(u) = \int_{t_0}^{T} f(t,(Au)(t))dt + \int_{t_0}^{T} \omega^0(t,(Au)(t),\upsilon^*(t))dt + \varphi((Au)(t_0),(Au)(t_1),\ldots,(Au)(t_k)),$$

получим, что $\Phi^*(0,-\upsilon^*) = J_1^*(\upsilon^*)$. Из соотношения (3) вытекает, что $J(\overline{u}) + J_1^*(\upsilon^*) = 0$. По определению $J_1^*(\upsilon^*) \geq \int_{t_0}^{T}(\overline{u}(t)|\upsilon^*(t))dt - J_1(\overline{u})$ и $J_1(\overline{u}) \leq \int_{t_0}^{T}(\overline{u}(t)|\upsilon^*(t))dt + J(\overline{u})$. Тогда из соотношения (4) вытекает, что

$$J_1^*(\upsilon^*) = \int_{t_0}^{T}(\overline{u}(t)|\upsilon^*(t))dt - J_1(\overline{u}), \qquad J_1(\overline{u}) = \int_{t_0}^{T}(\overline{u}(t)|\upsilon^*(t))dt + J(\overline{u}).$$

Из второго соотношения имеем, что

$$\int_{t_0}^{T} \omega^0(t,(A\overline{u})(t),\upsilon^*(t))dt = \int_{t_0}^{T}(\overline{u}(t)|\upsilon^*(t))dt + \int_{t_0}^{T} \omega(t,(A\overline{u})(t),\overline{u}(t))dt.$$

Отсюда, используя неравенства Юнга-Фенхеля получим

$$\omega^0(t,(A\overline{u})(t),\upsilon^*(t)) = (\overline{u}(t)|\upsilon^*(t)) + \omega(t,(A\overline{u})(t),\overline{u}(t)).$$

Из равенства $J_1^*(\upsilon^*) = \int_{t_0}^{T}(\overline{u}(t)|\upsilon^*(t))dt - J_1(\overline{u})$ вытекает, что $\upsilon^* \in \partial J_1(\overline{u})$.

Обозначив

$$I_1(u) = \int_{t_0}^{T} f(t,(Au)(t))dt, \quad I_2(u) = \int_{t_0}^{T} \omega^0(t,(A\overline{u})(t),\upsilon^*(t))dt,$$
$$I_3(u) = \varphi((Au)(t_0),(Au)(t_1),\ldots,(Au)(t_k))$$

имеем, что $J_1(u) = I_1(u) + I_2(u) + I_3(u)$. Так как существует число $\lambda > 0$ такое, что $\|F(t,x)\| \leq \lambda(1+|x|)$ при $x \in R^n$, $\text{gr } F_t$ выпукло и замкнуто, то имеем, что

$$\omega^0(t,y,\upsilon^*(t)) = \inf_{u \in R^n}\{(u|\upsilon^*(t)) + \omega(t,y,u)\} = \inf_{u \in R^n}\{(u|\upsilon^*(t)) : u \in F(t,y)\} \leq$$
$$\leq \inf_{u \in R^n}\{|u||\upsilon^*(t)| : u \in F(t,y)\} \leq \lambda(1+|y|)|\upsilon^*(t)|$$

при $y \in R^n$. Отсюда вытекает, что $I_2(u)$ непрерывен в точке $\overline{u}(\cdot)$. Тогда из теоремы Моро-Рокафеллара имеем, что $\partial J_1(\overline{u}) = \partial I_1(\overline{u}) + \partial I_2(\overline{u}) + \partial I_3(\overline{u})$. Так как $\upsilon^* \in \partial J_1(\overline{u})$, то существуют функции $\upsilon_1^* \in \partial I_1(\overline{u})$, $\upsilon_2^* \in \partial I_2(\overline{u})$, $\upsilon_3^* \in \partial I_3(\overline{u})$ такие, что $\upsilon^* = \upsilon_1^* + \upsilon_2^* + \upsilon_3^*$.

По теореме 2.4.1 $\upsilon_1^*$, где $\upsilon_1^*(\cdot) \in L_q^n[t_0,T]$, принадлежит в $\partial I_1(\overline{u})$ в том и только в том случае, когда существует функция $u_1^*(t) \in \partial f(t,\overline{x}(t))$, $u_1^*(\cdot) \in L_1^n[t_0,T]$, такая, что $\upsilon_1^*(\cdot) = (A^*u_1^*)(t)$.

По теореме 2.4.1 $\upsilon_2^*$, где $\upsilon_2^*(\cdot) \in L_q^n[t_0,T]$, принадлежит в $\partial I_2(\overline{u})$ в том и только в том случае, когда существует функция $u_2^*(t) \in \partial \omega^0(t,(A\overline{u})(t),\upsilon^*(t))$, $u_2^*(\cdot) \in L_1^n[t_0,T]$, такое, что $\upsilon_2^*(\cdot) = (A^*u_2^*)(t)$.

По теореме 2.5.1 $\upsilon_3^* \in \partial I_3(\overline{u})$, где $\upsilon_3^*(\cdot) \in L_q^n[t_0,T]$, в том и только в том случае, когда существует вектор $(-c_0,-c_1,\ldots,-c_k) \in \partial\varphi((A\overline{u})(t_0),(A\overline{u})(t_1),\ldots,(A\overline{u})(t_k))$, где $(c_0,c_1,\ldots c_k) \in R^{(k+1)n}$, такой, что

$$\int_{t_0}^{T}(\upsilon_3^*(t)|u(t))dt + \sum_{i=0}^{k}((Au)(t_i)|c_i) = 0$$

при $u(\cdot) \in L_p^n[t_0,T]$. Кроме того имеем, что $(A^*u_1^*)(t) + (A^*u_2^*)(t) + \upsilon_3^*(\cdot) = 0$. Поэтому $\upsilon_3^*(\cdot) = -(A^*u_1^*)(t) - (A^*u_2^*)(t)$. Теорема доказана.

Отметим, что если выполнено условие леммы 2, то условия 1)-5) также являются необходимыми.

**2. Оператор типа Вольтера.**

**Лемма 3.** Пусть $F:[t_0,T]\times R^n \to 2^{R^n}$, отображение $t \to F(t,x)$ измеримо на $[t_0,T]$, отображение $x \to F(t,x)$ замкнуто и выпукло почти при всех $t \in [t_0,T]$, т.е. $\mathrm{gr}F_t = \{(x,y): y \in F(t,x)\}$ замкнуто и выпукло почти при всех $t \in [t_0,T]$, существует такая суммируемая функция $\lambda(t)$, что $\|F(t,x)\| \le \lambda(t)(1+|x|)$ при $x \in R^n$, где $\|F(t,x)\| = \sup\{|y|: y \in F(t,x)\}$, $\|\varnothing\|=0$ и существует решение $u_0(t)$ задачи $u_0(t) \in F(t,(Au_0)(t))$ такое, что $x_0(t) = (Au_0)(t)$ принадлежит $\mathrm{dom}\,F_t = \{x: F(t,x) \ne \varnothing\}$ вместе с некоторой $\varepsilon$ трубкой, т.е. $\{x: |x_0(t)-x| \le \varepsilon\} \subset \mathrm{dom}\,F_t$, $f:[t_0,T]\times R^n \to (-\infty,+\infty]$ нормальный выпуклый интегрант, $\varphi: R^{(k+1)n} \to R_{+\infty}$ выпуклая функция, $A: L_1^n[t_0,T] \to C^n[t_0,T]$ линейный непрерывный оператор, $|(Au)(t)| \le L\int_{t_0}^t |u(s)|ds$ при $u(\cdot) \in L_1^n[t_0,T]$, где $L>0$. Допустим, что $\inf_{u \in L_1^n[t_0,T]} J(u)$ конечен, функция $f(t,(Au_0)(t)+y)$ суммируема при $y \in R^n$, $|y|<r$, где $r>0$, а функция $\varphi((Au)(t_0),(Au)(t_1),\ldots,(Au)(t_k))$ непрерывна в точке $u_0(\cdot) \in L_1^n[t_0,T]$. Тогда функция $h$ субдифференцируема в нуле, т.е. задача (1),(2) стабильна.

**Доказательство.** По условию имеем, что функционал
$$I(u) = \int_{t_0}^T f(t,(Au)(t))dt$$
непрерывен в точке $u_0(\cdot) \in L_1^n[t_0,T]$. Из непрерывности функционала $I(u)$ в точке $u_0(\cdot) \in L_1^n[t_0,T]$ получим, что существуют $\alpha_1 > 0$ и $M_1$ такие, что $I(u) \le M_1$ при $\|u(\cdot)-u_0(\cdot)\|_{L_1^n} < \alpha_1$. Аналогично, существуют $\alpha_2 > 0$ и $M_2$ такие, что $\varphi((Au)(t_0),(Au)(t_1),\ldots,(Au)(t_k)) \le M_2$ при $\|u(\cdot)-u_0(\cdot)\|_{L_1^n} < \alpha_2$. Обозначим $\alpha = \min\{\alpha_1,\alpha_2\}$.

Рассмотрим множество $S_t = \{x: |x_0(t)-x| \le \tfrac{\varepsilon}{2}\}$. Тогда из доказательства леммы 2.1.3 вытекает, что
$$\rho_x(F(t,y),F(t,x)) \le \frac{6\lambda(t)\left(1+|x_0(t)|+\dfrac{5\varepsilon}{6}\right)}{\varepsilon}|y-x|$$
при $x,y \in S_t$. Если $\upsilon(\cdot) \in L_1^n[t_0,T]$, то обозначив $M(t) = \dfrac{6\lambda(t)\left(1+|x_0(t)|+\dfrac{5\varepsilon}{6}\right)}{\varepsilon}$ имеем, что
$$\rho_x(F(t,x)+\upsilon(t),F(t,y)+\upsilon(t)) \le M(t)|x-y|$$
при $x,y \in S_t$. Тогда по теореме 2.1.1 и замечанию 2.1.1 существует решение $u_\upsilon(\cdot) \in L_1^n[t_0,T]$ задачи
$$u(t) \in F(t,(Au)(t))-\upsilon(t), \quad t \in [t_0,T], \quad u(\cdot) \in L_1^n[t_0,T],$$
такое, что
$$|(Au_\upsilon)(t)-(Au_0)(t)| \le L\int_{t_0}^t e^{m(t)-m(\tau)}|\upsilon(\tau)|d\tau, \quad |u_\upsilon(t)-u_0(t)| \le |\upsilon(t)|+LM(t)\int_{t_0}^t e^{m(t)-m(\tau)}|\upsilon(\tau)|d\tau$$
при $t \in [t_0,T]$, где $(Au_\upsilon)(\cdot) \in C^n[t_0,T]$, $m(t) = L\int_{t_0}^t M(s)ds$.

Поэтому
$$\int_{t_0}^T |u_\upsilon(t)-u_0(t)|dt \le \left(\int_{t_0}^T (|\upsilon(t)|+LM(t)\int_{t_0}^t e^{m(t)-m(\tau)}|\upsilon(\tau)|d\tau\right)dt \le$$
$$\le \left(\int_{t_0}^T |\upsilon(t)|dt + \int_{t_0}^T LM(t)\int_{t_0}^t e^{m(t)-m(\tau)}|\upsilon(\tau)|d\tau dt\right) \le \left(\int_{t_0}^T |\upsilon(t)|dt + Le^{m(T)}\int_{t_0}^T |\upsilon(t)|dt\int_{t_0}^T M(t)dt\right) \le$$

$$\leq (\int_{t_0}^T |\upsilon(t)|dt(1+Le^{m(T)}\int_{t_0}^T M(t)dt).$$

Если $\|\upsilon(\cdot)\|_{L_1^n} \leq \dfrac{\alpha}{(1+Le^{m(T)}\int_{t_0}^T M(t)dt)}$, то имеем, что $\|u_\upsilon(\cdot) - u_0(\cdot)\|_{L_1^n} \leq \alpha$.

Тогда получим, что
$$h(\upsilon) = \inf_{u \in L_1^n[t_0,T]} \Phi(u,\upsilon) \leq \Phi(u_\upsilon,\upsilon) \leq M_1 + M_2$$

при $\upsilon(\cdot) \in L_1^n[t_0,T]$, $\|\upsilon(\cdot)\|_{L_1^n} \leq \dfrac{\alpha}{(1+Le^{m(T)}\int_{t_0}^T M(t)dt)}$. Так как $h(0) = \inf_{u \in L_1^n[t_0,T]} J(u)$ конечен, тогда из

предложения 1.2.5 [26] следует, что функция h непрерывна в нуле. Тогда из предложения 1.5.2 [26] имеем, что h субдифференцируема в точке нуле. Лемма доказана.

**Замечание 1.** В лемме 2 и 3 считаем, что $|(Au)(t) - (Au_0)(t)| \leq \varepsilon$ при $\|u(\cdot) - u_0(\cdot)\|_{L_1^n} \leq \alpha$.

Отметим, что если выполнено условие леммы 3 при $u_0(\cdot) = \bar{u}(\cdot)$ и $p = 1$, то теорема 1 также верна.

### 2.8. Невыпуклая экстремальная задача для операторных включений

Пусть $f: [t_0,T] \times R^n \times R^n \to R_{+\infty}$ нормальный интегрант, $\varphi: R^n \times R^n \to R_{+\infty}$ функция, $t_0 < T$ и $F: [t_0,T] \times R^n \to \text{comp} R^n$ многозначное отображение.

Рассматривается задача минимизации функционала
$$J(u) = \varphi((Au)(t_0),(Au)(T)) + \int_{t_0}^T f(t, Au(t), u(t))dt, \tag{1}$$

при следующих ограничениях
$$u(t) \in F(t, Au(t)), \quad t \in [t_0,T], \ u(\cdot) \in L_p^n[t_0,T], \tag{2}$$

где $A: L_p^n[t_0,T] \to C^n[t_0,T]$.

Функцию $u(\cdot) \in L_p^n[t_0,T]$, удовлетворяющую (2) назовем решением задачи (2).

Пусть $\bar{u}(\cdot) \in L_p^n[t_0,T]$ является решением задачи (2). Будем говорить, что кривая $\bar{u}(t)$, $t_0 \leq t \leq T$, является решением задачи (1), (2), если $|J(\bar{u})| < +\infty$ и справедливо неравенство $J(u) \geq J(\bar{u})$ среди всех решений $u(\cdot) \in L_p^n[t_0,T]$ задачи (2).

Положим $\psi(s,x,y) = \inf\{|z-y| : z \in F(s,x)\}$ и рассмотрим минимизацию функционала
$$J_r(u) = \varphi((Au)(t_0),(Au)(T)) + \int_{t_0}^T f(t,(Au)(t),u(t))dt + r(\int_{t_0}^T \psi^p(t,(Au)(t),u(t))dt)^{\frac{1}{p}}$$

среди всех функций $u(\cdot) \in L_p^n[t_0,T]$.

Решение задачи (1), (2) обозначим через $\bar{u}(t)$. Обозначим
$$F_0(u) = \int_{t_0}^T \psi(t,(Au)(t),u(t))dt.$$

**Теорема 1.** Если $\bar{u}(\cdot) \in L_p^n[t_0,T]$ является решением задачи (1),(2), $F: [t_0,T] \times R^n \to \text{comp } R^n \cup \{\varnothing\}$ и $t \to F(t,x)$ измеримо по $t$, $A: L_p^n[t_0,T] \to C^n[t_0,T]$ оператор, существуют $k(\cdot) \in L_1[t_0,T]$, $M(\cdot) \in L_p[t_0,T]$, $k_1 > 0$, $k_2 > 0$ и $\alpha > 0$ такие, что $B((A\bar{u})(t),\alpha) \subset \text{dom } F_t = \{x \in R^n : F(t,x) \neq \varnothing\}$ при $t \in [t_0,T]$ и
$$|\varphi(z_2,u_2) - \varphi(z_1,u_1)| \leq k_2(|z_2-z_1|+|u_2-u_1|), \quad |f(t,x_1,y_1) - f(t,x_2,y_2)| \leq k(t)|x_1-x_2|+k_1|y_1-y_2|,$$
$$\rho_X(F(t,x_1),F(t,x_2)) \leq M(t)|x_1-x_2|$$

при $x_1, x_2 \in B((A\bar{u})(t), \alpha)$, $z_1, z_2 \in B((A\bar{u})(t_0), \alpha)$, $u_1, u_2 \in B((A\bar{u})(T), \alpha)$, $y_1, y_2 \in R^n$, существует $L > 0$ такое, что $|(Au)(t) - (A\upsilon)(t)| \leq L\|u - \upsilon\|_{L_p^n}$ при $u(\cdot), \upsilon(\cdot) \in L_p^n[t_0, T]$ и $L(\int_{t_0}^T M^p(s)ds)^{\frac{1}{p}} < 1$, то существует число $r_0 > 0$ такое, что $\bar{u}(t)$ минимизирует функционал $J_r(u)$ в $D$ при $r \geq r_0$, где $D = \left\{u(\cdot) \in L_p^n[t_0, T] : \|u(\cdot) - \bar{u}(\cdot)\|_{L_p^n[t_0, T]} \leq \frac{\alpha}{\beta}\right\}$,

$$\beta \geq \frac{L3^{\frac{2\tilde{p}+2}{p}}(L\int_{t_0}^T(M^p(s)ds)^{\frac{1}{p}} + 1)}{1 - L(\int_{t_0}^T M^p(s)ds)^{\frac{1}{p}}} + L + 1, \quad \tilde{p} = \begin{cases} p-1, \text{ если } p \in N, \\ p, \text{ если } p \notin N. \end{cases}$$

**Доказательство.** Положим

$$C = \left\{u(\cdot) \in L_p^n[t_0, T] : F_0(u) = \int_{t_0}^T \psi(t, (Au)(t), u(t))dt = 0\right\},$$

Обозначим $x(t) = (Au)(t)$. Покажем, что существует $\nu > 0$ такое, что

$$d_C(u) \leq \nu(\int_{t_0}^T \psi^p(t, x(t), u(t))dt)^{\frac{1}{p}}$$

при $u(\cdot) \in D$, где $d_C(u) = \inf\{\|\upsilon - u\|_{L_p^n} : \upsilon \in C\}$.

Пусть $u(\cdot) \in D$. Положив в теореме 2.2.2 $\rho(t) = \psi(t, x(t), u(t))$ имеем, что существует решение $u_0(\cdot) \in L_p^n[t_0, T]$ задачи $u(t) \in F(t, (Au)(t))$ такое, что

$$|u(t) - u_0(t)| \leq \rho(t) + \frac{L(\int_{t_0}^T \rho^p(s)ds)^{\frac{1}{p}} M(t)}{1 - L(\int_{t_0}^T M^p(s)ds)^{\frac{1}{p}}}, \qquad |(Au)(t) - (Au_0)(t)| \leq \frac{L(\int_{t_0}^T \rho^p(s)ds)^{\frac{1}{p}}}{1 - L(\int_{t_0}^T M^p(s)ds)^{\frac{1}{p}}}.$$

Покажем, что $\dfrac{L(\int_{t_0}^T \rho^p(s)ds)^{\frac{1}{p}}}{1 - L(\int_{t_0}^T M^p(s)ds)^{\frac{1}{p}}} \leq \alpha - \dfrac{\alpha}{\beta}$. Действительно

$$\frac{L(\int_{t_0}^T \rho^p(s)ds)^{\frac{1}{p}}}{1 - L(\int_{t_0}^T M^p(s)ds)^{\frac{1}{p}}} = \frac{L(\int_{t_0}^T \psi(s, x(s), u(s))^p ds)^{\frac{1}{p}}}{1 - L(\int_{t_0}^T M^p(s)ds)^{\frac{1}{p}}} = \frac{L(\int_{t_0}^T (\psi(\tau, x(\tau), u(\tau)) - \psi(\tau, \bar{x}(\tau), \bar{u}(\tau)))^p d\tau)^{\frac{1}{p}}}{1 - L(\int_{t_0}^T M^p(s)ds)^{\frac{1}{p}}} \leq$$

$$\leq \frac{L(\int_{t_0}^T (M(s)|(Au)(s) - (A\bar{u})(s)| + |u(s) - \bar{u}(s)|)^p ds)^{\frac{1}{p}}}{1 - L(\int_{t_0}^T M^p(s)ds)^{\frac{1}{p}}} \leq$$

$$\leq \frac{L(\int_{t_0}^T (M(s)L\|u - \bar{u}\|_{L_p^n} + |u(s) - \bar{u}(s)|)^p ds)^{\frac{1}{p}}}{1 - L(\int_{t_0}^T M^p(s)ds)^{\frac{1}{p}}} \leq$$

$$\leq \frac{3^{\frac{\tilde{p}+1}{p}} L \|u-\bar{u}\|_{L_p^n} (L(\int_{t_0}^T M^p(s)ds)^{\frac{1}{p}} + 1)}{1 - L(\int_{t_0}^T M^p(s)ds)^{\frac{1}{p}}} \leq \frac{L 3^{\frac{\tilde{p}+1}{p}} (L(\int_{t_0}^T M^p(s)ds)^{\frac{1}{p}} + 1)}{1 - L(\int_{t_0}^T M^p(s)ds)^{\frac{1}{p}}} \frac{\alpha}{\beta} \leq \alpha - \frac{\alpha}{\beta}.$$

Кроме того получим, что

$$(\int_{t_0}^T |u(t)-u_0(t)|^p dt)^{\frac{1}{p}} \leq (\int_{t_0}^T (\rho(t) + \frac{L(\int_{t_0}^T \rho^p(s)ds)^{\frac{1}{p}} M(t)}{1 - L(\int_{t_0}^T M^p(s)ds)^{\frac{1}{p}}})^p dt)^{\frac{1}{p}} \leq ((\int_{t_0}^T \rho^p(t)dt)^{\frac{1}{p}} +$$

$$+ \frac{L(\int_{t_0}^T \rho^p(s)ds)^{\frac{1}{p}}}{1 - L(\int_{t_0}^T M^p(s)ds)^{\frac{1}{p}}} (\int_{t_0}^T M^p(t)dt)^{\frac{1}{p}}) 3^{\frac{\tilde{p}+1}{p}} = (\int_{t_0}^T \rho^p(t)dt)^{\frac{1}{p}} (1 + \frac{L(\int_{t_0}^T M^p(t)dt)^{\frac{1}{p}}}{1 - L(\int_{t_0}^T M^p(s)ds)^{\frac{1}{p}}}) 3^{\frac{\tilde{p}+1}{p}} =$$

$$= (\int_{t_0}^T \psi^p(t,(Au)(t),u(t))dt)^{\frac{1}{p}} (1 + \frac{L(\int_{t_0}^T M^p(t)dt)^{\frac{1}{p}}}{1 - L(\int_{t_0}^T M^p(s)ds)^{\frac{1}{p}}}) 3^{\frac{\tilde{p}+1}{p}},$$

т.е.

$$(\int_{t_0}^T |u(t)-u_0(t)|^p dt)^{\frac{1}{p}} \leq (\int_{t_0}^T \psi^p(t,(Au)(t),u(t))dt)^{\frac{1}{p}} (1 + \frac{L(\int_{t_0}^T M^p(t)dt)^{\frac{1}{p}}}{1 - L(\int_{t_0}^T M^p(s)ds)^{\frac{1}{p}}}) 3^{\frac{\tilde{p}+1}{p}}.$$

Поэтому $d_c(u) \leq \|u_0(\cdot) - u(\cdot)\|_{L_p^n} \leq \nu (\int_{t_0}^T \psi^p(t,(Au)(t),u(t))dt)^{\frac{1}{p}}$ при $u(\cdot) \in D$, где

$$\nu = \frac{3^{\frac{\tilde{p}+1}{p}}}{1 - L(\int_{t_0}^T M^p(s)ds)^{\frac{1}{p}}}.$$

Пусть $u(\cdot) \in D$. Так как $\|u(\cdot) - \bar{u}(\cdot)\|_{L_p^n[t_0,T]} \leq \frac{\alpha}{\beta}$, то имеем, что

$$|(Au)(t) - (A\bar{u})(t)| \leq L\|u(\cdot) - \bar{u}(\cdot)\|_{L_p^n[t_0,T]} \leq \frac{L\alpha}{\beta} \leq \alpha.$$

Также имеем, что

$$|(Au)(t) - (Au_0)(t)| \leq L\|u(\cdot) - u_0(\cdot)\|_{L_p^n[t_0,T]} \leq$$

$$\leq L 3^{\frac{\tilde{p}+1}{p}} (\int_{t_0}^T \psi^p(t,(Au)(t),u(t))dt)^{\frac{1}{p}} \frac{1}{1 - L(\int_{t_0}^T M^p(s)ds)^{\frac{1}{p}}}) \leq$$

$$\leq L 3^{\frac{\tilde{p}+1}{p}} \int_{t_0}^T (M(t)|(Au)(t) - (A\bar{u})(t)| + |u(t) - \bar{u}(t)|)^p dt)^{\frac{1}{p}} \frac{1}{1 - L(\int_{t_0}^T M^p(s)ds)^{\frac{1}{p}}} \leq$$

$$\leq L3^{\frac{\tilde{p}+1}{p}} \int_{t_0}^{T}(M(t)L\|u-\overline{u}\|_{L_p^n}+|u(t)-\overline{u}(t)|)^p dt)^{\frac{1}{p}} \frac{1}{1-L(\int_{t_0}^{T}M^p(s)ds)^{\frac{1}{p}}} \leq$$

$$\leq L3^{\frac{\tilde{p}+1}{p}}\|u-\overline{u}\|_{L_p^n}(L(\int_{t_0}^{T}M^p(s)ds)^{\frac{1}{p}}+1)3^{\frac{\tilde{p}+1}{p}}\frac{1}{1-L(\int_{t_0}^{T}M^p(s)ds)^{\frac{1}{p}}} \leq$$

$$\leq L3^{\frac{2\tilde{p}+2}{p}}(L(\int_{t_0}^{T}M^p(s)ds)^{\frac{1}{p}}+1)\frac{1}{1-L(\int_{t_0}^{T}M^p(s)ds)^{\frac{1}{p}}})\frac{\alpha}{\beta} = \frac{L3^{\frac{2\tilde{p}+2}{p}}(L(\int_{t_0}^{T}M^p(s)ds)^{\frac{1}{p}}+1)}{1-L(\int_{t_0}^{T}M^p(s)ds)^{\frac{1}{p}}}\frac{\alpha}{\beta}.$$

Поэтому
$$|(A\overline{u})(t)-(Au_0)(t)| \leq |(A\overline{u})(t)-(Au)(t)|+|(Au)(t)-(Au_0)(t)| \leq$$

$$\leq \frac{L3^{\frac{2\tilde{p}+2}{p}}(L(\int_{t_0}^{T}M^p(s)ds)^{\frac{1}{p}}+1)}{1-L(\int_{t_0}^{T}M^p(s)ds)^{\frac{1}{p}}}\frac{\alpha}{\beta}+|(A\overline{u})(t)-(Au)(t)| \leq$$

$$\leq \frac{L3^{\frac{2\tilde{p}+2}{p}}(L(\int_{t_0}^{T}M^p(s)ds)^{\frac{1}{p}}+1)}{1-L(\int_{t_0}^{T}M^p(s)ds)^{\frac{1}{p}}}\frac{\alpha}{\beta}+\frac{L\alpha}{\beta} = (\frac{L3^{\frac{2\tilde{p}+2}{p}}(L(\int_{t_0}^{T}M^p(s)ds)^{\frac{1}{p}}+1)}{1-L(\int_{t_0}^{T}M^p(s)ds)^{\frac{1}{p}}}+L)\frac{\alpha}{\beta} \leq \alpha.$$

Получим $x_0(t)=(Au_0)(t) \in B(\overline{x}(t),\alpha)$. Также имеем,

$$\|u(\cdot)-u_0(\cdot)\|_{L_p^n[t_0,T]} \leq \frac{3^{\frac{2\tilde{p}+2}{p}}(L(\int_{t_0}^{T}M^p(s)ds)^{\frac{1}{p}}+1)}{1-L(\int_{t_0}^{T}M^p(s)ds)^{\frac{1}{p}}}\frac{\alpha}{\beta}.$$ Тогда имеем, что

$$\|u_0(\cdot)-\overline{u}(\cdot)\|_{L_p^n[t_0,T]} \leq \|u_0(\cdot)-u(\cdot)\|+\|u(\cdot)-\overline{u}(\cdot)\| \leq \frac{3^{\frac{2\tilde{p}+2}{p}}(L(\int_{t_0}^{T}M^p(s)ds)^{\frac{1}{p}}+1)}{1-L(\int_{t_0}^{T}M^p(s)ds)^{\frac{1}{p}}}\frac{\alpha}{\beta}+\frac{\alpha}{\beta} \leq \alpha.$$

Пусть $u(\cdot), u_1(\cdot) \in L_p^n$, такие, что $(Au)(t), (Au_1)(t) \in B(\overline{x}(t),\alpha)$. Тогда имеем, что
$$|J(u)-J(u_1)| \leq |\varphi((Au)(t_0),(Au)(T))-\varphi((Au_1)(t_0),(Au_1)(T))|+$$

$$+\left|\int_{t_0}^{T}f(t,(Au)(t),u(t))dt-\int_{t_0}^{T}f(t,(Au_1)(t),u_1(t))dt\right| \leq$$

$$\leq k_2(|(Au)(t_0))-(Au_1)(t_0)|+|(Au)(T)-(Au_1)(T)|+$$

$$+\int_{t_0}^{T}(k(t)|(Au)(t)-(Au_1)(t)|+k_1|u(t)-u_1(t)|)dt \leq$$

$$\leq (L\int_{t_0}^{T}k(t)dt+k_1(T-t_0)^{\frac{1}{q}})\|u(\cdot)-u_1(\cdot)\|_{L_p^n}+2Lk_2\|u-u_1\|_{L_p^n}=$$

$$=(L\int_{t_0}^{T}k(t)dt+k_1(T-t_0)^{\frac{1}{q}}+2Lk_2)\|u(\cdot)-u_1(\cdot)\|_{L_p^n},$$

где $\dfrac{1}{p}+\dfrac{1}{q}=1$. Положив $d = L\int_{t_0}^{T} k(t)dt + k_1(T-t_0)^{\frac{1}{q}} + 2Lk_2$, имеем, что

$$|J(u) - J(u_1)| \leq d\|u(t) - u_1(t)\|_{L_p^n}.$$

Пусть $r \geq d3^{\frac{\tilde{p}+1}{p}} \dfrac{1}{1 - L(\int_{t_0}^{T} M^p(s)ds)^{\frac{1}{p}}}$ и существует $\tilde{u}(\cdot) \in D$, что $J_r(\tilde{u}) < J_r(\overline{u})$. Положим

$\tilde{\rho}(t) = \psi(t,(A\tilde{u})(t),\tilde{u}(t))$. Тогда существует решение $w(\cdot) \in L_p^n[t_0,T]$ задачи $u(t) \in F(t,(Au)(t))$ такое, что

$$(\int_{t_0}^{T} |\tilde{u}(t) - w(t)|^p dt)^{\frac{1}{p}} \leq (\int_{t_0}^{T} \psi^p(t,(A\tilde{u})(t),\tilde{u}(t))dt)^{\frac{1}{p}} \dfrac{3^{\frac{\tilde{p}+1}{p}}}{1 - L(\int_{t_0}^{T} M^p(s)ds)^{\frac{1}{p}}}.$$

Тогда получим, что

$$J(w) \leq J(\tilde{u}) + d\|w(t) - \tilde{u}(t)\|_{L_p^n} \leq J(\tilde{u}) + r(\int_{t_0}^{T} \psi^p(t,(A\tilde{u})(t),\tilde{u}(t))dt)^{\frac{1}{p}} = J_r(\tilde{u}) < J_r(\overline{u}) = J(\overline{u}).$$

Получим противоречие. Полученное противоречие означает, что $\overline{u}(t)$ минимизирует также функционал $J_r(u)$ в $D$ при $r \geq r_0$. Теорема доказана.

**Теорема 2.** Пусть выполняется условие теоремы 1 и функция $\overline{u}(t)$ среди всех решений задачи (2) минимизирует функционал (1), где $p = 1$, $A$ − линейный непрерывный компактный оператор. Тогда существуют функции $\upsilon^*(\cdot) \in L_1^n[t_0,T]$, $u^*(\cdot) \in L_1^n[t_0,T]$ и $c_1, c_2 \in R^n$ такие, что

1) $(u^*(t), -\upsilon^*(t)) \in \partial(f(t,(A\overline{u})(t),\overline{u}(t)) + r\psi(t,(A\overline{u})(t),\overline{u}(t)))$,
2) $(-c_1, -c_2) \in \partial\varphi((A\overline{u})(t_0),(A\overline{u})(T))$,
3) $\int_{t_0}^{T} (\upsilon^*(t) - (A^*u^*)(t)|u(t))dt + ((Au)(t_0)|c_1) + ((Au)(T)|c_2)) = 0$ при $u(\cdot) \in L_1^n[t_0,T].$

**Доказательство.** По теореме 1 существует число $r_0 > 0$ такое, что $\overline{u}(t)$ минимизирует функционал $J_r(u)$ в $D$ при $r \geq r_0$, т.е. $J_r(\overline{u}) \leq J_r(u)$ при $u \in D$. Поэтому $J_r^0(\overline{u};u) \geq J_r^+(\overline{u};u) \geq 0$ при $u \in L_1^n[t_0,T]$, где

$$J_r^0(\overline{u};u) = \overline{\lim_{\upsilon \to \overline{u}, \lambda \downarrow 0}} \dfrac{1}{\lambda}(J_r(\upsilon + \lambda u) - J_r(\upsilon)),$$

$$J_r^+(\overline{u};u) = \overline{\lim_{\lambda \downarrow 0}} \dfrac{1}{\lambda}(J_r(\overline{u} + \lambda u) - J_r(\overline{u})).$$

Обозначим $\overline{f}(t,x,y) = f(t,x,y) + r\psi(t,x,y)$. Отметим, что (см.[17]), если существует суммируемая функция $M(t) > 0$ такая, что $\rho_x(F(t,x),F(t,y)) \leq M(t)|x-y|$ при $x, y \in B(\overline{x}(t),\alpha)$, то $|\psi(t,x_1,y_1) - \psi(t,x_2,y_2))| \leq M(t)|x_1 - x_2| + |y_1 - y_2|$ при $x_1, x_2 \in B(\overline{x}(t),\alpha)$, $y_1, y_2 \in R^n$. Если выполняется условие теоремы 1, то применяя лемму Фату 2.4.6 [27] имеем, что

$$J_r^0(\overline{u};u) \leq \int_{t_0}^{T} \overline{f}^0(t,(A\overline{u})(t),\overline{u}(t);(Au)(t),u(t))dt + \varphi^0((A\overline{u})(t_0),(A\overline{u})(T);((Au)(t_0),(Au)(T)))$$

в $L_1^n[t_0,T]$. Отсюда следует, что $\overline{u}(t) = 0$ минимизирует функционал

$$J(u) = \int_{t_0}^{T} \overline{f}^0(t,(A\overline{u})(t),\overline{u}(t);(Au)(t),u(t))dt + \varphi^0((A\overline{u})(t_0),(A\overline{u})(T);((Au)(t_0),(Au)(T)))$$

в $L_1^n[t_0,T]$. Легко проверить, что для $I(u)$ также выполняется условие следствия 2.6.1 при $p = 1$. Поэтому существуют функции $\upsilon^*(\cdot) \in L_1^n[t_0, T]$, $u^*(\cdot) \in L_1^n[t_0, T]$ и $c_1, c_2 \in R^n$ такие, что

1) $(u^*(t), -\upsilon^*(t)) \in \partial \bar{f}(t, (A\bar{u})(t), \bar{u}(t))$,

2) $(-c_1, -c_2) \in \partial \varphi((A\bar{u})(t_0), (A\bar{u})(T))$,

3) $\int_{t_0}^{T} (\upsilon^*(t) - (A^*u^*)(t) | u(t)) dt + ((Au)(t_0) | c_1) + ((Au)(T) | c_2)) = 0$ при $u(\cdot) \in L_1^n[t_0, T]$.

Теорема доказана.

## Гл.3. ЭКСТРЕМАЛЬНАЯ ЗАДАЧА ДЛЯ МНОГОМЕРНОГО

# ОПЕРАТОРНОГО ВКЛЮЧЕНИЯ

В гл.3 получены необходимые и достаточные условия экстремума для вариационных задач типа Фредгольма и экстремальной задачи для операторного включения.

### 3.1. Непрерывная зависимость решений операторного включения

Пусть $R^n$ $n$-мерное евклидово пространство, $D$ компактное множество в $R^m$. Совокупность всех непустых компактных (выпуклых компактных) подмножеств $R^n$ обозначим через $\text{comp}R^n (\text{conv}R^n)$.

Если $A, B \in \text{comp}R^n$, то положим $\rho_x(A,B) = \max\left\{\sup_{y \in B} d(y,A), \sup_{x \in A} d(x,B)\right\}$, где

$$d(y,A) = \inf\{|x-y| : x \in A\}, \quad |x-y| = \left(\sum_{i=1}^{n}(x_i - y_i)^2\right)^{\frac{1}{2}}.$$

Пусть $A: L_p^n(D) \to C^n(D)$ линейный непрерывный оператор, где $L_p^n(D) = \{u(t) = (u^1(t),\ldots,u^n(t)) \in R^n : u^i(\cdot) \in L_p(D), i = \overline{1,n}\}$ пространство с нормой $\|u(\cdot)\|_{L_p^n} = \left(\int_D |u(s)|^p ds\right)^{\frac{1}{p}}$. Множество $\{(u, Au) \in L_p^n(D) \times C^n(D) : u \in L_p^n(D)\}$ относительно нормы $\|(u, Au)\| = \|u\|_{L_p^n(D)} + \|Au\|_{C^n(D)}$ является подпространство в $L_p^n(D) \times C^n(D)$.

Так как $\|u\|_{L_p^n(D)} \leq \|u\|_{L_p^n(D)} + \|Au\|_{C^n(D)} \leq (1+\|A\|)\|u\|_{L_p^n(D)}$, то имеем, что в пространстве $\{(u, Au) \in L_p^n(D) \times C^n(D) : u \in L_p^n(D)\}$ нормы $\|(u, Au)\|$ и $\|u\|_{L_p^n(D)}$ эквивалентны.

Пусть $F: D \times R^n \to \text{comp}R^n$ многозначное отображение, $1 \leq p < +\infty$.

Рассмотрим задачу для включения
$$u(t) \in F(t, (Au)(t)), \quad t \in D, \quad u(\cdot) \in L_p^n(D) \tag{1}$$

Функцию $u(\cdot) \in L_p^n(D)$, удовлетворяющую (1) назовем решением задачи (1).

**Теорема 1.** Пусть $F: D \times R^n \to \text{comp}R^n$ многозначное отображение, $t \to F(t,x)$ измеримо по $t$, существует функция $M(\cdot) \in L_p(D)$, где $M(t) > 0$ такая, что
$$\rho_x(F(t,x), F(t,x_1))) \leq M(t)|x - x_1|$$
при $x, x_1 \in R^n$. Кроме того, пусть $A: L_p^n(D) \to C^n(D)$ линейный непрерывный оператор и $\|A\|\left(\int_D M^p(s)ds\right)^{\frac{1}{p}} < 1$; $\rho(\cdot) \in L_p(D)$ и $\overline{u}(\cdot)$ такие, что
$$d(\overline{u}(t), F(t, (A\overline{u})(t))) \leq \rho(t), \quad t \in D.$$

Тогда существует такое решение $u(\cdot) \in L_p^n(D)$ задачи (1), что

$$|u(t) - \overline{u}(t)| \leq \rho(t) + \frac{\|A\|\left(\int_D \rho^p(s)ds\right)^{\frac{1}{p}} M(t)}{1 - \|A\|\left(\int_D M^p(s)ds\right)^{\frac{1}{p}}}, \qquad |(Au)(t) - (A\overline{u})(t)| \leq \frac{\|A\|\left(\int_D \rho^p(s)ds\right)^{\frac{1}{p}}}{1 - \|A\|\left(\int_D M^p(s)ds\right)^{\frac{1}{p}}}.$$

**Доказательство.** Построим последовательность $u_i(t)$ ($i = 0,1,2\ldots$), где $t \in [t_0, T]$, с помощью рекурентного соотношения
$$x_0(t) = \overline{x}(t) = (A\overline{u})(t), \qquad x_{i+1}(t) = (Au_i)(t), \quad i = 0,1,2,\ldots \tag{2}$$
где $u_i(s) \in F(s, x_i(s))$ при $i \geq 0$, $u_i(s)$ измеримы и
$$|u_0(s) - \overline{u}(s)| = d(\overline{u}(s), F(s, x_0(s))), \quad u_0(s) \in F(s, x_0(s)),$$
$$|u_i(s) - u_{i-1}(s)| = d(u_{i-1}(s), F(s, (Au_{i-1})(s))), \quad i = 1,2,\ldots$$

при $s \in D$. По лемме 2.1.4 [25] такая функция $u_i(t)$ существует. По условию, при $i \geq 1$ имеем

$$|u_{i+1}(s) - u_i(s)| = d(u_i(s), F(s, x_{i+1}(s))) \leq \rho_x(F(s, x_i(s)), F(s, x_{i+1}(s))) \leq M(s)|x_{i+1}(s) - x_i(s)|$$

при $s \in D$. Ясно, что

$$|x_{i+1}(s) - x_i(s)| = |(Au_i)(s) - (Au_{i-1})(s)| \leq \|A\| (\int_D |u_i(s) - u_{i-1}(s)|^p ds)^{\frac{1}{p}}.$$

По условию, имеем, что $|u_0(s) - \bar{u}(s)| = d(\bar{u}(s), F(s, (A\bar{u})(s))) \leq \rho(s)$ при $s \in [t_0, T]$. Поэтому имеем, что

$$|x_1(t) - x_0(t)| = |(Au_0)(t) - (A\bar{u})(t)| = |A(u_0 - \bar{u})(t)| \leq \|A\| (\int_D |u_0(s) - \bar{u}(s)|^p ds)^{\frac{1}{p}} \leq \|A\| (\int_D |\rho(s)|^p ds)^{\frac{1}{p}}.$$

Тогда получим

$$|u_1(s) - u_0(s)| = d(u_0(s), F(s, x_1(s))) \leq \rho_x(F(s, x_0(s)), F(s, x_1(s))) \leq$$

$$\leq M(s)|x_1(s) - x_0(s)| \leq \|A\| M(s)(\int_D \rho^p(\tau)d\tau)^{\frac{1}{p}}$$

при $s \in D$. Ясно, что

$$|x_2(t) - x_1(t)| \leq \|A\| (\int_D |u_1(s) - u_0(s)|^p ds)^{\frac{1}{p}} \leq \|A\| (\int_D M^p(s)\|A\|^p \int_D \rho^p(v) dv ds)^{\frac{1}{p}} =$$

$$= \|A\|^2 (\int_D \rho^p(s)ds)^{\frac{1}{p}} (\int_D M^p(s)ds)^{\frac{1}{p}}$$

при $t \in D$. Поэтому

$$|u_2(s) - u_1(s)| \leq M(s)\|A\|^2 (\int_D \rho^p(s)ds)^{\frac{1}{p}} (\int_D M^p(s)ds)^{\frac{1}{p}}$$

при $s \in D$. Отсюда получим, что

$$|x_3(t) - x_2(t)| \leq \|A\|^3 (\int_{t_0}^T \rho^p(s)ds)^{\frac{1}{p}} (\int_{t_0}^T M^p(s)ds)^{\frac{2}{p}}$$

при $t \in D$. Тогда

$$|u_3(s) - u_2(s)| \leq M(s)\|A\|^3 (\int_{t_0}^T \rho^p(s)ds)^{\frac{1}{p}} (\int_{t_0}^T M^p(s)ds)^{\frac{2}{p}}$$

при $s \in D$. Поэтому

$$|x_4(t) - x_3(t)| \leq \|A\|^4 (\int_{t_0}^T \rho^p(s)ds)^{\frac{1}{p}} (\int_{t_0}^T M^p(s)ds)^{\frac{3}{p}}$$

при $t \in D$. Также имеем

$$|u_4(s) - u_3(s)| \leq M(s)\|A\|^4 (\int_{t_0}^T \rho^p(s)ds)^{\frac{1}{p}} (\int_{t_0}^T M^p(s)ds)^{\frac{3}{p}}$$

при $s \in D$. Поэтому

$$|x_5(t) - x_4(t)| \leq \|A\|^5 (\int_D \rho^p(s)ds)^{\frac{1}{p}} (\int_D M^p(s)ds)^{\frac{4}{p}}$$

при $t \in D$. Продолжая процесс получим

$$|u_{m+1}(s) - u_m(s)| \leq M(s)\|A\|^{m+1} (\int_D \rho^p(s)ds)^{\frac{1}{p}} (\int_D M^p(s)ds)^{\frac{m}{p}} \tag{3}$$

$$|x_{m+1}(s) - x_m(s)| \leq \|A\|^{m+1} (\int_D \rho^p(s)ds)^{\frac{1}{p}} (\int_D M^p(s)ds)^{\frac{m}{p}} \tag{4}$$

при $s \in D$. Тогда получим

$$|u_{m+1}(t) - \bar{u}(t)| \le |u_0(t) - \bar{u}(t)| + |u_1(t) - u_0(t)| + |u_2(t) - u_1(t)| + |u_3(t) - u_2(t)| + \ldots +$$
$$+ |u_{m+1}(t) - u_m(t)| \le \rho(t) + M(t)(\int_D \rho^p(s)ds)^{\frac{1}{p}} \sum_{i=0}^{m} \|A\|^{i+1} (\int_D M^p(s)ds)^{\frac{i}{p}}$$

при $s \in D$. Если $\|A\|(\int_D M^p(s)ds)^{\frac{1}{p}} < 1$, то имеем, что

$$|u_{m+1}(t) - \bar{u}(t)| \le \rho(t) + M(t)(\int_D \rho^p(s)ds)^{\frac{1}{p}} \frac{\|A\|}{1 - \|A\|(\int_D M^p(s)ds)^{\frac{1}{p}}}. \tag{5}$$

Также имеем, что
$$|x_{m+1}(t) - \bar{x}(t)| \le |x_1(t) - \bar{x}(t)| + |x_2(t) - x_1(t)| + |x_3(t) - x_2(t)| + \ldots + |x_{m+1}(t) - x_m(t)| =$$
$$\le (\int_D \rho^p(s)ds)^{\frac{1}{p}} \sum_{i=0}^{m} \|A\|^{i+1} (\int_D M^p(s)ds)^{\frac{i}{p}} \le (\int_D \rho^p(s)ds)^{\frac{1}{p}} \frac{\|A\|}{1 - \|A\|(\int_D M^p(s)ds)^{\frac{1}{p}}} \tag{6}$$

при $t \in D$. Из оценки (3)-(6) вытекает, что последовательности $x_m(t)$ и $u_m(t)$ сходятся соответственно к функциям $x(\cdot)$ и $u(\cdot) \in L_p^n(D)$. Из теоремы 1.2.19[4] и из (2) следует, что $x(t)$ непрерывно. Так как $u_i(t) \in F(t, x_i(t))$, то из теоремы 1.2.23 и 1.2.28 [3] имеем, что $u(t) \in F(t, x(t))$ при $t \in D$, т.е. из (2) следует, что $u(t) \in F(t, (Au)(t))$ при $t \in D$. Кроме того, из (5) и (6) получим, что

$$|u(t) - \bar{u}(t)| \le \rho(t) + M(t)(\int_D \rho^p(s)ds)^{\frac{1}{p}} \frac{\|A\|}{1 - \|A\|(\int_{t_0}^T M^p(s)ds)^{\frac{1}{p}}},$$

$$|x(t) - \bar{x}(t)| \le \frac{\|A\|(\int_D \rho^p(s)ds)^{\frac{1}{p}}}{1 - \|A\|(\int_D M^p(s)ds)^{\frac{1}{p}}}.$$

Теорема доказана.

В теореме 1 заменив неравенство $|(Au)(s) - (Au)(s)| \le \|A\|(\int_D |u(s) - u(s)|^p ds)^{\frac{1}{p}}$ неравенством $|(Au)(t) - (A\upsilon)(t)| \le L\|u - \upsilon\|_{L_p^n}$ аналогично теореме 1, доказывается следующая теорема.

**Теорема 2.** Пусть $F: D \times R^n \to comp R^n$ многозначное отображение, $t \to F(t, x)$ измеримо по $t$, существует функция $M(\cdot) \in L_p(D)$, где $M(t) > 0$ такая, что
$$\rho_x(F(t, x), F(t, x_1))) \le M(t)|x - x_1|$$
при $x, x_1 \in R^n$. Кроме того, пусть $A: L_p^n(D) \to C^n(D)$ оператор и существует $L > 0$ такое, что $|(Au)(t) - (A\upsilon)(t)| \le L\|u - \upsilon\|_{L_p^n}$ при $u(\cdot), \upsilon(\cdot) \in L_p^n(D)$; пусть $\rho(\cdot) \in L_p(D)$ и $\bar{u}(\cdot)$ такие, что $d(\bar{u}(t), F(t, (A\bar{u})(t))) \le \rho(t)$, $t \in D$, где $L(\int_D M^p(s)ds)^{\frac{1}{p}} < 1$. Тогда существует такое решение $u(\cdot) \in L_p^n(D)$ задачи (1), что

$$|u(t)-\overline{u}(t)| \leq \rho(t) + \frac{L(\int_D \rho^p(s)ds)^{\frac{1}{p}} M(t)}{1-L(\int_D M^p(s)ds)^{\frac{1}{p}}}, \qquad |(Au)(t)-(A\overline{u})(t)| \leq \frac{L(\int_D \rho^p(s)ds)^{\frac{1}{p}}}{1-L(\int_D M^p(s)ds)^{\frac{1}{p}}}.$$

Аналогично теореме 1, доказываются следующие теоремы.

**Теорема 3.** Пусть $F: D \times R^n \to \mathrm{comp}R^n$ многозначное отображение, $t \to F(t,x)$ измеримо по $t$, существует функция $M(\cdot) \in L_p(D)$, где $M(t) > 0$ такая, что
$$\rho_x(F(t,x), F(t,x_1))) \leq M(t)|x-x_1|$$
при $x, x_1 \in R^n$. Кроме того, пусть $A: L_p^n(D) \to L_\infty^n(D)$ линейный непрерывный оператор; $\rho(\cdot) \in L_p(D)$ и $\overline{u}(\cdot)$ такие, что
$$d(\overline{u}(t), F(t,(A\overline{u})(t))) \leq \rho(t), \quad t \in D,$$
где $\|A\|(\int_D M^p(s)ds)^{\frac{1}{p}} < 1$. Тогда существует такое решение $u(\cdot) \in L_p^n(D)$ задачи (1), что
$$|u(t)-\overline{u}(t)| \leq \rho(t) + \frac{\|A\|(\int_D \rho^p(s)ds)^{\frac{1}{p}} M(t)}{1-\|A\|(\int_D M^p(s)ds)^{\frac{1}{p}}}, \qquad |(Au)(t)-(A\overline{u})(t)| \leq \frac{\|A\|(\int_D \rho^p(s)ds)^{\frac{1}{p}}}{1-\|A\|(\int_D M^p(s)ds)^{\frac{1}{p}}}.$$

**Теорема 4.** Пусть $F:[t_0,T] \times R^n \to \mathrm{comp}R^n$ многозначное отображение, $t \to F(t,x)$ измеримо по $t$, существует функция $M(\cdot) \in L_p(D)$, где $M(t) > 0$ такая, что
$$\rho_x(F(t,x), F(t,x_1))) \leq M(t)|x-x_1|$$
при $x, x_1 \in R^n$. Кроме того, пусть $A: L_p^n(D) \to L_\infty^n(D)$ оператор и существует $L > 0$ такое, что $|(Au)(t)-(A\upsilon)(t)| \leq L\|u-\upsilon\|_{L_p^n}$ при $u(\cdot), \upsilon(\cdot) \in L_p^n(D)$; $\rho(\cdot) \in L_p(D)$ и $\overline{u}(\cdot)$ такие, что
$$d(\overline{u}(t), F(t,(A\overline{u})(t))) \leq \rho(t), \quad t \in D,$$
где $L(\int_D M^p(s)ds)^{\frac{1}{p}} < 1$. Тогда существует такое решение $u(\cdot) \in L_p^n(D)$ задачи (1), что
$$|u(t)-\overline{u}(t)| \leq \rho(t) + \frac{L(\int_D \rho^p(s)ds)^{\frac{1}{p}} M(t)}{1-L(\int_D M^p(s)ds)^{\frac{1}{p}}}, \qquad |(Au)(t)-(A\overline{u})(t)| \leq \frac{L(\int_D \rho^p(s)ds)^{\frac{1}{p}}}{1-L(\int_D M^p(s)ds)^{\frac{1}{p}}}.$$

Отметим, что в теореме 2 и 4 оператор $A$ в общем случае не является линейным.

Положив $\overline{u}(t) = 0$ и $(A\overline{u})(t) = 0$, из теоремы 1 имеем, что верно следующее следствие.

**Следствие 1.** Пусть $A: L_p^n(D) \to C^n(D)$ линейный непрерывный оператор, $F: D \times R^n \to \mathrm{comp}R^n$ многозначное отображение, $t \to F(t,x)$ измеримо по $t$, существует функция $M(\cdot) \in L_p(D)$, $M(t) > 0$ такая, что $\rho_x(F(t,x), F(t,x_1)) \leq M(t)|x-x_1|$ при $x, x_1 \in R^n$ и $\|F(t,0)\| \in L_p(D)$. Тогда существует решение $u(\cdot) \in L_p^n(D)$ задачи (1).

**Замечание 1.** Из доказательства теоремы 1 следует, что в теореме 1-4 условие $F: D \times R^n \to \mathrm{comp}R^n$ можно заменить условием: $F: D \times (\overline{x}(t) + \alpha B) \to \mathrm{comp}R^n$ при $t \in D$, где $\overline{x}(t) = (A\overline{u})(t)$, $\alpha \geq \dfrac{\|A\|(\int_D \rho^p(s)ds)^{\frac{1}{p}}}{1-\|A\|(\int_D M^p(s)ds)^{\frac{1}{p}}}$ или $\alpha \geq \dfrac{L(\int_D \rho^p(s)ds)^{\frac{1}{p}}}{1-L(\int_D M^p(s)ds)^{\frac{1}{p}}}$ соответственно,

$B = \{z \in R^n : \|z\| \le 1\}$ - единичный шар в $R^n$.

**Лемма 2.** Пусть $A : L_p^n(D) \to L_\infty^n(D)$ линейный непрерывный оператор, $F : D \times (x_0(t) + \alpha B) \to \text{comp} R^n$, где $x_0(t) = (Au_0)(t)$, $\alpha > 0$, функция $u_0(\cdot) \in L_p^n(D)$ - решение задачи $u(t) \in F(t,(Au)(t))$ и $t \to F(t,x)$ измеримо по $t$, существует функция $M(\cdot) \in L_p(D)$, где $M(t) > 0$ такая, что

$$\rho_x(F(t,x), F(t,x_1))) \le M(t)|x - x_1|$$

при $x, x_1 \in x_0(t) + \alpha B$, $\|A\|(\int_D M^p(s)ds)^{\frac{1}{p}} < 1$. Тогда существует такое $\delta > 0$, что при $s(\cdot) \in L_p^n(D)$, $\|s(\cdot)\|_{L_p^n} \le \delta$ найдется решение $u_s(t)$ задачи $u(t) \in F(t,(Au)(t)) + s(t)$, что $|x_0(t) - x_s(t)| \le \alpha$ при $t \in D$, где $x_s(t) = (Au_s)(t)$ и $\|x_s(\cdot) - x_0(\cdot)\|_{L_\infty^n} \to 0$ при $\|s(\cdot)\|_{L_p^n} \to 0$.

**Доказательство.** Ясно, что
$$\rho_x(F(t,x) + s(t), F(t,x') + s(t)) \le M(t)|x - x'|$$
при $|x - x_0(t)| \le \alpha$, $|x' - x_0(t)| \le \alpha$ и $d(u_0(t), F(t,x_0(t)) + s(t)) \le |s(t)|$. По теореме 3 существует решение $u_s(t)$ задачи $u_s(t) \in F(t,x_s(t)) + s(t)$, такое, что

$$|x_0(t) - x_s(t)| \le \frac{\|A\|(\int_D s^p(t)dt)^{\frac{1}{p}}}{1 - \|A\|(\int_D M^p(s)ds)^{\frac{1}{p}}}. \quad (7)$$

Если определять $\delta$ из неравенства $\dfrac{\|A\|\delta}{1 - \|A\|(\int_D M^p(s)ds)^{\frac{1}{p}}} \le \alpha$, то получим, что верна первая часть утверждения леммы.

Из (7) получим, что $\|x_s(\cdot) - x_0(\cdot)\|_{L_\infty^n(D)} \le \dfrac{\|A\|(\int_D s^p(t)dt)^{\frac{1}{p}}}{1 - \|A\|(\int_D M^p(s)ds)^{\frac{1}{p}}}$, т.е. $\|x_s(\cdot) - x_0(\cdot)\|_{L_\infty^n(D)} \to 0$ при $\|s(\cdot)\|_{L_p^n} \to 0$. Лемма доказана.

**Лемма 3.** Пусть $A : L_p^n(D) \to L_\infty^n(D)$ линейный непрерывный оператор, $F : D \times R^n \to \text{comp} R^n \cup \{\varnothing\}$, отображение $t \to F(t,x)$ измеримо на $D$, отображение $x \to F(t,x)$ замкнуто и выпукло почти при всех $t \in D$, т.е. $\text{gr} F_t = \{(x,y) : y \in F(t,x)\}$ замкнуто и выпукло почти при всех $t \in D$. Пусть существует функция $\lambda(\cdot) \in L_p^n(D)$ такая, что

$$\|F(t,x)\| \le \lambda(t)(1 + |x|)$$

при $x \in R^n$, где $\|F(t,x)\| = \sup\{|y| : y \in F(t,x)\}$, $\|\varnothing\| = 0$ и существует решение $u_0(t)$ задачи $u_0(t) \in F(t,x_0(t))$, что $x_0(t) = (Au_0)(t)$ принадлежит $\text{dom } F_t = \{x : F(t,x) \ne \varnothing\}$ вместе с некоторой $\varepsilon$ трубкой, т.е. $\{x : |x_0(t) - x| \le \varepsilon\} \subset \text{dom } F_t$, $\|A\|(\int_D M^p(s)ds)^{\frac{1}{p}} < 1$, где $M(t) = \dfrac{6\lambda(t)\left(1 + |x_0(t)| + \dfrac{5\varepsilon}{6}\right)}{\varepsilon}$. Тогда существуют такие $\delta > 0$ и решение $u_s(t)$ задачи $u(t) \in F(t,(Au)(t)) + s(t)$, где $s(\cdot) \in L_p^n(D)$, $\|s(\cdot)\|_{L_p^n} \le \delta$, что $|x_0(t) - x_s(t)| \le \varepsilon$ при $t \in D$, где $x_s(t) = (Au_s)(t)$, и

$$\|x_s(\cdot) - x_0(\cdot)\|_{L_\infty^n(D)} \to 0 \text{ при } \|s(\cdot)\|_{L_p^n(D)} \to 0.$$

**Доказательство.** Рассмотрим множество $S_t = \{x : |x_0(t) - x| \le \frac{\varepsilon}{2}\}$ и обозначим $\delta = \frac{\varepsilon}{3}$. Пусть $x \in S_t$, тогда $x + \delta B \subset \operatorname{int} \operatorname{dom} F_t$, где $B$ единичный шар в $R^n$ с центром в нуле. Так как множество $\operatorname{gr} F_t$ выпукло, то аналогично доказательству теоремы 3.1.1 [13] можно показать, что

$$\rho_x(F(t,y), F(t,x)) \le \frac{2\lambda(t)(1+|x|+\delta)}{\delta}|y-x| \qquad (8)$$

для $y \in x + \delta B$ п.в. $t$. Ясно, что $S_t + \delta B \subset \operatorname{int} \operatorname{dom} F_t$. Поэтому (8) верно для всех $x, y \in S_t$ и почти всех $t$, т.е.

$$\rho_x(F(t,y), F(t,x)) \le \frac{2\lambda(t)\left(1+|x_0(t)|+\frac{\varepsilon}{2}+\delta\right)}{\delta}|y-x|.$$

при $x, y \in S_t$. Поэтому доказательство леммы 3 вытекает из леммы 2. Лемма доказана.

Отметим, что используя теоремы 1, 2 или 4 можно получить ряд аналог леммы 2 и 3.

### 3.2. О субдифференцируемости интегрального функционала с оператором типа Фредгольма

Пусть $D$ компактное множество в $R^m$, $f: D \times R^n \to R_{+\infty}$ нормальный выпуклый интегрант, т.е. $f$ -такая функция из $D \times R^n$ в $R_{+\infty}$, что $\operatorname{ep} f_t = \{(z, \alpha) \in R^n \times R : f(t,z) \le \alpha\}$ - выпуклое замкнутое множество и $t \to \operatorname{ep} f_t$ измеримо.

Рассмотрим функционал

$$J(u(\cdot)) = \int_D f(t, (Au)(t))dt,$$

где $u \in L_p^n(D)$, $1 \le p < +\infty$.

**Лемма 1.** Пусть $f: D \times R^n \to R_{+\infty}$ нормальный выпуклый интегрант, функция $x \to f(t,x)$ не убывает, $A: L_p^n(D) \to L_\infty^n(D)$ выпуклый оператор, то $J: L_p^n(D) \to R_{+\infty}$ выпуклый функционал.

**Доказательство.** Пусть $u_1(\cdot), u_2(\cdot) \in L_p^n(D)$ и $\alpha \in [0,1]$. Тогда получим, что

$$J(\alpha u_1 + (1-\alpha)u_2) = \int_D f(t, A(\alpha u_1 + (1-\alpha)u_2)(t))dt \le$$

$$\le \int_D f(t, \alpha A(u_1)(t) + (1-\alpha)A(u_2)(t))dt \le$$

$$\le \alpha \int_D f(t, A(u_1)(t))dt + (1-\alpha)\int_D f(t, A(u_2)(t))dt.$$

Лемма доказана.

Пусть $A: L_p^n(D) \to C^n(D)$ линейный непрерывный оператор. Рассмотрим субдифференцируемость интегрального функционала

$$J(u(\cdot)) = \int_D f(t, (Au)(t))dt$$

в $L_p^n(D)$, где $1 \le p < +\infty$.

Отметим, что если $f: D \times R^n \to R_{+\infty}$ нормальный выпуклый интегрант и $A: L_p^n(D) \to C^n(D)$ линейный оператор, то $J: L_p^n(D) \to R_{+\infty}$ выпуклый функционал.

Субградиентами $J$ в точке $u_0(\cdot) \in L_p^n(D)$ являются по определению, элементы $\upsilon^* \in L_q^n(D)$, $\frac{1}{p} + \frac{1}{q} = 1$, для которых

$$J(u(\cdot)) - J(u_0(\cdot)) \geq \langle \upsilon^*(\cdot), u(\cdot) - u_0(\cdot) \rangle$$

при всех $u(\cdot) \in L_p^n(D)$ или, что равносильно,

$$J^*(\upsilon^*) + J(u_0(\cdot)) = \langle \upsilon^*(\cdot), u_0(\cdot) \rangle$$

где $\langle \upsilon^*(\cdot), u_0(\cdot) \rangle = \int_D (\upsilon^*(t) | u_0(t)) dt$.

Множество всех таких субградиентов обозначается через $\partial J(u_0(\cdot))$ и называется субдифференциалом функционала $J$ в точке $u_0(\cdot)$.

В этом параграфе устанавливается связь между $\partial J(u_0(\cdot))$ и $\partial f_t((Au_0)(t))$.

**Лемма 2.** Пусть $A: L_p^n(D) \to C^n(D)$ линейный непрерывный компактный оператор, $f: D \times R^n \to R_{+\infty}$ нормальный интегрант, $x \to f(t, x)$ сублинейная функция в $R^n$ и для всякого $x \in R^n$ функция $f(t, x)$ суммируема и $y^*(\cdot) \in L_q^n(D)$, где $\frac{1}{p} + \frac{1}{q} = 1$. Тогда

$$J^*(y^*(\cdot)) = \begin{cases} 0, & y^*(\cdot) = (A^* u^*)(t),\ u^*(t) \in \partial f(t, 0),\ u^*(\cdot) \in L_1^n(D), \\ +\infty, & \text{в других случаях}. \end{cases}$$

**Доказательство.** По условию леммы 2 имеем, что $J(u(\cdot)) = \int_D f(t, (Au)(t)) dt$ сублинейная функция в $L_p^n(D)$. Так как для всякого $x \in R^n$ функция $f(t, x)$ суммируема, то функционал $x(\cdot) \to \int_D f(t, x(t)) dt$ непрерывен в $C^n(D)$ (см.[14]).

По условию оператор $A: L_p^n(D) \to C^n(D)$ непрерывен. Так как $x(\cdot) \to \int_D f(t, x(t)) dt$ непрерывен в $C^n(D)$, то существует $\delta > 0$ такое, что функционал $\int_D f(t, x(t)) dt$ ограничен в множестве $\{x(\cdot) \in C^n(D) : \|x(\cdot)\| \leq \delta\}$. Так как $A: L_p^n(D) \to C^n(D)$ непрерывен, то существует $\alpha > 0$ такое, что $\|Au\|_{C^n} \leq \delta$ при $\|u\|_{L_p^n} \leq \alpha$. Поэтому получим, что функционал $J(u(\cdot)) = \int_D f(t, (Au)(t)) dt$ ограничен в множестве $\{u(\cdot) \in L_p^n(D) : \|u(\cdot)\|_{L_p^n} \leq \alpha\}$. Так как $J(u)$ выпуклый функционал, то имеем, что $u \to J(u)$ непрерывный функционал в $L_p^n(D)$.

Ясно, что $I_0(y) = \int_D f(t, y(t)) dt$ сублинейная функция в $L_\infty^n(D)$ и по следствию 2A[14] непрерывный функционал в $L_\infty^n(D)$. Поэтому по предложению 3 ([6], стр.210) и по теореме 3 ([6], стр.362) $\partial I_0(0)$ слабо* компактно и $\partial I_0(0) \subset L_1^n(D)$. Из теоремы 6.4.8[8] следует, что

$$I_0(y(\cdot)) = \max_{u^* \in \partial I_0(0)} \langle u^*, y \rangle = \max_{u^* \in \partial I_0(0)} \int_D (u^*(t) | y(t)) dt$$

при $y(\cdot) \in L_\infty^n(D)$. Известно, что $u^* \in \partial I_0(0)$ в том и только в том случае, когда $u^*(t) \in \partial f(t, 0)$ и $u^*(\cdot) \in L_1^n(D)$.

Ясно, что $J(u(\cdot)) = \max_{u^* \in \partial I_0(0)} \int_D (u^*(t) | (Au)(t)) dt$ при $u(\cdot) \in L_p^n(D)$. По условию $A: L_p^n(D) \to C^n(D)$

линейный непрерывный компактный оператор. Если $y^*(\cdot) \in L_q^n(D)$, где $\frac{1}{p}+\frac{1}{q}=1$, то применяя теоремы 6.2.7[12] имеем

$$J^*(y^*(\cdot)) = \sup_{u(\cdot) \in L_p^n(D)} \{\int_D (y^*(t)|u(t))dt - \int_D f(t,(Au)(t))dt\} =$$

$$= \sup_{u(\cdot) \in L_p^n(D)} \{\int_D (y^*(t)|u(t))dt - \sup_{u^* \in \partial I_0(0)} \int_D (u^*(t)|(Au)(t))dt\} =$$

$$= \sup_{u(\cdot) \in L_p^n(D)} \inf_{u^* \in \partial I_0(0)} \{\int_D (y^*(t)-(A^*u^*)(t)|u(t))dt =$$

$$= \inf_{u^* \in \partial I_0(0)} \sup_{u(\cdot) \in L_p^n(D)} \{\int_D (y^*(t)-(A^*u^*)(t)|u(t))dt =$$

$$= \begin{cases} 0, & y^*(\cdot)=(A^*u^*)(t), \ u^*(t) \in \partial f(t,0), \ u^*(\cdot) \in L_1^n(D), \\ +\infty, & \text{в других случаях.} \end{cases}$$

Лемма доказана.

**Следствие 1.** Если удовлетворяется условие леммы 2, то $y^*(\cdot) \in \partial J(0)$ в том и только в том случае, когда существует функция $u^*(t) \in \partial f(t,0)$, $u^*(\cdot) \in L_1^n(D)$, такая, что $y^*(\cdot)=(A^*u^*)(t)$.

**Теорема 1.** Пусть $f:D \times R^n \to R_{+\infty}$ нормальный выпуклый интегрант, существует $\delta > 0$ такое, что $f(t,\bar{x}(t)+x)$ суммируема при $x \in R^n$, $|x| \leq \delta$, где $\bar{x}(t)=(A\bar{u})(t)$, $\bar{u}(\cdot) \in L_p^n(D)$, $A: L_p^n(D) \to C^n(D)$ линейный непрерывный компактный оператор. Тогда $\partial J(\bar{u})$ непусто и $y^*(\cdot) \in L_q^n(D)$ принадлежит $\partial J(\bar{u})$ в том и только в том случае, когда существует функция $u^*(t) \in \partial f(t,\bar{x}(t))$, $u^*(\cdot) \in L_1^n(D)$, такая, что $y^*(\cdot)=(A^*u^*)(t)$.

**Доказательство.** Ясно, что

$$J'(\bar{u}(\cdot);u(\cdot)) = \lim_{\lambda \downarrow 0} \frac{1}{\lambda}(J(\bar{u}(\cdot)+\lambda u(\cdot))-J(\bar{u}(\cdot)) =$$

$$= \lim_{\lambda \downarrow 0} \frac{1}{\lambda} \int_{t_0}^T (f(t,(A\bar{u})(t)+\lambda(Au)(t))-f(t,(A\bar{u})(t)))dt.$$

Из доказательства предложения 4.1[11] имеем

$$f(t,(A\bar{u})(t))-f(t,(A\bar{u})(t)-(Au)(t)) \leq \frac{1}{\lambda}(f(t,(A\bar{u})(t)+\lambda(Au)(t))-f(t,(A\bar{u})(t))) \leq$$

$$\leq (f(t,(A\bar{u})(t)+(Au)(t))-f(t,(A\bar{u})(t))).$$

при $\lambda \in (0,1)$. Пусть $\alpha > 0$ такое, что $\|Au\|_{C^n} \leq \delta$ при $\|u\|_{L_p^n} \leq \alpha$. Так как $f(t,(A\bar{u})(t)-(Au)(t))$ и $f(t,(A\bar{u})(t)+(Au)(t))$ суммируемы при $u(\cdot) \in L_p^n(D)$, $\|u\|_{L_p^n} \leq \alpha$, то используя теорему Лебега получим

$$J'(\bar{u}(\cdot);u(\cdot)) = \int_D f'(t,(A\bar{u})(t);(Au)(t))dt.$$

Если учесть, что $\partial J(\bar{u}(\cdot)) = \partial J'(\bar{u}(\cdot);0)$, $\partial f'(t,(A\bar{u})(t)) = \partial f'(t,(A\bar{u})(t);0)$, то из следствия 1 имеем, что $y^*(\cdot) \in L_q^n(D)$ принадлежит $\partial J(\bar{u})$ в том и только в том случае, когда существует функция $u^*(t) \in \partial f(t,\bar{x}(t))$, $u^*(\cdot) \in L_1^n(D)$, такая, что $y^*(\cdot)=(A^*u^*)(t)$. Теорема доказана.

Пусть $f: D \times R^n \to R_{+\infty}$ нормальный выпуклый интегрант и $A: L_p^n(D) \to C^n(D)$ линейный непрерывный сюръективный оператор. Тогда по лемме о правом обратном отображение (см.[1], стр.128) существует отображении $M: C^n(D) \to L_p^n(D)$ (вообще говоря, разрывное и нелинейное), удовлетворяющее условиям $A \circ M = I_{C^n}$, $\|M(x(\cdot))\| \leq \alpha \|x(\cdot)\|_{C^n}$ при $x(\cdot) \in C^n(D)$, где $\alpha > 0$.

Субградиентами $J$ в точке $u_0(\cdot) \in L_p^n(D)$ являются по определению, элементы $\upsilon^* \in L_q^n(D)$, $\frac{1}{p} + \frac{1}{q} = 1$, для которых

$$J(u(\cdot)) - J(u_0(\cdot)) \geq \langle \upsilon^*(\cdot), u(\cdot) - u_0(\cdot) \rangle$$

при всех $u(\cdot) \in L_p^n(D)$, т.е.

$$\int_D f(t, (Au)(t)) dt - \int_D f(t, (Au_0)(t)) dt \geq \int_D (\upsilon^*(t) | u(t) - u_0(t)) dt$$

при всех $u(\cdot) \in L_p^n(D)$. Обозначив $x(t) = (Au)(t)$ и $x_0(t) = (Au_0)(t)$ имеем, что

$$\int_D f(t, x(t)) dt - \int_D f(t, x_0(t)) dt \geq \int_D (\upsilon^*(t) | (Mx)(t) - (Mx_0)(t)) dt$$

при $x(\cdot) \in C^n(D)$. Если $M : C^n(D) \to L_p^n(D)$ линейное отображение, то отсюда следует, что

$$\int_D f(t, x(t)) dt - \int_D f(t, x_0(t)) dt \geq \int_D (M^*(\upsilon^*)(t) | x(t) - x_0(t)) dt$$

при $x(\cdot) \in C^n(D)$. Ясно, что $I(y) = \int_D f(t, y(t)) dt$ выпуклый функционал в $C^n(D)$. Если существует $\delta > 0$ такое, что $f(t, x_0(t) + x)$ суммируема при $x \in R^n, |x| \leq \delta$, то по следствию 4В [14] или по следствию 2.2.1 [17] имеем, что существует функция $u^*(\cdot) \in L_1^n(D)$ такая, что $u^*(t) \in \partial f(t, \bar{x}(t))$ и $M^*(\upsilon^*)(t) = u^*(t)$. Поэтому $\upsilon^*(t) = (A^* u^*)(t)$, т.е. функционал $\upsilon^*(\cdot) \in L_q^n(D)$ принадлежит $\partial J(u_0)$, то существует функция $u^*(\cdot) \in L_1^n(D)$ такая, что $\upsilon^*(t) = (A^* u^*)(t)$, где $u^*(t) \in \partial f(t, \bar{x}(t))$, $u^*(\cdot) \in L_1^n(D)$.

Если $A : L_p^n(D) \to W_{p,1}^n(D)$ линейный непрерывный компактный оператор, где $D$ ограниченная область и $p > n$, или $A : L_p^n(D) \to L_\infty^n(D)$ линейный непрерывный компактный оператор, то теорема 1 остается также верной.

### 3.3. О субдифференцируемости граничного функционала

Пусть $D$ компактное множество в $R^m$ и $\operatorname{int} D \neq \varnothing$, $\varphi : \partial D \times R^n \to (-\infty, +\infty]$ нормальный выпуклый интегрант, $A : L_p^n(D) \to C^n(D)$ линейный непрерывный оператор.

Рассмотрим субдифференцируемость функционала

$$F(u(\cdot)) = \int_{\partial D} \varphi(s, (Au)(s)) ds$$

в $L_p^n(D)$, где $1 \leq p < +\infty$. По определению $F^* : L_q^n(D) \to R_{+\infty}$ и

$$F^*(\upsilon(\cdot)) = \sup_{u \in L_p^n} \{ \int_{t_0}^T (\upsilon(t) | u(t)) dt - \int_{\partial D} \varphi(s, (Au)(s)) ds \},$$

где $\upsilon(\cdot) \in L_q^n(D)$, $\frac{1}{p} + \frac{1}{q} = 1$.

Пусть $\alpha \in [0,1]$. Ясно, что

$$F(\alpha u_1(\cdot) + (1-\alpha)u_2(\cdot)) = \int_{\partial D} \varphi(s, A(\alpha u_1 + (1-\alpha)u_2)(s))ds =$$

$$= \int_{\partial D} \varphi(s, \alpha(Au_1)(s) + (1-\alpha)(Au_2)(s))ds \leq$$

$$\leq \alpha \int_{\partial D} \varphi(s, (Au_1)(s))ds + (1-\alpha)\int_{\partial D} \varphi(s, (Au_2)(s))ds = \alpha F(u_1(\cdot)) + (1-\alpha)F(u_2(\cdot)).$$

Получим, что $F: L_p^n(D) \to R_{+\infty}$ выпуклый функционал.

Положим $Q = \{y(\cdot) \in L_\infty^n(\partial D) : \int_{\partial D} \varphi(s, y(s))ds < +\infty\}$.

**Теорема 1.** Если $\varphi: \partial D \times R^n \to R_{+\infty}$ нормальный выпуклый интегрант, существуют функция $x_0(\cdot) = (Au_0)(\cdot)$, где $u_0(\cdot) \in L_p^n(\partial D)$, и $r > 0$ такие, что $\varphi(s, x_0(s) + x)$ суммируемы при $|x| \leq r$ в $\partial D$, то $\upsilon(\cdot) \in \partial F(\bar{u}(\cdot))$, где $\bar{u}(\cdot) \in L_p^n(\partial D)$, $\upsilon(\cdot) \in L_q^n(D)$, в том и только в том случае, когда существует функционал $x^* \in L_\infty^n(\partial D)^*$, где $x^* = x_a^* + x_s^*$, $x_a^*(y) = \int_{\partial D}(y^*(s)|y(s))ds$ при $y(\cdot) \in L_\infty^n(\partial D)$, где $y^*(\cdot) \in L_1^n(\partial D)$, такой, что $y^*(s) \in \partial\varphi(s, (A\bar{u})(s))$, $x_s^*((A\bar{u})) = \sup_{y(\cdot) \in Q} x_s^*(y(\cdot))$ и $\int_D (\upsilon(t)|u(t))dt - x^*(Au) = 0$ при $u(\cdot) \in L_p^n(D)$.

**Доказательство.** Так как $\varphi(s, x_0(s) + x)$ суммируемы при $|x| \leq r$ в $\partial D$, то функционал $I(y(\cdot)) = \int_{\partial D} \varphi(s, y(s))ds$ непрерывен в точке $x_0(\cdot)$ относительно пространства $L_\infty^n(\partial D)$.

Положив $C = \{(u(\cdot), (Au)(\cdot)) \in L_p^n(D) \times C^n(\partial D) : u(\cdot) \in L_p^n(D)\}$, по теореме 3.4.1 [6] существует функционал $x^* \in L_\infty^n(\partial D)^*$, где $x^* = x_a^* + x_s^*$, такой, что

$$F^*(\upsilon(\cdot)) = \sup_{u(\cdot) \in L_p^n(D)} \{\int_D (\upsilon(t)|u(t))dt - \int_{\partial D} \varphi(s, (Au)(s))ds\} =$$

$$= \sup_{(u(\cdot), y(\cdot)) \in L_p^n(D) \times L_\infty^n(\partial D)} \{\int_D (\upsilon(t)|u(t))dt - \int_{\partial D} \varphi(s, y(s))ds - \delta_C(u(\cdot), y(\cdot))\} =$$

$$= \sup_{(u(\cdot), y(\cdot)) \in L_p^n(D) \times L_\infty^n(\partial D)} \{\int_D (\upsilon(t)|u(t))dt - x^*(y) - \delta_C(u(\cdot), y(\cdot))\} +$$

$$+ \sup_{y(\cdot) \in L_\infty^n(\partial D)} \{x^*(y) - \int_{\partial D} \varphi(s, y(s))ds\} = \sup_{u \in L_p^n(D)} \{\int_D (\upsilon(t)|u(t))dt - x^*(Au)\} +$$

$$+ \sup_{y(\cdot) \in L_\infty^n(\partial D)} \{x_a^*(y) + x_s^*(y) - \int_{\partial D} \varphi(s, y(s))ds\} =$$

$$= \begin{cases} \int_{\partial D} \varphi^*(s, y^*(s))ds + \sup_{y \in Q} x_s^*(y) : \int_D (\upsilon(t)|u(t))dt - x^*(Au) = 0 \text{ при } u(\cdot) \in L_p^n[t_0, T], \\ +\infty : \text{в других случаях.} \end{cases}$$

По определению $\upsilon(\cdot) \in \partial F(\bar{u}(\cdot))$ в том и только в том случае, когда $F(\bar{u}(\cdot)) + F^*(\upsilon(\cdot)) = \langle \upsilon(\cdot), \bar{u}(\cdot) \rangle$. Так как $\int_D (\upsilon(t)|u(t))dt - x^*(Au) = 0$ при $u(\cdot) \in L_p^n[t_0, T]$, то имеем, что $\int_D (\upsilon(t)|u(t))dt = x^*(Au)$. Отсюда следует, что

$$\int_{\partial D} \varphi(s, (A\bar{u})(s))ds + \int_{\partial D} \varphi^*(s, y^*(s))ds + \sup_{y \in Q} x_s^*(y) = x^*(A\bar{u}),$$

т.е.

$$\int_{\partial D} \varphi(s, (A\bar{u})(s))ds + \int_{\partial D} \varphi^*(s, y^*(s))ds + \sup_{y \in Q} x_s^*(y) = \int_{\partial D} (y^*(s)|(A\bar{u})(s))ds + x_s^*(A\bar{u}).$$

Отсюда следует, что

$$\int_{\partial D}\varphi(s,(A\overline{u})(s))ds + \int_{\partial D}\varphi^*(s,y^*(s))ds = \int_{\partial D}(y^*(s)|(A\overline{u})(s))ds, \qquad \sup_{y\in Q}x_s^*(y) = x_s^*(A\overline{u}).$$

Поэтому $y^*(s)\in \partial\varphi((A\overline{u})(s))$ и $\sup_{y\in Q}x_s^*(y) = x_s^*(A\overline{u})$. Теорема доказана.

**Следствие 1.** Если $\varphi:\partial D\times R^n \to (-\infty,+\infty]$ нормальный выпуклый интегрант, существуют функция $x_0(\cdot) = (Au_0)(\cdot)$, где $u_0(\cdot)\in L_p^n(\partial D)$, и $r>0$ такие, что $\varphi(s,x_0(s)+x)$ суммируемы при $|x|\leq r$ в $\partial D$, то $\upsilon(\cdot)\in \partial F(\overline{u}(\cdot))$, где $\overline{u}(\cdot)\in L_p^n(\partial D)$, $\upsilon(\cdot)\in L_q^n(D)$ в том и только в том случае, когда существует функционал $x^*\in L_\infty^n(\partial D)^*$, где $x^* = x_a^* + x_s^*$, $x_a^*(y) = \int_{\partial D}(y^*(s)|y(s))ds$ при $y(\cdot)\in L_\infty^n(\partial D)$, где $y^*(\cdot)\in L_1^n(\partial D)$, такой, что

$y^*(s)\in \partial\varphi(s,(A\overline{u})(s))$, $x_s^*((A\overline{u})(\cdot)) = \sup_{y(\cdot)\in Q}x_s^*(y(\cdot))$ и $\int_D (\upsilon(t)|u(t))dt - x^*(Au) = 0$ при $u(\cdot)\in L_p^n(D)$.

**Следствие 2.** Если при $u_0(\cdot) = \overline{u}(\cdot)$ выполняется условие теоремы 1, то функционал $x^*\in L_\infty^n(\partial D)^*$ «абсолютно непрерывен», т.е $x_s^* = 0$.

Из следствия 1 непосредственно следует следующее следствие.

**Следствие 3.** Если $\varphi:\partial D\times R^n \to (-\infty,+\infty]$ нормальный выпуклый интегрант, существуют функция $x_0(\cdot) = (Au_0)(\cdot)$, где $u_0(\cdot)\in L_p^n(\partial D)$, и $r>0$ такие, что $\varphi(s,x_0(s)+x)$ суммируемы при $|x|\leq r$ в $\partial D$, то $\upsilon(\cdot)\in \partial F(\overline{u}(\cdot))$, где $\overline{u}(\cdot)\in L_p^n(\partial D)$, $\upsilon(\cdot)\in L_q^n(D)$ в том и только в том случае, когда существует функционал $x^*\in C^n(\partial D)^*$, где $x^* = x_a^* + x_s^*$, $x_a^*(y) = \int_{\partial D}(y^*(s)|y(s))ds$ при $y(\cdot)\in C^n(\partial D)$, где $y^*(\cdot)\in L_1^n(\partial D)$, такой, что

$y^*(s)\in \partial\varphi(s,(A\overline{u})(s))$, $x_s^*((A\overline{u})(\cdot)) = \sup_{y(\cdot)\in Q_1}x_s^*(y(\cdot))$ и $\int_D (\upsilon(t)|u(t))dt - x^*(Au) = 0$ при $u(\cdot)\in L_p^n(D)$.

**Теорема 2.** Если $\varphi:\partial D\times R^n \to R_{+\infty}$ нормальный выпуклый интегрант, $c>0$ и $a(\cdot)\in L_1(\partial D)$ такие, что $|\varphi(s,x)|\leq a(s) + c|x|^\alpha$ при $x\in R^n$, $s\in \partial D$, то $\upsilon(\cdot)\in \partial F(\overline{u}(\cdot))$, где $\upsilon(\cdot)\in L_q^n(D)$, в том и только в том случае, когда существует функционал $y^*(\cdot)\in L_{\alpha'}^n(\partial D)$, где $1\leq \alpha < +\infty$, $\dfrac{1}{\alpha}+\dfrac{1}{\alpha'}=1$, такой, что $y^*(s)\in \partial\varphi(s,(A\overline{u})(s))$ и

$$\int_D (\upsilon(t)|u(t))dt - \int_{\partial D}(y^*(s)|(A\overline{u})(s)) = 0 \text{ при } u(\cdot)\in L_p^n(D).$$

**Доказательство.** Так как $|\varphi(s,x)|\leq a(s) + c|x|^\alpha$ при $x\in R^n$, $s\in \partial D$, то функционал $I(y(\cdot)) = \int_{\partial D}\varphi(s,y(s))ds$ непрерывен в $L_\alpha^n(\partial D)$.

Положим $C = \{(u(\cdot),(Au)(\cdot))\in L_p^n(D)\times C^n(\partial D): u(\cdot)\in L_p^n(D)\}$. Тогда по теореме 3.4.1[6] существует функционал $y^*(\cdot)\in L_{\alpha'}^n(\partial D)$ такой, что

$$F^*(\upsilon(\cdot)) = \sup_{u(\cdot)\in L_p^n(D)}\{\int_D (\upsilon(t)|u(t))dt - \int_{\partial D}\varphi(s,(Au)(s))ds\} =$$
$$= \sup_{(u(\cdot),y(\cdot))\in L_p^n(D)\times L_\alpha^n(\partial D)}\{\int_D (\upsilon(t)|u(t))dt - \int_{\partial D}\varphi(s,y(s))ds - \delta_C(u(\cdot),y(\cdot))\} =$$
$$= \sup_{(u(\cdot),y(\cdot))\in L_p^n(D)\times L_\alpha^n(\partial D)}\{\int_D (\upsilon(t)|u(t))dt - \int_{\partial D}(y^*(s)|y(s))ds - \delta_C(u(\cdot),y(\cdot))\} +$$

$$+ \sup_{y(\cdot) \in L_\alpha^n(\partial D)} \{\int_{\partial D} (y^*(s)|y(s))ds - \int_{\partial D} \varphi(s,y(s))ds\} =$$

$$= \sup_{u \in L_p^n(D)} \{\int_D (\upsilon(t)|u(t))dt - \int_{\partial D} (y^*(s)|(Au)(s))ds\} + \int_{\partial D} \varphi^*(s,y^*(s))ds =$$

$$= \begin{cases} \int_{\partial D} \varphi^*(s,y^*(s))ds : \int_D (\upsilon(t)|u(t))dt - \int_{\partial D} (y^*(s)|(Au)(s))ds = 0 \text{ при } u(\cdot) \in L_p^n(D), \\ +\infty : \text{в других случаях.} \end{cases}$$

По определению $\upsilon(\cdot) \in \partial F(\bar{u}(\cdot))$ в том и только в том случае, когда

$F(\bar{u}(\cdot)) + F^*(\upsilon(\cdot)) = \langle \upsilon(\cdot), \bar{u}(\cdot) \rangle$. Так как $\int_D (\upsilon(t)|u(t))dt - \int_{\partial D} (y^*(s)|(Au)(s))ds = 0$ при $u(\cdot) \in L_p^n(D)$,

то имеем, что $\int_D (\upsilon(t)|u(t))dt = \int_{\partial D} (y^*(s)|(Au)(s))ds$ при $u(\cdot) \in L_p^n(D)$. Отсюда следует, что

$$\int_{\partial D} \varphi(s,(A\bar{u})(s))ds + \int_{\partial D} \varphi^*(s,y^*(s))ds = \int_{\partial D} (y^*(s)|(A\bar{u})(s))ds.$$

Поэтому $y^*(s) \in \partial \varphi((A\bar{u})(s))$. Теорема доказана.

### 3.4. О минимизации вариационных задач типа Фредгольма

Пусть $D$ компактное множество в $R^m$ и $\text{int } D \neq \varnothing$, $A : L_p^n(D) \to C^n(D)$ линейный непрерывный компактный оператор. В дальнейшем будем предполагать, что $f : D \times R^n \times R^n \to (-\infty, +\infty]$ и $\varphi : \partial D \times R^n \to (-\infty, +\infty]$ нормальные выпуклые интегранты.

Рассматривается задача минимизации функционала

$$J(u) = \int_{\partial D} \varphi(s,(Au)(s))ds + \int_D f(t,(Au)(t),u(t))dt \tag{1}$$

в пространстве $L_p^n(D)$.

Будем говорить, что кривая $\bar{u}(t)$, $t \in D$, является решением задачи (1), если $|\Phi_0(\bar{u})| < +\infty$ и справедливо неравенство $\Phi_0(u) \geq \Phi_0(\bar{u})$ при любой $u(\cdot) \in L_p^n(D)$.

Рассмотрим функционал

$$\Phi(u,\upsilon) = \int_{\partial D} \varphi(s,(Au)(s))ds + \int_D f(t,(Au)(t),u(t)+\upsilon(t))dt,$$

где $\upsilon(\cdot) \in L_p^n(D)$. Ясно, что $\Phi : L_p^n(D) \times L_p^n(D) \to R \cup \{+\infty\}$ и $\Phi(u,0) = J(u)$. Положим $h(\upsilon) = \inf_{u \in L_p^n(D)} \Phi(u,\upsilon)$. Из предложения 2.5 [11] вытекает, что $h$ выпуклая функция.

Задача (1) называется стабильной, если $h(0)$ конечно и $h$ субдифференцируема в нуле.

**Лемма 1.** Допустим, что $\inf_{u \in L_p^n(D)} J(u)$ конечен и существует такая функция $u_0(\cdot) \in L_p^n(D)$, что функции $f(t,(Au_0)(t)+y,u_0(t))$ и $\varphi(s,(Au_0)(s)+y)$ суммируемы при $y \in R^n$, $|y| < r$ и для некоторого $r > 0$. Тогда функция $h$ субдифференцируема в нуле, т.е. задача (1) стабильна (см. [26]).

**Доказательство.** Рассмотрим функционал

$$I(u) = \int_D f(t,(Au_0)(t)+(Au)(t),u_0(t))dt$$

в $L_p^n(D)$. Обозначим $x_0(t) = (Au_0)(t)$. Так как $I(x) = \int_D f(t,x_0(t)+x(t),u_0(t))dt$ непрерывен в точке $x(t) = 0$ относительно нормированной топологии пространства $C^n(D)$, то существуют такие $\alpha_1 > 0$ и $d_1$, что $I(x) \leq d_1$ при $x \in \{z(\cdot) \in C^n(D) : \|z(\cdot)\|_C \leq \alpha_1\}$. Те же свойства функционала $I_1(x) = \int_{\partial D} \varphi(s,x_0(s)+x(s))ds$ непрерывен в точке $x(s) = 0$ относительно нормированной топологии пространства $C^n(\partial D)$. Поэтому существуют $\alpha_2 > 0$ и $d_2$, что $I_1(x) \leq d_2$ при

$x \in \{z(\cdot) \in C^n(\partial D): \|z(\cdot)\|_{C^n(\partial D)} \le \alpha_2 \}$.

Обозначив $u_\upsilon(\cdot) = u_0(\cdot) - \upsilon(\cdot)$, $\alpha = \min\{\alpha_1, \alpha_2\}$, $x_\upsilon(t) = x_0(t) - (A\upsilon)(t)$ получим

$$h(\upsilon) = \inf_{u \in L_p^n(D)} \Phi(u,\upsilon) \le \Phi(u_\upsilon, \upsilon) \le d_1 + d_2$$

при $\upsilon(\cdot) \in L_p^n(D)$, $\|\upsilon\| \le \alpha$. Так как $h(0) = \inf_{u \in L_p^n(D)} J(u)$ конечен, тогда из предложения 1.5.2 [26] вытекает, что h субдифференцируема в точке нуле. Лемма доказана.

Пусть $\upsilon, w \in R^n$. Положим $f^0(t,z,\upsilon) = \inf_w \{(w|\upsilon) + f(t,z,w)\}$.

**Лемма 2.** Пусть f нормальный выпуклый интегрант на $D \times (R^n \times R^n)$, $u_0(\cdot) \in L_p^n(D)$ и функция $f(t, (Au_0)(t)+y, u_0(t))$ суммируема при $y \in R^n$, $|y| < r$ и для некоторого $r > 0$, $\upsilon^*(\cdot) \in L_q^n(D)$ и функционал $I(u) = \int_D f^0(t, (Au)(t), \upsilon^*(t))dt$ собственный, то функционал I(u) непрерывен в точке $u_0(\cdot)$.

**Доказательство.** Покажем, что для любого $\upsilon^*(\cdot) \in L_q^n(D)$ функция $f^0(t, (Au_0)(t)+y, \upsilon^*(t))$ суммируема при $y \in R^n$, $|y| < r$. Ясно, что

$$f^0(t, (Au_0)(t)+y, \upsilon^*(t)) = \inf_{\omega \in R^n} \{(\omega|\upsilon^*(t)) + f(t, (Au_0)(t)+y, \omega)\} \le$$
$$\le (u_0(t)|\upsilon^*(t)) + f(t, (Au_0)(t)+y, u_0(t)).$$

По условию отсюда вытекает, что функция $f^0(t, (Au)(t)+y, \upsilon^*(t))$ суммируема при $y \in R^n$, $|y| < r$. Тогда из предложения 1.2.5[26] следует, что функционал $u(\cdot) \to \int_D f^0(t, (Au)(t), \upsilon^*(t))dt$ непрерывен в точке $u_0(\cdot)$. Лемма доказана.

**Теорема 1.** Пусть $A: L_p^n(D) \to C^n(D)$ линейный непрерывный оператор, $f: D \times (R^n \times R^n) \to R_{+\infty}$ и $\varphi: \partial D \times R^n \to R_{+\infty}$ нормальные выпуклые интегранты. Для того, чтобы функция $\bar{u}(t)$ среди всех функций $u(\cdot) \in L_p^n(D)$ минимизировала функционал (1), достаточно, чтобы нашлись $\upsilon^*(\cdot) \in L_q^n(D)$, $y^*(\cdot) \in L_1^n(\partial D)$, $u^*(\cdot) \in L_1^n(D)$ такие, что

1) $u^*(t) \in \partial f^0(t, (A\bar{u})(t), \upsilon^*(t))$,
2) $y^*(s) \in \partial \varphi(s, (A\bar{u})(s))$,
3) $f^0(t, (A\bar{u})(t), \upsilon^*(t)) = f(t, (A\bar{u})(t), \bar{u}(t)) + (\bar{u}(t)|\upsilon^*(t))$,
4) $\int_D (\upsilon^*(t) - (A^*u^*)(t)|u(t))dt - \int_{\partial D} (y^*(s)|(Au)(s))ds = 0$ при $u(\cdot) \in L_p^n(D)$,

а если выполнено условие леммы 1 при $u_0(t) = \bar{u}(t)$ и $A: L_p^n(D) \to C^n(D)$ линейный непрерывный компактный оператор, то условия 1)-4) и являются необходимыми.

**Доказательство.** Достаточность теоремы непосредственно проверяется.

**Необходимость.** Из леммы 1 вытекает, что h субдифференцируема в точке нуле. Поэтому из замечания 3.2.3 и из предложения 3.2.4 [26] вытекает, что все решения $\bar{u}(\cdot)$ задачи $\inf\{J(u): u \in L_1^n(D)\}$ и все решения $\upsilon^*(\cdot)$ задачи $\sup_{\upsilon^* \in L_q^n(D)} \{-\Phi^*(0, \upsilon^*)\}$ связаны экстремальным соотношением

$$\Phi(\bar{u}, 0) + \Phi^*(0, -\upsilon^*) = 0 \qquad (2)$$

По определению

$$\Phi^*(0, -\upsilon^*) = \sup_{u, \upsilon \in L_p^n(D)} \{-\int_D (\upsilon(t)|\upsilon^*(t))dt - \int_D f(t, (Au)(t), u(t)+\upsilon(t))dt - \qquad (3)$$
$$- \int_{\partial D} \varphi(s, (Au)(s))ds\} = \sup_{u, \upsilon \in L_p^n(D)} \{\int_D (u(t)|\upsilon^*(t))dt - \int_D (u(t)+\upsilon(t)|\upsilon^*(t))dt -$$

$$-\int_D f(t,(Au)(t),u(t)+\upsilon(t))dt - \int_{\partial D}\varphi(s,(Au)(s))ds\} =$$

$$= \sup_{u\in L_p^n(D)}\left\{\int_D (u(t)|\upsilon^*(t))dt - \int_D f^0(t,(Au)(t),\upsilon^*(t))dt - \int_{\partial D}\varphi(s,(Au)(s))ds\right\}.$$

Обозначим $J_1(u) = \int_D f^0(t,(Au)(t),\upsilon^*(t))dt$, $J_2(u) = \int_{\partial D}\varphi(s,(Au)(s))ds$. Из (2) и (3) вытекает, что $J_1$ и $J_2$ собственные функционалы. При условии теоремы 1 из леммы 2 вытекает, что функционал $J_1$ непрерывен в точке $\bar{u}(\cdot)$. Аналогично имеем, что $J_2(u(\cdot))$ непрерывен в точке $\bar{u}(\cdot)$.

По соотношению (2) имеем
$$\int_{\partial D}\varphi(s,(A\bar{u})(s))ds + \int_D f(t,(A\bar{u})(t),\bar{u}(t))dt + \Phi^*(0,-\upsilon^*) = 0. \qquad (4)$$

Положив
$$S(u) = \int_D f^0(t,(A\bar{u})(t),\upsilon^*(t))dt + \int_{\partial D}\varphi(s,(A\bar{u})(s))ds$$

имеем, что $\Phi^*(0,-\upsilon^*) = S(\upsilon^*)$. Так как
$$S^*(\upsilon^*) \geq \int_D (\bar{u}(t)|\upsilon^*(t))dt - S(\bar{u}), \qquad S(\bar{u}) \leq \int_D (\bar{u}(t)|\upsilon^*(t))dt + J(\bar{u}).$$

Тогда из соотношения (2) следует, что
$$S^*(\upsilon^*) = \int_D (\bar{u}(t)|\upsilon^*(t))dt - S(\bar{u}), \qquad S(\bar{u}) = \int_D (\bar{u}(t)|\upsilon^*(t))dt + J(\bar{u}).$$

Из второго соотношения вытекает, что
$$\int_D f^0(t,(A\bar{u})(t),\upsilon^*(t))dt + \int_{\partial D}\varphi(s,(A\bar{u})(s))ds -$$
$$- \int_{\partial D}\varphi(s,(A\bar{u})(s))ds - \int_D f(t,(A\bar{u})(t),\bar{u}(t))dt = \int_D (\bar{u}(t)|\upsilon^*(t))dt,$$

т.е.
$$\int_D f^0(t,(A\bar{u})(t),\upsilon^*(t))dt - \int_D f(t,(A\bar{u})(t),\bar{u}(t))dt = \int_D (\bar{u}(t)|\upsilon^*(t))dt,$$

Отсюда используя неравенства Юнга-Фенхеля получим
$$f^0(t,(A\bar{u})(t),\upsilon^*(t)) - f(t,(A\bar{u})(t),\bar{u}(t)) = (\bar{u}(t)|\upsilon^*(t)),$$

т.е. $\qquad f^0(t,(A\bar{u})(t),\upsilon^*(t)) = f(t,(A\bar{u})(t),\bar{u}(t)) + (\bar{u}(t)|\upsilon^*(t)).$

Из равенства $S^*(\upsilon^*) = \int_D (\bar{u}(t)|\upsilon^*(t))dt - S(\bar{u})$ следует, что $\upsilon^* \in \partial S(\bar{u})$. Из теоремы Моро-Рокафеллара (см.[6]) имеем, что $\partial S(\bar{u}) = \partial J_1(\bar{u}) + \partial J_2(\bar{u})$. Поэтому существуют функции $\upsilon_1^* \in \partial J_1(\bar{u})$, $\upsilon_2^* \in \partial J_2(\bar{u})$ такие, что $\upsilon^* = \upsilon_1^* + \upsilon_2^*$. Так как $\upsilon_1^* \in \partial J_1(\bar{u})$, то из теоремы 3.2.1 следует, что $\upsilon_1^*(\cdot) \in L_q^n(D)$ принадлежит $\partial J_1(\bar{u})$ в том и только в том случае, когда существует функция $u^*(t) \in \partial f^0(t,(A\bar{u})(t),\bar{u}(t))$, $u^*(\cdot) \in L_1^n(D)$, такая, что $\upsilon_1^*(\cdot) = (A^*u^*)(t)$.

По теореме 3.3.1 $\upsilon_2^*(\cdot) \in L_q^n(D)$ принадлежит $\partial J_2(\bar{u})$ в том и только в том случае, когда существует функционал $x^* \in L_\infty^n(\partial D)^*$, где $x^* = x_a^*$, $x_a^*(y) = \int_{\partial D}(y^*(s)|y(s))ds$ при $y(\cdot) \in L_\infty^n(\partial D)$, где $y^*(\cdot) \in L_1^n(D)$, такой, что $y^*(s) \in \partial \varphi(s,(A\bar{u})(s))$ и $\int_D (\upsilon_2^*(s)|u(t))dt - \int_{\partial D}(y^*(s)|(Au)(s))ds = 0$ при $u(\cdot) \in L_p^n(D)$. Ясно, что $\upsilon_2^*(t) = \upsilon^*(t) - \upsilon_1^*(t) = \upsilon^*(t) - (A^*u^*)(t)$. Теорема доказана.

**Следствие 1.** Пусть $A: L_p^n(D) \to C^n(D)$ линейный непрерывный оператор, $\varphi: \partial D \times R^n \to R_{+\infty}$ и $f: D \times (R^n \times R^n) \to R_{+\infty}$ нормальные выпуклые интегранты. Для того, чтобы функция $\bar{u}(t)$ среди всех функций $u(\cdot) \in L_p^n(D)$ минимизировала функционал (1), достаточно, чтобы нашлись $\upsilon^*(\cdot) \in L_q^n(D)$, $y^*(\cdot) \in L_1^n(\partial D)$, $u^*(\cdot) \in L_1^n(D)$ такие, что

1) $(u^*(t), -\upsilon^*(t)) \in \partial f(t, (A\overline{u})(t), \overline{u}(t))$,
2) $y^*(s) \in \partial \varphi(s, (A\overline{u})(s))$,
3) $\int_D (\upsilon^*(t) - (A^*u^*)(t) | u(t)) dt - \int_{\partial D} (y^*(s) | (Au)(s)) ds = 0$ при $u(\cdot) \in L_p^n(D)$,

а если выполнено условие леммы 1 при $u_0(t) = \overline{u}(t)$ и $A : L_p^n(D) \to C^n(D)$ линейный непрерывный компактный оператор, то условия 1)-3) и являются необходимыми.

Пусть $A : L_p^n(D) \to C^n(D)$ линейный непрерывный оператор, $\varphi : \partial D \times R^n \to R_{+\infty}$ и $f : D \times (R^n \times R^n) \to R_{+\infty}$ нормальные выпуклые интегранты. Аналогично можно рассмотреть (см. п.1.5) минимизацию функционала

$$J(u) = \int_{\partial D} \varphi(s, (Au)(s)) ds + \int_D f(t, (Au)(t), u(t)) dt$$

при следующих ограничениях

$$u(t) \in F(t, (Au)(t)), \quad t \in D, \quad u(\cdot) \in L_p^n(D).$$

Пусть $A : L_p^n(D) \to L_\infty^n(D)$ линейный непрерывный оператор. Аналогично можно рассмотреть (см. п.1.5) минимизацию функционала

$$J(u) = \int_{\partial D} \varphi(s, (Au)(s)) ds + \int_D f(t, (Au)(t), u(t)) dt$$

в пространстве $u(\cdot) \in L_p^n(D)$.

### 3.5. Невыпуклая экстремальная задача для операторных включений

Пусть $D$ компактное множество в $R^m$ и $\text{int } D \neq \varnothing$, $A : L_p^n(D) \to C^n(D)$ линейный непрерывный компактный оператор. В дальнейшем будем предполагать, что $f : D \times R^n \times R^n \to R_{+\infty}$ и $\varphi : \partial D \times R^n \to R_{+\infty}$ нормальные интегранты.

Рассматривается задача минимизации функционала

$$J(u) = \int_{\partial D} \varphi(s, (Au)(s)) ds + \int_D f(t, (Au)(t), u(t)) dt \qquad (1)$$

при следующих ограничениях

$$u(t) \in F(t, (Au)(t)), \quad t \in D, \quad u(\cdot) \in L_p^n(D). \qquad (2)$$

Положим $\psi(t, x, y) = \inf\{|z - y| : z \in F(t, x)\}$ и

$$J_r(u) = \int_{\partial D} \varphi(s, (Au)(s)) ds + \int_D f(t, (Au)(t), u(t)) dt + r(\int_D \psi^p(t, (Au)(t), u(t)) dt)^{\frac{1}{p}}.$$

Решение задачи (1),(2) обозначим через $\overline{u}(t)$. Обозначим

$$F_0(u) = (\int_D \psi^p(t, (Au)(t), u(t)) dt)^{\frac{1}{p}}.$$

**Теорема 1.** Пусть $F : D \times R^n \to \text{comp} R^n$ - многозначное отображение, $F(t, x)$ измеримо по $t$, $f : D \times R^n \times R^n \to R_{+\infty}$ и $\varphi : \partial D \times R^n \to R_{+\infty}$ нормальные интегранты, существует функция $M(\cdot) \in L_p(D)$, где $M(t) > 0$ такая, что

$$\rho_x(F(t, x), F(t, x_1))) \leq M(t) |x - x_1|$$

при $x, x_1 \in R^n$. Кроме того, пусть $A : L_p^n(D) \to C^n(D)$ линейный непрерывный оператор и $\|A\| (\int_D M^p(s) ds)^{\frac{1}{p}} < 1$, существуют неотрицательные функции $k(\cdot) \in L_1(D)$, $k_1(\cdot) \in L_q(D)$, где

$\frac{1}{p} + \frac{1}{q} = 1$, и $k_2(\cdot) \in L_1(\partial D)$ такие, что

$$|f(t,x_1,y_1) - f(t,x_2,y_2)| \le k(t)|x_1 - x_2| + k_1(t)|y_1 - y_2|$$

при $x_1, y_1, x_2, y_2 \in R^n$,

$$|\varphi(s,x) - \varphi(s,y)| \le k_2(s)|x - y|$$

при $x, y \in R^n$. Тогда если $\bar{u}(\cdot) \in L_p^n(D)$ является решением задачи (1),(2), то существует число $r_0 > 0$ такое, что $\bar{u}(t)$ минимизирует функционал $J_r(u)$ в пространстве $L_p^n(D)$ при $r \ge r_0$.

**Доказательство.** Пусть $\tilde{u}(\cdot) \in L_p^n(D)$, $\tilde{x}(t) = (A\tilde{u})(t)$. Положим $\rho(t) = \psi(t,(A\tilde{u})(t),\tilde{u}(t))$. По теореме 3.3.1 существует решение $u(\cdot) \in L_p^n(D)$ задачи $u(t) \in F(t,(Au)(t))$ такое, что

$$|u(t) - \tilde{u}(t)| \le \rho(t) + \frac{\|A\|(\int_D \rho^p(s)ds)^{\frac{1}{p}} M(t)}{1 - \|A\|(\int_D M^p(s)ds)^{\frac{1}{p}}}, \qquad |(Au)(t) - (A\tilde{u})(t)| \le \frac{\|A\|(\int_D \rho^p(s)ds)^{\frac{1}{p}}}{1 - \|A\|(\int_D M^p(s)ds)^{\frac{1}{p}}}.$$

Ясно, что

$$|J(\tilde{u}) - J(u)| \le \left| \int_{\partial D} \varphi(s,(A\tilde{u})(s))ds - \int_{\partial D} \varphi(s,(Au)(s))ds \right| +$$

$$+ \left| \int_D f(t,(A\tilde{u})(t),\tilde{u}(t))dt - \int_D f(t,(Au)(t),u(t))dt \right| \le$$

$$\le \int_{\partial D} k_2(s)|A(\tilde{u}(s) - u(s))|ds + \int_D (k(t)|A(\tilde{u}(t) - u(t))| + k_1(t)|\tilde{u}(t) - u(t)|)dt \le$$

$$\le \int_{\partial D} k_2(s)ds \; \frac{\|A\|(\int_D \rho^p(s)ds)^{\frac{1}{p}}}{1 - \|A\|(\int_D M^p(s)ds)^{\frac{1}{p}}} + \int_D k(t)dt \; \frac{\|A\|(\int_D \rho^p(s)ds)^{\frac{1}{p}}}{1 - \|A\|(\int_D M^p(s)ds)^{\frac{1}{p}}} +$$

$$+ (\int_D k_1^q(t)dt)^{\frac{1}{q}} (\int_D \rho^p(t)dt)^{\frac{1}{p}} + \frac{\|A\|(\int_D \rho^p(s)ds)^{\frac{1}{p}} \int_D k_1(t)M(t)dt}{1 - \|A\|(\int_D M^p(s)ds)^{\frac{1}{p}}} \le$$

$$\le (\int_D \rho^p(s)ds)^{\frac{1}{p}} (\frac{\|A\|}{1 - \|A\|(\int_D M^p(s)ds)^{\frac{1}{p}}} (\int_D k_2(s)ds + \int_D k(t)dt + \int_D k_1(t)M(t)dt) + (\int_D k_1^q(t)dt)^{\frac{1}{q}}).$$

Положив

$$r_0 = \frac{\|A\|}{1 - \|A\|(\int_D M^p(s)ds)^{\frac{1}{p}}} (\int_D k_2(s)ds + \int_D k(t)dt + \int_D k_1(t)M(t)dt) + (\int_D k_1^q(t)dt)^{\frac{1}{q}}$$

имеем, что $|J(\tilde{u}) - J(u)| \le r_0 F_0(\tilde{u})$.

Покажем, что $\bar{u}(t)$ минимизирует также функционал $J_r(u)$ в пространстве $L_p^n(D)$ при $r \ge r_0$. Предположим противное. Пусть существует функция $\upsilon(\cdot) \in L_p^n(D)$ такая, что $J_r(\upsilon) < J(\bar{u})$. Так как для $\upsilon(\cdot)$ существует решение $\upsilon_0(\cdot)$ задачи и (2) такое, что $|J(\upsilon) - J(\upsilon_0)| \le rF_0(\upsilon)$ при $r \ge r_0$, то

$$J(\upsilon_0) \le J(\upsilon) + rF_0(\upsilon) = J_r(\upsilon) < J(\bar{u}).$$

Получим противоречие. Полученное противоречие означает, что $\bar{u}(t)$ минимизирует также функционал $J_r(u)$ в пространстве $L_p^n(D)$ при $r \geq r_0$. Теорема доказана.

**Теорема 2.** Пусть выполняется условие теоремы 1 при $p=1$, $A: L_1^n(D) \to C^n(D)$ линейный непрерывный компактный оператор и функция $\bar{u}(t)$ среди всех решений задачи (2) минимизирует функционал (1). Тогда существуют функции $\upsilon^*(\cdot) \in L_\infty^n(D)$, $y^*(\cdot) \in L_1^n(\partial D)$, $u^*(\cdot) \in L_1^n(D)$ такие, что

1) $(u^*(t), -\upsilon^*(t)) \in \partial(f(t, (A\bar{u})(t), \bar{u}(t)) + r\psi(t, (Au)(t), u(t)))$ при $r \geq r_0$,

2) $y^*(s) \in \partial\varphi(s, (A\bar{u})(s))$,

3) $\int_D (\upsilon^*(t) - (A^*u^*)(t) | u(t)) dt - \int_{\partial D} (y^*(s) | (Au)(s)) ds = 0$ при $u(\cdot) \in L_1^n(D)$.

**Доказательство.** По теореме 1 существует число $r_0 > 0$ такое, что $\bar{u}(t)$ минимизирует функционал $J_r(u)$ в пространстве $L_1^n(D)$ при $r \geq r_0$, т.е. $J_r(\bar{u}) \leq J_r(u)$ при $u \in L_1^n(D)$. Поэтому $J_r^0(\bar{u}; u) \geq J_r^+(\bar{u}; u) \geq 0$ при $u \in L_1^n(D)$, где

$$J_r^0(\bar{u}; u) = \overline{\lim_{\upsilon \to \bar{u}, \lambda \downarrow 0}} \frac{1}{\lambda}(J_r(\upsilon + \lambda u) - J_r(\upsilon)), \qquad J_r^+(\bar{u}; u) = \overline{\lim_{\lambda \downarrow 0}} \frac{1}{\lambda}(J_r(\bar{u} + \lambda u) - J_r(\bar{u})).$$

Обозначим $\bar{f}(t, x, y) = f(t, x, y) + r\psi(t, x, y)$. Так как существует суммируемая функция $M(t) > 0$ такая, что $\rho_X(F(t, x), F(t, y)) \leq M(t)|x - y|$ при $x, y \in R^n$, то

$$|\psi(t, x_1, y_1) - \psi(t, x_2, y_2)| \leq M(t)|x_1 - x_2| + |y_1 - y_2|$$

при $x_1, y_1, x_2, y_2 \in R^n$ (см. [17]). Поэтому

$$|\bar{f}(t, x_1, y_1) - \bar{f}(t, x_2, y_2)| \leq (k(t) + M(t))|x_1 - x_2| + (k_1(t) + 1)|y_1 - y_2|$$

при $x_1, y_1, x_2, y_2 \in R^n$. Тогда применяя лемму Фату (см. [27]) имеем, что

$$J_r^0(\bar{u}; u) \leq \int_D \bar{f}^0(t, (A\bar{u})(t), \bar{u}(t); (Au)(t), u(t)) dt + \int_D \varphi^0(s, (A\bar{u})(s); (Au)(s)) ds.$$

Положив

$$I(u) = \int_D \bar{f}^0(t, (A\bar{u})(t), \bar{u}(t); (Au)(t), u(t)) dt + \int_D \varphi^0(s, (A\bar{u})(s); (Au)(s)) ds$$

имеем, что функция $u(t) = 0$ минимизирует функционал $I(u)$ в пространстве $L_1^n(D)$. Легко проверяется, что для функционала $I(u)$ выполняется условие следствия 3.4.1. Поэтому по следствию 3.4.1 существуют функции $\upsilon^*(\cdot) \in L_\infty^n(D)$, $y^*(\cdot) \in L_1^n(\partial D)$, $u^*(\cdot) \in L_1^n(D)$ такие, что выполняется соотношение 1)-3) теоремы 2. Теорема доказана.

Отметим, что если $D$ компактное метрическое пространство, $(D, \Sigma, \mu)$ пространство с конечной положительной мерой, то все результаты гл.3 остаются также верными.

### 3.6. О субдифференцируемости интегрального функционала с оператором

Пусть $D$ компактное метрическое пространство, $(D, \Sigma, \mu)$ пространство с конечной положительной мерой, $f: D \times R^n \to R_{+\infty}$ нормальный выпуклый интегрант, т.е. $f$ -такая функция из $D \times R^n$ в $R_{+\infty}$, что $\mathrm{epf}_t = \{(z, \alpha) \in R^n \times R : f(t, z) \leq \alpha\}$ -выпуклое замкнутое множество и $t \to \mathrm{epf}_t$ измеримо.

Пусть $A: L_{p_1}^n(D) \to L_{p_2}^n(D)$ линейный непрерывный оператор. Рассмотрим субдифференцируемость интегрального функционала $J(u(\cdot)) = \int_D f(t, (Au)(t)) dt$ в $u \in L_{p_1}^n(D)$, $1 \leq p_1 < +\infty, 1 \leq p_2 < +\infty$.

Отметим, что если $f: D \times R^n \to R_{+\infty}$ нормальный выпуклый интегрант и $A: L_{p_1}^n(D) \to L_{p_2}^n(D)$ линейный оператор, то $J: L_{p_1}^n(D) \to R_{+\infty}$ выпуклый функционал.

Субградиентами $J$ в точке $u_0(\cdot) \in L_{p_1}^n(D)$ являются по определению, элементы $\upsilon^* \in L_{q_1}^n(D)$, $\dfrac{1}{p_1} + \dfrac{1}{q_1} = 1$, для которых

$$J(u(\cdot)) - J(u_0(\cdot)) \geq \langle \upsilon^*(\cdot), u(\cdot) - u_0(\cdot) \rangle$$

при всех $u(\cdot) \in L_{p_1}^n(D)$ или, что равносильно,

$$J^*(\upsilon^*) + J(u_0(\cdot)) = \langle \upsilon^*(\cdot), u_0(\cdot) \rangle,$$

где $\langle \upsilon^*(\cdot), u_0(\cdot) \rangle = \int_D (\upsilon^*(t) | u_0(t)) dt$.

Множество всех таких субградиентов обозначается через $\partial J(u_0(\cdot))$ и называется субдифференциалом функционала $J$ в точке $u_0(\cdot)$.

В этом параграфе устанавливается связь между $\partial J(u_0(\cdot))$ и $\partial f_t((Au_0)(t))$.

**Лемма 2.** Пусть $A: L_{p_1}^n(D) \to L_{p_2}^n(D)$ линейный непрерывный компактный оператор, $f: D \times R^n \to R_{+\infty}$ нормальный интегрант, $x \to f(t,x)$ сублинейная функция в $R^n$, существуют функция $a(\cdot) \in L_1^n(D)$ и $c > 0$ такие, что $|f(t,x)| \leq a(x) + c|x|^{p_2}$ при $x \in R^n$ и $y^*(\cdot) \in L_{q_1}^n(D)$, где $\dfrac{1}{p_1} + \dfrac{1}{q_1} = 1$. Тогда

$$J^*(y^*(\cdot)) = \begin{cases} 0, & y^*(\cdot) = (A^*u^*)(t), \ u^*(t) \in \partial f(t,0), \ u^*(\cdot) \in L_1^n(D), \\ +\infty, & \text{в других случаях.} \end{cases}$$

**Доказательство.** По условию леммы 2 имеем, что $J(u(\cdot)) = \int_D f(t,(Au)(t))dt$ сублинейная функция в $L_{p_1}^n(D)$. Так как $|f(t,x)| \leq a(x) + c|x|^{p_2}$ при $x \in R^n$, то по теореме 3.2.1[6] функционал $x(\cdot) \to \int_D f(t,x(t))dt$ непрерывен в $L_{p_2}^n(D)$.

Ясно, что $I_0(y) = \int_D f(t,y(t))dt$ сублинейная непрерывная функция в $L_{p_2}^n(D)$. Поэтому по предложению 3 ([6], стр.210) и по теореме 3 ([6], стр.362) $\partial I_0(0)$ слабо$^*$ компактно и $\partial I_0(0) \subset L_{q_2}^n(D)$. Из теоремы 6.4.8[8] следует, что

$$I_0(y(\cdot)) = \max_{u^* \in \partial I_0(0)} \langle u^*, y \rangle = \max_{u^* \in \partial I_0(0)} \int_D (u^*(t)|y(t))dt$$

при $y(\cdot) \in L_{p_2}^n(D)$. Известно, что (см. [6]) $u^* \in \partial I_0(0)$ в том и только в том случае, когда $u^*(t) \in \partial f(t,0)$ и $u^*(\cdot) \in L_{q_2}^n(D)$.

Ясно, что $J(u(\cdot)) = \max_{u^* \in \partial I_0(0)} \int_D (u^*(t)|(Au)(t))dt$ при $u(\cdot) \in L_{p_1}^n(D)$. По условию $A: L_{p_1}^n(D) \to L_{p_2}^n(D)$ линейный непрерывный компактный оператор. Если $y^*(\cdot) \in L_{q_1}^n(D)$, где $\dfrac{1}{p_1} + \dfrac{1}{q_1} = 1$, то применяя теорему 6.2.7[12] имеем

$$J^*(y^*(\cdot)) = \sup_{u(\cdot) \in L_{p_1}^n(D)} \{\int_D (y^*(t)|u(t))dt - \int_D f(t,(Au)(t))dt\} =$$

$$= \sup_{u(\cdot) \in L_{p_1}^n(D)} \{\int_D (y^*(t)|u(t))dt - \sup_{u^* \in \partial I_0(0)} \int_D (u^*(t)|(Au)(t))dt\} =$$

$$= \sup_{u(\cdot) \in L_{p_1}^n(D)} \inf_{u^* \in \partial I_0(0)} \{ \int_D (y^*(t) - (A^*u^*)(t) | u(t)) dt =$$

$$= \inf_{u^* \in \partial I_0(0)} \sup_{u(\cdot) \in L_{p_1}^n(D)} \{ \int_D (y^*(t) - (A^*u^*)(t) | u(t)) dt =$$

$$= \begin{cases} 0, & y^*(\cdot) = (A^*u^*)(t), \ u^*(t) \in \partial f(t,0), \ u^*(\cdot) \in L_{q_2}^n(D), \\ +\infty, & \text{в других случаях.} \end{cases}$$

Лемма доказана.

**Следствие 1.** Если удовлетворяется условие леммы 2, то $y^*(\cdot) \in \partial J(0)$ в том и только в том случае, когда существует функция $u^*(t) \in \partial f(t,0)$, $u^*(\cdot) \in L_{q_2}^n(D)$, такая, что $y^*(\cdot) = (A^*u^*)(t)$.

**Теорема 1.** Пусть $f: D \times R^n \to R_{+\infty}$ нормальный выпуклый интегрант, существуют функция $a(\cdot) \in L_1^n(D)$ и $c > 0$ такие, что $|f(t,x)| \leq a(x) + c|x|^{p_2}$ при $x \in R^n$, $A: L_{p_1}^n(D) \to L_{p_2}^n(D)$ линейный непрерывный компактный оператор, $\bar{u}(\cdot) \in L_{p_1}^n(D)$, $\bar{x}(t) = (A\bar{u})(t)$. Тогда $\partial J(\bar{u})$ непусто и $y^*(\cdot) \in L_{q_1}^n(D)$ принадлежит $\partial J(\bar{u})$ в том и только в том случае, когда существует функция $u^*(t) \in \partial f(t, \bar{x}(t))$, $u^*(\cdot) \in L_{q_2}^n(D)$, такая, что $y^*(\cdot) = (A^*u^*)(t)$.

**Доказательство.** Ясно, что

$$J'(\bar{u}(\cdot); u(\cdot)) = \lim_{\lambda \downarrow 0} \frac{1}{\lambda} (J(\bar{u}(\cdot) + \lambda u(\cdot)) - J(\bar{u}(\cdot)) =$$

$$= \lim_{\lambda \downarrow 0} \frac{1}{\lambda} \int_{t_0}^{T} (f(t, (A\bar{u})(t) + \lambda(Au)(t)) - f(t, (A\bar{u})(t))) dt.$$

Из доказательства предложения 4.1[11] имеем

$$f(t, (A\bar{u})(t)) - f(t, (A\bar{u})(t) - (Au)(t)) \leq \frac{1}{\lambda} (f(t, (A\bar{u})(t) + \lambda(Au)(t)) - f(t, (A\bar{u})(t))) \leq$$

$$\leq f(t, (A\bar{u})(t) + (Au)(t)) - f(t, (A\bar{u})(t))$$

при $\lambda \in (0,1)$. Так как

$$|f(t, (A\bar{u})(t)) - f(t, (A\bar{u})(t) - (Au)(t))| \leq |f(t, (A\bar{u})(t))| + |f(t, (A\bar{u})(t) - (Au)(t))| \leq$$

$$\leq 2a(x) + c|A\bar{u})(t)|^{p_2} + c|(A\bar{u})(t) - (Au)(t)|^{p_2},$$

$$|f(t, (A\bar{u})(t) + (Au)(t)) - f(t, (A\bar{u})(t))| \leq |f(t, (A\bar{u})(t))| + |f(t, (A\bar{u})(t) + (Au)(t))| \leq$$

$$\leq 2a(x) + c|(A\bar{u})(t))|^{p_2} + c|(A\bar{u})(t) + (Au)(t)|^{p_2},$$

то $f(t, (A\bar{u})(t)) - f(t, (A\bar{u})(t) - (Au)(t))$ и $f(t, (A\bar{u})(t) + (Au)(t)) - f(t, (A\bar{u})(t))$ суммируемы при $u(\cdot) \in L_{p_1}^n(D)$. Тогда используя теорему Лебега получим

$$J'(\bar{u}(\cdot); u(\cdot)) = \int_D f'(t, (A\bar{u})(t); (Au)(t)) dt.$$

Если учесть, что $\partial J(\bar{u}(\cdot) = \partial J'(\bar{u}(\cdot); 0)$, $\partial f'(t, (A\bar{u})(t)) = \partial f'(t, (A\bar{u})(t); 0)$, то из следствия 1 имеем, что $y^*(\cdot) \in L_{q_1}^n(D)$ принадлежит $\partial J(\bar{u})$ в том и только в том случае, когда существует функция $u^*(t) \in \partial f(t, \bar{x}(t))$, $u^*(\cdot) \in L_{q_2}^n(D)$, такая, что $y^*(t) = (A^*u^*)(t)$. Теорема доказана.

Пусть $D \subset R^m$ ограниченная область. Если $A: W_{p_1,1}^n(D) \to W_{p_2,1}^n(D)$ линейный непрерывный компактный оператор, то аналогично можно изучить субдифференцируемость интегрального функционала.

**Приложение**

### 1. Процессы, описываемые системой первого порядка с частными производными с обобщенным управлением

Пусть $U \subset R^r$ непустое компактное множество, $f^i : [0,S] \times [0,T] \times R^n \times U \to R$, $i = 0,1,\ldots,n$, $x(s,t) = (x^1(s,t),\ldots,x^n(s,t))$, $\varphi^i : [0,T] \to R$, $g^i : [0,S] \to R$.

Рассмотрим минимизацию функционала

$$J(u) = \int_0^S \int_0^T f^0(s,t,x(s,t),u(s,t))ds dt \qquad (1)$$

при условиях

$$\begin{aligned} x_s^i(s,t) &= f^i(s,t,\ x(s,t),\ u(s,t)), & i &= 1,2,\ldots,m, \\ x_t^i(s,t) &= f^i(s,t,\ x(s,t),\ u(s,t)), & i &= m+1,\ldots,n, \end{aligned} \qquad (2)$$

$$\begin{aligned} x^i(0,t) &= \varphi^i(t), & i &= 1,2,\ldots,m, \\ x^i(s,0) &= g^i(s), & i &= m+1,\ldots,n, \end{aligned} \qquad (3)$$

где $u(s,t) \in U$ измеримая функция, $\varphi^i(\cdot) \in L_2[0,T]$, $g^i(\cdot) \in L_2[0,S]$.

Положим $D = [0,S] \times [0,T]$, $V = \{u : u(s,t) \in U, u(s,t) \text{ измеримы}\}$.

Далее считаем, что функции $f^i(s,t,x,u)$, $f_x^i(s,t,x,u)$, $i = 0,1,2,\ldots,n$, в множестве $(s,t,x,u) \in D \times R^n \times U$ непрерывны и по переменной x удовлетворяют условию Липшица с постоянной L.

По условию имеем

$$\left|f^i(s,t,x,u)\right| \leq \left|f^i(s,t,x,u) - f^i(s,t,0,u)\right| + \left|f^i(s,t,0,u)\right| \leq L|x| + \sup_{(s,t,u) \in D \times U}\left|f^i(s,t,0,u)\right| \qquad (4)$$

при всех $(s,t,x,u) \in D \times R^n \times U$, откуда следует, что $f^i(s,t,x(s,t),u(s,t)) \in L_2(D)$ при любых $x(\cdot) \in L_2(D)$ и $u(\cdot) \in V$, где $i = 0,1,2,\ldots,n$. Аналогично имеем, что $f_x^i(s,t,x(s,t),u(s,t)) \in L_2(D)$ при любых $x(\cdot) \in L_2(D)$ и $u(\cdot) \in V$, где $i = 0,1,2,\ldots,n$.

Из леммы 3.2.1[17] следует, что при выполнении перечисленных условий решение задачи (2),(3) при измеримых функции $u(s,t) \in U$ существует и единственно.

Обозначим $\Im = \{\sigma(s,t) \in \text{rpm}(U) : \sigma(s,t)(U) = 1 \text{ п.в. } (s,t) \in D,\ \sigma(s,t) \text{ измеримы}\}$, где $\text{rpm}(U)$ множество всех вероятностных мер Радона в U (см.[4], стр. 294).

Рассмотрим обобщенную задачу для задачи (1)-(3), т.е. рассмотрим минимизацию функционала

$$J(u) = \int_0^S \int_0^T f^0(s,t,x(s,t),\sigma(s,t))ds dt \qquad (5)$$

при условиях

$$\begin{aligned} x^i(s,t) &= \varphi^i(t) + \int_0^s f^i(\nu,t,\ x(\nu,t),\ \sigma(\nu,t))d\nu, & i &= 1,2,\ldots,m, \\ x^i(s,t) &= g^i(s) + \int_0^t f^i(s,\tau,\ x(s,\tau),\ \sigma(s,\tau))d\tau, & i &= m+1,\ldots,n. \end{aligned} \qquad (6)$$

Под решением задачи (6) соответствующим управлению $\sigma \in \Im$ будем понимать вектор-функцию $x(s,t) = (x^1(s,t),\ldots,x^n(s,t)) \in L_2^n(D)$, имеющие производную $(x_s^1(s,t),\ldots,x_s^m(s,t), x_t^{m+1}(s,t),\ldots,x_t^n(s,t)) \in L_2^n(D)$, которая удовлетворяет (6).

Обозначим

$$Y = \{x(s,t) = (x^1(s,t),\ldots,x^n(s,t)) \in L_2^n(D) : (x_s^1(s,t),\ldots,x_s^m(s,t), x_t^{m+1}(s,t),\ldots,x_t^n(s,t)) \in L_2^n(D)\},$$

$$F(x,\sigma) = (\varphi^1(t) + \int_0^s f^1(\nu,t, x(\nu,t), \sigma(\nu,t))d\nu,\ldots,\varphi^m(t) + \int_0^s f^m(\nu,t, x(\nu,t), \sigma(\nu,t))d\nu,$$

$$g^{m+1}(s) + \int_0^t f^{m+1}(s,\tau, x(s,\tau), \sigma(s,\tau))d\tau,\ldots, g^n(s) + \int_0^t f^n(s,\tau, x(s,\tau), \sigma(s,\tau))d\tau),$$

$$g_0(x,\sigma) = \int_0^S \int_0^T f^0(s,t,x(s,t),\sigma(s,t))dsdt.$$

Отметим, что

$$(F_x(x,\sigma)\Delta x)(s,t) = (\int_0^s f_x^1(\nu,t, x(\nu,t), \sigma(\nu,t))\Delta x(\nu,t)d\nu,\ldots,\int_0^s f_x^m(\nu,t, x(\nu,t), \sigma(\nu,t))\Delta x(\nu,t)d\nu,$$

$$\int_0^t f_x^{m+1}(s,\tau, x(s,\tau), \sigma(s,\tau))\Delta x(s,\tau)d\tau,\ldots, \int_0^t f_x^n(s,\tau, x(s,\tau), \sigma(s,\tau))\Delta x(s,\tau)d\tau),$$

$$g_{0x}(\overline{x},\overline{\sigma})\Delta x = \int_0^S \int_0^T f_x^0(s,t,\overline{x}(s,t),\overline{\sigma}(s,t))\Delta x(s,t)dsdt.$$

Обозначим $h^1 = (f^1,\ldots,f^m)$, $h^2 = (f^{m+1},\ldots,f^n)$.

Покажем, что условие теоремы 5.3.2[4] выполняется. В теореме 5.2.3[4] положив $m = 0$ покажем, что отображения

$$(x,\theta) \to F(x,\overline{\sigma}+\theta(\sigma-\overline{\sigma})): Y \times [0,1] \to Y \quad \text{и} \quad (x,\theta) \to g_0(x,\overline{\sigma}+\theta(\sigma-\overline{\sigma})): Y \times [0,1] \to Y$$

непрерывны и имеют непрерывные производные в точке $(\overline{x},0)$. Пусть $(x_n,\theta_n) \in Y \times [0,1]$ и $(x_n,\theta_n) \to (x,\theta)$. Тогда из оценки (4) и теоремы 4.2.8[4] имеем, что

$$g_0(x_n, \overline{\sigma} + \theta_n(\sigma - \overline{\sigma})) = \int_0^S \int_0^T f^0(s,t,x_n(s,t), \overline{\sigma}(s,t) + \theta_n(\sigma(s,t) - \overline{\sigma}(s,t))) ds dt =$$

$$= \int_0^S \int_0^T \int_U f^0(s,t,x_n(s,t), r)((1-\theta_n)\overline{\sigma}(s,t) + \theta_n \sigma(s,t))(dr) ds dt =$$

$$= (1-\theta_n) \int_0^S \int_0^T \int_U f^0(s,t,x_n(s,t), r)\overline{\sigma}(s,t)(dr) ds dt + \theta_n \int_0^S \int_0^T \int_U f^0(s,t,x_n(s,t),r)\sigma(s,t))(dr) ds dt =$$

$$= (1-\theta_n)(\int_0^S \int_0^T \int_U f^0(s,t,x_n(s,t), r)\overline{\sigma}(s,t)(dr) ds dt - \int_0^S \int_0^T \int_U f^0(s,t,x(s,t), r)\overline{\sigma}(s,t)(dr) ds dt +$$

$$+ \theta_n (\int_0^S \int_0^T \int_U f^0(s,t,x_n(s,t),r)\sigma(s,t))(dr) ds dt - \int_0^S \int_0^T \int_U f^0(s,t,x(s,t),r)\sigma(s,t))(dr) ds dt +$$

$$+ (1-\theta_n) \int_0^S \int_0^T \int_U f^0(s,t,x(s,t), r)\overline{\sigma}(s,t)(dr) ds dt + \theta_n \int_0^S \int_0^T \int_U f^0(s,t,x(s,t),r)\sigma(s,t))(dr) ds dt =$$

$$= (1-\theta) \int_0^S \int_0^T \int_U f^0(s,t,x(s,t), r)\overline{\sigma}(s,t)(dr) ds dt + \theta \int_0^S \int_0^T \int_U f^0(s,t,x(s,t),r)\sigma(s,t))(dr) ds dt =$$

$$= \int_0^S \int_0^T f^0(s,t,x(s,t), \overline{\sigma}(s,t) + \theta(\sigma(s,t) - \overline{\sigma}(s,t))) ds dt,$$

т.е. функционал $(x, \theta) \to g_0(x, \overline{\sigma} + \theta(\sigma - \overline{\sigma})) : Y \times [0,1] \to Y$ непрерывен.

Аналогично имеем, что (см. доказательство теоремы 4.2.9[4]) оператор $(x, \theta) \to F(x, \overline{\sigma} + \theta(\sigma - \overline{\sigma})) : Y \times [0,1] \to Y$ непрерывен.

По условию, легко проверяется, что отображения

$(x, \theta) \to F(x, \overline{\sigma} + \theta(\sigma - \overline{\sigma})) : Y \times [0,1] \to Y$ и $(x, \theta) \to g_0(x, \overline{\sigma} + \theta(\sigma - \overline{\sigma})) : Y \times [0,1] \to Y$ имеют непрерывные производные в точке $(\overline{x}, 0)$ (см. также теоремы 4.2.8[4]).

По условию уравнение $y(s,t) - F_x(\overline{x}, \overline{\sigma})y(s,t) = z(s,t)$, где $z(\cdot) \in Y$, имеет единственное решение и решение непрерывно зависит от правой части, т.е. отображение $I - F_x(\overline{x}, \overline{\sigma})$ является гомеоморфизм из $Y$ в $Y$. Поэтому условие теоремы 5.2.3[4] выполняется.

Если положить $g_1(x, \sigma) \equiv 0$, $g_2(x, \sigma) \equiv 0$, $C_1 = \{0\}$, $C_2 = \{0\}$, то из теоремы 5.2.3[4] имеем следующее следствие 1.

**Следствие 1.** Если $(\overline{x}, \overline{\sigma}) \in Y \times \Im$ является решением задачи (4),(5), то

$$g_{0x}(\overline{x}, \overline{\sigma}) \circ (I - (F_x(\overline{x}, \overline{\sigma}))^{-1} D_2 F(\overline{x}, \overline{\sigma}; \sigma - \overline{\sigma}) + D_2 g_0(\overline{x}, \overline{\sigma}; \sigma - \overline{\sigma}) \geq 0$$

при $\sigma \in \Im$.

Из следствие 1 следует следующее следствие 2.

**Следствие 2.** Если $(\overline{x}, \overline{\sigma}) \in Y \times \Im$ является решением задачи (5),(6), то

$$\int_0^S \int_0^T f_x^0(s,t,\overline{x}(s,t), \overline{\sigma}(s,t))(I - (F_x(\overline{x}, \overline{\sigma}))^{-1} (\int_0^s h^1(\nu, t, \overline{x}(\nu, t), \sigma(\nu, t) - \overline{\sigma}(\nu, t)) d\nu,$$

$$\int_0^t h^2(s, \tau, \overline{x}(s, \tau), \sigma(s, \tau) - \overline{\sigma}(s, \tau)) d\tau) ds dt + \int_0^S \int_0^T f^0(s,t,\overline{x}(s,t), \sigma(s,t) - \overline{\sigma}(s,t)) ds dt \geq 0$$

при $\sigma \in \Im$.

Пусть $\sigma \in \Im$. Положим

$$p(s,t) = (I - (F_x(\overline{x},\overline{\sigma}))^{-1} D_2 F(\overline{x},\overline{\sigma};\sigma-\overline{\sigma})(s,t).$$

Тогда $(I - (F_x(\overline{x},\overline{\sigma}))p(s,t) = D_2 F(\overline{x},\overline{\sigma};\sigma-\overline{\sigma})(s,t)$. Ясно, что

$$p(s,t) - F_x(\overline{x},\overline{\sigma})p(s,t) = p(s,t) - (\int_0^s h_x^1(v,t,\overline{x}(v,t),\overline{\sigma}(v,t))p(v,t)dv,$$

$$\int_0^t h_x^2(s,\tau,\overline{x}(s,\tau),\overline{\sigma}(v,t))\, p(s,\tau)d\tau).$$

Так как $(I - (F_x(\overline{x},\overline{\sigma}))p(s,t) = D_2 F(\overline{x},\overline{\sigma};\sigma-\overline{\sigma})(s,t)$ и

$$D_2 F(\overline{x},\overline{\sigma};\sigma-\overline{\sigma})(s,t) = (\int_0^s h^1(v,t,\overline{x}(v,t),\sigma(v,t)-\overline{\sigma}(v,t))dv, \int_0^t h^2(s,\tau,\overline{x}(s,\tau),\sigma(s,\tau)-\overline{\sigma}(s,\tau))d\tau),$$

то имеем, что

$$p(s,t) - (\int_0^s h_x^1(v,t,\overline{x}(v,t),\overline{\sigma}(v,t))p(v,t)dv, \int_0^t h_x^2(s,\tau,\overline{x}(s,\tau),\overline{\sigma}(v,t))\, p(s,\tau)d\tau) =$$

$$= (\int_0^s h^1(v,t,\overline{x}(v,t),\sigma(v,t)-\overline{\sigma}(v,t))dv, \int_0^t h^2(s,\tau,\overline{x}(s,\tau),\sigma(s,\tau)-\overline{\sigma}(s,\tau))d\tau).$$

Положив $\overline{p} = (p^1,\ldots,p^m)$, $\widetilde{p} = (p^{m+1},\ldots,p^n)$,

$$m_1(s,t) = h^1(s,t,\overline{x}(s,t),\sigma(s,t)-\overline{\sigma}(s,t)), \qquad m_2(s,t) = h^2(s,t,\overline{x}(s,t),\sigma(s,t)-\overline{\sigma}(s,t))$$

имеем, что

$$\begin{aligned}\overline{p}_s(s,t) - h_x^1(s,t,\overline{x}(s,t),\overline{\sigma}(s,t))p(s,t) &= m_1(s,t), \\ \widetilde{p}_t(s,t) - h_x^2(s,t,\overline{x}(s,t),\overline{\sigma}(s,t))p(s,t) &= m_2(s,t),\end{aligned} \qquad (7)$$

где $\overline{p}(0,t) = 0$, $\widetilde{p}(s,0) = 0$.

Из следствия 2 имеем, что

$$\int_0^S \int_0^T (f_x^0(s,t,\overline{x}(s,t),\overline{\sigma}(s,t))p(s,t) + f^0(s,t,\overline{x}(s,t),\sigma(s,t)-\overline{\sigma}(s,t))dsdt \geq 0 \qquad (8)$$

при $\sigma \in \mathfrak{I}$.

Положим $H(s,t,x,u,\psi) = -f^0(s,t,x,u) + \sum_{i=1}^n f^i(s,t,x,u)\psi^i$ и выпишем сопряженную задачу для $\psi = \psi(s,t) = (\psi^1(s,t),\ldots,\psi^n(s,t))$:

$$\begin{aligned}\psi_s^i &= -H_{x^i} = f_{x^i}^0(s,t,\overline{x}(s,t),\overline{\sigma}(s,t)) - \sum_{j=1}^n f_{x^i}^j(s,t,\overline{x}(s,t),\overline{\sigma}(s,t))\psi^j, \quad \psi^i(S,t) = 0, \ i=1,2,\ldots,m \\ \psi_t^i &= -H_{x^i} = f_{x^i}^0(s,t,\overline{x}(s,t),\overline{\sigma}(s,t)) - \sum_{j=1}^n f_{x^i}^j(s,t,\overline{x}(s,t),\overline{\sigma}(s,t))\psi^j, \quad \psi^i(s,T) = 0, \ i=m+1,\ldots,n.\end{aligned} \qquad (9)$$

Обозначим $\overline{\psi} = (\psi^1,\ldots,\psi^m)$, $\widetilde{\psi} = (\psi^{m+1},\ldots,\psi^n)$,

$$A = \begin{pmatrix} f_{x_1}^1(s,t,\overline{x}(s,t),\overline{\sigma}(s,t)) & \ldots & f_{x_n}^1(s,t,\overline{x}(s,t),\overline{\sigma}(s,t)) \\ \ldots & \ldots & \ldots \\ f_{x_1}^n(s,t,\overline{x}(s,t),\overline{\sigma}(s,t)) & \ldots & f_{x_n}^n(s,t,\overline{x}(s,t),\overline{\sigma}(s,t)) \end{pmatrix}.$$

Из (7) и (9) имеем

$$\begin{pmatrix} \overline{p}_s(s,t) \\ \widetilde{p}_t(s,t) \end{pmatrix} = Ap(s,t) + m(s,t), \tag{10}$$

$$\begin{pmatrix} \overline{\psi}_s(s,t) \\ \widetilde{\psi}_t(s,t) \end{pmatrix} = -A^*\psi(s,t) + f_x^0(s,t,\overline{x}(s,t),\overline{\sigma}(s,t)), \tag{11}$$

где $m(s,t) = (m_1(s,t), m_2(s,t))$. Из (11) получим

$$f_x^0(s,t,\overline{x}(s,t),\overline{\sigma}(s,t))p(s,t) = A^*\psi(s,t)p(s,t) + \begin{pmatrix} \overline{\psi}_s(s,t) \\ \widetilde{\psi}_t(s,t) \end{pmatrix} p(s,t).$$

Тогда из (8) следует, что

$$\int_0^S \int_0^T (A^*\psi(s,t)p(s,t) + \sum_{i=1}^m \psi_s^i(s,t)p^i(s,t) + \sum_{i=m+1}^n \psi_t^i(s,t)p^i(s,t) +$$
$$+ f^0(s,t,\overline{x}(s,t), \sigma(s,t) - \overline{\sigma}(s,t)) ds dt \geq 0$$

при $\sigma \in \Im$. Отсюда следует

$$\int_0^S \int_0^T (A^*\psi(s,t)p(s,t) - \sum_{i=1}^m \psi^i(s,t)p_s^i(s,t) - \sum_{i=m+1}^n \psi^i(s,t)p_t^i(s,t) +$$
$$+ f^0(s,t,\overline{x}(s,t), \sigma(s,t) - \overline{\sigma}(s,t))) ds dt \geq 0$$

при $\sigma \in \Im$. Поэтому

$$\int_0^S \int_0^T (-\psi(s,t)m(s,t) + f^0(s,t,\overline{x}(s,t),\sigma(s,t) - \overline{\sigma}(s,t))) ds dt \geq 0$$

при $\sigma \in \Im$. Тогда

$$\int_0^S \int_0^T (-\sum_{i=1}^n f^i(s,t,\overline{x}(s,t), \sigma(s,t) - \overline{\sigma}(s,t))\psi^i(s,t) + f^0(s,t,\overline{x}(s,t),\sigma(s,t) - \overline{\sigma}(s,t))) ds dt \geq 0$$

при $\sigma \in \Im$. Отсюда имеем

$$\int_0^S \int_0^T H(s,t,\overline{x}(s,t),\sigma(s,t) - \overline{\sigma}(s,t),\psi(s,t)) ds dt \leq 0$$

при $\sigma \in \Im$.

Пусть $\{\rho_1, \rho_2,...\}$ такое семейство измеримых функций из D в $R^r$, что $\{\rho_1(s,t), \rho_2(s,t),...\}$ плотно в U п.в. $(s,t) \in D$. Если положить $\sigma(s,t) = \delta_{\rho_j(s,t)}$ при $(s,t) \in E$, где $\delta_{\rho_j(s,t)}$ мера Дирака, $\sigma(s,t) = \overline{\sigma}(s,t)$ при $(s,t) \notin E$ для каждого $j \in N$ и каждого измеримого подмножества $E \subset D$, то по условию имеем, что

$$\iint_E H(s,t,\overline{x}(s,t),\delta_{\rho_j(s,t)} - \overline{\sigma}(s,t),\psi(s,t)) ds dt \leq 0$$

Отсюда по условию 1.4.37(7)[4] следует, что

$$H(s,t,\overline{x}(s,t),\delta_{\rho_j(s,t)} - \overline{\sigma}(s,t),\psi(s,t)) \leq 0$$

п.в. $(s,t) \in D$. Поэтому

$$H(s,t,\overline{x}(s,t),\overline{\sigma}(s,t),\psi(s,t)) \geq \max_{r \in U} H(s,t,\overline{x}(s,t),r,\psi(s,t)).$$

С другой стороны, так как мера $\overline{\sigma}(s,t)$ сосредоточена на $U$ почти всюду в $D$, то

$$H(s,t,\overline{x}(s,t),\overline{\sigma}(s,t),\psi(s,t)) = \int_U H(s,t,\overline{x}(s,t),r,\psi(s,t))\overline{\sigma}(s,t)(dr) \leq \max_{r \in U} H(s,t,\overline{x}(s,t),r,\psi(s,t))$$

п.в. $(s,t) \in D$. Поэтому

$$H(s,t,\overline{x}(s,t),\overline{\sigma}(s,t),\psi(s,t)) = \max_{r \in U} H(s,t,\overline{x}(s,t),r,\psi(s,t)).$$

Таким образом доказана следующая теорема.

**Теорема 1.** Если $(\overline{x},\overline{\sigma}) \in Y \times \Im$ является решением задачи (5),(6), то

$$H(s,t,\overline{x}(s,t),\overline{\sigma}(s,t),\psi(s,t)) = \max_{r \in U} H(s,t,\overline{x}(s,t),r,\psi(s,t)).$$

Отметим, что ряд утверждений о непрерывности функционалов с обобщенным управлением имеются в [17] (см. леммы 4.1.1-4.1.5) и [4] (см. теоремы 4.2.8).

Если решение из класса $x(s,t) = (x^1(s,t),\ldots,x^n(s,t)) \in L_p^n(D)$, где $(x_s^1(s,t),\ldots,x_s^m(s,t),\, x_t^{m+1}(s,t),\ldots,x_t^n(s,t)) \in L_p^n(D)$, $1 \leq p \leq +\infty$, то для задачи (5),(6) аналогично можно получить необходимое условие оптимальности.

## 2. Оптимальное управление процессами, описываемыми системой первого порядка с частными производными

Пусть $U \subset R^r$ $U_1 \subset R^{r_1}$, $U_2 \subset R^{r_2}$ непустые компактные множества, $f:[0,T] \times [0,S] \times R^{2k} \times R^{2n} \to R$, $f_1:[0,T] \times [0,S] \times R^k \times R^n \times U \to R^k$, $f_2:[0,T] \times [0,S] \times R^k \times R^n \times U \to R^n$, $\varphi_1^0:[0,S] \times R^k \times R^k \to R$, $\varphi_2^0:[0,T] \times R^n \times R^n \to R$, $\varphi_1:[0,S] \times U_1 \to R^k$, $\varphi_2:[0,T] \times U_2 \to R^n$.

Положим $D = [0,T] \times [0,S]$, $V_t^p = \{u \in L_p^k(D): u_t \in L_p^k(D)\}$, $V_s^p = \{\upsilon \in L_p^n(D): \upsilon_s \in L_p^n(D)\}$.

Отметим, что $V^p = V_t^p \times V_s^p$ является нормированное пространство относительно нормы:

$$\|(u,v)\|_{V^2} = \left\{\int_0^T\int_0^S |u(t,s)|^p dtds\right\}^{\frac{1}{p}} + \left\{\int_0^T\int_0^S |u_t(t,s)|^p dtds\right\}^{\frac{1}{p}} +$$

$$+ \left\{\int_0^T\int_0^S |\upsilon(t,s)|^p dtds\right\}^{\frac{1}{p}} + \left\{\int_0^T\int_0^S |\upsilon_s(t,s)|^p dtds\right\}^{\frac{1}{p}}$$

Пусть $h_0(\cdot):[0,T] \times [0,S] \to U$, $h_1(\cdot):[0,S] \to U_1$, $h_2(\cdot):[0,T] \to U_2$ измеримые функции. Далее считаем, что включение и равенства удовлетворяются почти всюду.

Функция $(u,\upsilon) \in V^p$ удовлетворяющая системы

$$\begin{aligned} u_t(t,s) &= f_1(t,s,\, u(t,s),\, \upsilon(t,s),\, h_0(t,s)), \\ \upsilon_s(t,s) &= f_2(t,s,\, u(t,s),\, \upsilon(t,s),\, h_0(t,s)) \end{aligned} \quad (1)$$

и условии

$$u(0,s) = \varphi_1(s,h_1(s)), \qquad \upsilon(t,0) = \varphi_2(t,h_2(t)), \qquad (2)$$

называется решением задачи (1), (2).

Рассмотрим задачу минимизацию функционала

$$J(u(t,s), \upsilon(t,s)) = \int_0^T\int_0^S f_0(t,s, u(t,s), u_t(t,s), \upsilon(t,s), \upsilon_s(t,s))dtds +$$
$$+ \int_0^S \varphi_1^0(s, u(0,s), u(T,s))ds + \int_0^T \varphi_2^0(t, \upsilon(t,0), \upsilon(t,S))dt \qquad (3)$$

среди всех решений задачи (1),(2).

Обозначим $a_1(t,s,u,\upsilon) = f_1(t,s, u, \upsilon, U)$, $a_2(t,s,u,\upsilon) = f_2(t,s, u, \upsilon, U)$, $M_1(s) = \varphi_1(s, U_1)$, $M_2(t) = \varphi_2(t, U_2)$.

Рассмотрим минимизации функционала (3) среди всех решений задачи
$$u_t(t,s) \in a_1(t,s, u(t,s), \upsilon(t,s)),$$
$$\upsilon_s(t,s) \in a_2(t,s, u(t,s), \upsilon(t,s)) \qquad (4)$$

при условии
$$u(0,s) \in M_1(s), \quad \upsilon(t,0) \in M_2(t). \qquad (5)$$

Пусть выполняются следующие условия:

1) Функции $f_0(t,s,u,u_1,\upsilon,\upsilon_1)$ и $f_i(t,s,u,\upsilon,h_0)$, $i=1,2$, удовлетворяют условию Каратеодори и удовлетворяют условия Липшица относительно переменных $(u,u_1,\upsilon,\upsilon_1)$ и $(u,\upsilon)$ соответственно с постоянной $L$.

2) Функции $\varphi_1^0(s,z)$, $\varphi_2^0(t,y)$, $\varphi_1(s,h_1)$ и $\varphi_2(t,h_2)$ удовлетворяют условию Каратеодори,

3) функции $\varphi_1^0(s,z)$, $\varphi_2^0(t,y)$ удовлетворяют условия Липшица относительно переменных $z$ и $y$ соответственно с постоянной $L$, где $z \in R^{2k}, y \in R^{2n}$.

Отметим, что аналогично [4], стр. 403, можно показать, что задачи (1)-(3) и (3)-(5) эквивалентны.

Положим $D = [0,S]\times[0,T]$, $V = \{u : u(s,t) \in U, u(s,t)$ измеримы $\}$,

$\Im = \{\sigma(s,t) \in rpm(U) : \sigma(s,t)(U) = 1$ п.в. $(s,t) \in D, \sigma(s,t)$ измеримы $\}$,

$\Im_1 = \{\sigma_1(s) \in rpm(U_1) : \sigma_1(s)(U_1) = 1$ п.в. $s \in [0,S], \sigma_1(s)$ измеримы $\}$,

$\Im_2 = \{\sigma_2(t) \in rpm(U_2) : \sigma_2(t)(U_2) = 1$ п.в. $t \in [0,T], \sigma_2(t)$ измеримы $\}$,

где $rpm(U)$, $rpm(U_1)$ и $rpm(U_2)$ множество всех вероятностных мер Радона в $U$, $U_1$ и $U_2$ соответственно.

Рассмотрим обобщенную задачу для задачи (1)-(3), т.е. рассмотрим задачу минимизации функционала

$$J(u(t,s), \upsilon(t,s)) = \int_0^T\int_0^S f_0(t,s, u(t,s), u_t(t,s), \upsilon(t,s), \upsilon_s(t,s))dtds +$$
$$+ \int_0^S \varphi_1^0(s, u(0,s), u(T,s))ds + \int_0^T \varphi_2^0(t, \upsilon(t,0), \upsilon(t,S))dt$$

среди функций, удовлетворяющих системе
$$u_t(t,s) = f_1(t,s, u(t,s), \upsilon(t,s), \sigma(t,s)),$$
$$\upsilon_s(t,s) = f_2(t,s, u(t,s), \upsilon(t,s), \sigma(t,s)) \qquad (6)$$

и условию
$$u(0,s) = \varphi_1(s, \sigma_1(s)), \qquad \upsilon(t,0) = \varphi_2(t, \sigma_2(t)), \qquad (7)$$

при $\sigma \in \Im$, $\sigma_1 \in \Im_1$, $\sigma_2 \in \Im_2$.

Положим

$$P_1(t,s,u,\upsilon) = \{f_1(t,s, u, \upsilon, \mu): \mu \in prm(U)\}, \qquad P_3(s) = \{\varphi_1(s,\mu_1): \mu_1 \in prm(U_1)\},$$
$$P_2(t,s,u,\upsilon) = f_2(t,s, u, \upsilon, \mu): \mu \in prm(U)\}, \qquad P_4(t) = \{\varphi_2(t,\mu_2): \mu_2 \in prm(U_2)\}.$$

Из теоремы 1.6.13 и 1.6.14[4] следует, что

$$P_1(t,s,u,\upsilon) = co\, f_1(t,s, u, \upsilon, U), \qquad P_3(s) = co\, \varphi_1(s,U_1),$$
$$P_2(t,s,u,\upsilon) = co\, f_2(t,s, u, \upsilon, U), \qquad P_4(t) = co\, \varphi_2(t,U_2).$$

Рассмотрим минимизацию функционала (3) среди всех решений задачи
$$u_t(t,s) \in P_1(t,s, u(t,s), \upsilon(t,s)),$$
$$\upsilon_s(t,s) \in P_2(t,s, u(t,s), \upsilon(t,s)) \tag{8}$$
при условии
$$u(0,s) \in P_1(s), \quad \upsilon(t,0) \in P_2(t). \tag{9}$$

Отметим, что задачи (3),(6),(7) и (3),(8),(9) эквивалентны.

Из леммы 5.7.1 [23] следует следующая лемма 1.

**Лемма 1.** Если $U \subset R^r$ компактное множество, $f_1:[0,T]\times[0,S]\times R^k \times R^n \times U \to R^k$ удовлетворяет условию Каратеодори и условию Липшица относительно переменных $(u,\upsilon)$ с постоянной $L$, то
$$\rho_X(a_1(t,s,u_1,\upsilon_1), a_1(t,s,u,\upsilon)) \le L|(u_1,\upsilon_1) - (u,\upsilon)|,$$
$$\rho_X(P_1(t,s,u_1,\upsilon_1), P_1(t,s,u,\upsilon)) \le L|(u_1,\upsilon_1) - (u,\upsilon)|$$
при $(u_1,\upsilon_1), (u,\upsilon) \in R^k \times R^n$.

Положим $\psi_1(t,s,u,\upsilon,z_1) = \inf\{|z_1 - y_1|: y_1 \in a_1(t,s,u,\upsilon)\},$
$$\psi_2(t,s,u,\upsilon,z_2) = \inf\{|z_2 - y_2|: y_2 \in a_2(t,s,u,\upsilon)\},$$
$q_1(s,z_1) = \inf\{|z_1 - y_1|: y_1 \in M_1(s)\}, \quad q_2(t,z_2) = \inf\{|z_2 - y_2|: y_2 \in M_2(t)\}.$

(или $\psi_1(t,s,u,\upsilon,z_1) = \inf\{|z_1 - y_1|: y_1 \in P_1(t,s,u,\upsilon)\},$
$$\psi_2(t,s,u,\upsilon,z_2) = \inf\{|z_2 - y_2|: y_2 \in P_2(t,s,u,\upsilon)\},$$
$q_1(s,z_1) = \inf\{|z_1 - y_1|: y_1 \in P_3(s)\}, \quad q_2(t,z_2) = \inf\{|z_2 - y_2|: y_2 \in P_4(t)\}$).

Из теоремы 3.2.1[17] следует, что верна следующая теорема 1.

**Теорема 1.** Если $U \subset R^r$, $U_1 \subset R^{r_1}$ и $U_2 \subset R^{r_2}$ компактные множества, функции $f_1:[0,T]\times[0,S]\times R^k \times R^n \times U \to R^k$, $f_2:[0,T]\times[0,S]\times R^k \times R^n \times U \to R^n$ удовлетворяют условию Каратеодори и условию Липшица относительно переменных $(u,\upsilon)$ с постоянной $L$, функция $f:[0,T]\times[0,S]\times R^{2k} \times R^{2n} \to R$ измерима по $(t,s)$ и удовлетворяет условию Липшица относительно переменных $(u,u_1,\upsilon,\upsilon_1)$ с постоянной $L$, функции $\varphi_1^0:[0,S]\times R^k \times R^k \to R$ и $\varphi_2^0:[0,T]\times R^n \times R^n \to R$ измеримы по $s$ и $t$ соответственно и удовлетворяют условию Липшица относительно переменных $(u,u_1)$ и $(\upsilon,\upsilon_1)$ соответственно с постоянной $L$, функции $\varphi_1:[0,S]\times R^{r_1} \to R^k$, $\varphi_2:[0,T]\times R^{r_2} \to R^n$ удовлетворяют условию Каратеодори и $(\bar{u}, \bar{\upsilon}) \in V_t^p \times V_s^p, (1 \le p < \infty)$ являются решением задачи (1),(2), то существуют число $\lambda > 0$ и функции $\bar{v} \in L_q^k(D)$, $\bar{\omega} \in L_q^n(D)$, где $\bar{v}_t \in L_q^k(D)$, $\bar{\omega}_s \in L_q^n(D)$, $p+q = pq$, такие, что

1) $(\bar{v}_t(t,s), \bar{v}(t,s), \bar{\omega}_s(t,s), \bar{\omega}(t,s)) \in \partial(f(t,s,\cdot) + \lambda(\psi_1(t,s,\cdot) +$
$+ \psi_2(t,s,\cdot))(\bar{u}(t,s), \bar{u}_t(t,s), \bar{\upsilon}(t,s), \bar{\upsilon}_s(t,s)),$

2) $(\bar{v}(0,s), -\bar{v}(T,s)) \in \partial(\varphi_1^0(s,\cdot) + \lambda q_1(s,\cdot))(\bar{u}(0,s), \bar{u}(T,s)),$

3) $(\overline{\omega}(t,0), -\overline{\omega}(t,S)) \in \partial(\varphi_2^0(t,\cdot) + \lambda q_2(t,\cdot))(\overline{\upsilon}(t,0), \overline{\upsilon}(t,S))$.

Отметим, что аналогично теореме 5.7.2 [23] из теоремы 1 можно получить необходимое условие оптимальности для задачи (1)-(3) или (3),(6),(7) в виде принципа максимума.

## 3. Градиент в одной дискретной задаче оптимального управления в бесконечном интервале

Рассмотрим следующую задачу оптимального управления с дискретном временем: минимизировать функцию

$$J(u) = \sum_{i=0}^{\infty} \Phi_i(x_i, u_i) \qquad (1)$$

при условиях

$$x_{i+1} = f_i(x_i, u_i), \quad x_0 = x^0, \quad i = 0,1,2,\ldots, \qquad (2)$$

$$u = (u_0, u_1, u_2, \ldots) \in \ell_2^r, \quad x = (x_0, x_1, x_2, \ldots) \in \ell_2^n, \qquad (3)$$

где $u_i = (u_i^1, u_i^2, \ldots, u_i^r) \in R^r$, $x_i = (x_i^1, x_i^2, \ldots, x_i^n) \in R^n$, функции $f_i = (f_i^1, \ldots, f_i^n)$ и $\Phi_i$ предполагается известными, $f_i : R^n \times R^r \to R^n$, $\Phi_i : R^n \times R^r \to R$.

Будем пользоваться следующими обозначениями:

$$f_{ix_i} = \begin{pmatrix} f_{ix_i^1}^1, \ldots, f_{ix_i^n}^1 \\ \cdots \\ f_{ix_i^1}^n, \ldots, f_{ix_i^n}^n \end{pmatrix}, \qquad f_{iu_i} = \begin{pmatrix} f_{iu_i^1}^1, \ldots, f_{iu_i^r}^1 \\ \cdots \\ f_{iu_i^1}^n, \ldots, f_{iu_i^r}^n \end{pmatrix}.$$

Через $\ell_2^r$ обозначено гильбертово пространство, состоящее из всевозможных последовательностей $u = (u_0, u_1, u_2, \ldots) \in \ell_2^r \subset (R^r)^\infty$ с нормой

$$\|u\| = \left(\sum_{i=0}^{\infty} \sum_{j=1}^{r} |u_i^j|^2\right)^{\frac{1}{2}} < +\infty.$$

Считаем, что существуют числа $c > 0$ и $a = (a_0, a_1, a_2, \ldots) \in \ell_1$ такие, что

$$|\Phi_i(x_i, u_i)| \le c(|x_i|_{R^n}^2 + |u_i|_{R^r}^2) + a_i.$$

Так как $\ell_2^r$ гильбертово пространство, то $(\ell_2^r)^* \in \ell_2^r$, т.е. каждый линейный непрерывный функционал $u^* \in (\ell_2^r)^*$ имеет вид

$$\langle u^*, u \rangle = \sum_{i=0}^{\infty} \langle u_i^*, u_i \rangle,$$

где $u = (u_0^*, u_1^*, u_2^*, \ldots) \in \ell_2^r$. Положим

$$H_i(x_i, u_i, \psi_i) = \Phi_i(x_i, u_i) + \langle \psi_i, f_i(x_i, u_i) \rangle,$$

где $\psi_{i-1} = H_{ix_i}(x_i, u_i, \psi_i)$, $\psi_i = (\psi_i^1, \ldots, \psi_i^n)$.

**Теорема 1**. Пусть функции $\Phi_i$ и $F_i$ имеют частные производные по переменным x и u при $x \in R^n$ и $u \in R^r$ и существуют числа $L > 0$, $L_1 > 0$, где $0 < L_1 < \dfrac{\sqrt{2}}{2}$, и $L_2 > 0$ такие, что

1) $|f_i(x+z, u+h) - f_i(x,u)| \leq L_1|z| + L_2|h|$,

2) $|f_{i_x}(x+z, u+h) - f_{i_x}(x,u)| \leq L(|z| + |h|)$,

3) $|f_{i_u}(x+z, u+h) - f_{i_u}(x,u)| \leq L(|z| + |h|)$,

4) $|\Phi_{i_x}(x+z) - \Phi_{i_x}(x)| \leq L|z|$,  5) $|\Phi_{i_u}(x+z) - \Phi_{i_u}(x)| \leq L|z|$

при $x, z \in R^n$, $u, h \in R^r$, кроме того $(H_{iu_i}(x_i, u_i, \psi_i))_{i=0}^\infty \in \ell_2^r$, $(H_{ix_i}(x_i, u_i, \psi_i))_{i=0}^\infty \in \ell_2^n$ и $(f_i(x_i, u_i))_{i=0}^\infty \in \ell_2^n$ при $u = (u_0, u_1, u_2, \ldots) \in \ell_2^r$, $x = (x_0, x_1, x_2, \ldots) \in \ell_2^n$ и $\psi = (\psi_0, \psi_1, \psi_2, \ldots) \in \ell_2^n$. Тогда функция (1) при условиях (2) непрерывно дифференцируема в норме $\ell_2^r$, причем ее градиент

$$J'(u) = (H_{iu_i}(x_i, u_i, \psi_i))_{i=0}^\infty \in \ell_2^r,$$

где $\psi_{i-1} = H_{ix_i}(x_i, u_i, \psi_i)$, $H_{iu_i}(x_i, u_i, \psi_i) = (H_{iu_i^1}(x_i, u_i, \psi_i), \ldots, H_{iu_i^r}(x_i, u_i, \psi_i))$.

**Доказательство.** Пусть $u$ и $u + h \in \ell_2^r$ и пусть $x$ и $x + \Delta x \in \ell_2^n$ соответствующие этим управлениям дискретные траектории задачи (2), а $J(u)$ и $J(u+h) = J(u) + \Delta J$ соответствующие значения функции (1). Из (2) следует, что приращение $\Delta x$ удовлетворяет условиям

$$\Delta x_{i+1} = f_i(x_i + \Delta x_i, u_i + h_i) - f_i(x_i, u_i), \qquad \Delta x_0 = 0, \quad i = 0, 1, 2, \ldots,$$

Так как

$$\Delta \Phi_i = \Phi_i(x_i + \Delta x_i, u_i + h_i) - \Phi_i(x_i, u_i),$$

то из (1) получим, что

$$\Delta J = \sum_{i=0}^\infty (\Phi_i(x_i + \Delta x_i, u_i + h_i) - \Phi_i(x_i, u_i)).$$

Положим

$$H_i(x_i, u_i, \psi_i) = \Phi_i(x_i, u_i) + \langle \psi_i, f_i(x_i, u_i) \rangle,$$

$$H_{iu_i}(x_i, u_i, \psi_i) = (H_{iu_i^1}(x_i, u_i, \psi_i), \ldots, H_{iu_i^r}(x_i, u_i, \psi_i)).$$

Ясно, что

$$\Delta J = \sum_{i=0}^{\infty}(\Phi_i(x_i+\Delta x_i, u_i+h_i)-\Phi_i(x_i,u_i)+\langle \psi_i, f_i(x_i+\Delta x_i, u_i+h_i)-f_i(x_i,u_i)\rangle -$$
$$-\langle \psi_i, \Delta x_{i+1}\rangle = \sum_{i=0}^{\infty}(H_i(x_i+\Delta x_i, u_i+h_i, \psi_i)-H_i(x_i,u_i,\psi_i))-\langle \psi_i, \Delta x_{i+1}\rangle.$$

Из формулы конечных приращений следует, что

$$H_i(x_i+\Delta x_i, u_i+h_i, \psi_i) = H_i(x_i,u_i,\psi_i)+\langle H_{ix_i}(x_i+\theta_i \Delta x_i, u_i+\theta_i h_i, \psi_i), \Delta x_i\rangle +$$
$$+\langle H_{iu_i}(x_i+\theta_i \Delta x_i, u_i+\theta_i h_i, \psi_i), h_i\rangle, \qquad 0\leq \theta_i \leq 1, \ i=0,1,2,\ldots.$$

Тогда имеем

$$\Delta J = \sum_{i=0}^{\infty}(\langle H_{ix_i}(x_i+\theta_i \Delta x_i, u_i+\theta_i h_i, \psi_i), \Delta x_i\rangle + \langle H_{iu_i}(x_i+\theta_i \Delta x_i, u_i+\theta_i h_i, \psi_i), h_i\rangle)-$$
$$-\langle \psi_i, \Delta x_{i+1}\rangle = \sum_{i=0}^{\infty}(\langle H_{ix_i}(x_i+\theta_i \Delta x_i, u_i+\theta_i h_i, \psi_i)-H_{ix_i}(x_i,u_i,\psi_i), \Delta x_i\rangle +$$
$$+\langle H_{iu_i}(x_i+\theta_i \Delta x_i, u_i+\theta_i h_i, \psi_i)-H_{iu_i}(x_i,u_i,\psi_i), h_i\rangle + \langle H_{ix_i}(x_i,u_i,\psi_i)-\psi_{i-1}, \Delta x_i\rangle +$$
$$+\langle H_{iu_i}(x_i,u_i,\psi_i), h_i\rangle + \langle \psi_{i-1}, \Delta x_i\rangle - \langle \psi_i, \Delta x_{i+1}\rangle).$$

Обозначим

$$R_1 = \sum_{i=0}^{\infty}\langle H_{ix_i}(x_i+\theta_i \Delta x_i, u_i+\theta_i h_i, \psi_i)-H_{ix_i}(x_i,u_i,\psi_i), \Delta x_i\rangle,$$
$$R_2 = \sum_{i=0}^{\infty}\langle H_{iu_i}(x_i+\theta_i \Delta x_i, u_i+\theta_i h_i, \psi_i)-H_{iu_i}(x_i,u_i,\psi_i), h_i\rangle,$$
$$R_3 = \sum_{i=0}^{\infty}(\langle \psi_{i-1}, \Delta x_i\rangle - \langle \psi_i, \Delta x_{i+1}\rangle) = -\lim_{i\to\infty}\langle \psi_i, \Delta x_{i+1}\rangle.$$

Тогда имеем, что

$$\Delta J = \sum_{i=0}^{\infty}(\langle H_{iu_i}(x_i,u_i,\psi_i), h_i\rangle)+R_1+R_2+R_3,$$

где $\psi_{i-1} = H_{ix_i}(x_i,u_i,\psi_i)$, $\lim_{i\to\infty}\psi_i = 0$.

Так как $\Delta x_{i+1} = f_i(x_i+\Delta x_i, u_i+h_i)-f_i(x_i,u_i)$, $\Delta x_0 = 0$, $i=0,1,2,\ldots$, то по условию $|\Delta x_{i+1}| \leq L_1|\Delta x_i|+L_2|h_i|$. Поэтому

$$|\Delta x_{i+1}|^2 \leq 2L_1^2|\Delta x_i|^2 + 2L_2^2|h_i|^2$$

при $\Delta x \in \ell_2^n$, $h \in \ell_2^r$. Отсюда следует, что

$$\sum_{i=0}^{\infty}|\Delta x_{i+1}|^2 \leq 2L_1^2\sum_{i=0}^{\infty}|\Delta x_i|^2 + 2L_2^2\sum_{i=0}^{\infty}|h_i|^2.$$

Тогда получим

$$\sum_{i=0}^{\infty}|\Delta x_{i+1}|^2 - 2L_1^2\sum_{i=0}^{\infty}|\Delta x_i|^2 \leq 2L_2^2\sum_{i=0}^{\infty}|h_i|^2.$$

Если $1 - 2L_1^2 > 0$, то имеем, что $0 < L_1 < \dfrac{\sqrt{2}}{2}$. Поэтому

$$\sum_{i=0}^{\infty} |\Delta x_{i+1}|^2 = \sum_{i=0}^{\infty} |\Delta x_i|^2 \leq \dfrac{2L_2^2}{1 - 2L_1^2} \sum_{i=0}^{\infty} |h_i|^2.$$

По условию существует число $c > 0$ такое, что $|\psi_i| \leq c$. Поэтому

$$|R_1| \leq \sum_{i=0}^{\infty} \left| \left\langle H_{i x_i}(x_i + \theta_i \Delta x_i, u_i + \theta_i h_i, \psi_i) - H_{i x_i}(x_i, u_i, \psi_i), \Delta x_i \right\rangle \right| \leq L \sum_{i=0}^{\infty} (|\Delta x_i| + |h_i|)|\Delta x_i| +$$

$$+ \sum_{i=0}^{\infty} |\psi_i| L(|\Delta x_i| + |h_i|)|\Delta x_i| \leq (L + cL) \sum_{i=0}^{\infty} 1{,}5 |\Delta x_i|^2 + 0{,}5(L + cL) \sum_{i=0}^{\infty} |h_i|^2) \leq$$

$$\leq 1{,}5(L + cL) \dfrac{2 L_2^2}{1 - 2 L_1^2} \sum_{i=0}^{\infty} |h_i|^2 + 0{,}5(L + cL) \sum_{i=0}^{\infty} |h_i|^2 = \left((L + cL)\left(\dfrac{3 L_2^2}{1 - 2 L_1^2} + 0{,}5\right)\right) \sum_{i=0}^{\infty} |h_i|^2.$$

Аналогично существует число $K > 0$ такое, что

$$|R_2| = \left| \sum_{i=0}^{\infty} \left\langle H_{i u_i}(x_i + \theta_i \Delta x_i, u_i + \theta_i h_i, \psi_i) - H_{i u_i}(x_i, u_i, \psi_i), h_i \right\rangle \right| \leq K \sum_{i=0}^{\infty} |h_i|^2.$$

Ясно, что

$$|R_3| = \left| \sum_{i=0}^{\infty} (\langle \psi_{i-1}, \Delta x_i \rangle - \langle \psi_i, \Delta x_{i+1} \rangle) \right| = \left| \lim_{i \to \infty} \langle \psi_i, \Delta x_{i+1} \rangle \right| = \lim_{i \to \infty} c |\Delta x_{i+1}| = 0.$$

Если $(H_{i u_i}(x_i, u_i, \psi_i))_{i=0}^{\infty} \in \ell_2^r$, то имеем, что $J'(u) = (H_{i u_i}(x_i, u_i, \psi_i))_{i=0}^{\infty} \in \ell_2^r$. Теорема доказана.

Отметим, что если существуют числа $c > 0$ и $a = (a_0, a_1, a_2, \ldots) \in \ell_2$ такие, что

$$\left| \Phi_{i_{x_i}}(x_i, u_i) \right| \leq c(|x_i|_{R^n} + |u_i|_{R^r}) + a_i, \quad \left| \Phi_{i_{u_i}}(x_i, u_i) \right| \leq c(|x_i|_{R^n} + |u_i|_{R^r}) + a_i,$$

$$\left| f_{i_{x_i}}(x_i, u_i) \right| \leq c(|x_i|_{R^n} + |u_i|_{R^r}) + a_i, \quad \left| f_{i_{x_i}}(x_i, u_i) \right| \leq c(|x_i|_{R^n} + |u_i|_{R^r}) + a_i$$

при $x_i \in R^n$ и $u_i \in R^r$, то $(\Phi_{0_{x_0}}(x_0, u_0), \Phi_{1_{x_1}}(x_1, u_1), \Phi_{2_{x_2}}(x_2, u_2), \ldots) \in \ell_2^n$,

$(\Phi_{0_{u_0}}(x_0, u_0), \Phi_{1_{u_1}}(x_1, u_1), \Phi_{2_{u_2}}(x_2, u_2), \ldots) \in \ell_2^r$, $(f_{0_{x_0}}(x_0, u_0), f_{1_{x_1}}(x_1, u_1), f_{2_{x_2}}(x_2, u_2), \ldots) \in \ell_2^n$

$(f_{0_{u_0}}(x_0, u_0), f_{1_{u_1}}(x_1, u_1), f_{2_{u_2}}(x_2, u_2), \ldots) \in \ell_2^r$ при $(x_0, x_1, x_2, \ldots) \in \ell_2^n$ и $(u_0, u_1, u_2, \ldots) \in \ell_2^r$.

Поэтому $(H_{i u_i}(x_i, u_i, \psi_i))_{i=0}^{\infty} \in \ell_2^r$ и $(H_{i x_i}(x_i, u_i, \psi_i))_{i=0}^{\infty} \in \ell_2^r$ при $u = (u_0, u_1, u_2, \ldots) \in \ell_2^r$, $x = (x_0, x_1, x_2, \ldots) \in \ell_2^n$ и $\psi = (\psi_0, \psi_1, \psi_2, \ldots) \in \ell_2^n$.

Если выполняется условие теоремы 1 и $(\Phi_{i_{x_i}}(\bar{x}_i, \bar{u}_i))_{i=0}^{\infty} \in \ell_2^n$, где $(\bar{x}_0, \bar{x}_1, \bar{x}_2, \ldots) \in \ell_2^n$ и $(\bar{u}_0, \bar{u}_1, \ldots) \in \ell_2^r$, то $(\Phi_{i_{x_i}}(x_i, u_i))_{i=0}^{\infty} \in \ell_2^n$ при $(x_0, x_1, x_2, \ldots) \in \ell_2^n$ и $(u_0, u_1, u_2, \ldots) \in \ell_2^r$.

ЛИТЕРАТУРА

1. Алексеев В.М., Тихомиров В.М., Фомин С.В. Оптимальное управление, М.: Наука, 1979, 429 с.
2. Бедельбаев А.А. О субдифференциалах второго порядка и их приложения в вариационном исчислении.-Автореферат канд. дис. Алма-Ата, 1986.
3. Борисович Ю.Г., Гельман Б.Д. и др. Введение в терию многозначных отображений. Воронеж, 1986, 103 с.
4. Варга Дж. Оптимальное управление дифренциальными и функциональными уравнениями. М. : Наука, 1977, 623 с.
5. Данфорд Н., Шварц Дж.Т. Линейные операторы. Общая теория. М.: Наука, 1962, 895р.
6. Иоффе А.Д., Тихомиров В.М. Теория экстремальных задач. М.: Наука,1974, 479с.
7. Кларк Ф. Оптимизация и негладкий анализ. М.: Наука, 1988, 280 с.
8. Лоран П.-Ж. Аппроксимация и оптимизация. М.: Мир, 1975, 496 с.
9. Михайлов В.П. Дифференциальные уравнения в частных производных. М.:Наука, 1983, 424 с.
10. Мордухович Б.Ш. Методы аппроксимаций в задачах оптимизации и управления. М.: Наука, 1988, 350 с.
11. Обен Ж.П. Нелинейный анализ и его экономические приложения. М.: Мир, 1988, 264с.
12. Обен Ж.П., Экланд И. Прикладной нелинейный анализ. М.: Мир, 1988, 510 с.
13. Пшеничный Б.Н. Выпуклый анализ и экстремальные эадачи. М.:Наука,1980, 319с.
14. Рокафеллар Р. Интегралы, являющиеся выпуклыми функционалами, II. В кн. Математическая экономика. –М.: Мир, 1974, с. 170-204.
15. Рудин У. Функциональный анализ. М.: Мир, 1975, 443 с.
16. Садыгов М.А. Необходимое условие экстремума высших порядков для негладких функций. Изв. АН Азерб. ССР, сер.физ.-техн. и матем. наук, 1989, №6, с.33-47.
17. Садыгов М.А. Экстремальные задачи для негладких систем. Баку, 1996, Элм, 148 с.

18. Садыгов М.А. Негладкий анализ и его приложения к экстремальной задаче для включения типа Гурса-Дарбу. Баку, Элм-1999, 135 с.
19. Садыгов М.А. Исследование негладких оптимизационных задач. Баку, Элм, 2002, 125с.
20. Садыгов М.А. Исследование субдифференциала первого и второго порядков. Баку, Наука, 2007, 224с.
21. Садыгов М.А. Экстремальные задачи для трехмерных дифференциальных включений. Препринт №1, Баку, 2012, 88р.
22. Садыгов М.А. Экстремальные задача для интегрального включения. Препринт №1, Баку, 2013, 129 с.
23. Садыгов М.А. Субдифференциал высшего порядка и оптимизация. LAP LAMBERT Academic Publishing. Saarbrücken, Deutschland, 2014, 359 p.
24. Садыгов М.А. Экстремальные задачи для включений в частных производных. LAP LAMBERT Academic Publishing. Saarbrücken, Deutschland, 2015, 390 p.
25. Толстоногов А.А. Дифференциальные включения в банаховом пространстве. Новосибирск: Наука, 1986, 296 с.
26. Экланд И., Темам Р. Выпуклый анализ и вариационные проблемы. М.: Мир, 1979, 309 с.
27. Федерер Г. Геометрическая теория меры. М.: Наука, 1989, 760 с.
28. Ioffe A., Milosz T. On characterization of $C^{1,1}$ functions. Системный анализ. 2002, №3, с.3-13.
29. Clarke F. Functional Analysis, Calculus of Variations and Optimal Control. Springer-Verlag, London, 2013, 591 p.
30. Rockafellar R.T., Wets R. J-B. Variational Analysis. Springer, 2009, 733p.
31. Sadygov M.A. Higher order conditions in nondiferentiable programming problems. Proc.of IMM of NAS of Azerb., 2017, p.

ОГЛАВЛЕНИЕ